# Jacobi's Bound *

## *Jacobi's results translated in Kőnig's, Egerváry's and Ritt's mathematical languages*

François Ollivier (CNRS)

March 2022

***Dedicated to the memory of Harold W. Kuhn (1925-2014)***

*דער ניגון האָט באַקומען תּיקון*
*און דער למדן אויך!*
י. ל. פּרץ

*The melody got out of Purgatory*
*and the erudite too!*
I. L. Peretz

## רעזומע

יאַקאָביס רעזולטאַטן וועגן דעם חשבון פֿון דער אָרדענונג און פֿון די נאָרמאַלע פֿאָרמעס פֿון אַ דיפֿערענציאַלן סיסטעם ווערן איבערגעזעצט אין דעם פֿאָרמאַליזם פֿון דיפֿערענציאַלער אַלגעברע. מע גיט גאַנצע דערווײַזונגען נאָך יאַקאָביס אַרגומענטן, אין דעם כּמו־רעגולערן פֿאַל. דער עצם־טעאָרעם איז *יאַקאָביס גרענעץ*: די אָרדענונג פֿון אַ דיפֿערענציאַלן סיסטעם $P_1, \ldots, P_n$ איז ניט גרעסער ווי דער מאַקסימום $\mathcal{O}$ פֿון אַלע סך־הכּלען $\sum_{i=1}^{n} a_{i,\sigma(i)}$ פֿאַר אַלע אינדעקסן סובסטיטוציעס $\sigma$, ווּ $a_{i,j} := \mathrm{ord}_{x_j} P_i$, ד״ה *דעם טראָפּישן דעטערמינענט פֿון דער מאַטריצע* $(a_{i,j})$. עס בלײַבט נאָך היפּאָטעטיש אין דעם אַלגעמיינעם פֿאַל. די אָרדענונג איז פּונקט גלייך צו $\mathcal{O}$ אויב און נאָר אויב יאַקאָביס *אָפּגעשניטער דעטערמינאַנט* איז נישט נול.

יאַקאָבי האָט אויך געפֿינען אַ פּאָלינאָמישע צײַט אַלגאָריטם צו חשבונען $\mathcal{O}$, ענלעך צו קונס „אונגערישן מעטאָד“ און אַ מין קירצערער וועג אַלגאָריטמען, שייך צו דעם חשבון פֿון די קלענסטע פּאָזיטיווע גאַנצצאָלן $\ell_i$, אַזוי אַז מען קען רעכענען אַ נאָרמאַלע פֿאָרמע דורך דיפֿערענצירן די גלײַכונג $P_i$ בלויז $\ell_i$ מאָל, אין דעם גענערישן פֿאַל.

מע דערקלערט אויך אַ פּאָר יסודותדיקע רעזולטאַטן וועגן סדר ענערונגען און וועגן די פֿאַרשיידענע נאָרמאַלע פֿאָרמעס וואָס אַ סיסטעם קען זיי נעמען, ווי למשל דיפֿערענציעלע רעזאָלוואַנטן.

## Abstract

Jacobi's results on the computation of the order and of the normal forms of a differential system are translated in the formalism of differential algebra. In the quasi-regular case, we give complete proofs according to Jacobi's arguments. The main result is *Jacobi's bound*, still conjectural in the general case: the order of a differential system $P_1, \ldots, P_n$ is not greater than the maximum $\mathcal{O}$ of the sums $\sum_{i=1}^{n} a_{i,\sigma(i)}$, for all permutations $\sigma$ of the indices, where $a_{i,j} := \mathrm{ord}_{x_j} P_i$, *viz.* the *tropical determinant of the matrix* $(a_{i,j})$. The order is precisely equal to $\mathcal{O}$ iff Jacobi's *truncated determinant* does not vanish.

Jacobi also gave a polynomial time algorithm to compute $\mathcal{O}$, similar to Kuhn's "Hungarian method" and some variants of shortest path algorithms, related to the computation of integers $\ell_i$ such that a normal form may be obtained, in the generic case, by differentiating $\ell_i$ times equation $P_i$.

Fundamental results about changes of orderings and the various normal forms a system may have, including differential resolvents, are also provided.

François Ollivier (CNRS)
LIX UMR 7161
CNRS–École polytechnique
F-91128 Palaiseau CEDEX (France)
Email francois.ollivier@lix.polytechnique.fr

# Contents

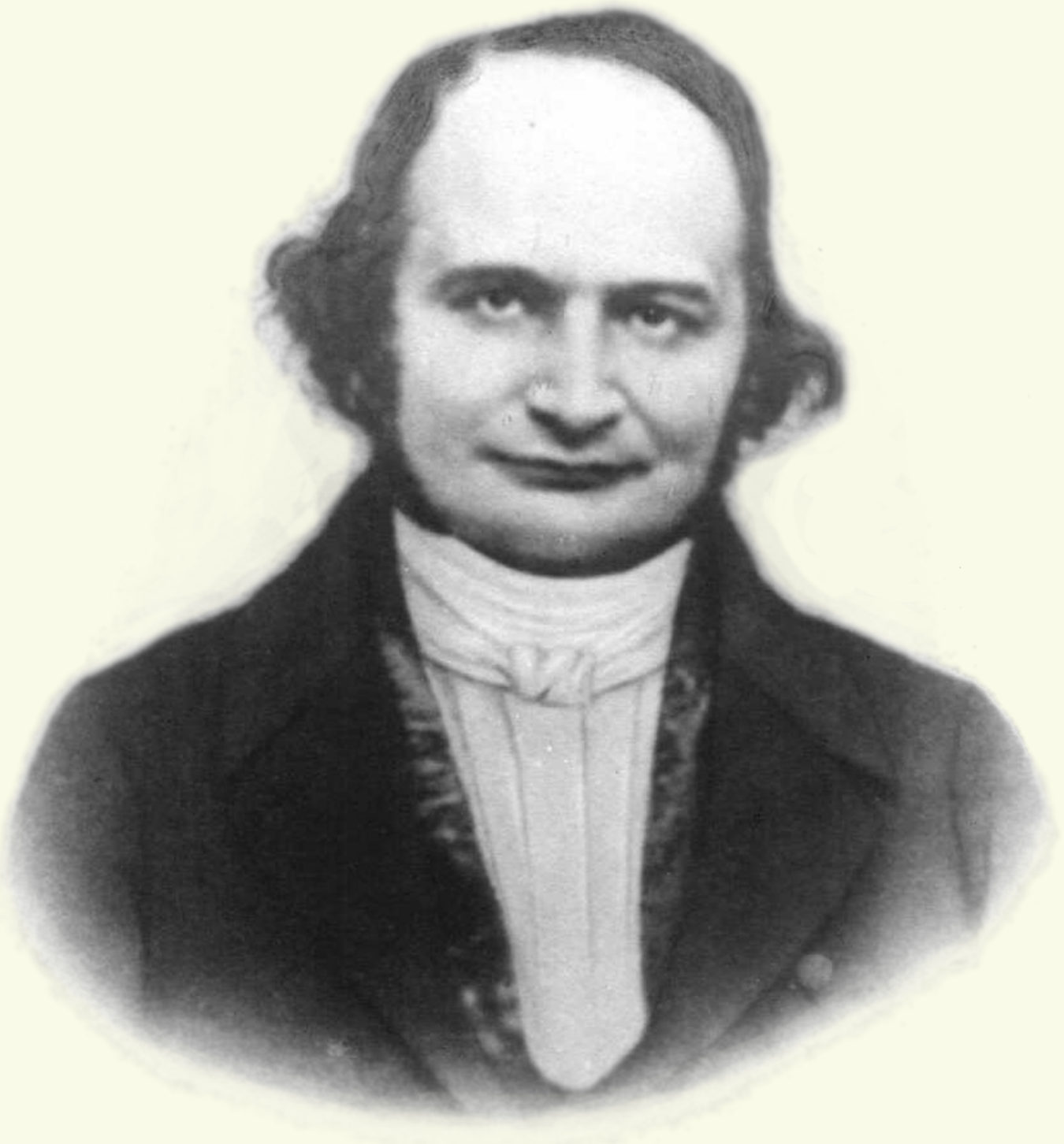

C. G. J. Jacobi.

# Introduction

## History

*Words in index appear in margins* · *Borchardt* · *Cohn, Sigismundus* · *Clebsch* · *order of a prime component* · *Jacobi's number, bound* · *transversal sum* · *Ritt* · *Kondratieva* · *Johnson, Joseph* · *normal form* · *Kuhn* · *Egerváry* · *Martello* · *Cohn, Richard*

In 1865 appeared in Crelle's journal a posthumous paper of Jacobi, edited by Borchardt [39] after a transcription by S. Cohn [II/13 a)], followed by a second one in the volume *Vorlesungen über Dynamik*, edited by Clebsch in 1866 [40]. These two papers contain the following main result, that remains conjectural in the general case: *the order of an ordinary differential system of n equations $P_i$ in n variables $x_{\sigma(i)}$ is, at most, the maximum $\mathcal{O}$ of the "transversal" sums $\sum_{i=1}^{n} \mathrm{ord}_{x_{\sigma(i)}} P_i$ for all permutations $\sigma \in S_n$*. Known as *Jacobi's bound*, it mostly survived during the xx$^{\mathrm{th}}$ century in the differential algebra community, thanks to J.F. Ritt, who gave a first complete proof in the linear case [85]. It was extended by Kondratieva *et al.* [60] to systems satisfying Johnson's regularity hypothesis [53] in the ordinary case and also in the partial differential case [61]. But some important aspects were completely forgotten, such as a simplest normal form reduction, bounds on the order of differentiations required for computing normal forms, including differential resolvent, and a first polynomial time algorithm to solve the assignment problem, *i.e.* in our setting computing the bound, faster than by trying the $n!$ permutations. A similar algorithm was rediscovered by Kuhn [64] in 1955, using previous results of Egerváry [23] (see Martello [72] for more details). R.M. Cohn [13] was the first, in 1983, to mention Jacobi's contribution.

*Jacobi's bound* · *Saltykow* · *Volevich* · *shortest reduction* · *Shaleninov* · *Pryce* · *Jacobi's number* · *tropical determinant* · *Dijkstra; Ford; Bellman*

Jacobi's bound was still mentioned by Saltykow [88] in 1950. In 1960, it was rediscovered in Moscow by Volevich [96] for differential operators and his simplest normal form reduction by Shaleninov [91] in 1990 and Pryce [81] in 2001 for the resolution of implicit DAE's. One may also mention that in modern vocabulary the expression of Jacobi's number $\mathcal{O}$ is known as the *tropical determinant* [71]. Two algorithms introduced by Jacobi to compute his *minimal canon* may be regarded as precursors of Dijkstra's [19] or Ford [25] and Bellman's [4] shortest paths algorithms.

*isoperimetric equations* · *Moritz Jacobi* · *normal form*

It is difficult to know precisely when the manuscripts were written related to Jacobi's bound. Jacobi did not use to date his writings. It seems to be a byproduct of his work on isoperimetric systems, mentioned in letters to his brother Moritz in 1836 [47, 1836 IX. 17., XII. 20., 1837 III. 5.]. The second part of his paper on the last multiplier [41], which appeared in 1845, contains a section devoted to these systems, where he promised to publish later a method for computing normal forms.

*isoperimetric equations*

Proofs are often omitted in Jacobi's manuscripts. The style of some passages may suggest a mathematical cookbook, providing computational methods without justifications, but no examples of precise differential systems are given, only general abstract families of systems like *isoperimetric equations*. It is clear that the efficiency is a constant preoccupation, even if it is not formalized. This work is closely related to Jacobi's interest in mechanics; there was at that time a strong need for fast computational tools, mostly for astronomical ephemerides [32][1].

In his 1840 letter to the Académie des Sciences de Paris [42], Jacobi said that he was working for some years on a publication that included his last multiplier method. One may guess that the various unpublished fragments were intended to take part in this never achieved ambitious project entitled *Phoronomia*. As Ritt [85] guessed, the bound may have been suggested to him by his method for computing normal forms. We refer to our survey [80] or Königsberger [58] for more historical details and also to Saltykow [87] that gives interesting precisions about datings and Jacobi's unachieved great project.

*last multiplier*
*Phoronomia*
*Ritt*
*Königsberger*
*Saltykow*

## Aims of this paper

We present Jacobi's main results related to the order and normal forms of differential systems, using the formalism of differential algebra. We prove them under hypotheses that could have been implicit in Jacobi's work and using, as far as possible, methods suggested in his work. Two particular aspects require attention.

Jacobi gives no detail about the nature of the functions he considers and he does not describe the tools to be used to perform the required eliminations, although he had worked on algebraic elimination (see *e.g.* [38]). We restrict here to polynomial equations. It seems implicit that Jacobi's attention was focused on physical equations, generating prime differential ideals. However, we tried to consider the case of systems defining many components, whenever the extra work remained little. Jacobi's results related to *normal forms* of differential systems will be translated using *characteristic sets* of differential ideals.

*normal form*
*characteristic set*

Jacobi often considers implicit genericity conditions and sometimes gives first a "generic" theorem (*i.e.* a proposition that holds in some Zariski open set) followed by a second theorem describing the cases where the first assertion fails to be true. We will try to provide explicitly such conditions, most of the time expressed by the non vanishing of some Jacobian determinant.

Keeping in mind such particularities of the XIX$^{\text{th}}$ century mathematical style, we recommend the reading of Jacobi's original papers, this text being only a partial commentary, completed with some technical parentheses.

The computation of the tropical determinant occupies a large part of Jacobi's manuscripts and of this paper too. In contrast with those related to differential systems, Jacobi [39] gave very precise proofs of his combinatorial results. This may have dispensed us from longer comments, but a careful study shows that complexity issues require some more attention, as well as the relations between Jacobi's *canons* and Egerváry's *covers*, a notion that allows us to make a link between Jacobi's shortest reduction and the choice of a corresponding ranking on derivatives used in differential algebra algorithms. The implicit, but pioneering,

*canon; Egerváry cover*

[1]Jacobi himself had an experience in practical computing, on a smaller scale and in a different field, when he published his *Canon arithmeticus* [46]. The revision of the half million numbers it contains required the help of friends and relatives [47], including Dirichlet's wife and mother!

*Dirichlet*

introduction of basic concepts and problems of graph theory must also be underlined.

## Content

Section 1 introduces Jacobi's bound in the context of applying his last multiplier method to isoperimetric equations. We limit ourselves here to an informal presentation of the genesis of the results. The next section 2 details Jacobi's algorithm, extended to rectangular matrices and makes a preliminary study of its complexity, considering the special case of maximal matching and the relations with Kuhn's Hungarian algorithm, using Egerváry's covers. Then, in section 3, deeper complexity results are provided, including a $O(sn^2)$ version of Jacobi's algorithm (subsec. 3.2) and a description of some related algorithms. Hopcroft and Karp's algorithm is reinterpreted in Jacobi's setting(3.1) and Jacobi's methods for computing a minimal canon when a canon (subsec. 3.3) or a maximal sum (subsec. 3.4) is known are shown to be equivalent to shortest paths algorithms.

It is followed by a short combinatorial parenthesis, sec. 4, about the "strong bound" (subsec. 4) and reduction to order one (subsec. 4.2), completed with algorithmic hints to get block decompositions (subsec. 4.3). An algebraic parenthesis, sec. 5, is devoted to *quasi-regularity*, a key implicit assumption in Jacobi's proof, and to *"Lazard's lemma"*, that plays a central part in establishing the results on shortest reduction that characterizes some quasi-regular components.

Jacobi's bound is proved in section 6, together with the necessary and sufficient conditions for the bound to be reached, expressed by the system's *truncated determinant* $\nabla$. The shortest normal form reduction is presented in section 7, followed in section 8 by a study of the various possible normal forms of a given system, including a complete description of the possible structures for zero dimensional differential ideals in two variables. Section 9 gives a method for computing a characteristic set for some ordering, knowing one for some other ordering. The last section 10 is devoted to the special case of resolvent computations.

This work is a *translation*, using contemporary mathematical concepts, of Jacobi's findings. In this regard, it is impossible to fix a precise border between our interpretation and Jacobi's ideas. We may only hope to encourage our readers to discover the original texts. To fix the ideas, most of Jacobi's combinatorial results appear with their original proofs, but everything that concerns covers, generalization to rectangular matrices, the strong bound, allowing entries $-\infty$ elements, and computational complexity are our comments.

For only two theorems about differential systems, we were able to find sketches of proofs in Jacobi's writings. The first is Jacobi's bound itself. Some steps are clearly indicated, others are interpretations of incomplete passages in the manuscript that may be contested. All algebraic technicalities, prime or radical ideals, are imposed by our formalism but entirely absent from Jacobi's manuscripts and from the mathematics of his time. This also includes considerations on various notions of "regularity". About the different normal forms of a system, prop. 149 and prop. 152 are our contribution. The generalization in sec. 7 of the shortest reduction idea to what we have called *Egerváry orderings* is ours, but some allusions by Jacobi suggest he may have noticed this possibility. *Egerváry ordering*

Concerning the resolvent computation, which is the second case where we can rely on Jacobi's indications for a proof, we give a modern treatment, proving in th. 175 Jacobi's bound for the resolvent polynomial and characteritic set, with the strong bound that he was certainly not considering. We are able to skip some steps in Jacobi's quite precise sketch, but in lemma 169, everything that does not concern the strong bound, that is mostly part ii) follows Jacobi's construction. *resolvent* *resolvent*

## Notations and conventions

We will often omit "differential" when notions such as differential characteristic sets are considered, and use "algebraic" to stress that a char. set is non differential.

An index of notations is provided p. 104. We denote by $[1, n]$ the set of integers $\{1, 2, \ldots, n\}$, and by $S_n$ its group of permutations. We will consider here equations $P_1, \ldots, P_m$ in the differential polynomials algebra $\mathcal{F}\{x_1, \ldots, x_n\}$, where $\mathcal{F}$ is a differential field of characteristic 0. The perfect differential ideal $\{P\}$ is equal to the intersection of prime components $\bigcap_{i=1}^{r} \mathcal{P}_i$. Jacobi's bound is denoted by $\mathcal{O}$, the $a_{i,j}$ are the entries of a matrix $A$, the notations $\lambda_i$ is introduced in algo. 9, $\nabla$, $J_P$ in def. 76, $S_{s,n}$ in def. 1.

If $\mathcal{A}$ is the characteristic set of a differential or algebraic ideal, we denote by $H_{\mathcal{A}}$ the product of initials $In_i$ and separants $Sep_i$ of its elements and by $In_{\mathcal{A}}$ or $Sep_{\mathcal{A}}$ the products of its initials or separants only. The notation $x_i \ll x_j$ means that $x_i$ and its derivatives are smaller that $x_j$ and its derivatives. To help comparisons with Jacobi's presentation, the elements of char. sets are not assumed to be listed in increasing order, as in Ritt [86, chap. 1].

We write $F(n_1, \ldots, n_p) = O(G(n_1, \ldots, n_p))$, with $F, G : \mathbf{N}^p \mapsto \mathbf{N}$ if there exist constants $A$ and $B$ such that $F \leq AG + B$.

# 1 The last multiplier and isoperimetric systems

## 1.1 The last multiplier

The last multiplier method is first mentioned by Jacobi in a short paper published in French in 1842, entitled "On a new principle of analytical mechanics" [43], followed by a second one in Italian in 1844: "On the principle of the last multiplier and its use as a new general principle of mechanics" [44]. It is not the place here to give details on the subject and we shall limit ourselves to a few hints in order to help understand the link with the genesis of Jacobi's bound. The reader will find illuminating illustrations on classical examples in Nucci and Leach's papers [77, 78]. *last multiplier* *Jacobi's bound* *Nucci and Leach*

Euler's multiplier Jacobi presents his last multiplier as a generalization of Euler's multiplier. If
Lagrange system one has a Lagrange system in two variables:

$$\frac{\mathrm{d}x_1}{f_1(x_1, x_2)} = \frac{\mathrm{d}x_2}{f_2(x_1, x_2)}, \tag{1}$$

Euler's multiplier Euler's multiplier $\mu$ is defined by the property $\mathrm{d}(\mu(f_2\,\mathrm{d}x_1 - f_1\,\mathrm{d}x_2)) = 0$. Knowing the exact differential $\mu(f_2\,\mathrm{d}x_1 - f_1\,\mathrm{d}x_2)$, finding a first integral for the system (1), which is a solution of $f_1\frac{\partial\omega}{\partial x_1} + f_2\frac{\partial\omega}{\partial x_2} = 0$ is reduced to integrations.

In the case of a Lagrange system in $n$ variables,

$$\frac{\mathrm{d}x_1}{f_1(x)} = \cdots = \frac{\mathrm{d}x_n}{f_n(x)}, \tag{2}$$

last multiplier the last multiplier may be defined in the following way. Let $\omega_i$, $1 \le i < n$, be first integrals for (2), any first integral $\omega$ is a solution of

$$\begin{vmatrix} \frac{\partial\omega}{\partial x_1} & \cdots & \frac{\partial\omega}{\partial x_n} \\ \frac{\partial\omega_1}{\partial x_1} & \cdots & \frac{\partial\omega_1}{\partial x_n} \\ \vdots & & \vdots \\ \frac{\partial\omega_{n-1}}{\partial x_1} & \cdots & \frac{\partial\omega_{n-1}}{\partial x_n} \end{vmatrix} = 0.$$

Let us denote by $D_i$ the Jacobian determinant

$$\begin{vmatrix} \frac{\partial\omega_1}{\partial x_1} & \cdots & \frac{\partial\omega_1}{\partial x_{i-1}} & \frac{\partial\omega_1}{\partial x_{i+1}} & \cdots & \frac{\partial\omega_1}{\partial x_n} \\ \vdots & & \vdots & \vdots & & \vdots \\ \frac{\partial\omega_{n-1}}{\partial x_1} & \cdots & \frac{\partial\omega_{n-1}}{\partial x_{i-1}} & \frac{\partial\omega_{n-1}}{\partial x_{i+1}} & \cdots & \frac{\partial\omega_{n-1}}{\partial x_n} \end{vmatrix}.$$

last multiplier A last multiplier $\mu$ is defined by

$$\mu\sum_{i=1}^{n} f_i\frac{\partial}{\partial x_i} = \sum_{i=1}^{n} \pm D_i\frac{\partial}{\partial x_i},$$

or also by

$$\sum_{i=1}^{n}\frac{\partial(\mu f_i)}{\partial x_i} = 0, \quad \textit{i.e.} \quad \mathrm{d}\left(\mu\sum_{i=1}^{n} \pm f_i \bigwedge_{j\neq i} \mathrm{d}x_j\right) = 0,$$

which is for $n = 2$ the definition of Euler multiplier. The quotient of two last multipliers is a first integral, possibly trivial, *i.e.* constant. And the product of a last multiplier by a first integral is a last multiplier.

Jacobi's goal is explicitly exposed in 1842 [43]: having first remarked that, when knowing a multiplier for a system in two variables, the computation of solutions only requires integrations, he claims that his last multiplier method allows to generalize this result to any system of ordinary differential equations in $n$ variables, provided that one already knows $n-2$ first integrals. A new first integral $\omega_{n-1}$ is then easily obtained by quadrature, as in the case of 2 variables. *first integral*

Knowing $n-2$ independent first integrals seems of course a very unlikely circumstance, but in 1840 Jacobi [42] insisted on the importance of a remark of Poisson [45], providing a method to compute a sequence of new first integrals, for any conservative mechanical system, that already possesses two first integrals, independently of energy. *Poisson*

The definition of a last multiplier obviously depends on the choice of a family of $n-2$ first integrals or of the coordinate functions $x_i$. For a system $x_i' = f_i(x)$, the multiplier is indeed given by the explicit formula:

$$\mu = e^{-\int \sum_{i=1}^{n} \frac{\partial f_i}{\partial x_i}\, \mathrm{d}t},$$

which may be interpreted as the inverse of a Wrońskian, expressing the variation of a volume form along a trajectory, so that the last multiplier $\mu_2$ associated to new coordinates $y_i$ must satisfy *Wrońskian*

$$\mu = \left|\frac{\partial y_i}{\partial x_j}\right| \mu_2,$$

a relation that appears (up to a logarithm) in manuscript [II 23 a)] [40, formula (7) p. 40].

The application of the last multiplier method thus requires the knowledge of a normal form for a system of equations and the result will depend on the chosen normal form, casting some light on Jacobi's interest for the various normal forms a given system may possess and for differential elimination.

## 1.2 Isoperimetric equations

In 1844 and 1845, Jacobi published in two parts a 135-page paper [41], describing his last multiplier method for the integration of differential systems. Among the examples of applications he gives, stands the *isoperimetric problem.* *isoperimetric equations*

"Let $U$ be a given function of the independent variable $t$, the dependent ones $x$, $y$, $z$ etc. and their derivatives $x'$, $x''$, etc., $y'$, $y''$, etc., $z'$, $z''$, etc. etc. If we propose the problem of determining the functions $x$, $y$, $z$ in such a way that the integral

$$\int U \;\mathrm{d}t$$

be *maximal or minimal* or more generally that the differential of this integral vanishes, it is known that the solution of the problem depends on the integration of the system of differential equations:

$$\begin{aligned}
0 &= \frac{\partial U}{\partial x} - \frac{\mathrm{d}\frac{\partial U}{\partial x'}}{\mathrm{d}t} + \frac{\mathrm{d}^2\frac{\partial U}{\partial x''}}{\mathrm{d}t^2} - \text{etc.},\\
0 &= \frac{\partial U}{\partial y} - \frac{\mathrm{d}\frac{\partial U}{\partial y'}}{\mathrm{d}t} + \frac{\mathrm{d}^2\frac{\partial U}{\partial y''}}{\mathrm{d}t^2} - \text{etc.},\\
0 &= \frac{\partial U}{\partial z} - \frac{\mathrm{d}\frac{\partial U}{\partial z'}}{\mathrm{d}t} + \frac{\mathrm{d}^2\frac{\partial U}{\partial z''}}{\mathrm{d}t^2} - \text{etc. etc.},
\end{aligned}$$

I will call these in the following *isoperimetric differential equations ...*"
[GW IV, p. 495]

For simplicity, we write $x_1$, ..., $x_n$, instead of $x$, $y$, $z$, etc. and denote by $P_i = 0$ the $i^{\text{th}}$ isoperimetric equation. Jacobi noticed the difficulty of applying his last multiplier method if he could not first reduce the system to a normal form (see also [39, first section]). If the highest order derivative of $x_i$ in $U$ is $x_i^{(e_i)}$, the order of $x_j$ in the $i^{\text{th}}$ isoperimetric equation is at most $e_i + e_j$. If the $e_i$ are not all equal to their maximum $e$, then we cannot compute a normal form without using *auxiliary equations* obtained by differentiating the $i^{\text{th}}$ isoperimetric equation $\lambda_i$ times, and a first problem is to determine minimal suitable values for the $\lambda_i$. In 1845, Jacobi had clearly in mind a thorough study of normal form computations for he wrote: "I will expose in another paper the various ways by which this operation may be done, for this question requires many remarkable theorems that necessitate a longer exposition." [GW IV, p. 502]

Jacobi's method for computing a normal form may be sketched in the following way. Assume that the Hessian matrix $\left(\partial^2 U/\partial x_i^{(e_i)}\partial x_j^{(e_j)}\right)$ has a non-zero determinant. We may further assume, up to a change of indices, that the sequence $e_i$ is non decreasing and that the principal minors of the Hessian are non-zero.

From the first isoperimetric equation $P_1$, as

$$\frac{\partial P_1}{\partial x_1^{(2e_1)}} = \pm\frac{\partial^2 U}{\partial\left(x_1^{(e_1)}\right)^2} \neq 0,$$

one will deduce on some open set, using the implicit function theorem, an expression

$$x_1^{(2e_1)} := F_1(x_1, \ldots, x_1^{(2e_1-1)}, x_2, \ldots, x_2^{(e_1+e_2)}, \ldots, x_n, \ldots, x_n^{(e_1+e_n)}).$$

Using the first equation and its derivatives up to the order $e_2 - e_1$, together with the second equation, one may invoke again the implicit function theorem, using the fact that the Jacobian matrix of $P_2$ and $P_1^{(e_2-e_1)}$, with respect to the derivatives $x_1^{(e_1+e_2)}$ and $x_2^{(2e_2)}$, is equal to the second principal minor of the Hessian of $U$, which is assumed not to vanish. One deduces an expression

$$x_2^{(2e_2)} := F_2(x_1, \ldots, x_1^{(2e_1-1)}, x_2, \ldots, x_2^{(2e_2-1)}, x_3, \ldots, x_3^{(e_2+e_3)}, \ldots, x_n, \ldots, x_n^{(e_2+e_n)}).$$

Iterating the process, we get at the end a last expression

$$x_n^{(2e_n)} := F_n(x_1, \ldots, x_1^{(2e_1-1)}, \ldots, x_n, \ldots, x_n^{(2e_n-1)}),$$

that may be obtained using each isoperimetric equation $P_i = 0$ and its derivatives up to order $\lambda_i := e_n - e_i$.

In this normal form, each variable $x_i$ appears with the order $2e_i$, so that the order of the system is $2\sum_{i=1}^n e_i$. This appears to be both a special case of Jacobi's bound (see sec. 6) and of Jacobi's algorithm for computing normal forms (sec. 7), using the minimal number of derivatives of the initial equations, provided that the "system determinant" or "truncated determinant" $\nabla$ (see def. 76), here equal to the Hessian of $U$, does not vanish. In case of arbitrary equations $P_i$, for which $a_{i,j} := \mathrm{ord}_{x_j} P_i$ can take any value, things become more complicated, starting with the computation of the bound $\max_\sigma \sum_{i=1}^n a_{i,\sigma(i)}$, that is the subject of our next section. But we see how this particular simple example may have suggested the whole theory.

*Jacobi's bound* *Jacobi's algorithm* *truncated determinant* *system determinant*

In section 2. of [39], we have restored a passage of [II/13 b), f° 2200] that quotes the isoperimetric equations as an example for which all the transversal sums have the same value.

*transversal sum*

# 2 Computing the bound. Jacobi's algorithm

*In algorithms, we will assume that matrices are represented by some array structure, so that one may get or change the value of some entry $a_{i,j}$ with constant cost.*

## 2.1 Preliminaries

The *assignment problem* has been first considered by Monge in 1781 [75], in the special case of the transportation problem (moving things from initial places to new places, minimizing the sum of the distances) and in a continuous setting (digging excavations somewhere in order to create some embankment somewhere else). Before stating the discrete version, Jacobi [39, § 3] has written that it was "also worth to be considered for itself…" a premonitory intuition. Such kind of problems reappeared indeed in the middle of the xx$^{\text{th}}$ Century—when Jacobi's method was forgotten—in the following form: $n$ workers must be assigned to $n$ tasks; assuming that the worker $i$ has a productivity $a_{i,j}$ when affected at task $j$, how can we find an affectation $j = \sigma(i)$ that maximizes the sum of productivity indices?

*assignment problem; Monge*

At a meeting of the American Psychological Association in 1950, a participant described the following reaction: *"[he] said that from the point of view of a mathematician there was no problem. Since the number of permutations was finite, one had only to try them all and chose the best. […] This is really cold comfort for the psychologist, however, when one considers that only ten men and ten jobs mean over three and a half million of permutations."*[2] quoted by Schrijver [93, p. 8]. Jacobi

[2] From some optimistic standpoint, it could have been a way to escape ethical issues raised by the use of psychology in management.

did not consider the brute force method as a solution… and he gave a polynomial time algorithm!

The assignment problem also appears as a weighted generalization of the marriage problem marriage or maximal bipartite matching problem: a graph describing couples of bipartite matching compatible boys and girls is represented by a $s \times n$ matrix of zeros and ones. The problem of computing the maximal number of compatible couples between these $s$ boys and $n$ girls amounts to computing a maximal transversal sum, according to def. 1.

Kuhn Kuhn's [64] and Jacobi's algorithms are quite similar. The main difference is the following. Jacobi remarks that if the columns of the matrix admit maxima placed in different rows, then their sum is the maximum to be found. He will then add minimal constants $\lambda_i$ to the rows in order to get a matrix with this property. Kuhn considers integers $\mu_i$ and $\nu_j$, such that $a_{i,j} \leq \mu_i + \nu_j$, with $\sum_{i=1}^{n}(\mu_i + \nu_i)$ cover minimal; this is called a *minimal cover*. He then uses Egerváry's theorem [23, 93]: Egerváry $\sum_{i=1}^{n}(\mu_i + \nu_i) = \max_{\sigma\in S_n} \sum_{i=1}^{n} a_{i,\sigma(i)}$. The relations between canons and minimal covers will be investigated in subsec. 2.4. On the topic of relations between Jacobi's algorithm and Kuhn's Hungarian method, I cannot do better than referring Kuhn to Kuhn's excellent—and moving—presentation [65]. See also rem. 21.

Some of Jacobi's results could be extended with no extra work to the case of underdetermined system underdetermined systems. This is why we will expose his algorithm in the case of an $s \times n$ matrix $A$, with $s \leq n$. Entries are assumed to belong to a totally ordered additive commutative group $\mathrm{M}$, *i.e.* a commutative group with a total order such that $x \geq y \iff x - y \geq 0$. Some definitions and results in this subsection may make sense on a weaker structure, *e.g.* a totally ordered additive monoïd, but as it is, Jacobi's algorithm requires subtractions. Ex. 8 shows that an exponential complexity is unavoidable with a partial order. The special case of $-\infty$ entries strong bound (the "strong bound") will be considered in subsec. 4.1.

**Definition 1.** — *Let $s \leq n$ be two integers, we denote by $S_{s,n}$ the set of injections $\sigma : [1, s] \mapsto [1, n]$.*

*Let $A$ be a $s \times n$ matrix of elements in $\mathrm{M}$, a totally ordered commutative monoïd,* Jacobi's number *the* Jacobi number *of $A$ is defined by the formula*

$$\max_{\sigma\in S_{s,n}} \sum_{i=1}^{s} a_{i,\sigma(i)}$$

*and is denoted by $\mathcal{O}_A$. If $s > n$, we define $\mathcal{O}_A := \mathcal{O}_{A^t}$.*

maximal element *Without further specification, a* maximum *$a_{i,j}$ in $A$, is understood as being a* maximum transversal maxima *maximal element in its column, i.e. such that $\forall 1 \leq i' \leq s$, $a_{i,j} \geq a_{i',j}$. We call* transversal maxima *a set of maxima placed in pairwise different rows and columns. It is said to be a* maximal set of transversal maxima *if there is no set of transversal maxima with more elements in $A$.*

canon *Let $\ell$ be a vector, we denote by $A + \ell$ the matrix $(a_{i,j} + \ell_i)$. We call a* canon *a matrix $A + \ell$, with $\ell_i \geq 0$, $1 \leq i \leq s$, that possesses $s$ transversal maxima and also the vector $\ell$ itself.*

*We will use on canons the partial order defined by $\ell \le \ell'$ if $\ell_i \le \ell'_i$, for all $1 \le i \le s$.*

***Remarks.* — 2)** In the case of a preorder or a partial order, a maximum is such that no element is strictly greater, and many maxima may then exist. Example 8 shows that, using a partial order, one may have to consider the $n!$ sums to find a Jacobi number for a square $n \times n$ matrix. *partial order*

**3)** In the case of a square matrix, with elements in a totally ordered commutative group, if $\ell$ is a canon, then $a_{i,\sigma(i)} + \ell_i$ is a maximal set of transversal maxima iff $\sum_{i=1}^{n} a_{i,\sigma(i)}$ is the *maximal transversal sum* we are looking for.

For $s < n$, one may complete $A$ with $n - s$ rows of zeros, which reduces the problem to the case of a square matrix.

**4)** Any totally ordered monoïd M, which is furthermore such that $\forall(a, b, c) \in \mathrm{M}^3$, $a < b$ implies $a + c < b + c$, may be seen as a submonoïd of the ordered group defined as the quotient of $\mathrm{M}^2$ by the congruence associated to the equivalence relation $(a, b) \equiv (c, d)$ if $a + d = b + c$.

Working with an arbitrary monoïd, the knowledge of a canon is not enough to compute a maximal transversal sum, as shown by the next example.

***Example* 5.** — Let M be the totally ordered commutative monoïd represented by $\mathbf{Z} \times \infty \cup \mathbf{Z}$ with the convention $a\infty + b = a\infty$ (meaning that $b \ne 0\infty + b = 0\infty$) and an order such that, for all $b \in \mathbf{Z}$, $0\infty < b < 1\infty$. Let

$$A = \begin{pmatrix} \mathbf{1} & -2\infty & 0 \\ \infty & \mathbf{1} & \infty \\ 0 & -2\infty & \mathbf{1} \end{pmatrix},$$

then, the unique maximal transversal sum is the diagonal $1+1+1$ and the unique minimal canon is $\lambda = (\infty, 0, \infty)$. However, the maximal transversal sum $(0 + \infty) + 1 + (0 + \infty) = 2\infty$ of $A + \lambda$, does not correspond to a transversal sum of $A$.

*rectangular matrix* In the group case, Jacobi's algorithm (algo. 9), applied to a rectangular matrix, returns the minimal canon[3] $\lambda$, that will be used in section 7 to compute the shortest reduction to normal form, $\lambda_i$ being the minimal number of times one needs to differentiate $P_i$ in order to compute a normal form (under some genericity hypotheses). But, when $s < n$, the sum of the corresponding maxima, and so the order of this normal form, may fail to be equal to $\mathcal{O}_A$. *Jacobi's algorithm* *normal form* *shortest reduction*

***Example* 6.** — Consider the matrix

$$\begin{pmatrix} \mathbf{1} & 0 & \mathbf{3} & \mathbf{4} \\ 0 & \mathbf{1} & 2 & 0 \end{pmatrix}.$$

The minimal canon is $\lambda = (0, 0)$; however the sum of the corresponding 3 pairs of transversal maxima are 2, 4 and 5 whereas the maximal transversal sum is 6. To

[3]The existence of a canon is a consequence of algo. 9, unicity is shown in prop. 7 below.

find it, we may add two rows of zeros:

$$\begin{pmatrix} 1 & 0 & 3 & 4 \\ 0 & 1 & 2 & 0 \\ 0 & 0 & 0 & 0 \\ 0 & 0 & 0 & 0 \end{pmatrix}.$$

Then the minimal canon is $\ell = (0, 1, 2, 2)$:

$$\begin{pmatrix} 1 & 0 & 3 & \mathbf{4} \\ 1 & 2 & \mathbf{3} & 1 \\ 2 & \mathbf{2} & 2 & 2 \\ \mathbf{2} & 2 & 2 & 2 \end{pmatrix}.$$

The following proof of the unicity of the minimal (simplest) canon, assuming canons do exist, is due to Jacobi [39, th. IV, § 2].

*minimal canon* **Proposition 7.** — *Let $A$ be a $s \times n$ matrix of elements in a totally ordered monoïd* $\mathrm{M}$ *and $\ell$, $\ell'$ be two canons:*

*i) The $s$-tuple $\ell''$ defined by $\ell''_i := \min(\ell_i, \ell'_i)$ is a canon for $A$.*

*ii) If a minimal canon exists for the ordering defined by $\ell \le \ell'$ if $\ell_i \le \ell'_i$, for all $1 \le i \le s$, it is unique.*

Proof. — Let $I := \{i \in [1, s] | \ell_i < \ell'_i\}$ and $\bar{I} := [1, s] \setminus I$. Let $\sigma$ and $\sigma'$ be the elements of $S_{s,n}$ corresponding to maximal sets of transversal maxima for the canons $A + \ell$ and $A + \ell'$. We define $\sigma''(i) = \sigma(i)$ if $i \in I$ and $\sigma''(i) = \sigma'(i)$ if not, so that $\sigma''(i) \neq \sigma''(i')$ if $i \neq i'$ are both in $I$ or both in $\bar{I}$. Furthermore, if $i \in I$ and $i' \in \bar{I}$, then $a_{i',\sigma'(i')} + \ell'_{i'} \ge a_{i,\sigma'(i')} + \ell'_i$ (as $a_{i',\sigma'(i')} + \ell'_{i'}$ is maximal in $A + \ell'$) and $a_{i,\sigma'(i')} + \ell'_i > a_{i,\sigma'(i')} + \ell_i$ (as $i \in I$), so that

$$a_{i',\sigma'(i')} + \ell'_i > a_{i,\sigma'(i')} + \ell_i.$$

In the same way, $a_{i,\sigma(i)} + \ell_i \ge a_{i',\sigma(i)} + \ell'_{i'}$, meaning that $a_{i,\sigma''(i)} + \ell''_i$ is maximal in its column in the matrix $A + \ell''$. Those inequalities also imply that $\sigma''(i) = \sigma(i) \neq \sigma'(i') = \sigma''(i')$, so that $\sigma$ is an injection. This completes the proof of i), of which ii) is a straightforward consequence. ∎

Before exposing the algorithm, a last example will show that its polynomial complexity cannot be achieved with a partial order.

***Example* 8.** — We work on a partially ordered monoïd $\mathrm{M}$, generated by elements $a_{i,j}$, $1 \le i, j \le n$ and $\infty$, with rules $a_{i,j} + a_{i,k} = \infty$, $a_{i,j} + a_{k,j} = \infty$. The element $0$ is minimal, $\infty$ maximal, no order relation is defined between non trivial transversal sums of less than $n$ elements $a_{i,j}$. Then, for the matrix $A := (a_{i,j})$, the order between transversal sums of $n$ elements can be chosen arbitrarily and there is no shorter way to find the maximal transversal sum than trying all $n!$ possibilities.

***In the sequel,*** $\mathrm{M}$ ***will always be a totally ordered commutative group.***

## 2.2 Jacobi's algorithm

See [39, § 3] for Jacobi's proof of the algorithm and [40, § 1] for a detailed example.

**Algorithm 9.** — **Input:** an $s{\times}n$ matrix $A$ of elements $a_{i,j}$ in some totally ordered commutative group M. We assume $s \le n$. The case $s = 1$ is trivial, so we assume $s \ge 2$. *Jacobi's algorithm*
**Output:** the minimal canon $\lambda$ of $A$.

***Step 1.*** **(Preparation process)** — Increase each row $i_0$ of the least element $\ell_{i_0} \in$ M such that one of its elements become maximal (in its column), *i.e.* $\ell_{i_0} := \min_{j=1}^{n} \max_{i=1}^{s} a_{i,j} - a_{i_0,j}$. It produces a new matrix $A^0 = A + \ell^0$ such that each row possesses a maximal element. The number of transversal maxima in $A^0$ is at least 2, that corresponds to the case where all elements in row $i$ and all elements in column $j$ are maximal (except perhaps the element $a_{i,j} + \ell_i^0$). We find 2 transversal maxima. If $s = 2$, the problem is solved. *This step requires $O(sn)$ operations.* *preparation process*

If $s > 2$, we enter step 2 with $A^0$, $\ell^0$ and a set of exactly $r := 2$ transversal maxima. If $s = 2$, we have finished and return $A^0$.

After $\kappa$ previous iterations of step 2), its input is $A^\kappa = A + \ell^\kappa$; in the algorithm, we denote its entries by $a_{i,j}$ for simplicity, according to the computer program convention $a_{i,j} := a_{i,j} + \ell^\kappa$.

***Step 2.*** — **a)** For readability, we may reorder the rows and columns, so that the transversal maxima in $A^\kappa$ are the elements $a_{i,i}$ for $1 \le i \le r < s$. *Left* (resp. *right*) columns are columns $j \le r$ (resp. $j > r$). *Upper* (resp. *lower*) rows are rows $i \le r$ (resp. $i > r$), as below. *column (left, right)* *row (upper, lower)*

We define the *starred elements* of $A^\kappa$ as being the left transversal maxima[4] $a_{i,i}$. *starred elements*

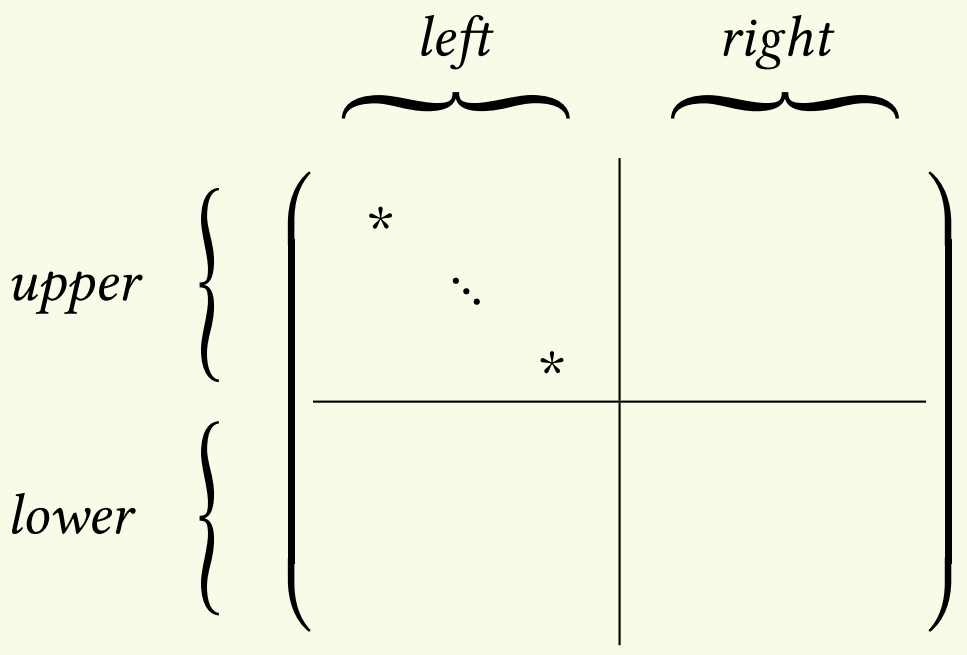

**b)** Assume that there is a maximal element located in a right column and a lower row. We can add it to the set of transversal maxima. If it now contains $s$ elements, the process is finished. If not, we repeat step 2.

[4] Jacobi defined also the maximal elements in right columns as "starred"; we prefer to reserve this denomination to left transversal maxima to underline the specific roles played by these two sets of maxima in the algorithm.

*path* **c)** We say that *there is a path*[5] from row $i$ to row $i'$ if there is a starred maximum in row $i$, equal to some element of row $i'$ located in the same column, or recursively if there is a path from row $i$ to row $i''$ and from row $i''$ to row $i'$. We also define *class (first, second, third)* *first class rows* as being upper rows with at least a right maximal element, or recursively rows to which there is a path from a first class row. The construction of the set of first class rows, together with paths to them from rows with a right maximum may be done in $O(sn)$ operations, using an array $F$ of Booleans with $F_i :=$ true if row $i$ belongs to the first class (we cannot afford looking into a list).

**d)** If there is no lower row of the first class, we go to substep e) p. 17.

Assume that there is a lower row of the first class, then there is a path to it from an upper row $i_0$ containing a right maximal element $a_{i_0,\beta}$. Let it be a path of length $p$, consisting of rows $i_0$ to $i_p$, so that for all $0 \le \alpha < p$ the element $a_{i_{\alpha+1},i_\alpha}$ is equal to the starred element $a^*_{i_\alpha,i_\alpha}$. We can then construct a set of $r+1$ transversal maxima by replacing $a^*_{i_\alpha,i_\alpha}$ by $a_{i_{\alpha+1},i_\alpha}$, for $0 \le \alpha < r$, and adding $a_{i_0,\beta}$ to the list, as illustrated in the figure below, where the new increased set of transversal maxima is written in red letters.

$$\left(\begin{array}{cccc|c}
a^*_{i_0,i_0} & \cdots & & & \mathbf{a_{i_0,\beta}} \\
\| & & & & \\
\mathbf{a_{i_1,i_0}} & \cdots & a^*_{i_1,i_1} & & \\
 & & \| & & \\
 & & \mathbf{a_{i_2,i_1}} & \cdots & \\
 & & & \vdots & \\
 & & & \cdots\ a^*_{i_{p-1},i_{p-1}} & \\
\hline
 & & & \| & \\
 & & & \mathbf{a_{i_p,i_{p-1}}} & \\
 & & & & 
\end{array}\right) \tag{3}$$

If $r+1 = s$, we have finished and return $A^\kappa$, if not, we repeat step 2 a), c) and d) until no lower row of the first class is found.

The next lemma is given by Jacobi in [39, sec. 3].

*transversal maxima* **Lemma 10.** — *The maximal number of transversal maxima in $A^\kappa$ is $r$ iff there is no lower right maximum, nor lower row of the first class.*

PROOF. — The substep d) of the algorithm above proves that the given condition is necessary. Let us assume that there is no lower right maximum, nor lower row of the first class, but that there exists a set $T$ of $r' > r$ transversal maxima $a_{p_i,q_i}$.

[5]This notion is closely related to that of *increasing path*, as defined in Hopcroft and Karp [33] (see 3.1), which explains the choice of that word to translate *transitum datur* in [39].

As there are only $r$ left columns, $h \geq r' - r$ of them, say $a_{p_i,q_i}$ $1 \leq i \leq h$, are upper right maxima. From their first class rows, one can build paths as above, starting from rows $p_{i_k}$, $1 \leq k \leq h$, and considering only maximal elements $a_{p_i,p_i}$ and $a_{p_i,q_i}$, $1 \leq i \leq r'$. The sets of rows in such paths are disjoint, for $T$ is a set of transversal maxima. As there is no lower first class row, these $h$ paths must end with some row $j$, such that there is no $q_i = j$, $1 \leq i \leq r'$, *i.e.* no maximum in $T$ located in the same column as $a_{j,j}$. So, we associate bijectively to any right maximum in $T$ a left column with no maximum in $T$ and there are in $T$ $h$ right elements and at most $r - h$ left elements: $r' \leq h + r - r = r$, a contradiction. ■

**Definition 11.** — *We define the rows of the* third class *as being the lower rows and all the rows from which there is a path to a lower row. The rows not in the first or third class form the* second class.

**e)** We increase the third class rows by the smallest value $\mu$ such that one of their elements $a_{i,j} + \mu$ become equal to a right maximum or to a starred element $a^*_{i,i}$ located in some first or second class row. *This may be done in at most $O(sn)$ operations.*

We then iterate step 2 with a new matrix $A^{\kappa+1} = A + \ell^{\kappa+1}$, such that $\ell^{\kappa+1} > \ell^{\kappa}$.

If this element $a_{i,j}+\mu$ is equal to a right maximum, then at step 2. b), the number of transversal maxima will be increased. Or else, if $a^*_{i,i}$ belongs to a second class row, this row will go to the third class and the cardinal of the second class will decrease. If it belongs to the first class, then at the next step there will be a lower right maximal element (if it is in a right column) or a first class lower row (if it belongs to a left one), so that the number of transversal maxima will increase. Let $r$ be the number of starred elements, there are at least 1 first class, 1 upper third class row and at most $r - 2$ second class rows, so that we need at most $r - 1$ iterations in order, first, to exhaust the second class, and then increase the number of transversal maxima, which can occur at most $s-2$ times. So step 2 is iterated at most $\sum_{r=2}^{s-1}(r-1) = (s-1)(s-2)/2$ times before the algorithm returns the requested result. This bound is sharp: see rem. 15.

As the complexity of step 2 is $O(sn)$, *the whole algorithm requires $O(s^3n)$ elementary operations. This gives for $s = n$ a $O(n^4)$ complexity*, which corresponds to that of Kuhn's original Hungarian algorithm (see Burkard *et al.* [7, Ch. 4.1 p. 77]), which was first proved in the presentation of Munkres [66].

We need here an easy lemma to evaluate in the non negative integer case the complexity of the algorithm, using by anticipation the minimality of the canon $A + \ell$, returned by the algorithm (lem. 14 below).

***Lemma 12.*** — *i) We have* $\min_{i=1}^{s} \ell_i = 0$.

*ii) If $a_{i,j} \geq 0$, for $1 \leq i \leq s$ and $1 \leq j \leq n$, then for all $1 \leq i \leq s$, $\ell_i \leq C :=$* $\max_{i=1}^{s} \max_{j=1}^{n} a_{i,j} \leq \mathcal{O}_A$.

Proof. — i) This is a consequence of the minimality of $\ell$, as $\ell - \min_{i=1}^{s} \ell_i$ is also a canon.

ii) Assume that $\mathcal{O}_A = \sum_{i=1}^{s} a_{i,\sigma(i)}$ and let $\ell_{i_0} = \min_{i=1}^{s} \ell_i = 0$. Then, for all $1 \le i \le s$, $a_{i_0,\sigma(i_0)} \ge a_{i,\sigma(i_0)} + \ell_i$, so that $\ell_i \le a_{i,\sigma(i)} - a_{i_0,\sigma(i_0)}$ This quantity is bounded by both $C$ and $\mathcal{O}_A$. Furthermore, if $C = a_{i_0,j_0}$, let $\sigma$ be an injection such that $\sigma(i_0) = j_0$, then $C \le \sum_{i=1}^{n} a_{i,\sigma(i)} \le \mathcal{O}_A$. ■

So, if the elements of $A$ are non negative integers bounded by $C$, then, by lem. 12, the elements of $A^\kappa$ encountered during iteration $\kappa + 1$ of step 2) are bounded by $2C$, so that each integer operation requires $O(\ln C)$ bit operations.

minimal canon **Theorem 13.** — *The above algorithm returns the minimal canon $\lambda$ of $A$ in at most $O(s^3 n)$ elementary operations. Assuming that the elements in the matrix are non negative integers $a_{i,j} \le C$, it requires at most $O(s^3 n \log C)$ bit operations.*

Proof. — The termination and complexity of the algorithm have already been proved. We only have to show that $\ell$ is the smallest canon $\lambda$.

The proof, that follows Jacobi's [39, p. 21], relies on the following lemma. The $a_{i,j}$ denote here the entries of the input matrix $A$.

***Lemma 14.*** — *Let $\lambda^\kappa$ be the minimal canon for the matrix $A^\kappa := A + \ell^\kappa$, encountered during iteration $\kappa + 1$ of step 2), assume that $a_{i,i} + \ell_i^\kappa$ $1 \le i \le r < s$ form the set of transversal maxima in $A^\kappa$, with respect to which the classes are defined at step e) and that there is no lower right maxima nor first class lower row. Then there is no* unchanged row *unchanged row of the third class in $A^\kappa$, i.e. a third class row of index i with $\lambda_i^\kappa = 0$.*

Proof of the lemma. — The $a_{i,i} + \ell_i^\kappa$ $1 \le i \le r$ form a maximal set of transversal maxima in $A^\kappa$, according to lem. 10. Let $a_{i,\sigma(i)} + \ell_i^\kappa + \lambda_i^\kappa$, $1 \le i \le s$, be a maximal set of transversal maxima in $A^\kappa + \lambda^\kappa$.

If row $i$ is an unchanged row of the third class, the element $a_{i,\sigma(i)} + \ell_i^\kappa$ is maximal (in its column) in $A^\kappa + \lambda^\kappa$, and so it is also maximal in $A^\kappa$. It cannot be an upper right element, for then the row $i$ would be of the first class, and it cannot be lower right. So, $1 \le \sigma(i) \le r$.

Let $H$ denote the set of integers $1 \le i \le r$ such that row $i$ is an unchanged row of the third class. For $i \in H$, the elements $a_{i,\sigma(i)} + \ell_i^\kappa$ and $a_{\sigma(i),\sigma(i)} + \ell_{\sigma(i)}^\kappa$ are both maximal elements of the column $\sigma(i)$. So, the row $\sigma(i)$ must be unchanged too and, as there is a path from it to row $i$, it belongs to the third class: $\sigma : H \mapsto H$ is a bijection. Hence, there is no unchanged lower row $i'$ of the third class, for we would have $\sigma^{-1}(i') \in H$.

Let the row $i_0$ be an unchanged row of the third class. Due to the third class definition, we can find a sequence of third class rows $i_\alpha$, $0 \le \alpha \le p$, such that:

i) $a_{i_{\alpha+1},i_\alpha} + \ell_{i_{\alpha+1}}^\kappa = a_{i_\alpha,i_\alpha} + \ell_{i_\alpha}^\kappa$, $0 \le \alpha < p$;

ii) rows $i_\alpha$, $0 \le \alpha < p$ are upper rows;

iii) row $i_p$ is lower.

Using i), we prove by recurrence that all rows $i_\alpha$, $0 \le \alpha \le p$ are unchanged. As row $i_p$ is lower, we arrive to a final contradiction, that concludes the proof of the *lemma.* ■

Each row of a canon must contain a maximal element. So $\lambda \geq \ell^0$, which is the canon produced by the preparation process. Recursively assuming that $\ell^\kappa \leq \lambda$, as there is no unchanged row of the third class, and as, during step 2) e) we increase third class rows by the minimal integer requested to change the class partition, we have $\ell^{\kappa+1} \leq \lambda$ and the canon $\ell$ returned by the algorithm must be the minimal canon $\lambda$, which is unique according to prop. 7. This concludes the proof of the theorem. ∎

***Remarks. — 15)*** Let the $n \times n$ integer matrix $A$ be defined by $a_{i,j} = (n-1)^2$ $-(i-1)(j-1)$, one shall apply step 2 precisely $(n-1)(n-2)/2$ times. E.g., for $n = 4$, the matrix is:

$$\begin{array}{c} \mathrm{I} \\ \mathrm{III} \\ \mathrm{III} \\ \mathrm{III} \end{array} \begin{pmatrix} 9 & \mathbf{9} & 9 & 9 \\ \mathbf{9} & 8 & 7 & 6 \\ 9 & 7 & 5 & 3 \\ 9 & 6 & 3 & 0 \end{pmatrix},$$

where we have indicated the classes of the rows on he left, the starred maxima being in **bold**. Step 2 shall be applied 3 times and here is the sequence of matrices it produces, with the increment of each row, the last matrix being the canon.

$$\begin{array}{c} \mathrm{I} \\ \mathrm{II} \\ \mathrm{III} \\ \mathrm{III} \end{array} \begin{pmatrix} 9 & 9 & \mathbf{9} & 9 \\ 10 & \mathbf{9} & 8 & 7 \\ \mathbf{10} & 8 & 6 & 4 \\ 10 & 7 & 4 & 1 \end{pmatrix} \begin{array}{c} 0 \\ 1 \\ 1 \\ 1 \end{array}, \quad \begin{array}{c} \mathrm{I} \\ \mathrm{III} \\ \mathrm{III} \\ \mathrm{III} \end{array} \begin{pmatrix} 9 & 9 & \mathbf{9} & 9 \\ 10 & \mathbf{9} & 8 & 7 \\ \mathbf{11} & 9 & 7 & 5 \\ 11 & 8 & 5 & 2 \end{pmatrix} \begin{array}{c} 0 \\ 1 \\ 2 \\ 2 \end{array}, \quad \begin{pmatrix} 9 & 9 & 9 & \mathbf{9} \\ 11 & 10 & \mathbf{9} & 8 \\ 12 & \mathbf{10} & 8 & 6 \\ \mathbf{12} & 9 & 6 & 3 \end{pmatrix} \begin{array}{c} 0 \\ 2 \\ 3 \\ 3 \end{array}$$

**16)** Jacobi gave the criterion of lemma 10 as a way to help finding a maximal set of transversal maxima, but seemed to assume that, most of the time, one will find them by inspection, as he did for the $10 \times 10$ matrix provided as an example in [40, fig. 1 p. 48]. Our presentation is a reinterpretation that does not fully reflect the spirit of a method intended for hand computation.

In his analysis of Jacobi's contribution [65], Kuhn made a distinction in his algorithm between a *Kőnig step*, 2) d) *i.e.* finding the maximal number of transversal maxima, and an *Egerváry step*, 2) e), *i.e.* increasing the number of transversal maxima. This is coherent with Jacobi's presentation and underlines its deep similarity with Kuhn's Hungarian method.

*Kuhn*
*Kőnig step*
*Egerváry step*
*Hungarian method*

Jacobi completed his work with a few more algorithms, allowing to compute the minimal canon, knowing an arbitrary canon or a maximal set of $s$ transversal maxima (see subsec. 3.3 and subsec. 3.4 below), which he did not use in his study of differential systems. They will be exposed in the next section 3, together with some complements about algorithms and complexity. We will conclude this section with the case of the maximal matching problem, followed by some properties of covers that will be needed in sections 6 and 7.

## 2.3 Maximal matching

*Egerváry* *Kondratieva* *Frobenius* *Kőnig's theorem* *Schrijver*

Egerváry's results were influenced by the following theorem of Kőnig [62, 63, 92] (see also Kondratieva [61, lemma 2]), which was in turn inspired by previous works of Frobenius [27, 28] (one may refer to Schrijver [93] for historical details). It is an easy consequence of Jacobi's criterion for characterizing maximal transversal families of maxima (lemma 10).

*maximal matching*

The problem is equivalent to *maximal matching*: a graph $G \subset [1, s] \times [1, n]$ being given, to find an injective partial function $\sigma : [1, s] \mapsto [1, n]$, with $(i, \sigma(i)) \in G$ and such that $\sharp\mathrm{Im}(\sigma)$ is maximal. Associating to $G$ the $s \times n$ matrix $A$ with $a_{i,j} = 1$ if $(i, j) \in G$ and $a_{i,j} = 0$ if not, we are reduced to the computation of a maximal transversal family of maxima.

*Kőnig's theorem*

**Theorem 17.** — *Let $A$ be a $s \times n$ matrix of zeros and ones, with $s \leq n$, $m$ be the smallest integer such that, for some $p \in \mathbf{N}$, the ones are all located in the union of $p$ rows and $m - p$ columns, then $m$ is the maximal diagonal sum $\mathcal{O}$ in $A$.*

Proof. — It is easily seen that $\mathcal{O} \leq m$: in any diagonal sum, at most $p$ "ones" belong to these $p$ rows, $m - p$ to these $m - p$ columns, and their total number is at most $(m - p) + p = m$.

To prove $\mathcal{O} \geq m$, we can use Jacobi's construction. Without loss of generality, one may discard rows and columns of zeros and reduce to the case where all columns contain a one. Assume that we have $r := \mathcal{O}$ diagonal starred ones, that we may assume to be $a_{1,1}, \ldots, a_{r,r}$. We can then use lemma 10 with the following change : the ones are the maximal elements, the zeros the non maximal elements. According to the lemma, there are no lower right ones. Let $p$ be the number of first class rows, that we may assume to be rows 1 to $p$. Rows $p + 1$ to $r$ do not belong to the first class and so they contain no ones located in columns 1 to $p$ nor $r + 1$ to $n$. Rows $r + 1$ to $s$ belong to the third class and, in the same way, cannot contain ones in columns 1 to $p$, nor $r + 1$ to $n$, as there are no lower right ones. So, all the ones belong to $p$ rows and $r - p$ columns, as illustrated by the figure below.

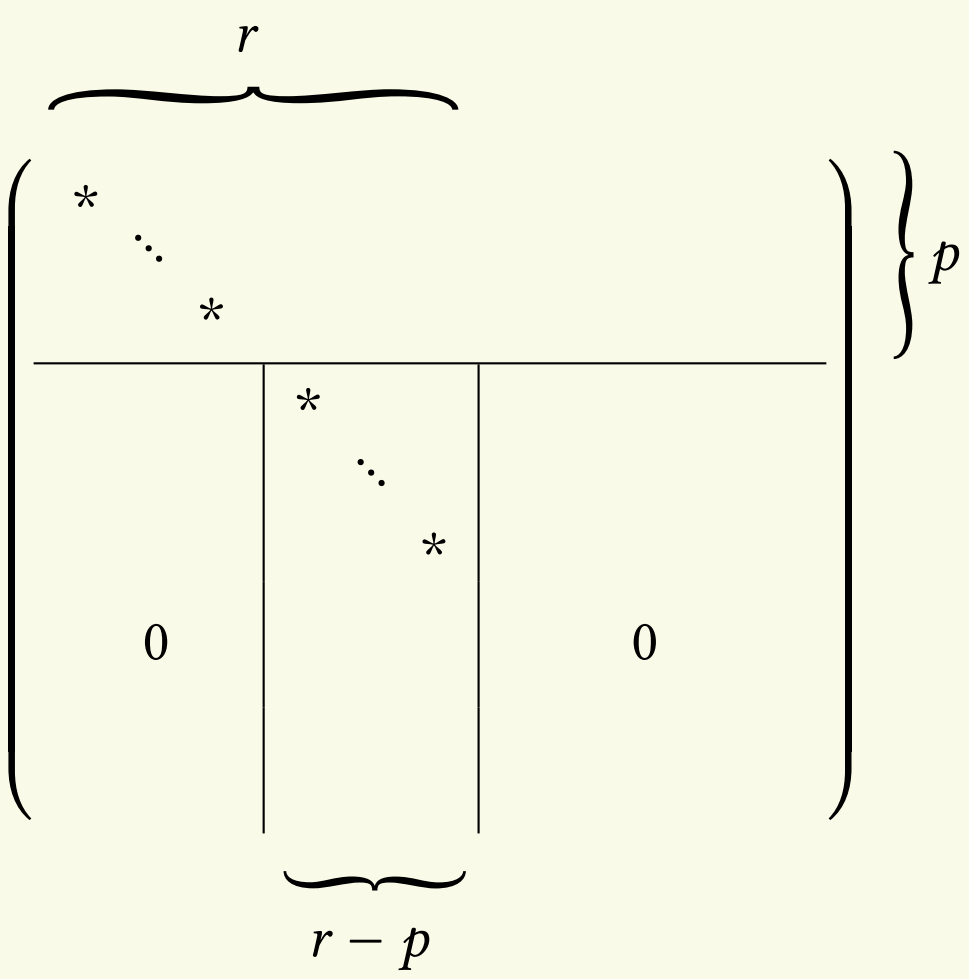

This concludes the proof. ■

*We can already get a naïve algorithm by adapting Jacobi's method to this special case, in order to avoid the computation of a canon.*

**Algorithm 18.** — **Input data:** $A$, an $s \times n$ matrix of zeros and ones.
**Output:** a maximal transversal sum in $A$.

We discard rows and columns of zeros and reduce to the case $s \leq n$ by transposition. Classes of rows will be constructed here, not with respect to maximal elements, but with respect to "ones".

To solve the problem, we only have to construct the set of first class rows, with a cost of $O(sn)$ operation, and to apply step 2) d) of algo. 9, which by lemma 10 can occur at most $r - 1$ times, where $r$ is the size of the matching; *hence a total cost of $O(s^2 n)$ operations for the whole algorithm.* This is the class of complexity of an improved version of Jacobi's algorithm (see below subsec. 3.2).

But it is possible to lower the complexity with a slight modification, due to Hopcroft and Karp [33]. See below 3.1. *Hopcroft and Karp*

## 2.4 Covers

*If not stated otherwise, we consider in this section only square matrices $A$ with elements in some totally ordered commutative group* $\mathrm{M}$.

Covers, at the basis of Ergeváry [23] and Kuhn's [64] approach do not appear in Jacobi's paper but it is interesting to investigate their relations with canons.

**Definition 19.** — *We call a* cover *for $A$ the data of two vectors $(\mu_1, \ldots, \mu_n) \in \mathrm{M}^n$ and $(\nu_1, \ldots, \nu_n) \in \mathrm{M}^n$, such that $a_{i,j} \leq \mu_i + \nu_j$. A cover $\mu, \nu$ is a* minimal cover *if the sum $\sum_{i=1}^n \mu_i + \nu_i$ is minimal.* *cover* *minimal cover*

*Let $\mu, \nu$ and $\mu', \nu'$ be two covers for $A$, then we say that they are equivalent if there exists $\gamma \in \mathrm{M}$ such that $\mu_i' = \mu_i + \gamma$ and $\nu_j' = \nu_j - \gamma$.*

**Proposition 20.** — *i) A cover $\mu, \nu$ of $A$ is minimal iff there exists a permutation $\sigma$ such that $a_{i,\sigma(i)} = \mu_i + \nu_{\sigma(i)}$. Then, $\mathcal{O}_A = \sum_{i=1}^n \mu_i + \nu_i$.*

*ii) Let $\ell$ be a canon for $A$, $\mu_i := (\max_{k=1}^n \ell_k) - \ell_i$ and $\nu_j := \max_{k=1}^n (a_{k,j} - \mu_k)$. The vectors $\mu, \nu$ form a minimal cover for $A$, that we define as* the cover associated to the canon $\ell$. *cover associated to a canon*

*iii) Let $\mu, \nu$ be a minimal cover for $A$, the $\ell_i := (M := \max_{k=1}^n \mu_k) - \mu_i$ form a canon for $A$, that will be called* the canon associated to the cover $\mu, \nu$. *canon associated to a cover*

*iv) a) If $\ell'$ is the canon associated to the cover associated to a canon $\ell$, then $\ell_i' = \ell_i - \min_{k=1}^n \ell_k$.*
*b) If $\mu', \nu'$ is the cover associated to the canon associated to a cover $\mu, \nu$, then $\mu', \nu'$ and $\mu, \nu$ are equivalent, with $\mu' = \mu - \min_{i=1}^n \mu_i$.*

Proof. — i) $\Longrightarrow$ We have, by hypothesis, $\sum_{i=1}^n \mu_i + \nu_i = \sum_{i=1}^n a_{i,\sigma(i)}$ and, by definition of a cover, $\sum_i \mu_i + \nu_i \geq \sum_i a_{i,\sigma(i)}$, hence the minimality of the cover $\mu, \nu$.

$\Longleftarrow$ Let us assume that there is no such permutation $\sigma$. Then, by th. 17 the entries $a_{i,j}$ with $a_{i,j} = \mu_i + \nu_j$ belong to $p$ rows and $m-p$ columns, with $m<n$, that we may suppose to be rows $1,\ldots,p$ and columns $1,\ldots,m-p$. Let

$$e := \min_{i=p+1}^{n} \min_{j=m-p+1}^{n} (\mu_i + \nu_j - a_{i,j}),$$

we have $e>0$ as $e=0$ would mean that $\mu_i + \nu_j = a_{i,j}$ for some $p<i\le n$ and $m-p<j\le n$, which is impossible. We define $\mu_i' := \mu_i$ if $1\le i\le p$ and $\mu_i' := \mu_i - e$ if $p<j\le n$, $\nu_j' := \nu_j + e$ if $1\le j\le m-p$ and $\nu_j' = \nu_j$ if $m-p<j\le n$: $\mu', \nu'$ is a cover, with $\sum_{i=1}^n(\mu_i' + \nu_i') = \sum_{i=1}^n(\mu_i + \nu_i) - (n-m)e$.

ii) By construction, $\mu_i + \nu_j \ge a_{i,j}$, so that $\mu, \nu$ is a cover. Minimality is a consequence of i), remarking that, if $a_{i,\sigma(i)} + \ell_i$ form a maximal transversal sum in $A+\ell$, then $a_{i,\sigma(i)} = \mu_i + \nu_i$.

iii) By i) there exists a permutation $\sigma$ such that $a_{i,\sigma(i)} = \mu_i + \nu_{\sigma(i)}$, so that $a_{i,\sigma(i)} + \ell_i = \mu_i + \nu_{\sigma(i)} + \ell_i = \nu_{\sigma(i)} + M = \mu_{i'} + \nu_{\sigma(i)} + \ell_{i'} \ge a_{i',\sigma(i)} + \ell_{i'}$ for all $1 \le i' \le n$, so that $\ell$ is a canon.

iv) a) We have

$$\ell_i' = \left(\max_{k=1}^{n} \mu_k\right) - \mu_i = \max_{k=1}^{n}\left[\left(\max_{h=1}^{n}\ell_h\right) - \ell_k\right] - \left(\max_{h=1}^{n}\ell_h\right) + \ell_i = \ell_i - \min_{k=1}^{n}\ell_k.$$

b) Similar computations show that $\mu_i' = \mu_i - \min_{i=1}^n \mu_i$ and $\nu_i' = \nu_i + \min_{i=1}^n \mu_i$. ■

***Remark* 21.** — The proof of i) relies on a recursive process to compute a minimal cover by successive reduction of the total sum $\sum_{i=1}^n \mu_i + \nu_i$. If the number $n-p$ of rows $i$ for which $\mu_i$ is decreased is chosen to be minimal, such an algorithm has *Kuhn* a polynomial complexity and is indeed very close to Kuhn's Hungarian method, *Hungarian method* but also to Jacobi's algorithm, as they will correspond to third class rows. Jacobi's algorithm may precisely be used to determine a set of $p$ rows and $m-p$ columns containing all $(i,j)$ with $a_{i,j} = \mu_i + \nu_j$, with $p$ maximal (see prop. 58).

With a different choice, some second class rows are also increased. Then, the complexity may be exponential with integer entries, and the process may never *Jüttner* end with entries in $\mathbf{Q}(\sqrt{5})$. See Jüttner [55].

**Definition 22.** — *The minimal cover associated to the minimal canon of A will be called the* Jacobi cover *or the* canonical cover, *denoted by* $\alpha, \beta$.

***Remarks.* — 23)** Knowing any canon, we can compute the associated cover in $O(n^2)$ operations.

**24)** If $A$ is a matrix of non negative elements of M, *i.e.* $a_{i,j} \ge 0$, then any cover is equivalent to a cover of non negative elements in M. It is easily seen that $\min_i \mu_i + \min_j \nu_j \ge \min_{i,j} a_{i,j}$, so that one just has to define $\mu_{i'}' := \mu_{i'} - \min_i \mu_i$ and $\nu_{j'}' := \nu_{j'} + \min_i \mu_i$ to be sure that $\mu_i' \ge 0$ and $\nu_j' \ge 0$.

**25)** In our definition, the terms of a cover are not assumed to be non negative *Egerváry* as in Egerváry's one [23]. We cannot restrict to this case for all our applications, because we will need to consider in sec. 4.1 entries equal to $-\infty$.

The matrix

$$\begin{pmatrix} 0 & 2 \\ -1 & 1 \end{pmatrix}$$

admits a maximal transversal sum of positive elements $0 + 1$, but no minimal cover of non negative integers: $2 - 1$ is also a maximal transversal sum, so that $\mu_2 + \nu_1 = -1$ and $\mu_2$ or $\nu_1$ must be negative.

**26)** For any integer matrix $A$ of zeros and ones, all integer minimal covers are equivalent to minimal covers that are vectors of zeros and ones, as well as their associated canons.

**27)** If $A$ is a canon, then $A^{\mathrm{t}}$ is not necessarily a canon, but we can easily compute the associated cover $\mu_i$, $\nu_j$ for $A$. Then, $\nu_j$, $\mu_i$ will be a cover for $A^{\mathrm{t}}$ and the associated canon $\ell_i := (\max_k \nu_k) - \nu_i$ will be computed in $O(n^2)$ operations. *transposition*

**28)** If $\alpha, \beta$ is the Jacobi cover of $A$, $\beta$, $\alpha$ is not in general the Jacobi cover of $A^{\mathrm{t}}$, which is not even in all cases equivalent to $\beta$, $\alpha$. For $A$ such that $a_{i,j} = \mu_i + \nu_j$, the canonical cover of $A$ (resp. $A^{\mathrm{t}}$) will be equivalent to $\mu, \nu$ (resp. $\nu$, $\mu$) and they will even be equal if some $\mu_{i_0} = 0$ (resp. some $\nu_{j_0} = 0$). But the canonical cover of the matrix

$$A := \begin{pmatrix} 1 & 1 & 1 & \mathbf{1} & 1 \\ 1 & 1 & \mathbf{1} & 1 & 1 \\ 1 & \mathbf{1} & 0 & 0 & 0 \\ \mathbf{1} & 0 & 0 & 0 & 0 \\ 1 & 0 & 0 & 0 & \mathbf{0} \end{pmatrix}$$

is $(1, 1, 1, 0, 0), (1, 0, 0, 0, 0)$, that of $A^{\mathrm{t}}$ is $(1, 1, 0, 0, 0), (1, 1, 0, 0, 0)$ and the two covers have no simple relation.

**29)** Let $A$ be a matrix. If $\ell$ is a canon for $A$, then for any $\bar{\ell} \in \mathrm{M}^n$, $\ell$ is a canon for $\left(a_{i,j} + \bar{\ell}_j | (i,j) \in [1,n]^2\right)$. This is in particular the case, when $\bar{\ell}$ is a canon for $A^{\mathrm{t}}$. So, denoting by $A^{\mathrm{c}}$ the minimal canon of $A$, $A^{\mathrm{ctct}} = A^{\mathrm{tctc}} = \left(a_{i,j} + \lambda_i \bar{\lambda}_j | (i,j) \in [1,n]^2\right)$, where $\lambda$ (resp. $\bar{\lambda}$) is the minimal canon for $A$ (resp. $A^{\mathrm{t}}$). It is the simplest matrix $B$ that is a canon for $A$ and such that $B^{\mathrm{t}}$ is a canon for $A^{\mathrm{t}}$, meaning that minimal quantities are added to the rows and columns of $A$ to obtain $B$ with such properties.

**30)** Let $A$ be a matrix of zeros and ones, and $\mu, \nu$ a cover of $A$. The non zero elements of $A$ are located in the rows $i$ with $\mu_i \neq 0$ and columns $j$ with $\nu_j \neq 0$. We recover Kőnig's theorem 17 when considering a minimal cover. Reciprocally, if $R$ and $C$ are two sets of rows and columns containing all the ones appearing in $A$, with $\sharp R + \sharp C$ minimal, then $\mu_i = 1$, if $i \in R$, and $\nu_j = 1$, if $j \in C$, defines a minimal cover for $A$. *Kőnig's theorem*

**31)** Let $A$ be a matrix of non negative integers and $a_{i,\sigma(i)}$ be transversal maxima for $A$, then this matrix has at most

$$\prod_{i=1}^{n} (a_{i,\sigma(i)} + 1)$$

covers of non negative integers and this number is reached if all elements of $A$ except these transversal values are 0. Assume that the $a_{i,j}$ belong to some ordered group $G$ where there is no infinite strictly decreasing sequence of positive elements, then any matrix admits a finite number of non negative covers.

On the other hand, if there are an infinite number of values $c \in G$ such that $\min_{j \neq \sigma(i_0)}(\mu_{i_0} + \nu_j - a_{i_0,j}) \geq c \geq 0$, there is an infinite number of non equivalent minimal covers $\mu'_{i_0} := \mu_{i_0} - c$, $\mu_{i \neq i_0} := \mu_i$; $\nu'_{\sigma(i_0)} := \nu_{\sigma(i_0)} + c$, $\nu'_{j \neq \sigma(i_0)} := \nu_j$.

elementary path **32)** Assume that the $a_{i,i}$ form a maximal transversal sum in the matrix $A$. For any cover $\mu, \nu$ of $A$, there is an elementary path (alg. 9 step 2 c) from row $i_0$ to row $i_1$ of the canon of $A$ associated to $\mu, \nu$ (prop. 20 iii) iff $a_{i_1,i_0} = \mu_{i_1} + \nu_{i_0}$.

The next proposition will help to clarify the situation and to compute, in case of need, non Jacobi covers and their canons.

minimal cover **Proposition 33.** — *Let $A$ be a matrix, such that the $a_{i,i}$ form a maximal transversal sum. Let $\mu, \nu$ be a minimal cover for $A$. Using remark 32, we will use the reflexive transitive closure of the path relation $\prec$, defined on rows of the associated canon of $A$ and the transposed relation $\prec^{\mathrm{t}}$ defined on the rows of the associated canon of $A^{\mathrm{t}}$, i.e. the columns of $A$. Rows and columns will be denoted by their indices. We repeat the elementary rules:*

*a) $i_1 \prec i_2$ if $a_{i_2,i_1} = \mu_{i_2} + \nu_{i_1}$;*

*b) $i_1 \prec^{\mathrm{t}} i_2$ if $a_{i_1,i_2} = \mu_{i_1} + \nu_{i_2}$.*

*i) For any integer $i_0$, the rules:*

$$\mu'_i := \mu_i + e \quad \text{and} \quad \nu'_i := \nu_i - e \quad \text{if} \quad i_0 \prec i$$

*and*

$$\mu'_i := \mu_i \quad \text{and} \quad \nu'_i := \nu_i \quad \text{if not,}$$

*where*

$$e \leq \min_{\substack{i_0 \prec i \\ i_0 \not\prec i'}} \left(\mu_{i'} + \nu_i - a_{i',i}\right), \tag{4}$$

*define a minimal cover for $A$.*

*i′) For any integer $i_0$, the rules*

$$\mu''_i := \mu_i - e \quad \text{and} \quad \nu''_i := \nu_i + e \quad \text{if} \quad i_0 \prec^{\mathrm{t}} i$$

*and*

$$\mu''_i := \mu_i \quad \text{and} \quad \nu''_i := \nu_i \quad \text{if not,}$$

*where*

$$e \leq \min_{\substack{i_0 \prec^{\mathrm{t}} i \\ i_0 \not\prec^{\mathrm{t}} i'}} \left(\mu_i + \nu_{i'} - a_{i,i'}\right),$$

*define a minimal cover for $A$.*

*ii) For any minimal cover $\mu', \nu'$ of $A$ and $i_0 \prec i$, we have $\mu'_i \geq \mu_i + \mu'_{i_0} - \mu_{i_0}$.*

*ii′) For any minimal cover $\mu', \nu'$ of $A$ and $i_0 \prec^{\mathrm{t}} i$, we have $\mu'_i \leq \mu_i + \mu'_{i_0} - \mu_{i_0}$.*

Proof. — i) As the transversal sum is unchanged, the minimality is granted. We only have to prove that we obtain a cover. If $i_0 \prec i$ and and $i_0 \prec j$ or $i_0 \not\prec i$ and $i_0 \not\prec j$, then $a_{i,j} \leq \mu'_i + \nu'_i = \mu_i + \nu_j$. If $i_0 \prec i$ and $i_0 \not\prec j$, $a_{i,j} \leq \mu_i + \nu_j \leq \mu'_i + \nu'_i = \mu_i + \nu_i + e$. If $i_0 \not\prec i$ and $i_0 \prec j$, $a_{i,j} \leq \mu'_i + \nu'_j = \mu_i + \nu_j - e$ by (4).

The proof of i′) is similar by transposition. *transposition*

ii) This is equivalent to $\mu'_{i_0} - \mu_{i_0} \leq \mu'_i - \mu_i$. Assume there is a path $i_k$, $0 \leq k \leq r$, with $i_r = i$. Then, $a_{i_k, i_k} = \mu_{i_k} + \nu_{i_k} = \mu'_{i_k} + \nu'_{i_k}$. By rem. 32, $a_{i_{k+1}, i_k} = \mu_{i_{k+1}} + \nu_{i_k} \leq \mu'_{i_{k+1}} + \nu'_{i_k}$, so that $\mu'_{i_{k+1}} - \mu_{i_{k+1}} \geq \nu_{i_k} - \nu'_{i_k} = \mu'_{i_k} - \mu_{i_k}$. We have then $\mu'_{i_0} - \mu_{i_0} \leq \mu'_{i_1} - \mu_{i_1} \leq \cdots \leq \mu'_i - \mu_i$.

The proof of ii′) is again similar by transposition. ∎

The following examples are easy illustrations of the last proposition.

***Examples.*** — **34)** A matrix $A$ with $a_{i,j} = \mu_i + \nu_j$, admits a single class of minimal covers: that of $\mu, \nu$.

**35)** We define the triangular matrix of zeros and ones:

$$A := \begin{pmatrix} 1 & 1 & 0 & \cdots & 0 \\ 0 & \ddots & \ddots & \ddots & \vdots \\ \vdots & \ddots & \ddots & \ddots & 0 \\ \vdots & & \ddots & \ddots & 1 \\ 0 & \cdots & \cdots & 0 & 1 \end{pmatrix} \quad \text{and} \quad B := \begin{pmatrix} 1 & 1 & 0 & \cdots & 0 \\ 1 & \ddots & \ddots & \ddots & \vdots \\ 0 & \ddots & \ddots & \ddots & 0 \\ \vdots & \ddots & \ddots & \ddots & 1 \\ 0 & \cdots & 0 & 1 & 1 \end{pmatrix}.$$

The matrix $A$ has exactly $n+1$ minimal covers of non negative integers indexed by $0 \leq k \leq n$, defined by

$$\begin{array}{llllll} \alpha_i = 1 & \text{if} & i \leq k & \text{and} & \alpha_i = 0 & \text{if} \quad i > k \\ \beta_j = 0 & \text{if} & j \leq k & \text{and} & \beta_j = 1 & \text{if} \quad j > k. \end{array}$$

There are only $n$ classes of integer covers up to equivalence: the covers obtained for $k = 0$ and $k = n$ are equivalent and correspond to the minimal canon.

The matrix $B$ has only 2 minimal covers of non negative integers: $\mu = (1, \ldots, 1)$, $\nu = (0, \ldots, 0)$ and $\mu = (0, \ldots, 0)$, $\nu = (1, \ldots, 1)$. All its minimal covers are equivalent.

***Remark* 36.** — We could have defined covers for rectangular matrices, but if $s < n$, *rectangular matrix* we need to impose some lower bound for the $\nu_j$ in order to ensure the existence of a minimal cover. Indeed, the covers $\mu_i + c$ and $\nu_j - c$ are equivalent, but

$$\sum_{i=1}^{s} (\mu_i + c) + \sum_{j=1}^{n} (\nu_j - c) = \left( \sum_{i=1}^{s} \mu_i + \sum_{j=1}^{n} \nu_j \right) - (n-s)c,$$

so that no minimal cover can be defined for $\nu_i \in \mathbf{Z}$ and $\mu_j \in \mathbf{Z}$. If we impose $\nu_j \geq 0$, for any matrix $A$ of non negative elements, a minimal cover $\mu, \nu$ is such that

$$\sum_{i=1}^{s} \mu_i + \sum_{j=1}^{n} \nu_j = \mathcal{O}.$$

so that Jacobi's bound (def. 1) could also be handled using covers in a direct way.

*minimal cover* Adding $n - s$ rows of 0, as in rem. 3, all $\mu_i$ for $s < i \leq n$ must be equal and we may choose their common value to be 0. Then, a cover $\mu, \nu$ is minimal iff there exists an injection $\sigma \in S_{s,n}$ such that $\mu_i + \nu_{\sigma(i)} = a_{i,\sigma(i)}$ and, for all $j \notin \mathrm{Im}(\sigma)$, $\nu_j = 0$. (See also rem. 57.)

# 3 Related algorithms and deeper complexity analysis

*complexity (of Jacobi's algorithm)* DISCUSSING the complexity of Jacobi's algorithm is an interesting subject, but we need to keep in mind that it is anachronistic to do it in the setting of modern computation models, when Jacobi's concern was to spare the work of useless rows rewriting, in a time when pen and paper remained the main computational tools [39, § 4, p. 32]. We will now provide some improvements that lead to a better complexity, in our contemporary formalism.

## 3.1 Finding a maximal set of transversal maxima. The bipartite matching problem and Hopcroft and Karp's algorithm in Jacobi's setting

*bipartite matching* We have encountered with Jacobi's algorithm the following special problem of finding a maximal set of transversal maxima. This amounts to solving the assignment problem with a matrix of zeros and ones, using Jacobi's characterization (see *transversal one* lem. 10). In what follows, all maximal values being 1, we will speak of *transversal* *starred one* *ones*, *starred ones* instead of transversal or starred maxima. This is known as the maximal bipartite matching, or marriage problem.

We could change the data structure and use the graph of the relation $a_{i,j} = 1$, smaller that the full matrix, as initial data, but for the sake of clarity, we will stay here in the dense setting.

*Hopcroft and Karp* We refer to Hopcroft and Karp [33] for a different and more detailed presentation. First, in Jacobi's setting, we may look repeatedly for lower right ones. This may be done in sequence, until no such "one" is found with a total cost $O(sn)$. *Kőnig's theorem* One may notice that Kőnig's theorem (th. 17) implies that this first step already produces at least $\lceil \mathcal{O}/2 \rceil$ transversal ones.

The elementary path relation can be constructed with cost $O(s^2)$; the cost of the first class construction is $O(sn)$. The main idea is to build a maximal set (in the sense that it is not strictly included in another such set) of disjoint paths of minimal length going to a lower first class row, before building a new path relation. So the main step of the algorithm is not to produce a single augmenting path, but, at each stage $k$, a maximal set of disjoint paths of the same length $\gamma$. In the sequel, we assume, as we may, that all rows and columns of zeros have been discarded and, up to transposition, that $s \leq n$.

*length (of a path)* **ALGORITHM 37.** — LENGTH. — **Input data:** a matrix $A$ and a transversal set of "ones", given by an injection $\sigma : [1, r] \mapsto [1, s] \times [1, n]$, that we will assume here to be $a_{i,i}$, $1 \le i \le r$.
**Outputs:** the minimal length $\gamma$ of a path from the set $L_0$ of first class rows with a right "one" to a lower row, or *"failed"* if no such path exists and the list of sets of rows $L_i$, $0 \le i \le \gamma$, to which there is a path of length $i$, and no shorter path, from a first class row with a right "one".

We start with the upper rows with right "ones", that form the set $L_0$. Let $M := L_0$.

At step 1, we define $L_1$ to be the set of elements not in $M$, such that there is a path of length 1 to them from some element of $L_0$. We increase $M$ with $L_1$. We then define $L_2$ to be the set of elements not in $M$ to which there is a path of length 1 from an element of $L_1$, etc.

We stop this process as soon as $L_\gamma$ is empty and we return then *"failed"*, or contains a lower line, that will be by construction a first class lower line. The integer $\gamma$ will be the minimal length of a path leading to a lower first class row, that we return.

*This process is achieved in $O(s^2)$ operations.*

To find a maximal set of disjoint paths, we may use the following recursive process. The maximal number of disjoint paths is bounded by the cardinal of $L_0$, as these elements are the possible origins of any of them. We define first a set $F$ of available rows, which is initialized with the set of rows not in $L_0$. In the following function, we assume that the length $\gamma$ has been computed using alg. 37.

**Algorithm 38.** — *Path.* **Input data:** an integer $\gamma_0$ and a path $[i_0, \ldots, i_j]$, where $i_0$ is a first class row with a right "one" and $i_\ell \in L_\ell \cap F$. *path*
**Output:** a path $[i_0, \ldots, i_j, \ldots, i_\gamma]$, with $i_\gamma$ a lower first class row and $i_\ell \in L_\ell \cap F$, for $j < \ell \le \gamma$, or *"failed"* if no such path exists.
**Global variables:** $F$, a set of elements to be used for building paths during the process, and $F_j := L_j \cap F$.

***Step 1)*** Let $C$ be the set of elements of $F_{j+1}$ such that there is a path from $i_j$ to them.

If $j = \gamma - 1$ and $C \ne \emptyset$, let $c \in C$, remove $c$ from $F$ and $F_{j+1}$ and return $[i_0, \ldots, i_j, c]$.

If $j < \gamma - 1$, go to step 2).

***Step 2)*** **For** $c \in C$ **do**: remove $c$ from $F$ and $F_{j+1}$; if $\mathrm{Path}(\gamma, [i_0, \ldots, i_j, c]) \ne$ "failed", then return $\mathrm{Path}(\gamma, [i_0, \ldots, i_j, c])$.
If $C$ is exhausted before a path could be found, then return "failed".

This process returns a path of length $\gamma$ to a lower first class row, if such a path exists. The elements in that path are removed from $F$, as well as the elements of $L_i$ from which no such path of length $\gamma - i$ has been found. So repeated call to that function will produce disjoint paths and the function Path can be applied only once to a given row.

*This implies that we can apply in sequence* Path *to the elements of $L_0$ to get a maximal set of disjoint paths of minimal length $\gamma$ in $O(s^2)$ operations.*

**Algorithm 39.** — *Increase.* — **Input data:** a list $T$ of transversal maxima, and an "increasing" path $[j_0, \ldots, j_\gamma]$ from a first class row with a right maximum to a lower first class row.
**Output:** an increased list of transversal maxima.

We proceed as in alg. 9 step 2) d) p. 16. *I.e.*, denoting by $1_{i,j}$ a one placed in row $i$ and column $j$, we have an an upper right "one" $1_{j_0,h}$ and "ones" $1_{j_{i+1},j_i}$, $0 \le i < \gamma$, according to the path relation, so that we may replace the stared "ones" $1_{j_i,j_i}$ by $1_{j_{i+1},j_i}$ and complete them with $1_{j_0,h}$.

*Hopcroft and Karp's algorithm* **Algorithm 40.** — *Hopcroft–Karp.* — **Input data:** a matrix $A$ of zeros and ones.
**Output:** the elements of a maximal transversal sum of $A$ and a minimal cover.

***Step 1)*** As stated above, we repeatedly look for lower right ones, producing first at least $\lceil \mathcal{O}/2 \rceil$ transversal ones $T$.

***Step 2)*** Let $(\gamma, L) := \mathit{Length}(A, T)$. If $\gamma =$*"failed"*, then return $T$.
If not, let $F := [1, s] \setminus L_0$ and **for** $i_0 \in L_0$ **do**: If $J := \mathit{Path}(\gamma, [i_0]) \neq$*"failed"* then $T :=\mathit{Increase}(T, J)$.
**Repeat** *Step 2)*.

*The total cost of steps 1) or 2) is $O(sn)$, so the key point in bounding the complexity is to evaluate how many times step 2) is performed, which is the goal of the next two lemmata 41 and 45.*

To keep the spirit of Jacobi's algorithm, we have distinguished in the algorithm the case of paths of length 0 (step 1) and the general case. In the sequel, step 1 will be step 0 and the $k^{th}$ iteration of step 2 will be considered as step $k$, according to the common length $\gamma_k$ of the increasing paths used at each step.

***Lemma* 41.** — *Let $\gamma_k$ be the length of the paths used at stage $k$, then the sequence $\gamma_k$ is strictly increasing.*

Proof. — Assume it is not the case and $\gamma_k \le \gamma_{k-1}$. We may assume that $k$ is minimal with that property. We call a changed row, a row that has been used in some path at stage $k-1$. In any path used at stage $k$, there must be a changed row. If not, either $\gamma_k = \gamma_{k-1}$ and this contradicts the fact that algorithm 38 produces a maximal set of disjoint paths of length $\gamma_{k-1}$, or $\gamma_k < \gamma_{k-1}$ and this contradicts the minimality of the length of paths produced by algorithm 37.

An injective function $\phi : [1, s] \mapsto [1, s]$ defines a unique set of disjoint paths and loops, the union of which is equal to the union of its image and its definition domain. If $\exists r\ \phi^r(i) = i$, then $i$ belongs to a loop, if not let $r_0 := \max\{r | \phi^{-r}(i)$ is defined$\}$ and $r_1 := \max\{r | \phi^r(i)$is defined$\}$, then $i$ belongs to the path $\phi^{-r_0}(i)$, …, $\phi^{r_1}(i)$. Reciprocally, any disjoint set of paths and loops defines a unique such function.

Let $\phi$ be the function defined by the $\tau$ paths of stage $k-1$ and $i_0$, …, $i_{\gamma_k}$ be a path of stage $k$: it must have some rows in common with the paths of stage $k-1$. Let them be $i_{h_1}$, …, $i_{h_r}$, $r \ge 1$. If $\phi^{-1}(i_{h_\kappa})$ is defined, we replace in the graph of $\phi$ the couples $(\phi^{-1}(i_{h_\kappa}), i_{h_\kappa})$ with $(\phi^{-1}(i_{h_\kappa}), i_{h_\kappa+1})$. We then add to the graph of $\phi$ the couples $(i_\zeta, i_{\zeta+1})$, $\zeta \notin \mathrm{Im}(h)$.

This construction is illustrated by the figure on the right. Elementary path relations at stage $k-1$ are indicated by ⋮ and at stage $k$ by |; after rearrangement by [. The starred ones of stage $k-1$ by *1* and those of stage $k$ by **1**. This defines an injection to which is associated a new set of $\tau+1$ paths and (possibly) loops.

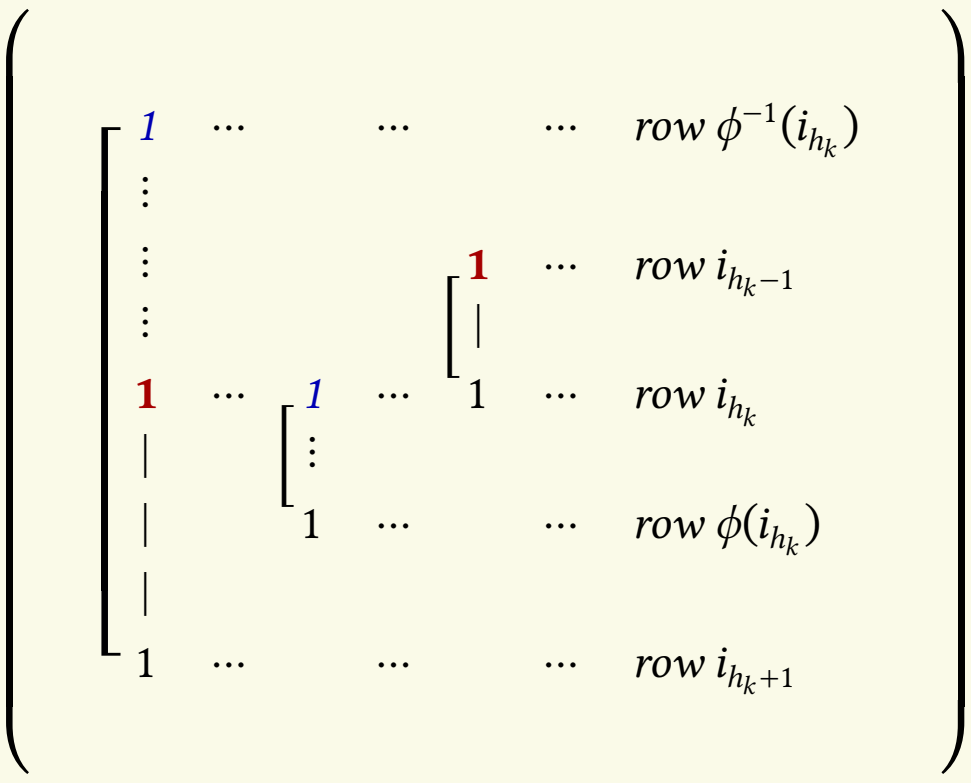

Then, the sum of their lengths is at most $\tau\gamma_{k-1} + \gamma_k - r$ (and strictly smaller iff loops do exists). So, as $\gamma_k \le \gamma_{k-1}$, one path must be of length strictly smaller that $\gamma_{k-1}$. This contradicts the minimality of $\gamma_{k-1}$. ■

***Remark* 42.** — Our paths of length $r$ correspond to paths of length $2r+1$ following the conventions of Hopcroft and Karp [33]. This is due to the fact that they define the path relations, not between rows but between the "ones" involved in the path relation. In their setting, starred elements appear with a minus sign and the others with a plus sign. So the process of reconstruction reduces to computing the sum of the two paths. E.G., denoting by $1_{i,j}$ a one placed in row $i$ and column $j$, we have: *path* *Hopcroft and Karp*

$$\begin{aligned} & (+1_{1,1} - 1^*_{1,3} + 1_{2,3} - 1^*_{2,4} + 1_{4,4}) + (+1_{2,2} - 1^*_{2,3} + 1_{3,3}) \\ = \ & (+1_{1,1} - 1^*_{1,3} + 1_{3,3}) + (+1_{2,2} - 1^*_{2,4} + 1_{4,4}). \end{aligned}$$

The element $1_{2,3}$ that appeared two times has vanished and the two paths of lengths 5 and 3 are replaced by two paths of length 3.

***Examples.* — 43)** In the following example, the stared ones of the first stage ($k-1$ in lem. 41) are in *blue italic*: *1* and those of the second stage ($k$ in lem. 41) in **red bold**: **1**. The first path includes rows 1, 2 and 4, the second rows 4 and 3. In Hopcroft and Karp's convention: $(1_{1,3} - 1_{1,2} + 1_{2,2} - 1_{2,1} + 1_{4,1}) + (1_{4,4} - 1_{4,1} + 1_{3,1}) = (1_{1,3} - 1_{1,2} + 1_{2,2} - 1_{2,1} + 1_{3,1}) + (1_{4,4})$. In Jacobi's setting: one path of length 2, formed of rows 1, 2 and 3, and one of length 0, *viz.* a lower right "one": $1_{4,4}$.

$$\begin{pmatrix} 0 & \mathit{1} & \mathbf{1} & 0 \\ \mathit{1} & \mathbf{1} & 0 & 0 \\ \mathbf{1} & 0 & 0 & 1 \\ 1 & 0 & 0 & 0 \end{pmatrix}$$

**44)** Using the same conventions, the first path includes rows 1, 2 and 3, the second rows 3, 1 and 4. In H. & K.'s convention: $(1_{1,3} - 1_{1,2} + 1_{2,2} - 1_{2,1} + 1_{3,1}) + (1_{3,4} - 1_{3,1} + 1_{1,1} - 1_{1,3} + 1_{4,3}) = (1_{1,1} - 1_{1,2} + 1_{2,2} - 1_{2,1}) + (1_{3,4}) + (1_{4,3})$. In Jacobi's: one loop, formed of rows 1 and 2, and two lower right "ones": $(1_{3,4})$ and $(1_{4,3})$.

$$\begin{pmatrix} 1 & \mathit{1} & \mathbf{1} & 0 \\ \mathit{1} & \mathbf{1} & 0 & 0 \\ \mathbf{1} & 0 & 0 & 1 \\ 0 & 0 & 1 & 0 \end{pmatrix}$$

**Lemma 45.** — *Let $(a_{i,j})_{(i,j)\in G_1}$ and $(a_{i,j})_{(i,j)\in G_2}$, where $G_1$, $G_2$ are the graphs of two injective functions $[1,s]\mapsto[1,n]$, be two families of $r_1 := \sharp G_1$ and $r_2 := \sharp G_2$ transversal ones of $A$. We assume $r_2 > r_1$. Lower right "ones" in the family $(a_{i,j})_{(i,j)\in G_2}$, i.e. elements that are not placed in the same rows or columns as the elements of $G_1$, will be considered to be paths of length $0$.*

*Then, using only the starred "ones" in $(a_{i,j})_{(i,j)\in G_1}$ and the "ones" in $(a_{i,j})_{(i,j)\in G_2}$ placed in the same columns, we define a path relation such that there exists a path of length at most $\lfloor r_1/(r_2-r_1)\rfloor$.*

Proof. — If some lower right "one" exists in $(a_{i,j})_{(i,j)\in G_2}$, then the result stands according to our convention. If not, we obtain possibly loops (if $G_1$ and $G_2$ have a element $(i,j)$ in common, then we consider it as a loop from row $i$ to itself) and $r_2-r_1$ open paths, as there are as many elements from $G_1$ and $G_2$ in loops and one more element in $G_2$ than in $G_1$ in some open path. As a path of length $m$ involves $m$ starred ones in $G_1$, the sum of the lengths of all paths is at most $r_1$ and there exists a path of length at most $\lfloor r_1/(r_2-r_1)\rfloor$. ∎

**Theorem 46.** — *Assume that the maximal number of transversal "ones" in $A$ is $s_0$, then the algorithm requires at most $2\lfloor\sqrt{s_0}\rfloor$ steps. Its complexity is $O(\sqrt{s_0}sn)$.*

Proof. — Let $\ell$ be the number of steps (at step $k$, paths of the same length $\gamma_k$ are considered). Let $G_1$ be the set of starred ones at step $\ell-\lfloor\sqrt{s_0}\rfloor$ and $G_2$ the maximal set of starred ones at step $\ell$. We have $\sharp G_1 \le s_0-\lfloor\sqrt{s_0}\rfloor$ and $\sharp G_2-\sharp G_1\ge\lfloor\sqrt{s_0}\rfloor$. Then, using lemma 45, the length of a path at step $\ell-\lfloor\sqrt{s_0}\rfloor$ is a most $(s_0-\lfloor\sqrt{s_0}\rfloor)/\lfloor\sqrt{s_0}\rfloor \le \lfloor\sqrt{s_0}\rfloor$, so that there are at most $\lfloor\sqrt{s_0}\rfloor$ steps before step $\ell-\lfloor\sqrt{s_0}\rfloor$ and $\ell\le 2\lfloor\sqrt{s_0}\rfloor$. ∎

*Frobenius* This problem was first investigated by Frobenius [27] in order to decide *a priori* if a matrix where non zero elements can appear at known places has an identically vanishing determinant. The best known asymptotic complexity for its resolution is bigger than that of a numerical determinant. We can only achieve the exponent of matrix multiplication with probabilistic algorithms using random numerical values. See Ibarra and Moran [37]. *Ibarra and Moran*

## 3.2 A $O(s^2n)$ version of Jacobi's algorithm

*complexity of Jacobi's algorithm* In order to improve the complexity of Jacobi's algorithm, we only have to remark that it is useless to reconstruct the whole path relation in order to reduce the number of second class rows or make some lower first class row appear, as the starred maxima will remain unchanged. We replace step 2 *e*) p. 17 with e′) and e″), defined as follows.

e′) We have already defined the first, second and third classes $C_{\rm I}$, $C_{\rm II}$ and $C_{\rm III}$; we moreover define the set $C_{\rm I,II}$ of first and second class rows. For each $i\in C_{\rm I,II}$, *starred maximum* we compute the minimal distance $d_i$ between its starred maximum $a_{i,i}$ and some third class row element in the same column and the minimal distance $d_0$ between some upper right maxima in a first class row and some third class row element in the same column: $d_0 := \min_{i\in C_{\mathrm{III}}} \min_{j=r+1}^{n} m_j - a_{i,j}$, where $m_j := \max_{i\in C_{\mathrm{I}}} a_{i,j}$, $r < j \leq n$ is the maximum in column $j$.

*All this is done with a cost at most $O(sn)$.*

e″) We will then increase every third class row by $d := \min(d_0, \min_{i\in C_{\mathrm{I,II}}} d_i)$. *third class*
If $d = d_0$, this creates a lower right maximum and we go to b) p. 15. If $d = d_{i_0}$ with $i_0 \in C_{\mathrm{I}}$, this creates a lower first class row, step e″) is finished and we go to substep d) p. 16. If $i_0 \in C_{\mathrm{II}}$, we remove $i_0$ from the sets $C_{\mathrm{II}}$ and $C_{\mathrm{I,II}}$, and add it to the third class. We redefine $d_i$, $i \in C_{\mathrm{I,II}}$ to be the minimum of $d_i - d$ and the distance between its starred maximum $a_{i,i}$ and $a_{i_0,i}$. We also redefine $d_0$ to be the minimum of $d_0 - d$ and $\min_{j=r+1}^{n} m_j - a_{i_0,j}$. We iterate e″) with these new values.

*The substep* e″) *is performed with a total cost of $O(n)$ operations and will be iterated at most $s-2$ times; so its total cost is $O(sn)$ operations. Hence the following theorem.*

**Theorem 47.** — *Using substep 2.* e′) *and* e″), *the complexity of Jacobi's algorithm is bounded by $O(s^2n)$.* ■

The improved complexity $O(n^3)$ was first obtained for square matrices by Dinic and Kronrod [20] in 1969, rediscovered independently by Tomizawa [95] in 1971 and then by Edmonds and Karp [22] in 1972. *Dinic and Kronrod; Tomizawa* *Edmonds and Karp*

***Remarks.* — 48)** We could obviously improve the "Kőnig step" b) c) by using alg. 40, but this does not change the exponents in the asymptotic complexity $O(s^2n)$ of th. 47. *Kőnig step*

**49)** As we have already seen (see rem. 15), we cannot avoid, in some cases, to repeat at least $s-2$ times step 2) of this improved version of Jacobi's algorithm. It could be possible to speed up the construction of elementary path relations, as in most cases they are unchanged or reversed. But a $O(s^2)$ complexity, when building the class partition, seems unavoidable. In this situation, we don't know how to construct in a single step a large set of augmenting paths, as for the maximal matching problem (see 2.3). *augmenting path* *maximal matching*

One may notice that pioneering aspects of Jacobi's work include reachability issues and computing the transitive closure of a directed graph. But this problem is not formalized and its solution is implicitly assumed to be achieved in some naïve way for small size data. However, some of his algorithms solve problems equivalent to instances of the shortest path problem. *directed graph*

## 3.3 A canon being given, to find the minimal one

In order to solve this problem, Jacobi [39, VII p. 23] proposes to compute first a maximal set of transversal maxima, which may be done using Hopcroft and Karp's algorithm 40 with complexity $O(n^{5/2})$ for a square matrix $A$. Knowing transversal maxima, we may then use the following method. *Hopcroft and Karp*

**ALGORITHM 50.** — ***Input:*** a square matrix $A$ of size $n$, a canon $\ell$ for $A$ and a maximal system of transversal maxima for $\ell$, that we assume for simplicity to be $a_{i,i} + \ell_i$.

***Output:*** the minimal canon of $A$.

**Step 1.** We decrease all the $\ell_i$ by $\min_{i=1}^n \ell_i$, so that at least one $\ell_i$ is 0.

**Step 2.** We build the path relation. Then, we establish the list $L_1$ of rows with $\ell_i = 0$, or from which there is a path to a row with $\ell_i = 0$.

By lemma 51 below, if $L_2 := [1, n] \setminus L_1 = \emptyset$, we have finished and return $\ell$. If not, we compute

$$b := \min(\min_{i \in L_2} \ell_i, \min_{i \in L_2} \min_{i' \in L_1} a_{i,i} + \ell_i - a_{i',i} - \ell_{i'}). \tag{5}$$

*Step 2) may be achieved with complexity $O(n^2)$.*

**Step 3)** We decrease all the $\ell_i$, $i \in L_2$ by $b$. In this way, some rows will go from $L_1$ to $L_2$. We then **repeat** step 2).

*The complexity for the whole process is $O(n^3)$.*

The next lemma shows that this algorithm actually returns the minimal canon.

*minimal canon* ***Lemma 51.*** — *i) The canon $\ell$ is the minimal canon $\lambda$ iff there is a path from all rows to a row $i$ with $\ell_i = 0$.*

*ii) If $\mathcal{O}$ is equal to $\sum_{i=1}^n a_{i,\sigma(i)}$, with $a_{i,\sigma(i)} \geq 0$ for $1 \leq i \leq n$, then we have:* $\max_{i=1}^n \lambda_i \leq (n-1) \max_{1 \leq i,j \leq n} a_{i,j}$.

PROOF. — i) $\Rightarrow$ — Assume that $\ell_{i_0} > \lambda_{i_0}$ where $\lambda$ is the minimal canon. We may choose such $i_0$ so that there is a path of minimal length from row $i_0$ to a row $i$ with *elementary path* $\ell_i = 0$. It means that there is an elementary path from row $i_0$ to some row $i_1$ with *unchanged row* $\ell_{i_1} = \lambda_{i_1}$ unchanged. Then, $a_{i_0,i_0} + \ell_{i_0} = a_{i_1,i_0} + \lambda_{i_1}$, so that $\lambda_{i_0} = \ell_{i_0}$; a contradiction.

$\Leftarrow$ — By the definition of $b$ in (5), we have a new canon $\tilde{\ell} < \ell$ defined by $\tilde{\ell}_i = \ell_i$ for $i \in L_1$ and $\tilde{\ell}_i = \ell_i - b$ for $i \in L_2$.

ii) There is a path $i_h$, $0 \leq h \leq p$ from any row $i_0$ to a row $i_p$ with $\lambda_{i_p} = 0$. So, $a_{i_{j+1},\sigma(i_j)} + \lambda_{i_{j+1}} = a_{i_j,\sigma(i_j)} + \lambda_{i_j}$, which implies $\lambda_{i_j} - \lambda_{i_{j+1}} \leq a_{i_{j+1},\sigma(i_j)}$, so that $\lambda_{i_0} \leq \sum_{j=0}^{p-1} a_{i_{j+1},\sigma(i_j)} \leq (n-1) \max_{1 \leq i,j \leq n} a_{i,j}$, as $p \leq n-1$. ∎

It is possible to turn the complexity to $O(n^2 \ln(n))$, using "balanced trees" or *Adelson-Velsky and Landis* "AVL trees", from the name of their inventors Adelson-Velsky and Landis [1]. See *Knuth* also Knuth [57, 6.2.3 p. 451]. This tree structure allows to maintain dynamically an ordered list of $p$ elements, allowing to insert, delete, search the order of an element or an element of a given order in $O(\ln p)$ operations.

**ALGORITHM 52.** — We use the same input, data and conventions as in alg. 50 above.

**Step 1) a)** Decrease all the $\ell_i$ by $\min_{i=1}^n \ell_i$. Create a list $L_1$ of the rows $i$ with $\ell_i = 0$ and a list $L_2$ of the remaining elements.

**b)** For $i \in L_2$, create a balanced tree $T_i$, containing for all the rows $i'$ of $L_1$ the pairs $(a_{i,i} - a_{i',i} - \ell_{i'}, i')$, sorted by lexicographical order.

*This may be achieved with total cost $O(n^2 \ln n)$.*

**Step 2)** Compute $d_i := \min_{i' \in L_1} a_{i,i} - a_{i',i} - \ell_{i'}$ and $b := \min_{i \in L_2} \min(d_i + \ell_i, \ell_i)$. Let $L_3 := \{i \in L_2 | \min(d_i + \ell_i, \ell_i) = b\}$. For all $i \in L_2$, decrease $\ell_i$ by $b$. Let $L_2 := L_2 \setminus L_3$ and $L_1 := L_1 \cup L_3$. For all $i \in L_2$ and all $i' \in L_3$, add $(a_{i,i} - a_{i',i} - \ell_{i'}, i')$ to the tree $T_i$.

*All this may be done with total cost $O(n \ln n)$.*

If $L_2$ is empty, we have finished and return $\ell$, if not we iterate step 2), *which will be performed at most $n$ times*, hence the next theorem.

**Theorem 53.** — *Knowing a canon for a square matrix of size $n \times n$ and a set of transversal maximal elements in this canon, algo. 52 computes the minimal canon with $O(n^2 \ln n)$ elementary operations.* *minimal canon*

Proof. — The complexity has already been proved. At step 1), $L_1$ contains rows with $\ell_i = 0$. An easy recurrence shows that $L_1$ contains rows $i$ with $\ell_i = 0$, or from which there is a path to such a row. Indeed, if $b = \ell_i$, then $\ell_i$ is set to 0 and, if $b = a_{i,i} + \ell_i - a_{i',i} - \ell_{i'}$, then $\ell_i$ is set to $a_{i',i} + \ell_{i'} - a_{i,i}$. So, $a_{i,i} + \ell_i = a_{i',i} + \ell_{i'}$ and there is a path to row $i' \in L_1$ from which there is a path to a row $i''$ with $\ell_{i''} = 0$. At the end of the algorithm, $L_1 = [1, n]$ and, by lem. 51, $\ell$ is the minimal canon. ■

**Proposition 54.** — *Let $A_{i,j}$ be a square matrix of size $n$ with at least a transversal set of $n$ maxima, that we may assume to be $a_{i,i}$. Then, the reflexive transitive closure of the path relation does not depend on the choice of this transversal set.* *path relation*

Proof. — Let $a_{i,i}$ and $a_{\sigma(i),i}$ be two transversal sets of maxima. We denote by $\prec_1$ (resp. $\prec_2$) the path preorder relation defined using the first family (resp. the second). Assume that there are elementary paths $i \prec_1 j$ and $\sigma(i) \prec_2 j$. Consider the cycle $i_0 = i$ and $i_{p+1} = \sigma(i_p)$. Let $i_r = i$, with $r$ minimal. According to the path definition, there is an elementary path $i_p \prec_1 i_{p+1}$ and $i_{p+1} \prec_2 i_p$. Using the cycle, we have

$$i \prec_2 i_{r-1} \prec_2 \cdots \prec_2 i_1 = \sigma(i) \prec_2 j \text{ and } \sigma(i) = i_1 \prec_1 i_2 \prec_1 \cdots \prec_1 i_r = i \prec_1 j,$$

so that the reflexive transitive closure of $\prec_1$ and $\prec_2$ are the same. ■

**Definition 55.** — *We will denote by $\pi_A$ the path relation associated with the minimal canon $A + \lambda$ of $A$.* *path relation*

***Remarks.*** — **56)** a) If there is an elementary path from row $i_0$ to row $i_1$, then to row $i_2$ ...up to row $i_r$ with $\lambda_{i_r} = 0$, then

$$\lambda_{i_0} = \sum_{k=1}^{r} a_{i_k, i_{k-1}} - a_{i_{k-1}, i_{k-1}}. \tag{6}$$

b) Generically, *i.e.* if $a_{i,i} + \lambda_i = a_{j,i} + \lambda_j$ for a single index $j \neq i$, $\pi_A$ defines a forest of rooted trees with $n$ labeled vertices, where the roots correspond to rows with $\lambda_i = 0$. By a variant of Cayley's formula there are then $(n+1)^{n-1}$ possibilities, and as much sets of formulae 6 for the values $\lambda_i$ in the minimal canon. *Cayley* *rectangular matrix* **57)** Assume that Jacobi's bound for some $s \times n$ matrix $A$, with $s < n$, is $\mathcal{O} = \sum_{i=1}^{s} a_{i,i}$. Then, completing $A$ with $n - s$ rows of 0, the minimal canon $\lambda$, is such that $\lambda_{i_0} = \max_{j=s+1}^{n} \max_{i=1}^{s} a_{i,j} + \lambda_i$, for $s < i_0 \leq n$. (See also rem. 36.)

Assume now that we want to compute the minimal canon for a $n \times s$ matrix $A$ completed with $n - s$ columns of 0. Let again the transversal maximal elements be $a_{i,i}$, $1 \leq i \leq n$. Then, it is easily seen that for $i' > s$, the minimal canon $\lambda$ is such that $\lambda_{i'} = \max_{i=1}^{s} \lambda_i$.

Before leaving this subject, we will emphasize the special case of matrices of zeros and ones, associated to maximum matching problems. In this case and with a canon of zeros and ones, the algorithm 50 will run in $O(n^2)$ steps, as we only need to apply step 2) one time. It provides an efficient way to maximize the *Kőnig's theorem* number of rows (or columns) in Kőnig's theorem (th. 17).

**Proposition 58.** — *Let $A$ be some $s \times n$ matrix of zeros and ones (possibly horizontally or vertically rectangular), $A'$ the $\max(s,n) \times \max(s,n)$ matrix obtained by adding $|n - s|$ columns or rows of 0 to $A$, $\alpha, \beta$ the associated canonical cover and $\mathcal{O}$ the maximal transversal sum of $A$ and $A'$.*

*i) Let $\lambda$ be the minimal canon of $A'$, then for any sets $R$ of rows and $C$ of columns, with $\sharp R + \sharp C$ minimal and containing all the 1 in $A$, $i_0 \in R$ implies $\lambda_{i_0} = 0$.*

*ii) In Kőnig's theorem, there exists a unique such couple of sets of rows $R$ and columns $C$ with $R$ maximal for inclusion (resp. with $C$ maximal for inclusion).*

Proof. — Let $\mu, \nu$ be the minimal cover defined by $\mu_i = 1$ if $i \in R$ and $\nu_j = 1$ if $j \in C$. If $R = \emptyset$, the assertion is true. If not, the canon associated to $\mu, \nu$ is such that $\ell_i = 0$ iff $\mu_i = 1$. The assertion i) is then an easy consequence of the minimality of $\lambda$ (prop. 7).

ii) The result is straightforward if the Jacobi number $\mathcal{O}$ of $A$ is $s$. Then $R = \{1, \ldots, s\}$ is the maximal set of rows.

If some $\lambda_i = 1$, then the result is a direct consequence of i), as the maximal set $R$ is defined by $\lambda_i = 0$ or equivalently $\alpha_i = 1$, where $\alpha, \beta$ is the minimal cover associated to $\lambda$: $\alpha_i = \max_{k=1}^{s} \lambda_k (= 1) - \lambda_i$.

If all the $\lambda_i$ are 0 and $\mathcal{O} < s$, one just has to consider the $(\max(s,n) + 1) \times (\max(s,n) + 1)$ matrix $A''$ obtained by adding to $A'$ a column of $\max(s,n)$ zeros and then a row of $\max(s,n) + 1$ ones. We get a maximal transversal sum of value $\mathcal{O} + 1$ for $A''$ by completing one for $A'$ with the 1 in column and row of index $\max(s,n) + 1$. As there must be some zero in any maximal transversal sum, some $\lambda_i$ must be 1 in the minimal canon of $A''$, so that we can now apply i).

*transposition* The statement for columns is obtained by considering the transpose matrix $A^{t}$. ■

**Definition 59.** — *We call these covers the* row maximal *and the* column maximal *minimal covers.*

**Algorithm 60.** — **Input data:** A $s \times n$ matrix $A$ of zeros, together with an injection $\sigma$ such that $\sum_{i=1}^{s} a_{i,\sigma(i)}$ is maximal.

**Output:** A row maximal minimal cover.

***Step 1.*** Compute the Jacobi number $\mathscr{O}$ of $A$. If $\mathscr{O} = s$, then the $s$ rows of $A$ form the row maximal cover.

*This requires* $O(\sqrt{\mathscr{O}}sn)$ *operations.*

***Step 2.*** Make a square matrix $A'$ by adding to $A$ $|n - s|$ rows or columns and add a column of $\max(s, n)$ zeros and a row of $\max(s, n) + 1$ ones to define a $(\max(s, n)+1)\times(\max(s, n)+1)$ matrix $A''$ as in the proof of prop. 58. Then compute the minimal canon of $A''$ using algorithm 50; the rows $1 \leq i \leq s$ with $\lambda_i = 0$ form the row maximal minimal cover of $A$.

*As* $\max(s, n) - \min(s, n)$ *rows share the same value for* $\ell_i$ *(rem. 57), this may be done in* $O(\min(s, n))^2$ *operations.*

*The asymptotic complexity of the algorithm is that of step. 1:* $O(\sqrt{\mathscr{O}}sn)$ *elementary operations.*

## 3.4 Transversal maxima being given, to find the minimal canon

If we don't have a canon but just know the places of transversal maxima in a canon $A + \ell$, the next algorithm computes the minimal canon (Jacobi [39, p. 25]). *transversal maxima*

**ALGORITHM 61.** — **Input data:** A $n \times n$ matrix $A$ and a transversal family, that we may assume to be $a_{i,i}$.
**Output:** The minimal canon of $A$.
**Repeat** the following *elementary step* **until** all rows remain unchanged.

ELEMENTARY STEP. For $1 \leq i \leq n$, increase row $i$ by $\max_{k=1}^{n} a_{k,i} - a_{i,i}$.

*The elementary step may be done in* $O(n^2)$ *operations.*

**PROPOSITION 62.** — *This algorithm produces the minimal canon in* $O(n^3)$ *operations.* *minimal canon*

PROOF. — The elementary step in the algorithm will be repeated at most $n$ times, the exact number being in the generic case $p+1$, where $p$ is the maximal distance from any row to a row with $\lambda_i = 0$, according to the path relation of def. 55. ■

***Remarks.*** — **63)** This algorithm may be easily modified to compute the path relation. Given any transversal family $a_{i,\sigma(i)}$, it may be used to test if it corresponds to a maximal sum, the stopping of the algorithm after $n$ steps being a necessary and sufficient condition.

If the algorithm does not stop after $n$ steps, it means that there is a loop $\sigma : I \subset [1, n] \mapsto I$, so that $\sum_{i\in I} a_{i,i} < \sum_{i\in I} a_{i,\sigma(i)}$.

**64)** The last example of [39, § 3] is the transpose of a canon. Then, this transpose is not a canon, but the terms of a maximal transversal sum are known and we can apply the above method. We may also compute a cover and deduce from it a canon (see rem. 27), allowing to use algo. 52 for better efficiency. In Jacobi's informal setting, the two methods have comparable complexities.

**65)** We will encounter in subsec. 7.1 prop. 115 the following problem: *"Assuming that a canon of $A$ exists such that the elements $a_{i_\kappa,j_\kappa}$, $1 \le \kappa \le r$, are maximal, to compute such a canon."* If we assume that the $a_{i,i}$ are transversal maxima, one may use a variant of alg. 61 where the elementary step is completed with: "For $1 \le \kappa \le r$, increase row $i_\kappa$ by $\max_{i=1}^{n} a_{i,j_\kappa} - a_{i_\kappa,j_\kappa}$."

*This process stops after at most n steps iff such a canon exists. So its complexity is also $O(n^3)$.*

The next proposition completes the preceding remark.

**Proposition 66.** — *Assuming that there exists a canon such that the $a_{i_\kappa,j_\kappa}$, $1 \le \kappa \le r$, are maximal, there exists a unique minimal such canon, that is given by the variant of alg. 61 in the last remark.*

Proof. — Let $\ell$ and $\ell'$ be two such canons. We know by prop. 7 that $\min(\ell_i, \ell'_i)$ is a canon. We also have $a_{j_\kappa,j_\kappa} + \ell_{j_\kappa} = a_{i_\kappa,j_\kappa} + \ell_{i_\kappa}$ and $a_{j_\kappa,j_\kappa} + \ell'_{j_\kappa} = a_{i_\kappa,j_\kappa} + \ell'_{i_\kappa}$ so that if $\min(\ell_{j_\kappa}, \ell'_{j_\kappa}) = \ell_{j_\kappa}$, then $\min(\ell_{i_\kappa}, \ell'_{i_\kappa}) = \ell_{i_\kappa}$, and reciprocally with $\ell'$. So, $a_{j_\kappa,j_\kappa} + \min(\ell_{j_\kappa}, \ell'_{j_\kappa}) = a_{i_\kappa,j_\kappa} + \min(\ell_{i_\kappa}, \ell'_{i_\kappa})$. The canon $\min(\ell_i, \ell'_i)$ satisfies the requested property and so a unique minimal such canon exists.

To conclude, we only have to remark that at each step of the algorithm, each row is increased by the minimal value in order to satisfy the requested inequalities, hence the minimality of the result. ■

## 3.5 Tropical geometry

*Maclagan and Sturmfels* We refer to Maclagan and Sturmfels [71] for more details on this topic. The basic idea is to replace products by sums and additions by max. It is then obvious that Jacobi's bound is a tropical determinant. *tropical geometry, determinant*

We will denote by $M \odot N$ the tropical matrix multiplication. One may wonder why the analogy with the determinant cannot be used in a straightforward way. One should remark that this analogy suffers important limitations: the analog of addition is "max", that has no inverse, and the tropical determinant of a tropical product of matrices is not in general the sum of their tropical determinants. Such a property holds only in special situations, *e.g.* $|A \odot B|_T = |A|_T \odot |B|_T$ if $B$ is a canon and $A$ the transpose of a canon. Moreover, the tropical determinant is also the tropical permanent... *canon* *transposition* *tropical permanent*

Assume that $a_{i,i}$ is a transversal family with a maximal sum. Then, reducing row $i$ by $a_{i,i}$, we get a new matrix $B$ with $b_{i,i} = 0$, the result of the last algorithm 61 is the tropical matrix product: $(a_{1,1}, \ldots, a_{n,n}) \odot B^n$. A $O(n^\alpha)$ algorithm for the tropical multiplication would produce a $O(n^\alpha \ln(n))$ algorithm for the problem of finding a minimal canon, knowing the elements of a maximal transversal sum.

***Remark* 67.** — Let $I \subset [1,n]$ be a subset of indices, we denote by $\bar{I}$ the complementary set $[1,n] \setminus I$. We denote by $A_{I,J}$ the submatrix of $A$ restricted to the rows of $I$ and the columns of $J$. For any subset $I$, we have the following formula, that mimics a classical property of the determinant:

$$|A|_T = \max_{J \subset [1,n],\ \sharp J = \sharp I} |A_{I,J}|_T + |A_{\bar{I},\bar{J}}|_T.$$

## 3.6 Minimal canons subject to inequalities

**Proposition 68.** — *Let $A$ be a $n \times n$ matrix, and $c_i$, $1 \le i \le n$ positive integers. Then there exists a unique minimal canon subject to the conditions $\ell_i \ge c_i$.*

Proof. — Let $A' := A + c$. Then, $\ell$ is a canon of $A$, subject to $\ell_i \ge c_i$ iff $\ell - c$ is a canon of $A'$, so the unique minimal canon $\lambda$ of $A'$ is such that $\lambda + c$ is the unique minimal canon of $A$, subject to $\ell_i \ge c_i$. ∎

This proof provides an algorithm to compute such minimal canons, that are used in sections 9 and 10 to bound the order of derivations of initial equations, *resolvent* necessary to perform changes of orderings or to compute resolvents. *change of ordering*

## 3.7 Minimal canons and shortest paths

Let $A + \ell$ be a canon for $A$; assume that the $a_{i,i}$ form a maximal transversal sum. Then, we define a weighted directed graph $G$ on the set $\{0, 1, \ldots, n\}$, by associating the weight $w_{i,j} := a_{i,i} + \ell_i - a_{j,i} - \ell_j \ge 0$ to the ordered pair $(i, j)$, and $w_{i,0} := \ell_i$ to the ordered pair $(i, 0)$.

Reciprocally, we may associate to any such directed graph with non negative weights a square matrix $A$ and a canon $A + \ell$, defined by $a_{i,i} = C$, for $C \ge 2\max_{(i,j)\in[0,n]^2} w_{i,j}$, $\ell_i = w_{i,0}$ and $a_{j,i} := C - w_{i,j} + \ell_i - \ell_j$.

**Proposition 69.** — *The vector of integers $\lambda$ is the minimal canon of $A$ iff there exists in $G$ a shortest path of length $\ell_i - \lambda_i$ from vertex $i$ to vertex 0.* *shortest path*

Proof. — To prove that $\lambda$ is a canon, assume that $a_{i,i} + \lambda_i < a_{j,i} + \lambda_j$. This is equivalent to $\ell_i - \lambda_i > \ell_j - \lambda_j + w_{i,j}$, meaning that there is a shortest path to vertex 0 passing by $j$, which contradicts the minimality of $\ell_i - \lambda_i$. It is then enough to remark that there exists such a shortest path from vertex $i$ to vertex 0 in $G$ iff there exists a path in $A + \lambda$, in the meaning of lemma 51, from row $i$ to a row $i_0$ with $\lambda_{i_0} = 0$. So, according to this lemma, iff $A + \lambda_i$ is the minimal canon of $A$. ∎

In the same way, let $A$ be a $n \times n$ square matrix. Define an oriented weighted graph on the set of vertices $\{0, 1, \ldots, n\}$ by setting $w_{i,0} := 0$ on edge $(i, 0)$ and $w_{i,j} := a_{i,i} - a_{j,i}$ on edge $(i, j)$. Reciprocally, define for any such weighted graph a matrix $A$ with $a_{i,i} := C := \max(0, \max_{i,j} w_{i,j})$ and $a_{j,i} := C - w_{i,j}$.

**Proposition 70.** — *i) The entries $a_{i,i}$ of $A$ form a maximal transversal sum iff $G$ admits no strictly negative cycle.*

*ii) Assuming the $a_{i,i}$ to form a maximal transversal sum, a vector $\lambda$ is the minimal canon of $A$ iff there is a shortest path of length $-\lambda_i$ from vertex $i$ to vertex 0 in $G$.* *shortest path*

**Proof.** — i) There is a cycle in $G$, with negative value $\gamma$ iff there exists a permutation $\sigma : [1, n] \mapsto [1, n]$, with $\sum_{i=1}^{n} a_{i,\sigma(i)} = -\gamma + \sum_{i=1}^{n} a_{i,i}$, so that the $a_{i,i}$ do no form a maximal sum.

ii) $\Rightarrow$. Minimality is granted by rem. 56 a). To see that $\lambda$ is a canon, it is enough to remark that if $a_{i,i} + \lambda_i < a_{j,i} + \lambda_j$, then there is a path from $i$ to 0 of length $-\lambda_j + w_{i,j} = -\lambda_j + a_{i,i} - a_{j,i} < -\lambda_i$. A contradiction.

$\Leftarrow$. By rem. 56 a), there is a path from $i$ to 0 of length $-\lambda_i$. Let be a path $\iota_1 = i, \ldots, \iota_r$. Then, $\sum_{h=1}^{r-1} a_{\iota_h,\iota_h} + \lambda_{\iota_h} - a_{\iota_{h+1},\iota_h} - \lambda_{\iota_{h+1}} = 0$, so that $\sum_{i=1}^{r-1} a_{\iota_h,\iota_h} - a_{\iota_{h+1},\iota_h} = \sum_{i=1}^{r-1} w_{\iota_h,\iota_{h+1}} = -\sum_{i=1}^{r-1} \lambda_h - \lambda_{h+1} = -\lambda_i + \lambda_{\iota_r} \geq -\lambda_i$, which proves that $-\lambda_i$ is minimal.
■

This means that the problems considered in sections 3.3 and 3.4 are equivalent to computing a shortest path, respectively for a directed graph with positive weights and a directed graph with arbitrary weights but no negative cycle. Although this contribution is not explicit, it shows that Jacobi deserves some mention as a pioneer of graph theory and shortest paths problems. See Schrijver's very interesting article [94] for more historical details and references. The complexity of Jacobi's original algorithms is $O(n^3)$, similar to that of the algorithms of Ford [25] and Bellman [4]. In the positive case, we are able to turn it to $O(n^2 \ln n)$, which is the complexity of Dijstra's algorith [19], using binary heaps. In the sparse setting, with $m$ edges, it may be turned to $O(\min(k(n^{1+1/k} + m), (n + m) \ln n))$, where $k \geq 2$ is any fixed integer that is the height of the heap (D. Johnson [49] ), or even $O(n \ln n + m)$ ( Fredman and Tarjan [26]), using Fibonacci heaps. Basically, the two methods given by Jacobi are respectively very close to the algorithms of Bellman–Ford and Dijkstra.

*shortest path* · *directed graph* · *graph theory* · *Schrijver* · *Ford and Bellman* · *Dijkstra* · *Johnson, Donald* · *Fredman and Tarjan*

## 3.8 Physical analogies

It is always a greatest help for mathematical intuition to rest on physical models; one may easily design devices computing minimal covers.

*minimal cover*

### 3.8.1 Mechanical computation of a minimal cover

E.g., one may consider a mechanical system consisting of $2n$ horizontal rods, $n$ standing for the rows and $n$ standing for the columns, crossing at right angles.

*mechanical system*

At the crossing of two rods $i$ and $j$, a cable passing to a pulley is attached to both of them, so that if the relative height $\mu_i$ and $\nu_j$ of the two rods, as well as the height $a_{i,j}$ of the pulley is defined to be 0 at rest, when the $a_{i,j}$ are increased to take new positive values, one has:

$$\frac{\mu_i + \nu_j}{2} \geq a_{i,j}, \text{ that becomes } \mu_i + \nu_j \geq a_{i,j},$$

by chosing a half scale for the pulley height. Under gravity, the total energy of the system, which for rods of equal masses is proportional to

$$\sum_{i=1}^{n} \mu_i + \nu_i,$$

will be minimal, so that this device will produce a minimal cover. Assuming that the weight of a rod is $M$, adding a little extra weight to those standing for the rows, say $0 < \epsilon < M/n$, the equilibrium point will be unique and will correspond to minimal values for the $\mu_i$, which corresponds to the minimal canon, provided that we impose $\mu_i \geq 0$, using some wedge. *minimal cover*

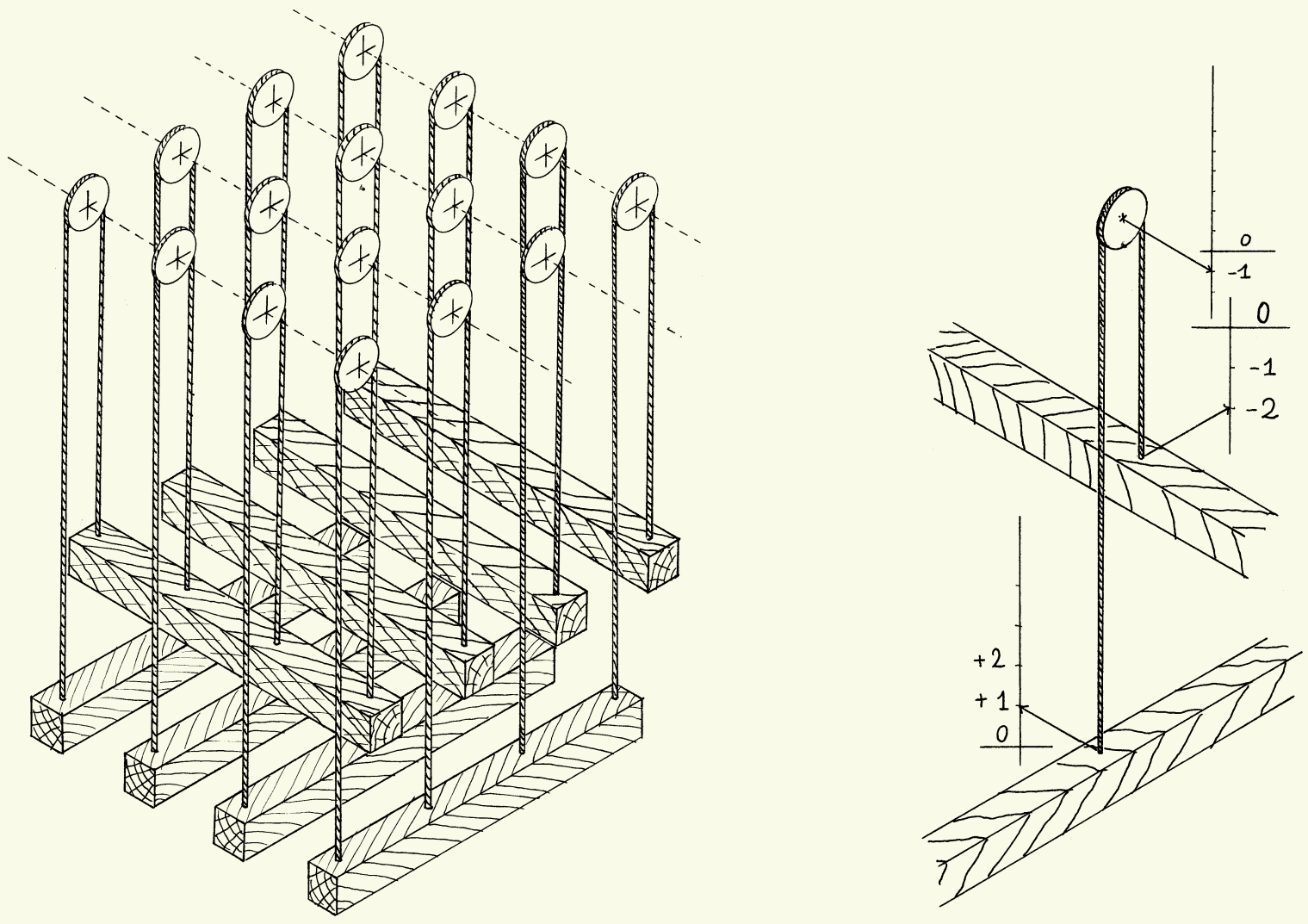

Using such a device, $-\infty$[6] entries can be modeled by suppressing the cable and pulley at some crossing. One can also allow rods to move down, so that negative values for the $\nu_j$ can be achieved too. If $\mathcal{O}$ is $-\infty$, then some rods will fall down... until they are stopped by the finite length of the cables or reach the ground.

### 3.8.2 Materialization of the path relation

*path relation* A second mechanical device may help visualize the graph of the path relation $\pi_A$ *mechanical system* (see def. 55) and rem. 56. Some vertical patterns reproduce the profile of each row of the matrix, *e.g.* below on the left row 3 of some $7 \times 7$ matrix. At the top of the part of each pattern $i$ corresponding to $a_{i,i}$, an orthogonal rod is fixed. The patterns are assumed to be able to move vertically, so that if some $a_{j,i}$ is greater than $a_{i,i}$, the rod of pattern $i$ will rest on pattern $j$. The lowest patterns rest on the floor, corresponding to $\ell_i = 0$. The drawing below on the right corresponds to the

[6] See subsec. 4.1.

minimal canon of

$$A=\begin{pmatrix}3&4&2\\1&3&4\\1&1&3\end{pmatrix},\quad\text{which is:}\quad A+\begin{pmatrix}0\\1\\2\end{pmatrix}=\begin{pmatrix}3&4&2\\2&4&5\\3&3&5\end{pmatrix}.$$

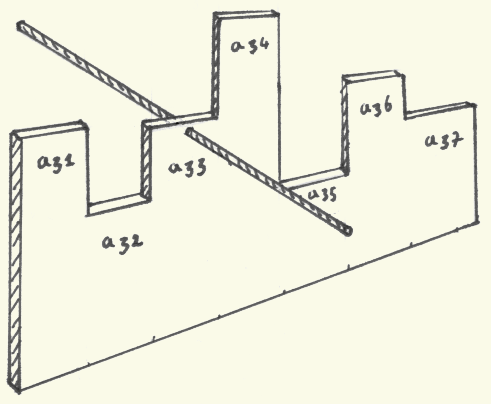

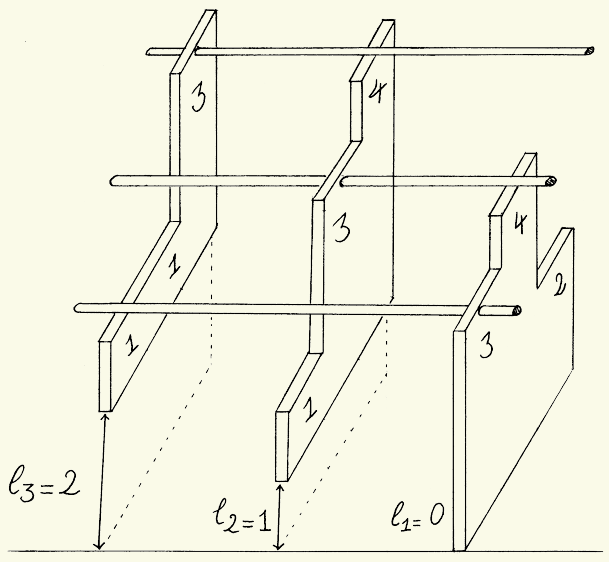

If we use some wedges to impose minimal values $c_i$ for the height of some profiles, one will obtain the minimal canon subject to the condition $\ell_i \geq c_i$.

### 3.8.3 Electrical computation of a minimal cover

*electric circuit* We finish with an electric circuit that may be used to compute a minimal cover. The voltmeters placed in the rows and columns of the circuit will measure quantities corresponding to the covers $\mu_i$ and $\nu_j$. Some adjustable voltage generators are connected at each crossing, providing a tension corresponding to $a_{i,j}$. The presence of a diode realizes the inequality: $\mu_i + \nu_j \geq a_{i,j}$. So, $\mu, \nu$ is a cover. We cannot with this device model entries $a_{i,j} < 0$; the absence of connection or a generator with a negative voltage are equivalent. If all the internal resistances of the voltmeters are equal to some value $R$, we need have $\sum_i \mu_i = \sum_j \nu_j$, as the intensity in and out of the circuit must be equal.

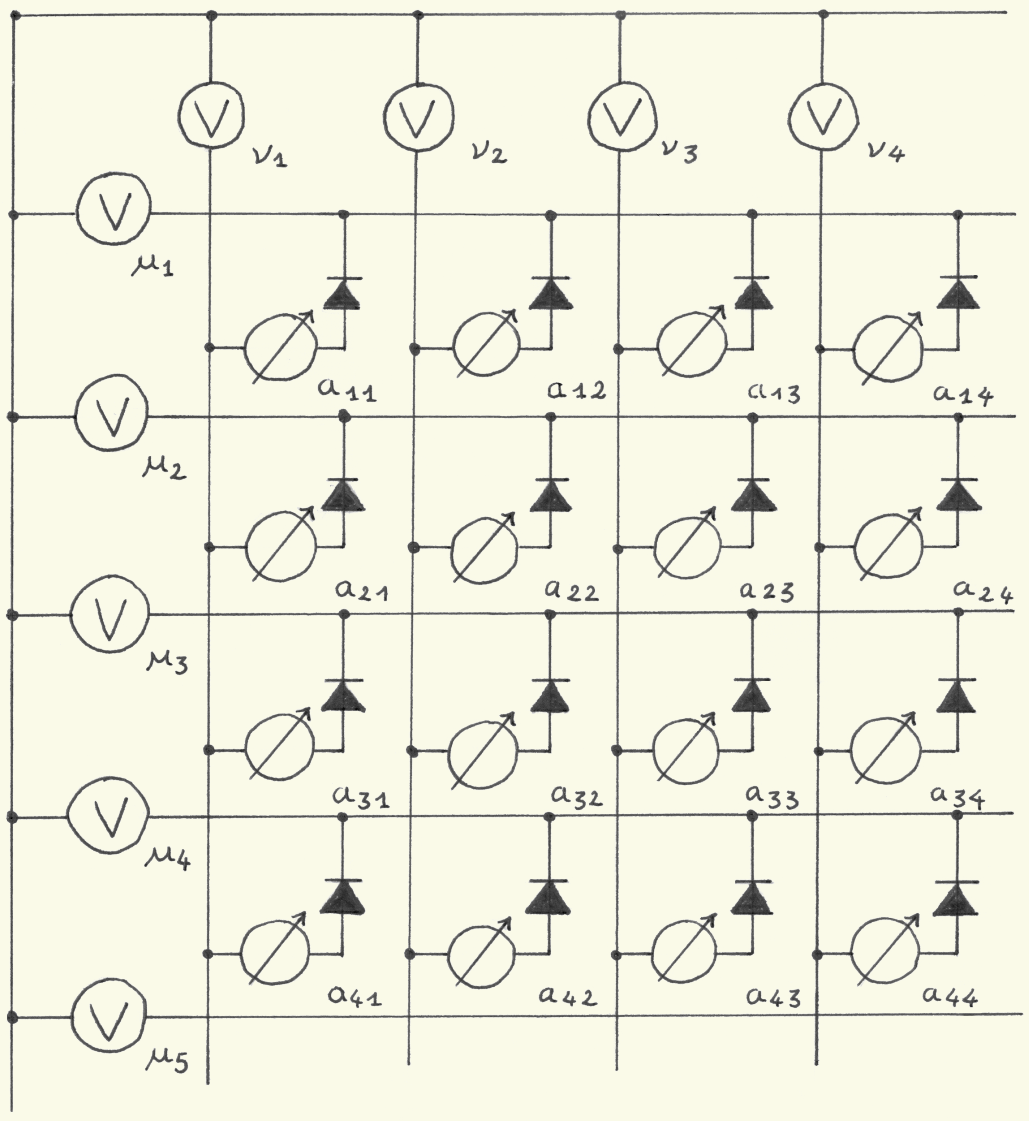

If we assume that the $a_{i,j}$ are 0, except for $a_{1,j} := 1$, $1 < j \leq n$ and $a_{i,1} := 1$, $1 < i \leq n$, one finds the solution $\mu_1 = \nu_1 = (n-1)/n$ and $\mu_i = \nu_j = 1/n$ for $1 < i, j \leq n$, which is not a minimal cover.

We need some extra assumption: let $|a| := \max_{i,j} a_{i,j}$ (assuming $|a| > 0$) it is then enough to replace the $a_{i,j}$ with $b_{i,j} := a_{i,j} + n|a|$ to have a working device.

Let indeed the $(i, j)$ with $b_{i,j} = \mu_i + \nu_i$ belong to the union of $r$ rows and $s$ columns, with $r + s$ minimal, in the spirit of Kőnig's theorem 17. For simplicity, let these $r$ rows be rows 1 to $r$ and these $s$ columns be columns 1 to $s$ and let $b_{i,j} = \mu_i + \nu_i$ for $i + j = r + s$. Let $I_1$ be the electrical intensity from the first $r$ rows to the last $n - s$ columns, $I_2$ from the last $n - r$ rows to the first $s$ columns and $J$ from the first $r$ rows to the first $s$ columns; by hypothesis, the intensity from the last $n - r$ rows to the last $n - s$ columns is 0. For $1 \leq i \leq r$, the minimal value of $\mu_i$ is at most $R(I_1 + J)/r$ and the minimal value of $\nu_j$ for $1 \leq j \leq s$ is $R(I_2 + J)/s$. In the same way, the minimal values $\mu_{i_0}$ of $\mu_i$ for $r < i \leq n$ and the minimal value $\nu_{j_0}$ of $\nu_j$ for $s < j \leq n$ are resp. at most $RI_2/(n-s)$ and $RI_1/(n-r)$.

As $\mu_{i_0} + \nu_j = b_{i_0,j_1}$ for $j = r + s - i_0 n$, we need have $R((I_1 + J)/r + I_1/(n-s)) \leq n|a|$, so that

$$I_1 \leq \frac{r(n-s)}{n-s+r}\frac{n|a|}{R}.$$

In the same way, we have

$$I_2 \leq \frac{s(n-r)}{n-r+s}\frac{n|a|}{R}.$$

This implies that $\mu_{i_0}$ is at most $sn|a|/(n-r+s)$ and $\nu_{i_0}$ is at most $rn|a|/(n-s+r)$.

We will show that $r + s = n$, so that $\mu, \nu$ is a minimal cover. If not, $r + s \leq n - 1$, $\mu_{i_0} \leq sn|a|/(2s+1) < n|a|/2$ and $\nu_{i_0} \leq rn|a|/(2r+1) < n|a|/2$, so that

$$\mu_{i_0} + \nu_{i_0} < n|a| \leq b_{i_0,j_0},$$

a contradiction. So, $r + s = n$ and $\mu$, $\nu$ is a minimal cover.

The conception of a better analog device for computing large tropical determinants may have some practical interest.

## 3.9 Conclusion of section 3

The best deterministic complexity bounds for the assignment problem with integer matrices rely on "scaling" methods, that is recursively replacing in $A$ $a_{i,i}$ by $\lfloor a_{i,i}/2 \rfloor$ to obtain an approximate maximum, as in Gabow and Tarjan [29], where a $O\left(n^{5/2}\ln(nC)\right)$ complexity is achieved (with $0 \leq a_{i,j} \leq C$). See Schrijver [93] or Burkard *et al.* [7] for more details. The basic idea is to use Hopcroft and Karp's algorithm, which is faster, to improve the approximation at each of the $\ln_2 C$ steps.

*Gabow and Tarjan* *Schrijver* *Hopcroft and Karp* *Burkard*

So, in the integer case, the best deterministic exponent for the tropical determinant is 2.5, which is greater than the best known exponent for matrix multiplication: $\omega \leq 2.3728639$ (*cf.* Le Gall [69]). The exponent $\omega$ may be reached using probabilistic algorithms. See *e.g.* Sankowski [89] and the references therein for

*Le Gall* *Sankowski*

*Duan and Petie* more details. Duan and Petie [21] have given a near linear time algorithm for solving the $(1-\epsilon)$-approximate weight matching problem in a graph with $m$ edges and $n$ vertices, with a $O(m\epsilon^{-2}\log^3 n)$ time.

For the best of my knowledge, the $O(n^3)$ complexity of the improved Jacobi's *Kuhn* or Kuhn's algorithms (see sec. 3.2) remains the best when working on an arbitrary *Fredman and Tarjan* monoïd M. In the sparse case, for a graph with $m$ edges, Fredman and Tarjan [26] obtain a $O(nm+n^2\ln n)$ complexity, using Fibonacci heaps.

We have seen that Jacobi's work contained the germs of important notions in combinatorial optimization and graph theory. The efficiency considerations in Jacobi's papers reflect his computational tools: pen and paper, but his algorithms and theoretical framework for the assignment problem may be easily adapted to express improved complexity bounds obtained in the early seventies.

# 4 A differential parenthesis. Various forms of the bound.

## 4.1 Ritt's strong and weak bound

Jacobi did not mention what should be done if some variable $x_j$ and its derivatives do not appear in some polynomial $P_i$. The easiest answer is to define $a_{i,j} := \mathrm{ord}_{x_j} P_i$ as 0, but a better choice in such a case is the convention introduced by Ritt [85] $\mathrm{ord}_{x_j} P_i = -\infty$. Lando [68] defined the first choice as the *weak bound*, and the second as the *strong bound*. Our definition also includes some minor modifications in order to extend the bound to under- or overdetermined systems. *Lando* *weak and strong bounds*

**Definition 71.** — *Let* $P = \{P_1, \ldots, P_s\} \subset \mathcal{F}\{x_1, \ldots, x_n\}$. *By convention,* $\mathrm{ord}_{x_j} P_i = -\infty$ *if* $x_j$ *and its derivatives do not appear in* $P_i$. *Let* $a_{i,j} := \mathrm{ord}_{x_j} P_i$, *we define the* order matrix *of* $P$, $A_P := (a_{i,j})$ *and recall that* $S_{s,n}$ *is the set of injections* $[1, s] \mapsto [1, n]$ *(def. 1). If* $s \le n$, *we define* Jacobi's number *as* $\mathcal{O}_P := \mathcal{O}_A = \max_{\sigma \in S_{s,n}} \sum_{i=1}^{n} a_{i,\sigma(i)}$. *If* $s > n$, $\mathcal{O}_P := \mathcal{O}_{A^{\mathrm{t}}}$. *order* *Jacobi's number*

***Remarks.*** — **72)** An easy consequence of remark 56 a) and of the proof of prop. 62 is the following. We assume that $A$ is a matrix of non negative integers and $-\infty$ elements, with $\max_{i,j} a_{i,j} = C$ and $\mathcal{O}_A \in \mathbf{N}$. Then, if for the minimal canon $\lambda$, the sequence $\lambda_i$ is non decreasing, we have $0 \le \lambda_i \le (i-1)C$, and the associated minimal cover $\alpha, \beta$ is such that $0 \le \alpha_i \le (n-i)C$, with $-(n-1)C \le \beta_j \le C$.

*Cf.* rem. 24 for the case of non negative integers, where $0 \le \alpha_i, \beta_j \le C$.

**73)** Of course, we described an algorithm for matrices with coefficients in an ordered abelian group, and $\mathbf{Z} \cup \{-\infty\}$ with the convention $\mathbf{Z} - \infty = \{-\infty\}$ has not such a structure. So, one may use the group $\{a\infty + b, (a,b) \in \mathbf{Z}^2\}$ with $(a_1\infty + b_1) + (a_2\infty + b_2) = (a_1 + a_2)\infty + (b_1 + b_2)$ and $(a_1\infty + b_1) < (a_2\infty + b_2) \iff a_1 < a_2$ or $a_1 = a_2$ and $b_1 < b_2$. One may also use first algo. 40 to check if $\mathcal{O}_A > -\infty$, *i.e.*, if there exists $\sigma \in S_{s,n}$ such that $a_{i,\sigma(i)} \in \mathbf{N}$, $1 \le i \le s$, by replacing elements in $\mathbf{N}$ by 1 and $-\infty$ elements by 0.

**74)** If $A$ is a $n \times n$ matrix of non negative integers and $-\infty$ values, we may use also any value $D < -(n-1)C$ instead of $-\infty$ to define a new matrix $\tilde{A}$ such that $\mathcal{O}(\tilde{A}) < 0 \iff \mathcal{O}_A = -\infty$. It may indeed be difficult to model in a suitable way $-\infty$ entries of a matrix, even in some computer algebra systems.

Moreover, the tropical determinant of $A$ is $\mathcal{O}_A = -a\infty + b$ iff the tropical determinant of $A'$ is $\mathcal{O}_{A'} = aD + b$ with $aD \le \mathcal{O}_{A'} < (a-1)D$. *tropical determinant*

**75)** If $s \ne n$ we can also complete the matrix of orders $A_P$ with $|n - s|$ rows or columns of zeros, in order to get a square matrix $A_P^{\boxminus} = (a_{i,j} | 1 \le i, j \le \max(s,n))$. Jacobi's bound is then equal to $\mathcal{O}_{A_P^{\boxminus}}$, an equivalent definition that allows to compute the bound, with Jacobi's algorithm. (*Cf.* rem. 36 and 57.) *rectangular matrix*

Jacobi also introduces a determinant $\nabla$, the non vanishing of which is a necessary and sufficient condition for the bound to be reached. In order to define it, he considers the matrix $(\partial P_i/\partial x_j^{(a_{i,j})})$ and forms its determinant. Then, he only keeps the products $\pm\prod_{i=1}^n \partial P_i/\partial x_{\sigma(i)}^{(a_{i,\sigma(i)})}$ such that $\sum_{i=1}^n a_{i,\sigma(i)} = \mathcal{O}$. this is why he calls this expression the *truncated determinant*[7] of the system. We may equivalently use the following definition.

truncated determinant

**Definition 76.** — *Let $v \in \mathbf{N}^n$, we will use the notations $\lambda_P := \lambda_{A_P}$, $\mathrm{ord}^v_{x_j} P := \mathrm{ord}_{x_j} P - v_j$ and $\mathrm{ord}^v P := \max_{j=1}^n \mathrm{ord}^v_{x_j} P$.*

Jacobi's cover · minimal canon · rectangular matrix

*The definition of Jacobi's cover is extended to the case $s \neq n$, using rem 75, by setting $\lambda^{\boxminus}_{P,i} := \lambda_{A^{\boxminus}_P,i} - \min_{h=1}^s \lambda_{A^{\boxminus}_P,h}$, so that $\min_{i=1}^s \lambda^{\boxminus}_{P,i} = 0$, $\alpha_{P,i} := \max_{h=1}^s \lambda^{\boxminus}_{P,h} - \lambda^{\boxminus}_{P,i} = \alpha_{A^{\boxminus}_P,i} - \min_{h=1}^s \alpha_{A^{\boxminus}_P,h}$, so that $\min_{i=1}^s \alpha_{P,i} = 0$, and $\beta_{P,j} := \max_{i=1}^s a_{i,j} - \alpha_{P,i} = \beta_{A^{\boxminus}_P,i} + \min_{h=1}^s \alpha_{A^{\boxminus}_P,h}$.*

*Let $\mu, \nu$ be a minimal cover of the order matrix $A^{\boxminus}_P$. We denote by $\mathrm{J}^{\mu,\nu}_P$ the matrix $(\partial P_i/\partial x_j^{(\mu_i+\nu_j)}; (i,j) \in [1,s] \times [1,n])$, by $\mathrm{J}_P$ the matrix $\mathrm{J}^{\alpha_P,\beta_P}_P$ and $\mathrm{ord}^{\mathrm{J}}$ denotes $\mathrm{ord}^{\beta_P}$.*

system determinant

*If $s = n$, we call the* system determinant *and denote by $\nabla_P$ the determinant of $\mathrm{J}_P$. For $s < n$ (resp. $s > n$), $\nabla_P := \{{}_J\nabla_P | J \subset [1,n], \sharp J = \min(s,n)\}$, where ${}_J\nabla_P$ denotes the determinant of the submatrix of $\mathrm{J}_P$ with columns (resp. rows) in $J$.*

***Remarks.* — 77)** When $s \neq n$, the minimal canon $\lambda^{\boxminus}_P$ of $A^{\boxminus}_P$ is in general different from the minimal canon $\lambda$ of $A_P$ returned by Jacobi's algorithm. Both are defined so that $\min_{i=1}^s \lambda^{\boxminus}_i = \min_{i=1}^s \lambda_i = 0$.

**78)** If $s < n$, $\mathcal{O}_A = \sum_{i=1}^s a_{i,\sigma_i}$, $J = \mathrm{Im}(\sigma)$ and ${}_JA_P$ is the restriction of $A_P$ to columns in $J$, then $\alpha_P$, $(\beta_{P,j} | j \in J)$ is Jacobi's cover for ${}_JA_P$.

***Examples.* — 79)** Let $P_1 = x_1^{(4)} + x_2'' + x_3'$ and $P_2 = x_1^{(5)} + x_2''' + x_3''$. Then $\alpha_{A^{\boxminus}_P} = (1,2,0)$, $\beta_{A^{\boxminus}_P} = (3,1,0)$ and $\alpha_P = (0,1)$, $\beta_P = (4,2,1)$.

**80)** Let $P$ be a characteristic set and the main variable of $P_i$ be $x_{\sigma(i)}$. Then, if $\sum_{i=1}^s a_{i,\sigma(i)} = \mathcal{O}_P$, we have: $\alpha_i = 0$, $1 \le i \le n$ and $\beta_j = \mathrm{ord}_{x_j} P$. This is always the case when $s = n$.

It is straightforward that this definition of $\nabla_P$ is equivalent to Jacobi's *truncated determinant* (See [39, prop. II p. 17]), which is to keep in the Jacobian determinant only the terms corresponding to maximal sums in the order matrix. For any cover $\mu, \nu$, partial derivatives $\partial P_i/\partial x_j^{(\mu_i+\nu_j)}$ are indeed non identically zero iff $\mathrm{ord}_{x_j} P_j = \mu_i + \nu_j$.

**Proposition 81.** — *Let $\mu_i$, $\nu_j$ be a cover of $A^{\boxminus}_P$ and $J \subset [1,n]$ with $\sharp J = s \le n$ (resp. $\sharp J = n < s$), then ${}_J\nabla_P = \left(\partial P_i/\partial x_j^{(\mu_i+\nu_j)} | (i,j) \in [1,s] \times J \text{ (resp. } J \times [1,n])\right)$.*

Proof. — For $s \le n$, we only keep in $|\partial P_i/\partial x_j^{(a_{i,j})}; (i,j) \in [1,s] \times J|$ the products corresponding to a bijection $\sigma : [1,s] \mapsto J$ with $\sum_{i=1}^n a_{i,\sigma(i)} = \sum_{i=1}^n \alpha_i + \sum_{j=1}^n \beta_j = \mathcal{O}_P = \sum_{i=1}^n \mu_i + \sum_{j=1}^n \nu_j$. The case $s > n$ is similar, by transposition. ■

[7] *Determinans mancum*, or *Determinans mutilatum*.

## 4.2 Reduction to order 1

*Lando* We proceed with the well known reduction to first-order equations. Lando [68] *first order reduction*
*underdetermined system* did prove Jacobi's bound for systems of order one, also considering underdetermined systems, but only with the weak bound. She remarks that the weak bound for the first order reduction may be greater than that of the original system, but that the strong bound remains the same. We can even prove that the truncated determinant is unchanged, up to its sign. *Jacobi's number* *weak bound*

We introduce new variables $u_{j,k}$, for $1 \le j \le n$ and $0 \le k < r_j := \max_{i=1}^{s} \mathrm{ord}_{x_j} P_i$ and replace, in the $P_i$, $x_j^{(k)}$ by $u_{j,k}$, for $0 \le k < r_j$ and by $u'_{j,r_j-1}$ for $k = r_j$, obtaining a new differential polynomial $Q_i$. We complete the system $Q = 0$ with the equations $W_{j,k} := u_{j,k} - u'_{j,k-1} = 0$, for $1 \le j \le n$ and $0 < k < r_j$.

***Lemma 82.*** — *Let $P$ be a system of $n$ differential polynomials in $\mathscr{F}\{x_1, \ldots, x_n\}$, and $Q, W$ the system of $\sum_{i=1}^{n} r_i$ equations in $\mathscr{F}\{y\}$ obtained by reduction to the first order, as defined above, then $\mathscr{O}_P = \mathscr{O}_{Q,W}$, $\nabla_P = \pm\nabla_{Q,W}$ and in the minimal canon $\bar\lambda$ for the order matrix $B$ of the system $Q, W$ the integers $\bar\lambda_i$, $1 \le i \le n$, corresponding to the $Q_i$, are the same as the $\lambda_i$ in the minimal canon $\lambda$ of the order matrix $A$ of $P$.*

Proof. — For the columns of the order matrix, we use the order on the $u_{j,k}$ defined by $u_{j,k} < u_{j',k'} \Leftrightarrow j < j'$ or $j = j'$ and $k > k'$. For the rows, we choose to put first in the system the $P_i$, in the same order, and then the $W_{j,k}$, using the same order as for the $u_{j,k}$.

The order matrix $B$ of $Q, W$ has the following shape: $(L_1 \cdots L_n)$ with $L_j :=$

$$
\left(\begin{array}{ccccccccc}
[{}^{1,j}_{r_j-1}] & & \cdots & & [{}^{1,j}_{k}] & & \cdots & & [{}^{1,j}_{0}] \\
\vdots & & & & \vdots & & & & \vdots \\
[{}^{i,j}_{r_j-1}] & & \cdots & & \mathbf{[{}^{i,j}_{k}]} & & \cdots & & [{}^{i,j}_{0}] \\
\vdots & & & & \vdots & & & & \vdots \\
[{}^{n,j}_{r_j-1}] & & \cdots & & [{}^{n,j}_{k}] & & \cdots & & [{}^{n,j}_{0}] \\
 & & & & & & & & \\
 & & \cdots & & \sum_{j<j} r_j \text{ rows} & & \cdots & & \\
 & & & & & & & & \\
\mathbf{0} & 1 & & & & & & & \\
 & \ddots & \ddots & & & & & & \\
 & & \mathbf{0} & 1 & & & & & \\
 & & & \mathbf{0} & 1 & & & & \\
 & & & & 0 & \mathbf{1} & & & \\
 & & & & & 0 & \mathbf{1} & & \\
 & & & & & & \ddots & \ddots & \\
 & & & & & & & 0 & \mathbf{1} \\
 & & & & & & & & \\
 & & \cdots & & \sum_{j>j} r_j \text{ rows} & & \cdots & & \\
 & & & & & & & &
\end{array}\right).
$$

For more readability, only terms possibly different from $-\infty$ are displayed. The terms $[{}^{i,j}_{r_j-1}]$ are 1 if $\mathrm{ord}_{x_j} P_i = r_j$, 0 if $\mathrm{ord}_{x_j} P_i = r_j - 1$ and $-\infty$ otherwise. The terms $[{}^{i,j}_{k}]$ for $0 \le k < r_j - 1$ are 0 if $x_j^{(k)}$ appears in $P_i$ and $-\infty$ otherwise. It is easily seen that $n$ transversal integer elements in the first $n$ rows may be completed to a maximal set of transversal integer elements in one and only one way. Indeed, once an element $\mathbf{[{}^{i,j}_{k}]}$ is chosen in $L_j$, there is a unique choice of integer transversal elements in the rows corresponding in $L_j$ to the $W_{j,k}$. Among these elements that appear in **bold** in $L_j$, there are exactly $k$ **1**, so that their sum is equal to $k$ if $k < r_j - 1$, or equal to $\mathrm{ord}_{x_j} P_i = [{}^{i,j}_{k}] + k$ if $k = r_j - 1$.

Hence, the maximum of these sums is precisely $\mathrm{ord}_{x_j} P_i$ and there is a one to one correspondence between maximal transversal sums in $A$ and $B$, so that both matrices have the same tropical determinant.

The structure of the system determinant $\nabla_{Q,W}$ is similar, so that its terms are the same as in $\nabla_P$, up to the sign. An easy computation relying on the number of inversions shows that

$$\nabla_{Q,W} = |J_{Q,W}| = (-1)^{\sum_{j=1}^{n}(n-j+1)(r_j-1)} \nabla_P. \tag{7}$$

Indeed, for any term in a determinant that corresponds to the permutation $\sigma$, we may define an inversion as the case when $i' < i$ and $\sigma(i') > \sigma(i)$. So the number of inversions is $\sum_i \iota_i$ with $\iota_i := \sharp\{i' \in \mathbf{N} | i' < i,\ \sigma(i') > \sigma(i)\}$. Then, the signature $e(\sigma)$ of a permutation being determined by the number $\rho$ of transpositions $\tau_j$ in any decomposition $\sigma = \tau_1 \circ \cdots \tau_\rho$, $e(\sigma) = (-1)^{\sum_i \iota_i}$.

Any permutation $\sigma$ of $A$ is associated to a single permutation of $B$ and the inversion numbers $\iota_i$ for the first $n$ rows are the same for $A$ and $B$. Then, one only has to consider the remaining rows of $J_{Q,W}$ that are decomposed in blocks $M_J :=$

$$\left(\begin{array}{ccccccccc}
 & \cdots & & \multicolumn{3}{c}{\sum_{j'<j} r_{j'} \text{ rows}} & & \cdots & \\
\mathbf{1} & -1 & & & & & & & \\
 & \ddots & \ddots & & & & & & \\
 & & \mathbf{1} & -1 & & & & & \\
 & & & \mathbf{1} & -1 & & & & \\
 & & & & 1 & \mathbf{-1} & & & \\
 & & & & & 1 & \mathbf{-1} & & \\
 & & & & & & \ddots & \ddots & \\
 & & & & & & & 1 & \mathbf{-1} \\
 & \cdots & & \multicolumn{3}{c}{\sum_{j'>j} r_{j'} \text{ rows}} & & \cdots &
\end{array}\right).$$

The presence of the $-1$ terms compensates the extra inversions, so that we may reduce to the case of the diagonal of 1. For each of them, the total inversion number is $n - j + 1$, corresponding to rows $j$ to $n$. This produces a total of $(n - j + 1)(r_j - 1)$, hence formula (7).

In order to show that $\bar\lambda_i = \lambda_i$, $1 \le i \le n$, we may consider first the conditions on the terms 1 and 0 in a maximal sum, that needs to become maximal elements in $B + \bar\lambda$. It is convenient to denote by $\bar\lambda_j^k$ the element of $\bar\lambda$ corresponding to $W_{j,k}$. Assuming that $\sum_j a_{j,j}$ is a maximal transversal sum for $A$, if $\bar\lambda$ is a canon of $B$, the transversal elements in the rows of the $W_{j,k}$ must be maximal, which is equivalent *transversal maxima* to

$$\begin{aligned} &\bar\lambda_j^k \ge \left[{}^{i,j}_{k}\right] + \bar\lambda_i, \text{ and } \bar\lambda_j^k \ge \bar\lambda_j^{k+1} + 1, \text{ for } 1 \le i \le n \text{ and } k > a_{j,j},\\ &\bar\lambda_j^k \ge \left[{}^{i,j}_{k}\right] + \bar\lambda_i - 1, \text{ and } \bar\lambda_j^k \ge \bar\lambda_j^{k-1} - 1, \text{ for } 1 \le i \le n \text{ and } k \le a_{j,j}. \end{aligned} \tag{8}$$

So, $\bar\lambda_j^k \ge \ell_j^k(\bar\lambda)$ with

$$\begin{aligned} \ell_j^k(\bar\lambda) \quad &:= \max_{h=k}^{r_k-1} \max(0, h - k + \max_{i=1}^{n} \left[{}^{i,j}_{h}\right] + \bar\lambda_i), &&\text{for } k > a_{j,j},\\ &:= \max_{h=0}^{k-1} \max(0, h - k - 1 + \max_{i=1}^{n} \left[{}^{i,j}_{h}\right] + \bar\lambda_i), &&\text{for } k \le a_{j,j}. \end{aligned} \tag{9}$$

The matrix $L_j + \bar\lambda$ has the following shape.

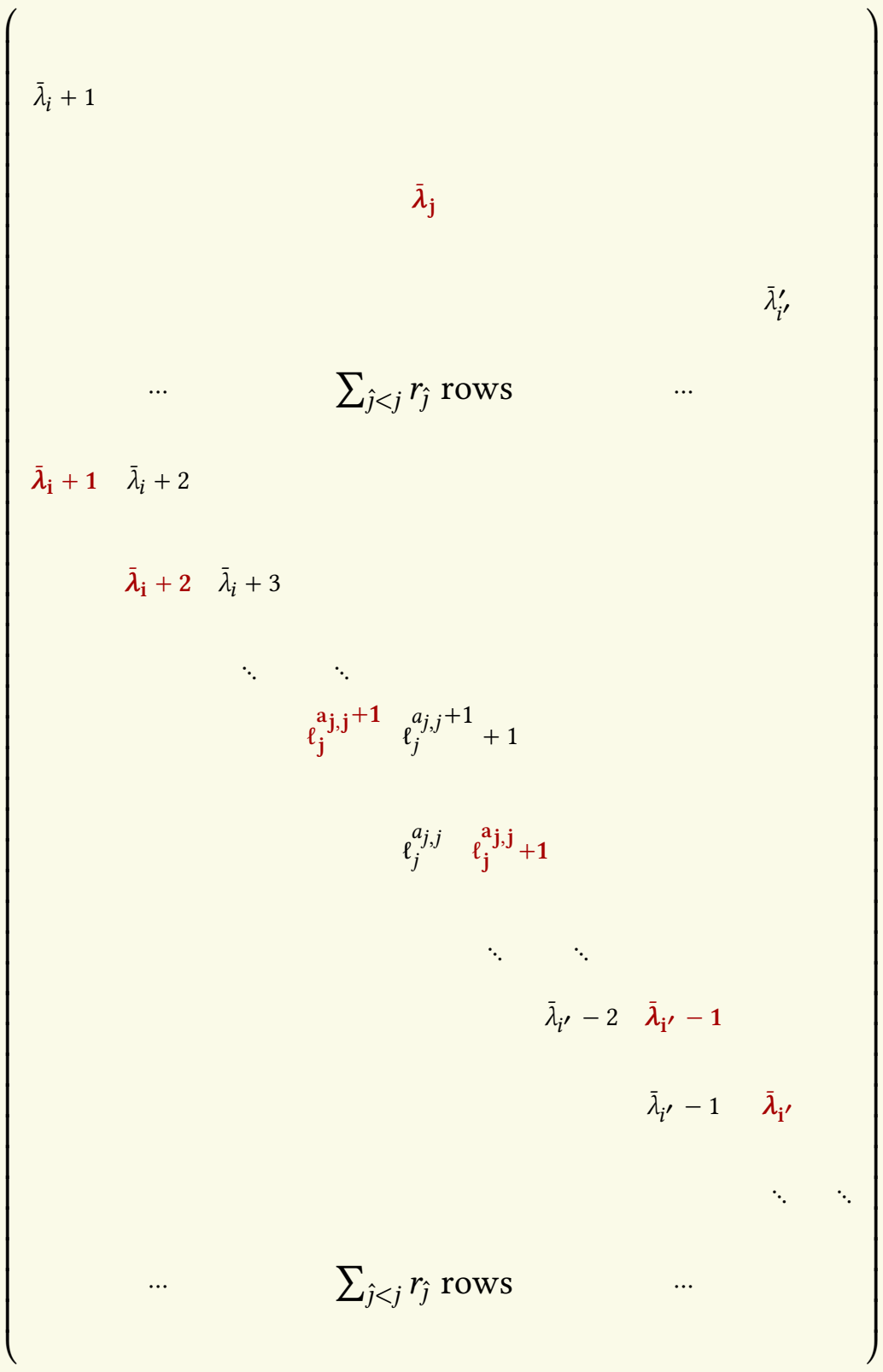

There are only two terms to consider in the column corresponding to $u_{j,a_{j,j}}$, the first from the row of $W_{j,a_{j,j}+1}$ is $\ell_j^{a_{j,j}+1}(\bar\lambda) + 1$, the second from the row of $W_{j,a_{j,j}}$ is

$\ell_j^{a_{j,j}}(\bar\lambda)$. So, the following inequality means that $a_{i,j}$ is maximal:

$$\bar\lambda_j \geq \max(\ell_j^{a_{j,j}+1}(\bar\lambda)+1, \ell_j^{a_{j,j}}(\bar\lambda)) = \max(0, \max_{i\neq j} a_{i,j} + \bar\lambda_i - a_{j,j}). \tag{10}$$

The vector $\lambda$ is the smallest positive solution of this system of inequalities, so that $\bar\lambda_i \geq \lambda_i$ for $1 \leq i \leq n$. But, replacing $\bar\lambda_i$ by $\lambda_i$ and $\bar\lambda_j^k$ by $\ell_j^k(\bar\lambda)$ in formula (9), we obtain a value for $\bar\lambda$ that satisfies both inequalities (10) and (8), meaning that the transversal terms with maximal sum are maximal, so that it is indeed the minimal canon for $B$ and $\bar\lambda_i = \lambda_i$ for $1 \leq i \leq n$. ■

## 4.3 Block decompositions

If the integer elements in the order matrix admit a minimal cover of $p$ rows and $q$ columns with $0 < p, q < n$, then the system $P$ admits a non trivial triangular block decomposition. In the case where $s = n$ and $p + q = n$, one may look for such a block decomposition using the reflexive transitive closure of the elementary path relation, as defined in subsection 3.3, but for a matrix $B$ with $b_{i,j} := 1$ if $x_j$, or some of its derivatives, appear in $P_i$ and $b_{i,j} := 0$ if not. One gets so a partial preorder that defines equivalence classes of rows $i, j$ with $i \prec j$ and $j \prec i$. Sorting the variables and equations according to this preorder produces a triangular block decomposition, the blocks corresponding to equivalence classes, that do not depend on the choice of a maximal transversal family, by prop. 54.

In the same spirit, considering the reflexive, transitive and symmetric closure provides a diagonal block decomposition. We will not develop these easy results, but they can be very helpful to clarify the structure of a system before any attempt to solve it, whenever its size makes difficult to find the requested form by simple inspection.

*block decomposition*

We may also consider the preorder defined by that path relation for a canon and the triangular block decomposition it defines. We have the next proposition.

**Proposition 83.** — *i) We define the equivalence relation $i \equiv i'$ if there is a path from $i$ to $i'$ and from $i'$ to $i$; it does not depend on the chosen canon.*

*ii) Let $E_p$, $1 \leq p \leq r$ be the equivalence classes defined by this relation, and $a_{i,\sigma(i)}$ a maximal transversal family, then the sets $\sigma(E_p)$, $1 \leq p \leq r$, do not depend on the chosen permutation $\sigma$ and $\nabla = \prod_{p=1}^s \nabla_p$, where $\nabla_p$ is the determinant of the submatrix of $\mathrm{J}$ corresponding to rows in $E_p$ and columns in $\sigma(E_p)$.*

Proof. — i) If there is a path from $i$ to $i'$ and from $i'$ to $i$, for some canon $\ell$ then by prop. 33, for any row $j$ in the defined cycle, and any canon $\ell'$, $\ell'_j - \ell_j = \ell'_i - \ell_i$, so that the same cycle exists in the canon $\ell'$.

ii) Independence with respect to the choice of $\sigma$ is a consequence of prop. 54. Then, the matrix $\mathrm{J}$ is block triangular and its determinant is the product of the determinants of the diagonal blocks, which are the $\nabla_p$. ■

# 5 An algebraic parenthesis. Quasi-regularity and "Lazard's lemma"

JACOBI considers functions without any precision about their nature. One may present his results in the framework of *diffiety theory*, provided that the equations are defined by $\mathscr{C}^\infty$ functions, satisfying some natural regularity hypotheses (see [79]). We use here the formalism of Ritt's differential algebra, that allows effective computations. Here, characteristic sets will be used instead of Jacobi's "normal forms", and Lazard's lemma will take the place of the implicit function theorem, the neighborhood of a point being then a Zariski open space, dense in a whole component.

*diffiety theory* · *Ritt* · *Lazard's lemma* · *normal form*

## 5.1 Quasi-regularity

As we will see, quasi-regularity is an implicit hypothesis, that plays a central role in Jacobi's proof of the bound. The informal meaning of this notion is that a differential system $P_i(x) = 0$, $1 \leq i \leq s$, "behaves like" the linearized system $\mathrm{d}P_i = 0$, *viz.* $\sum_{j=1}^{n} \sum_{k=0}^{a_{i,j}} \partial P_i / \partial x_j^{(k)}\, \mathrm{d}x_j^{(k)} = 0$, in the neighborhood of a generic point of some component of $\{P\}$.

*quasi-regularity* · *linearized system*

This idea was formalized by Johnson [50, 51, 52] who used it to prove Janet's conjecture [53]. It is also the key of the first partial proof of Jacobi's bound in the non linear case, given by Kondratieva *et al.* [60, 61]. Only in the case $s = n = 2$, was Ritt [86, chap. VII § 6. p. 136] able to prove the bound for general components, that is without the quasi-regularity hypothesis.

*Janet's conjecture* · *Johnson, Joseph* · *Jacobi's bound* · *Kondratieva* · *Ritt* · *quasi-regularity*

We will provide here a more general definition than Kondratieva's one [61, 79], in order to underline that the property used is weaker than the "independence" of Kähler differentials $\mathrm{d}P_i$ of which it is a consequence. Quasi-regular was chosen because this property is shared by some components of a differential equation, that the classical theory considers as "singular" (See Houtain [34] or Hubert [36]); it is indeed a property, not of the system alone, but of a component, with respect to the system.

*Kähler's differentials* · *Houtain; Hubert* · *singular component*

In the following, $\mathscr{F}$ will denote a differential field of characteristic 0. We refer to Ritt [86] and Kolchin [59] for more details about differential algebra, and to Boulier [6] for characteristic sets. It is natural here to state the definition for an arbitrary differential field $\mathscr{F}_\Delta$, with a finite set $\Delta := \{\delta_1, \ldots, \delta_m\}$ of commuting derivations, possibly empty. We will use Kähler's differentials for differential fields and rings extensions, denoted by $\Omega_{A/B}$, for which we refer to Johnson [51].

*differential algebra* · *Ritt and Kolchin* · *characteristic set* · *Boulier* · *Kähler's differentials* · *Johnson, Joseph*

**DEFINITION 84.** — *Let $\mathscr{G}_\mathscr{P}/\mathscr{F}$ denote the differential field extension defined by a prime differential ideal $\mathscr{P} \subset \mathscr{F}\{x_1, \ldots, x_n\}$, $\Delta$ denote the set of commuting derivations of the differential fields $\mathscr{F}$ and $\mathscr{G}_\mathscr{P}$ and $\Theta$ the commutative monoïd it generates. We denote by $\mathscr{D}_\mathscr{P} := \mathscr{G}_\mathscr{P}[\Delta]$ the non commutative ring of differential operators and by $\mathscr{M}_\mathscr{P}$ the $\mathscr{D}_\mathscr{P}$-module $\mathscr{G}_\mathscr{P} \otimes_{\mathscr{F}\{x\}} \Omega_{\mathscr{F}\{x\}/\mathscr{F}}$. For any $Q \in \mathscr{F}\{x\}$, $\mathrm{d}_\mathscr{P} Q \in \mathscr{M}_\mathscr{P}$ denotes the differential of $Q$.*

*differentials* *The module $\mathscr{M}_{\mathscr{P}}$ is a $\mathscr{D}_{\mathscr{P}}$-free module generated by the* $\mathrm{d}_{\mathscr{P}} x_i$, $1 \le i \le n$.

*Let $P_i$ $1 \le i \le s$ be differential polynomials in $\mathscr{F}\{x_1, \ldots, x_n\}$, and $\{P\} = \bigcap_{j=1}^{r} \mathscr{P}_j$, where the $\mathscr{P}_j$ are prime differential ideals such that $\mathscr{P}_i \subset \mathscr{P}_j$ implies*[8] *$i = j$. The prime ideals $\mathscr{P}_j$ are called the* components *of $\{P\}$ and by extension the components of $P$.*

*quasi-regularity* *i) We say that $\mathscr{P}_j$ is a* quasi-regular *component of the system $P = \{P_1, \ldots, P_s\}$ if* $\mathrm{d}_{\mathscr{P}_j} \mathscr{P}_j = (\mathrm{d}_{\mathscr{P}_j} P)_{\mathscr{M}_{\mathscr{P}_j}}$.

*strongly quasi-regular component* *ii) We say that $\mathscr{P}_j$ is* strongly quasi-regular *if the family $(\theta\, \mathrm{d}_{\mathscr{P}_j} P_i)$, $\theta \in \Theta$, is linearly independent over $\mathscr{G}_{\mathscr{P}_j}$.*

*regular component* *iii) In the case $\sharp\Delta = 1$, a) a prime component $\mathscr{P}$ of $\{P\}$ is said to be* regular[9] with respect to $\mathrm{ord}^{v}$ *if there exists $S \in \mathrm{R}_r := \mathscr{F}[x_j^{(k)} | 1 \le j \le n,\ 0 \le k \le r + v_j]$, with $r := \mathrm{ord}^{v} P$ such that $\mathscr{P} = \{P\} : S^{\infty}$, where $\mu, v$ is a minimal cover of $P$; b) it is said to be* regular *if it is regular with respect to* $\mathrm{ord}^{\mathrm{J}}$.

Strong quasi-regularity implies quasi-regularity (lem. 85 iv). It corresponds to
*Kondratieva, Mikhalev and Pankratiev* the *independence* assumption in Kontratieva, Mikhalev and Pankratiev [61] and
*Sadik* also to quasi-regularity as formerly defined with Sadik [79] in the framework of diffiety theory. This property may also be shared by some singular components, as defined in iii)[10]. In most practical situations, strong quasi-regularity is easier to prove *a priori*, but quasi-regularity is enough to conclude that the differential Hilbert function of the differential ideal $(\mathrm{d}_{\mathscr{P}} P)$ and of the $\mathscr{D}_{\mathscr{P}}$-module $\mathrm{d}_{\mathscr{P}_j} \mathscr{P}_j$ are the same.

Quasi-regularity is thus very useful, mostly combined with the properties of $\mathrm{d}_{\mathscr{P}_j} \mathscr{P}_j$ in the next lemma, as it allows to reduce to the linear case. In the sequel, we will omit the index "$_{\mathscr{P}_j}$" for the differential "d" operator or other notations when there is no ambiguity.

*linearized system* **Lemma 85. —** *i) If $\mathscr{P}$ is a quasi-regular component, we have: a) A characteristic set $\mathscr{A}$ exists for $\mathscr{P}$ for some ordering $\prec$ on derivatives with main derivatives $v_A$, $A \in \mathscr{A}$ iff a standard basis*[11] *exists for the $\mathscr{D}$-module $\mathrm{d}\mathscr{P}$, for the ordering induced by $\prec$ on differentials $\mathrm{dY}$ (where $\mathrm{Y} := \Theta X$ denotes the set of derivative of the $x_i$), with main derivatives $\mathrm{d}v_A$, $A \in \mathscr{A}$, where $v_A$ denotes the main derivative of $A$;*

*b) Let $\mathrm{Y} \subset \{x_1, \ldots, x_n\}$, $\mathscr{P} \cap \mathscr{F}\{Y\} \ne (0)$ iff $(\mathrm{d}_{\mathscr{P}} P) \cap (\mathrm{d}Y) \ne (0)$.*

*ii) The main component $\mathscr{P} = [\mathscr{A}] : H_{\mathscr{A}}^{\infty}$ of a system $\mathscr{A}$ that is a characteristic set*
*quasi-regularity* *of some prime differential ideal $\mathscr{P}$ is quasi-regular and strongly quasi-regular in the ordinary case;*

*iii) a) Any quasi-regular component $\mathscr{P}$ of $P$ has diff. codimension at most $s = \sharp P$.*

---

[8]Such a decomposition is known to be unique, see Ritt [86, ch. II § 19] or Kolchin [59, th. 1, p. 14].

[9]Ritt [86, chap. II § 17 p. 32] defines singular components only in the case of an ideal generated by a single polynomial.

[10]The meaning of this classical notion of regularity is of a different nature and will be investigated in subsec. 7.4.

*Castro-Jiménez* [11]We refer to Castro-Jiménez [9, 10] for more details on standard (or Gröbner) bases of $\mathscr{D}$- *standard or Gröbner bases*
modules.

*b) A component $\mathscr{P}$ of $P$ is a strongly quasi-regular component iff it is quasi-regular and of differential codimension $s$.*

Proof. — i) a) $\Longrightarrow$. Assume that, for some ordering $\prec$, $\mathscr{A}$ is a characteristic set of $\mathscr{P}$. Then any $Q \in \mathscr{P}$ is reducible by $\mathscr{A}$, so that $\mathrm{d}Q$ is also reducible by $\mathrm{d}\mathscr{A}$ for the ordering induced by $\prec$. This implies that $\mathrm{d}\mathscr{A}$ is a standard basis for the $\mathscr{D}$-module $\mathrm{d}\mathscr{P}$ for this same ordering *standard basis*

$\Longleftarrow$. If $G$ is a standard basis of $\mathrm{d}\mathscr{P}$ for some ordering, consider a characteristic set $\mathscr{A}$ of $\mathscr{P}$ for the corresponding ordering. By what precedes, $\mathrm{d}\mathscr{A}$ is also a standard basis for the same ordering and $G$ and $\mathrm{d}\mathscr{A}$ have the same leading terms, hence the result.

b) Using a), it is enough to consider an ordering $\prec$ that eliminates letters not in $Y$.

ii) It is a straightforward consequence of i).

iii) a) If $\mathscr{P}$ is quasi-regular, then its differential codimension is the dimension $s = \sharp P$ of the free $\mathscr{D}$-module generated by the $\mathrm{d}P_i$, minus the dimension of the module of relations between the $\mathrm{d}P_i$.

b) $\Longrightarrow$. So, if the differential codimension is $s$, the $\mathrm{d}P_i$ and their derivatives must be independent, meaning that $\mathscr{P}$ is strongly quasi-regular.

$\Longleftarrow$. Let $\mathscr{A}$ be a differential characteristic set of $\mathscr{P}$. It may be extracted from an algebraic characteristic set of a prime component $\mathscr{I}$ of $\sqrt{(P, \Delta P, \ldots, \Delta^r P)}$ for $r$ great enough. The codimension of $\mathscr{I}$ is at most $s\binom{r+m}{m}$, which is the number of its generators. As the derivatives of the $\mathrm{d}P$ are independent, the codimension is exactly $s\binom{r+m}{m}$, so that $\mathrm{d}\mathscr{I} = \langle \mathrm{d}P, \mathrm{d}\Delta P, \ldots, \mathrm{d}\Delta^r P\rangle$ and $\mathrm{d}\mathscr{A} \in (\mathrm{d}P)$, which implies that $\mathscr{P}$ is quasi-regular. The differential codimension of $\mathscr{P}$ is then the dimension of the $\mathscr{D}$-module $(\mathrm{d}P)$, which is $s$. ■

***Remark* 86.** — Proving that the differential codimension of a system of $s$ equations is at most $s$, is easy in the quasi-regular case, but the general case is a difficult conjecture: the *dimensional conjecture*. R. Cohn [13] has shown that it would be implied by Jacobi's bound (even weak) for arbitrary systems. *Cohn, Richard* *dimensional conjecture*

***Examples.* — 87)** The component $x$ of $(x')^p - p^p x^{p-1}$, $1 < p$, is singular. It is strongly quasi-regular only when $p$ is equal to 1 or 2.

**88)** If $\mathscr{A}$ is a characteristic set of a prime ideal $\mathscr{P}$, this ideal is strongly quasi-regular and regular with respect to $\mathrm{ord}^\nu$ for any $\nu$. *characteristic set*

In the "algebraic" case, that is when $\Delta = \emptyset$, Lazard's lemma provides a simple criterion for quasi-regularity.

## 5.2 Lazard's lemma

Many proofs of this folkloric result are already available in the differential algebra literature (see *e.g.* Morrison [74] or Boulier *et al.* [6]). The motivation of the following one is to make a link with the implicit function theorem by using Newton's method. *Boulier* *Morrison*

**Theorem 89.** — *Let $P_1, \ldots, P_s$ be polynomials in $k[x_1, \ldots, x_n]$ with $s \le n$ and $\mathrm{J} := (\partial P_i/\partial x_j | 1 \le i, j \le s)$. If $\mathcal{I} := (P) : |\mathrm{J}|^\infty \ne (1)$, then*

*i) Any component $\mathcal{P}$ of $\mathcal{I}$ is a strongly quasi-regular component of the system $P$, of codimension $s$;*

*ii) For any component $\mathcal{P}$ of $\mathcal{I}$, $\mathcal{P} \cap k[x_{s+1}, \ldots, x_n] = (0)$, so $\mathcal{I} \cap k[x_{s+1}, \ldots, x_n] = (0)$;*

*iii) $\mathcal{I}$ is radical.*

Proof. — i) We notice that $\mathrm{d}P_i = \sum_{j=1}^n \partial P_i/\partial x_j \,\mathrm{d}x_j$, so that the differentials $\mathrm{d}P_i$, $1 \le i \le s$, are linearly independent, $\mathcal{P}$ is a strongly quasi-regular component of the system $P$ and $\operatorname{codim}\mathcal{P} = s$ by lem. 85 iii).

ii) It is a straightforward consequence of i) and lemma 85 i).

iii) Using Rabinowitsch's trick [82], $\bar{\mathcal{I}} := (P, u|\mathrm{J}|)$ is such that $\bar{\mathcal{I}} \cap k[x] = \mathcal{I}$, so that it is enough to prove that $\bar{\mathcal{I}}$ is radical. We will work in the module $\mathcal{M} := k[x,u]/\mathcal{Q} \otimes_{k[x,u]} \Omega_{k[x,u]/k}$, where $\mathcal{Q} = \sqrt{\bar{\mathcal{I}}}$. Let $Q \in \mathcal{Q}$, we will show that $Q \in \bar{\mathcal{I}}$. We need first the following lemma.

***Lemma* 90.** — *With the above hypotheses, for any $Q \in \mathcal{Q}$, $\mathrm{d}Q (\in \mathcal{M}) = 0$ iff $Q \in \mathcal{Q}^2$.*

Proof of the lemma. — $\Longrightarrow$. Immediate. $\Longleftarrow$. The condition $|\mathrm{J}| \ne 0$ implies that any affine solution $\eta_1$ of $\mathcal{Q}$ is regular and that $\mathcal{Q}$ defines components of codimension $s$. In the projective space, their parts at infinity are of codimension $s+1$. Then, one may find an affine space $L_i(x) = 0$, $1 \le i \le s$, of dimension $s$ that contains $\eta_1$ and no point solution of $Q$ at infinity, so that, for any finite set $G$ that generates $\mathcal{Q}$ as an ideal, $G, L$ defines $d$ isolated points $\eta_k$, $1 \le k \le d$, where $d$ is the sum of the degrees of the prime components $\mathcal{P}_i$ of $\mathcal{Q}$. Then, there exists a generic homogeneous linear form $L_0$, such that the values $L_0(\eta_k)$, $1 \le k \le d$ are paiwise different. We may work in the new coordinates $y_0 = u$, $y_i = x_i$ for $1 \le i < s$ and $y_{s+i} = L_i$ for $0 \le i \le n-s$.

In such coordinates, with the ordering $y_i > y_j$ if $i < j$, $\mathcal{Q}$ admits a char. set $\mathcal{A}$ with $A_s(y_s, \ldots, y_n)$ of degree $d$, with coefficient in $y_s^d$ equal to 1. If the discriminant of $A_s$ does not vanish, the $d$ values of $y_s$ solution of $A_s$ are all distinct, so that, for $1 \le i < s$, one may choose $A_i$ linear in $y_i$ and such that its initial $In_i$ does not vanish if the discriminant of $A_s$ does not vanish. (See Giusti, Lecerf and Salvy [30, lem. 1 § 3] for more details.) So, for any solution $\eta$ of $\mathcal{Q}$, one may find a char. set $\mathcal{A}_\eta \subset \mathcal{Q}$, such that its product of initials and separants $H_\eta$ does not vanish at $\eta$.

*Giusti, Lecerf and Salvy*

Using the reduction process, we have, for some $q \in \mathbf{N}$, $H_\eta^q Q = \sum_{i=1}^s M_{\eta,i} A_{\eta,i}$. So, if $\mathrm{d}Q = 0$, then for all $1 \le j \le s$, $\sum_{i=1}^s M_{\eta,i} \partial A_{\eta,i}/\partial y_j \in \mathcal{Q}$. As $|\partial A_{\eta,i}/\partial y_j| = \prod_{i=1}^s Sep_i$, that divides $H_\eta$, we have $H_\eta M_{\eta,i} \in \mathcal{Q}$ and $H_\eta^{q+1} Q \in \mathcal{Q}^2$. By Noetherianity, there is a finite family of solutions of $\mathcal{Q}$, $\eta_\ell$, $\ell \in [1, p]$, such that $1 \in (H_{\eta_1}, \ldots, H_{\eta_p}) + \mathcal{Q}$, so that $Q \in \mathcal{Q}^2$. This concludes the proof of the *lemma*. ∎

Let $\tilde{\mathrm{J}}$ denote the adjugate matrix of $\mathrm{J}$, $\psi_0(Q) := Q$, $\phi_0(Q) := 0$ and

$$\psi_1(Q) := u(\partial Q/\partial x_1, \ldots, \partial Q/\partial x_s)\tilde{\mathrm{J}}(P_1, \ldots, P_s)^{\mathrm{t}}, \ \phi_1(Q) = Q - \psi_1(Q).$$

We have $d\phi_1(Q) = 0 \in \mathscr{M}$, so that, by lem. 90, $\phi_1(Q) \in \mathscr{Q}^2$. We may write

$$\phi_1(Q) = \sum_{(g_1,g_2)\in G^2} M_g(Q) g_1 g_2,$$

where $G$ generates $\mathscr{Q}$ as an ideal. We recursively define $\psi_r(Q)$ and $\phi_r(Q)$, for $r > 1$, by the formulae

$$\psi_{r+1}(Q) := \psi_1(Q) + \sum_{(g_1,g_2)\in G^2} M_g(Q)[\psi_r(g_1)\psi_r(g_2) + \psi_r(g_1)\phi_r(g_2) + \phi_r(g_1)\psi_r(g_2)]$$

and

$$\phi_{r+1}(Q) := \sum_{(g_1,g_2)\in G^2} M_g(Q)\phi_r(g_1)\phi_r(g_2).$$

We easily prove by recurrence that $\psi_r(Q) + \phi_r(Q) = Q$, $r \in \mathbf{N}$. The result is true for $r = 0$ and $r = 1$. Assume it stands for $r$. Then, $\psi_{r+1}(Q) + \phi_{r+1}(Q) =$

$$\begin{aligned} & \psi_1(Q) + \textstyle\sum_{(g_1,g_2)\in G^2} M_{g_1,g_2}(Q)[\psi_r(g_1) + \phi_r(g_1)][\psi_r(g_2) + \phi_r(g_2)] \\ = \ & \psi_1(Q) + \textstyle\sum_{(g_1,g_2)\in G^2} M_{g_1,g_2}(Q) g_1 g_2 \\ = \ & Q. \end{aligned}$$

We prove in the same way that $\psi_r(Q) \in \mathscr{I}$ and $\phi_r(Q) \in \mathscr{Q}^{2^r}$. By noetherianity, there exists $r_0$ such that $\mathscr{Q}^{2^{r_0}} \subset \mathscr{I}$, so $Q = \psi_{r_0}(Q)[\in \mathscr{I}] + \phi_{r_0}(Q)[\in \mathscr{Q}^{2^{r_0}} \subset \mathscr{I}] \in \mathscr{I}$. ∎

***Example* 91.** — Let $P_1 := [(x_2 + x_1)^2 - 2x_3^2](x_1 + x_2 - x_3)(x_1 + x_2 + x_3)^2$, $P_2 := x_2^2 + (x_1 + x_2)^2 - 2x_3^2$ in the polynomial ring $\mathbf{Q}[x_1, x_2, x_3]$. The ideal $[P] : \mathrm{J}^\infty$ admits 2 prime components with characteristic sets $\{x_1 - 2x_3, x_2 + x_3\}$ and $\{x_1, x_2 - x_3\}$.

We now give a flexible corollary that is a real Swiss knife to secure sharp bounds in many situations, mostly when computing characteristic sets with elimination orderings.

**Corollary 92.** — *i) Assume that $P_{s_0+1}, \ldots, P_s$ belong to $k[x_{s_0+1}, \ldots, x_n]$. Then, $(P) : |\mathrm{J}|^\infty \cap k[x_{s_0+1}, \ldots, x_n] \supset (P_{s_0+1}, \ldots, P_s) : |\mathrm{J}|^\infty \supset (P_{s_0+1}, \ldots, P_s) : |\mathrm{J}_0|^\infty$, where $\mathrm{J}_0$ is the Jacobian matrix of the polynomials $P_{s_0+1}, \ldots, P_s$ with respect to the variables $x_{s_0+1}, \ldots, x_s$.*

*ii) Moreover, if a) $\sqrt{(\mathrm{J})} = \sqrt{(\mathrm{J}_0)}$ and b) the $P_i$, $1 \le i \le s_0$ are linear in $x_j$, $1 \le j \le s_0$, $(P) : |\mathrm{J}|^\infty \cap k[x_{s_0+1}, \ldots, x_n] = (P_{s_0+1}, \ldots, P_s) : |\mathrm{J}_0|^\infty$.*

Proof. — i) Let $\mathscr{P}$ be a component of $[P] : |\mathrm{J}|^\infty$. As $|\mathrm{J}| = |\mathrm{J}_1||\mathrm{J}_0|$, where $\mathrm{J}_1 := (\partial P_i/\partial x_j | 1 \le i, j \le s_0)$, $|\mathrm{J}| \notin \mathscr{P}$ and $|\mathrm{J}_1| \notin \mathscr{P}$, so $\langle dP_1, \ldots, dP_{s_0}\rangle \cap \langle dx_{s_0+1}, \ldots, dx_s\rangle$ is equal to 0. This implies that $\langle dP\rangle \cap dx_{s_0+1}, \ldots, dx_s\rangle$ is equal to $\langle dP_{s_0+1}, \ldots, dP_s\rangle$ and $\mathscr{P} \cap k[x_{s_0+1}, \ldots, x_s]$ is a component of $(P_{s_0+1}, \ldots, P_s) : |\mathrm{J}|^\infty$, which obviously contains $(P_{s_0+1}, \ldots, P_s) : |\mathrm{J}_0|^\infty$, as $|\mathrm{J}_0|$ divides $|\mathrm{J}|$.

ii) It is easily seen that any solution $(x_{s_0+1}, \ldots, x_n)$ of $P_i$, $s_0 < i \le s$, such that $\mathrm{J}_0$ does not vanish may be prolongated to get a solution $(x_{s_0+1}, \ldots, x_n)$ of $P_i$, $1 \le i \le s$ by solving a linear system with a determinant equal to $|J|/|J_0|$, that does not vanish if $|J_0|$ does not vanish. ■

The well known *Jacobian conjecture* (Keller [56], see also Bass *et al.* [3]) would imply that condition b) may be omitted in ii).

## 6 Jacobi's proof of the bound

In manuscript [II/13 b)] (*cf.* [39, prop. 1 p. 16 and prop. 2 p. 17]), Jacobi gives two different versions of this result. In the first, he writes that the order[12] is $\mathcal{O}$, but in the second, he claims that the order is $\mathcal{O}$ iff the truncated determinant vanishes. A modern reader may be surprised by this way of giving in a first theorem a generic result, and then describing more precisely, in a second theorem, the possible exceptions to the first one. Jacobi's proof also contains paradoxical arguments that led Ritt to conclude it was whimsical.

It seems however possible to save the proof and get the second version of *Kondratieva* Jacobi's theorem, more precise in the ordinary case than that of Kondratieva *et al.* [61], for it provides also a necessary and sufficient condition of equality of the order to Jacobi's bound.

Another easy improvement is to replace Jacobi's bound by the following def- *diffiety theory* inition, related to a precise component of a system, following the diffiety theory version, given with Sadik [79, def. 4.2].

*order w.r.t. a component* **Definition 93.** — *Let $\mathcal{P} \in \mathcal{F}\{x\}$ be a prime ideal. For any $P_i \in \mathcal{F}\{x\}$, we denote*[13] *by ${}^{\mathcal{P}}a_{i,j}$ or $\mathrm{ord}_{\mathcal{P}, x_j} P_i$ the least integer $k$ such that $\partial P_i/\partial x_j^{(k)} \notin \mathcal{P}$, or $-\infty$ if $\partial P_i/\partial x_j^{(k)} \in \mathcal{P}$ for all integer $k$.*

*Jacobi's number w.r.t. a component* *For any system $P \subset \mathcal{F}\{x\}$, the order matrix $A_{\mathcal{P}} = ({}^{\mathcal{P}}a_{i,j})$, the corresponding Jacobi* *Jacobi's cover w.r.t. a component* *bound $\mathcal{O}_{\mathcal{P}, P}$, Jacobi's cover $\alpha_{\mathcal{P}, P}, \beta_{\mathcal{P}, P}$, the canon $\lambda_{\mathcal{P}, P}$, $\mathrm{J}_{\mathcal{P}, P} := (\partial P_i/\partial x_j^{(a_{\mathcal{P},i,j})})$,* *system determinant w.r.t. a component* *$\nabla_{\mathcal{P}} := |\mathrm{J}_{\mathcal{P}}|$ and the twisted order $\mathrm{ord}^{\mathrm{J}_{\mathcal{P}}}$ are defined accordingly for any component $\mathcal{P}$ of $\{P\}$.*

Obviously, $\mathcal{O}_P \geq \mathcal{O}_{\mathcal{P}}$ for any component $\mathcal{P}$ of $P$. This new definition provides theoretical sharper bounds, without any change in the scheme of the proof of the theorem, that we express first in its original form, that is assuming $s = n$, before *underdetermined system* giving an easy generalization to under- or overdetermined systems.

*Jacobi's bound (proof)* **Theorem 94.** — *Let $\mathcal{P}$ be a strongly quasi-regular component $\mathcal{P}$ of the system $P_i$, $1 \leq i \leq n$ of differential polynomials in $\mathcal{F}[x_1, \ldots, x_n]$.*

*i) The order of $\mathcal{P}$ is at most $\mathcal{O}_{\mathcal{P}}$.*

*ii) The order of $\mathcal{P}$ is equal to $\mathcal{O}_{\mathcal{P}}$ iff $\nabla_{\mathcal{P}} \notin \mathcal{P}$.*

---

*Jacobi's number* [12] Jacobi used many different notations for his bound. In our translation [39], we used $H$ for *Borchardt* consistency, following Borchardt's posthumous edition [Crelle 64]. The bound is actually denoted by $\mu$ in the manuscript [II/13 b)]. See [39, fig. p. 32]. We prefer here the notation $\mathcal{O}$, used in [II/23 b)], as $H$ is now a standard notation in differential algebra for the product of initials and separants.

[13] We allow ourselves to write ${}^{\mathcal{P}}a_{i,j}$ instead of $a_{\mathcal{P},i,j}$ for more readability.

Proof. — Before considering Jacobi's arguments, we need first the following lemma[14]. Indices $\mathscr{P}$ will be omitted when there is no ambiguity.

**Lemma 95.** — *If $\mathscr{P}$ is a strongly quasi-regular component of $P$, $\mathscr{O}_{\mathscr{P}} \in \mathbf{N}$.*

Proof. — The variant of Kőnig's theorem (see above th. 17), that is also stated by R. Cohn [13], shows that, if $\mathscr{O}_{\mathscr{P}} = -\infty$, then one may find $a$ rows and $b$ columns in $A'_{\mathscr{P}}$, containing all the elements in $\mathbf{N}$, with $a+b<n$. So that $n-a>b$ equations in $\mathrm{d}P$ must depend on $b$ differentials $\mathrm{d}x_j$, which contradicts strong quasi-regularity. Hence the *lemma*. ∎ *Cohn, Richard*

a) First argument *Linearization.* — Jacobi first claims that one may reduce the problem to the case of a linear system. This, of course, cannot stand in all cases: we need the quasi-regularity hypothesis. We can assume that such assumptions were implicit in the physical situations that were considered by Jacobi: for proving that the order of $P$ corresponds to the dimension of the space of solutions of $\mathrm{d}P$, Jacobi used the fact that, if the general solution $y(\gamma, t)$ of $P$ depends on parameters $\gamma_i$, that may correspond to initial conditions, then $\partial y/\partial \gamma_i$ is a solution of the linearized system $\mathrm{d}P = 0$. *linearization* *general solution* *initial conditions* *parameters* *linearized system*

In our algebraic setting, the order of the differential field extension $\mathscr{G}_{\mathscr{P}}/\mathscr{F}$, defined by the prime component $\mathscr{P}$, is, by lemma 85 i), the dimension of the quotient module $\Omega_{\mathscr{G}_{\mathscr{P}}/\mathscr{F}} = \mathscr{M}_{\mathscr{P}}/(\mathrm{d}P)_{\mathscr{M}_{\mathscr{P}}}$ as a vector space.

b) Second argument *Stationary systems* — Jacobi claims then that one can assume the linearized system $\mathrm{d}P$ to have constant coefficients. This affirmation seems paradoxical, but it is dubious that he could have written it without a precise idea in mind. One may also notice that S. Cohn and Borchardt made no remark on that point. We have a single indication, the start of an argument that have been ruled out by Jacobi: *In explorando ordine systematis cum tantum altissima differentialia respiciuntur, in æquationibus differentialibus linearibus, ad quas proposita revocata sunt, supponere licet Coëfficientes esse constantes. Nam æquationibus 3) iteratis vicibus differentiatis, ut novæ obtinentur æquationes [...]* [15]. *stationary system* *constant coefficients* *Borchardt* *Cohn, Sigismundus*

We propose the following interpretation, expressed in the framework of the theory of standard bases of $\mathscr{D}$-module[16], which agrees with Jacobi's idea of looking at highest derivatives in the linearized system $\mathrm{d}P$. *standard basis*

**Definition 96.** — *We denote by $K$ the field $\mathscr{G}$ equipped with the derivation $\delta_0$, with $\delta_0 c = 0$ $\forall c \in K$, and by $\mathscr{M}_0$ the free $K[\delta_0]$-module generated by the $\mathrm{d}x_i$.*

*Let $m \in \mathscr{M}$ or $m \in \mathscr{M}_0$, $m = \sum_{v\in\Upsilon} c_v v$, where $\Upsilon$ denotes the set of derivatives of the $x_i$. We extend the definition 76 of* $\mathrm{ord}^{J_{\mathscr{P}}}$ *to $\mathscr{M}$ and $\mathscr{M}_0$ and define the* head *of $m$ to be $\kappa m := \sum_{\mathrm{ord}^{J_{\mathscr{P}}} v = \mathrm{ord}^{J_{\mathscr{P}}} m} c_v v$, the sum of terms of greatest order.*

---

[14] In the linear case, the result has been proved by Ritt [85].

[15] Looking for the order of the system, as one only considers the highest derivatives in the linear equations to which the proposed ones are reduced, one may assume [their] coefficients to be constants. Because, differentiating the equations 3) iterated times in order to obtain new equations [...]

[16] *Cf.* Castro-Jiménez [9, 10]. *Castro-Jiménez*

***Lemma* 97.** — *If $\nabla_{\mathscr{P}} \notin \mathscr{P}$, then $\mathscr{P}$ is a strongly quasi-regular component of $\{P\}$ and* $\kappa\,(\mathrm{d}P)_{\mathscr{M}} = \kappa\,(\kappa\ \mathrm{d}P)_{\mathscr{M}} \cong \kappa\,(\mathrm{d}P)_{\mathscr{M}_0} = \kappa\,(\kappa\ \mathrm{d}P)_{\mathscr{M}_0}$.

Proof. — If $\nabla_{\mathscr{P}} \notin \mathscr{P}$, the matrix $\mathrm{J} = (\partial P_i/\partial x_j^{(\alpha_{\mathscr{P},i}+\beta_{\mathscr{P},j})})$ is invertible in $\mathscr{G}$, so that the family $(\kappa\,\mathrm{d}P_i^{(k)})$, $1 \le i \le n$, $k \in \mathbf{N}$ is independent, as

$$\kappa\,\mathrm{d}P_i = \sum_{j=1}^{n} \partial P_i/\partial x_j^{(\alpha_{\mathscr{P},i}+\beta_{\mathscr{P},j})}\ \mathrm{d}x_j^{(\alpha_{\mathscr{P},i}+\beta_{\mathscr{P},j})}.$$

This implies that $\kappa(\mathrm{d}P) = (\kappa\,\mathrm{d}P)$ in $\mathscr{M}$ and $\mathscr{M}_0$. There is then a unique homomorphism of $\Theta$-monoïdeals, from $\kappa(\mathrm{d}P)_{\mathscr{M}}$ to $\kappa(\mathrm{d}P)_{\mathscr{M}_0}$, that send $\kappa\,\mathrm{d}P_i$ to $\kappa\,\mathrm{d}P_i$. Hence the *lemma*. ∎

So, if $\nabla_{\mathscr{P}} \notin \mathscr{P}$, it is indeed enough to prove the bound for some constant coefficient linear system. Assume that $\nabla_{\mathscr{P}} \in \mathscr{P}$ and that the $P_i$ are ordered with non increasing $\lambda_{\mathscr{P},i}$, let then $i_0$ be the smallest integer such that the first $i_0$ rows of $\mathrm{J}_{\mathscr{P}}$ are dependent. We may find some $c_i \in \mathscr{G}$, $1 \le i \le i_0$ with $c_{i_0} = 1$ such that $\sum_{i=1}^{i_0} c_i \delta_0^{(\lambda_i-\lambda_{i_0})} \kappa\ \mathrm{d}P_i = 0$. So, $(\mathrm{d}P)$ is generated by the family $\mathrm{d}P_1, \ldots, \mathrm{d}P_{i_0-1}$, $\sum_{i=1}^{i_0} c_i \delta^{(\lambda_i-\lambda_{i_0})}\,\mathrm{d}P_i$, $\mathrm{d}P_{i_0+1}, \ldots, \mathrm{d}P_n$. Jacobi's bound for this new linear system will be strictly smaller than $\mathscr{O}_{\mathscr{P}}$, as orders have decreased in equation $i$.

We may iterate the process until we find a free linear system $m_i$, $1 \le i \le n$, generating $(\mathrm{d}P)$, with a non vanishing system determinant $\nabla_m$. This must happen, for $\mathscr{O}_m \ge 0$, as $\mathscr{P}$ is strongly quasi-regular, and if $\mathscr{O}_m = 0$, then $(\mathrm{d}P) = (m) = (\kappa m)$, so that $\nabla_m$ cannot vanish, using again strong quasi-regularity.

c) Third argument *Determinant degree.* Assume that we have a linear system with constant coefficients $m_i = 0$, $1 \le i \le n$. We may represent it as a matrix of differential operators $M(\delta_0) = (m_{i,j})$ with $m_{i,j} = \sum_{p=0}^{a_{\mathscr{P},i,j}} c_{i,j,p}\delta_0^p$, where the $a_{\mathscr{P},i,j}$ are the elements of the order matrix $A_{\mathscr{P}}$ related to the component $\mathscr{P}$. The linearized system is equivalent to $M(\delta_0)(\mathrm{d}x_1, \ldots, \mathrm{d}x_n)^{\mathrm{t}} = 0$. The number of independent solutions of such a system is the number of roots $\xi$ of $|M(y)| = 0$, counted with multiplicities. Jacobi did only consider the simple case of pairwise distinct roots. *Chrystal* The general situation was later investigated by Chrystal [11].

This equation has degree at most $\mathscr{O}_m := \max_{\sigma\in S_n} \sum_{i=1}^{n} a_{i,\sigma(i)}$, and the coefficient of $\delta_0^{\mathscr{O}_m}$ is equal to $\nabla_m = |c_{i,j,a_{i,j}}|$, so that the order of the system is exactly $\mathscr{O}_{\mathscr{P}}$ iff $\nabla_{\mathscr{P}} \notin \mathscr{P}$. This concludes the proof of the theorem.

The next classical example shows how the bound can actually vary, according to the different components.

***Examples.*** — **98)** Let $P := (x')^2 - 4x \in \mathbf{Q}\{x\}$. This is a special case of ex. 87; it is well known that $P$ defines two components: the main regular component $\mathscr{P} := [P, x'' - 2]$ and the singular component $[x]$. We have $\mathscr{O}_{\mathscr{P}} = 1 = \mathscr{O}_P$, which is indeed the order of $\mathscr{P}$. The singular component $[x]$ is strongly quasi-regular and $\mathscr{O}_{[x]} = 0$, which is also again the order of $[x]$.

**99)** The system $x_1'' + x_2' = 0$, $x_1 = 0$ has a non vanishing $\nabla = -1$ and its order is indeed exactly $\mathcal{O} = 1$.

**100)** The system $x_1'' + x_2' + x_1 = 0$, $x_1' + x_2 = 0$ has a vanishing $\nabla$ and $\mathcal{O} = 2$. Indeed, $P_1 - P_2' = x_1$. The new system $P_2, x_1$ has a non vanishing $\nabla$ and a strictly smaller Jacobi bound $\mathcal{O} = 0$, that is the order of the system.

It is now easy to extend Jacobi's bound to underdetermined or (non strongly) quasi-regular systems. We need first to define the order of such a system, by analogy with the degree of an algebraic system, as done for diffieties with Sadik [79]. *underdetermined system* *Sadik* *diffieties*

**Definition 101.** — *Let $\mathcal{P}$ be a prime differential ideal of $\mathcal{F}\{x\}$ of differential dimension $d$. The order of $\mathcal{P}$ is the maximal order of quasi-regular components of differential dimension $0$ of the ideals $\mathcal{P} + [L_1, \ldots, L_m]$, where the $L_i := c_{i_0} + \sum_{j=1}^{n} c_{i,j} x_j$ are generic linear equations of order $0$, with coefficients in the differential transcendental extension $\mathcal{F}\langle c_{i,j} | (i,j) \in [1,d] \times [0,n] \rangle / \mathcal{F}$.* *order of a prime component*

*Jacobi's bound (generalization)*

**Corollary 102.** — *i) For any quasi-regular component $\mathcal{P}$ of $P_i \in \mathcal{F}\{x_1, \ldots, x_n\}$, $1 \le i \le s$, the order of $\mathcal{P}$ is at most $\mathcal{O}_P$.*

*ii) If $s \le n$, the order is equal to $\mathcal{O}_{\mathcal{P}} \in \mathbf{N}$ iff the matrix $\mathrm{J}_P$ (def. 76) has rank $s$ in the differential field extension $\mathcal{G}$ defined by $\mathcal{P}$.*

**Proof.** — i) We may reduce to the strongly quasi-regular case by using lem. 85 i). Indeed, we can choose $\hat{P} \subset P$ so that $\mathrm{d}\hat{P}$ is a maximal differentially independent family: $\mathrm{ord}\mathcal{P} \le \mathcal{O}_{\hat{P}} \le \mathcal{O}_P$. This includes the case when the differential codimension of $\mathcal{P}$ is strictly less than $s$ and we have $\sharp\hat{P} \le n$.

When $s < n$, we have seen that $\mathcal{O}_{\mathcal{P}}$ is obtained by completing matrix $A$ with $n-s$ rows of zeros that correspond to the orders of generic linear equations $L_i$, of order $0$. If $\mathcal{P}$ is a strongly quasi-regular component of $P$, one may find a generic system of $n-s$ equations $L_i$, of order $0$, such that $\mathcal{P} + [L]$ is a strongly quasi-regular component of $P, L$. So, the theorem, applied to the system $P, L$ implies that the order of $\mathcal{P}$ is bounded by $\mathcal{O}_{\mathcal{P}}$.

ii) For a generic system $L_i$, $s < i \le n$, $\nabla_{\mathcal{P}+L}$ does not vanish iff $\mathrm{J}_{\mathcal{P}}$ has full rank. So, using the theorem again, the order is equal to $\mathcal{O}$ iff $\mathrm{J}_{\mathcal{P}}$ has full rank. ∎

***Examples.*** — **103)** Let $P_1 := (x_1''' + x_2'')^3 x_2''' + x_1' + x_2$, $P_2 := x_1'' + x_2'$. It defines a single prime component $\mathcal{P}$, which is $[x_1' + x_2]$, of differential dimension 1, so that it is not strongly quasi-regular. However, $\mathrm{d}P = (\mathrm{d}x_1' + \mathrm{d}x_2)$ generates $\mathrm{d}\mathcal{P}$, so that it is quasi-regular. Its order is 1, that is indeed bounded by $\mathcal{O}_P = 5$.

**104)** Consider now the system defined by $P_1$ alone. It has a main regular prime component $\mathcal{P}_1 = [P_1] : Q^\infty$, with $Q := x_1''' + x_2''$, for which $P_1$ is a characteristic set for an orderly or a Jacobi ordering. For it, $\mathrm{J}_{P_1}$ is of full rank, so that it has order $3 = \mathcal{O}_{P_1}$. As $P_1'' = Q(S := (3QQ'' + 6Q'^2)x_2''' + 6QQ'x_2^{(4)} + Q^2 x_2^{(5)} + 1)$, $\mathcal{P}_2 := [P] : S^\infty = [x_1' + x_2]$ is another prime component, that is singular (see th. 139 iv) below) and no other component exists. We may notice that $\mathrm{J}_{\mathcal{P}_2}$ has full rank, so that the order is $1 = \mathcal{O}_{\mathcal{P}_2, P_1}$.

The order of $\mathscr{P}$ is strictly less than $\mathscr{O}_{\mathscr{P}}$ whenever there exists $\hat{P} \subset P$ with $\mathscr{O}_{\mathscr{P},\hat{P}} < \mathscr{O}_{\mathscr{P},P}$ and $\nabla_{\mathscr{P},\hat{P}} \notin \mathscr{P}$.

## 7 Shortest normal form reduction

*normal form* *characteristic set* We consider here one of Jacobi's results that may have the greatest consequences for improving the resolution of differential systems in most practical cases. Jacobi describes a method that, *generically*, *i.e.* when his system determinant $\nabla$ does not vanish, allows to compute a normal form or a characteristic set, using as few derivatives as possible of the system equations: indeed, it is then enough to differentiate $P_i$ up to order $\lambda_i$, where $\lambda$ is the minimal canon of $A_P$ and, *under stronger genericity hypotheses*, it is impossible to compute a normal form by differentiating one of the $P_i$ fewer times (see subsec. 7.3).

*canon* *truncated determinant* *Egerváry ordering* In fact, except for minimality, Jacobi's results stand for any canon. One may guess that Jacobi was aware of this fact although he did not state it explicitly. In [40, § 3 p. 58], he claims indeed that, if a normal form can be computed using the equations $P_1, \ldots, P_n$ and a minimal (for inclusion) set of derivatives of these equations, we can deduce from them a canon, and from this canon the minimal one. Such a property is not general, but is valid if the truncated determinant does not vanish and the normal form is associated to an *Egerváry ordering* defined by a non minimal canon $\mathfrak{l}$: this generalization of the Jacobi ordering is given in def. 107 below. Another implicit appearance of Egerváry orderings in Jacobi's manuscripts will be described in § 10.

*shortest reduction* *change of ordering; Pardi !* *Boulier; Lemaire* *Moreno Maza* The shortest reduction method could be suggested as a default strategy in computer algebra systems, when it is requested to compute a characteristic set, without specifying a precise ordering. It may also be used as a first step in methods using a change of orderings, such as *Pardi !*, designed by Boulier, Lemaire and Moreno Maza [5].

*Shaleninov* *Pryce* *DAE* Shaleninov [91] and Pryce [81] proposed strategies for the integration of implicit DAE that turn to be equivalent to Jacobi's shortest reduction. It appears that, in many practical situations, the system determinant $\nabla$ actually does not vanish, so that this method can be efficiently used.

*underdetermined system* Jacobi only considers the case when there are as many equations as variables. The generalization to underdetermined systems is easy but may be done in different ways. The rule of the game is to produce a bound $r_i$ on the order of successive derivatives of each equations $P_i$, requested to compute a characteristic set $\mathscr{A}$ of a given component $\mathscr{P}$. More precisely, one should be able to extract $\mathscr{A}$ of the characteristic set $\mathscr{B}$ of a prime component of the algebraic ideal $\sqrt{\left(P_i^{(k)}|1 \le i \le s,\ 0 \le k \le r_i\right)} : \nabla_P^\infty$, that turns to be equal to

$$\left(P_i^{(k)}|1 \le i \le s,\ 0 \le k \le r_i\right) : \nabla_P^\infty.$$

We must notice that bounding the order of derivatives of the $P_i$ involved in computing a characteristic set of a *singular* component—even strongly quasi-regular—is a challenge and that, using the minimal canon $\lambda_{\mathcal{P}}$ related to this component, instead of $\lambda_P$, does not help (see def. 93). One may think of example 98 and its singular component $[x]$. We need to differentiate $P$ one time to get a splitting $2x'(x''-2)$, so that the singular component may appear, but $\lambda_{\mathcal{P}} = \lambda = 0$. Such considerations are strongly related to the difficult Ritt problem, that is to decide inclusion of two components given by characteristic sets.

*singular component*

*Ritt problem*

*characteristic set*

## 7.1 Differential dimension zero

We first need some more definitions and preliminary results.

**Definition 105.** — *Let $\mathcal{I}$ and $\mathcal{J}$ be two ideals of some ring $A$, we denote by $\mathcal{I}:\mathcal{J}^\infty$ the ring $\{a\in\mathcal{I}|\forall b\in\mathcal{J}\ \exists q\in\mathbf{N}\ ab^q\in\mathcal{I}\}$.*

*Let $Q$ be any subset of $A$, then $\mathcal{I}:Q^\infty$ denotes $\mathcal{I}:(Q)^\infty$.*

The following proposition states a few easy properties; part iii) is a folkloric avatar of Rabinowitsch trick [82].

*Rabinowitsch*

**Proposition 106.** — *Under the hypotheses of def. 105, we have:*

*i) a) $\mathcal{I}:\mathcal{J}^\infty=\mathcal{I}:\sqrt{\mathcal{J}}$, b) if $\mathcal{I}_i\subset\mathcal{I}_{i+1}$, $\left(\bigcup_{i\in\mathbf{N}}\mathcal{I}_i\right):\mathcal{J}^\infty=\bigcup_{i\in\mathbf{N}}(\mathcal{I}_i:\mathcal{J}^\infty)$, c) $\sqrt{\mathcal{I}}:\mathcal{J}=\sqrt{\mathcal{I}:\mathcal{J}^\infty}$, d) if $\mathcal{J}_1\subset\mathcal{J}_2$, then $\mathcal{I}:\mathcal{J}_1^\infty=(\mathcal{I}:\mathcal{J}_2^\infty):\mathcal{J}_1^\infty\supset\mathcal{I}:\mathcal{J}_2^\infty$;*

*ii) If $\mathcal{J}=(Q_i|1\le i\le r)$, $\mathcal{I}:\mathcal{J}^\infty=\bigcap_{i=1}^r(\mathcal{I}:Q_i^\infty)$;*

*iii) If $A=k[x_1,\ldots,x_n]$ and $\mathcal{J}=(Q_i|1\le i\le r)$,*
*then $\mathcal{I}:\mathcal{J}^\infty=(\mathcal{I};\sum_{i=1}^r u_iQ_i-1)_{k[x;u]}\cap k[x]$.*

Proof. — i) a), b), c) $\subset$ and d) are immediate. i) c) $\supset$. If $a\in\sqrt{\mathcal{I}:\mathcal{J}^\infty}$, then $\forall b\in\mathcal{J}$ $\exists(p,q)\in\mathbf{N}^2$ $a^pb^q\in\mathcal{I}$ and so $(ab)^{\max(p,q)}\in\mathcal{I}$, $ab\in\sqrt{\mathcal{I}}$ and $a\in\sqrt{\mathcal{I}}:\mathcal{J}$.

ii) $\subset$ is immediate. $\supset$. Assume $P\in\bigcap_{i\in I}\mathcal{I}:Q_i^\infty$, and let $Q\in\mathcal{J}$. There exists a finite sum $Q=\sum_{i=1}^r M_iQ_i$. Assuming $PQ_i^{q_i}\in\mathcal{I}$, with $q=\max_{i=1}^r q_i$, we have $P(\sum_{i=1}^r M_iQ_i)^{rq}=PQ^{rq}\in\mathcal{I}$.

iii) $\subset$. Assume $PQ_i^{q_i}\in\mathcal{I}$, with $q:=\max_{i=1}^r q_i$, we have $P(\sum_{i=1}^r u_iQ_i)^{rq}\in\mathcal{I}$, so that $P\in(\mathcal{I};\sum_{i=1}^r u_iQ_i-1)$. $\supset$. Assume $P\in(\mathcal{I};\sum_{i=1}^r u_iQ_i-1)\cap k[x]$. Let $Q=\sum_{i=1}^r M_iQ_i\in\mathcal{J}$. As $P\in k[x]$, we may replace all the $u_i$ by $uM_i$ in the expression of $P$ as an element of $\mathcal{I}+[\sum_{i\in I}u_iQ_i-1]$. Then, $\sum_{i\in I}u_iQ_i-1$ becomes $uQ-1$. Let $d$ be the maximal degree of this expression in $u$. If $u$ is replaced by $1/Q$, we may get rid of denominators by multiplying by $Q^d$, so that $PQ^d\in\mathcal{I}$. ■

**Definition 107.** — *Let $A$ be the order matrix of a differential system $P$ and $\mu$, $\nu$ a cover for $A$, we say that an ordering $\prec$ on derivatives is* an Egerváry ordering *associated to the cover $\mu,\nu$ of $P$, if for any two derivatives $v_1$ and $v_2$,* $\mathrm{ord}^\nu v_1<\mathrm{ord}^\nu v_2$ *implies $v_1\prec v_2$.*

An Egerváry ordering $\prec$ is a Jacobi ordering if $\mu,\nu$ is equivalent to the Jacobi or canonical cover $\alpha,\beta$ (as defined in def. 22).

***Remark* 108.** — Considering general Egerváry orderings may prove to be useful even if they require a greater number of derivation. E.g. if one needs to compute a characteristic set for an ordering that is by chance of this kind, then it may be hopefully easier to compute it directly, rather than computing first a characteristic
change of ordering set for a Jacobi ordering and performing a change of ordering using a package such
Pardi ! as *Pardi !* However, there is no extra work to expose this more general case that will be used in section 8 for proving th. 143 and is implicit in Jacobi's computation of the resolvent in section 10; see rem. 171.

system determinant We have seen with prop. 81 that the system determinant $\nabla$ (see def. 76), that is equal to $|J^{\mu,\nu}|$, does not depend on the choice of the cover. This is obviously not the case of minors of $J^{\mu,\nu}$. We recall that, according to prop. 20 iv) and def. 76, if
cover $\mu, \nu$ is the cover associated to a canon $\ell$ and $\ell'$ the canon associated to $\mu$, $\nu$, then canon
$\min_{i=1}^n \mu_i = 0$ and $\min_{i=1}^n \ell'_i = 0$.

***Remark* 109.** — In the next theorem and the remaining of the subsection, we will consider a system of equations with coefficients in a differential ring of polynomials $\mathscr{R} := \mathscr{K}\{y_1, \ldots, y_p\}$, with fraction field $\mathscr{F}$, which includes the case of a differential field for $p = 0$.

Egerváry reduction **Theorem 110.** — *With the hypotheses of the last def. 107 and rem. 109, let $\ell$ be the canon of the order matrix $A_P$, associated to the cover $\mu, \nu$. We assume the $P_i$ to be ordered with non increasing $\ell_i$—and so non decreasing $\mu_i$.*

*We define the sets $F_k := \{i \in [1, n] | \mu_i \le k\}$. Let $\sigma \in S_n$ be a permutation and $\Delta_{\sigma,k}$ denote the minor of $J^{\mu,\nu}$ corresponding to the rows $F_k$ and the columns $\sigma(F_k)$. Let $\nabla_{\sigma,k} := \prod_{\kappa=\mu_1}^{k} \Delta_{\sigma,\kappa}$ and $\nabla_\sigma := \nabla_{\sigma,\mu_n}$.*

*Let $\mathrm{R}_k := \mathscr{R}[x_j^{(\kappa+\nu_j)} \mid 1 \le j \le n,\ -\nu_j \le \kappa \le k]$ and*

$$\mathscr{I}_{\sigma,k} := (P_i^{(\kappa)} | 1 \le i \le n,\ 0 \le \kappa \le k - \mu_i)_{\mathrm{R}_k} : \nabla_{\sigma,k}^\infty,$$

*and $\mathscr{I}_\sigma := \mathscr{I}_{\sigma,\mu_n}$.*

*i) With these assumptions, the ideals $\mathscr{I}_{\sigma,k}$, $k \in \mathbf{N}$, and $[P] : \nabla_\sigma^\infty$ are radical, with*

$$\mathscr{I}_{\sigma,k_1} \subset \mathscr{I}_{\sigma,k_2} \cap \mathrm{R}_{k_1},\ \textit{for}\ k_2 \ge k_1\ \textit{with equality if}\ k_1 \ge \mu_n \qquad (11)$$

*and*

$$\mathscr{I}_{\sigma,k} = [P] : \nabla_\sigma^\infty \cap \mathrm{R}_k, \quad \textit{for}\ k \ge \mu_n. \qquad (12)$$

Egerváry ordering *ii) Each prime components $\mathscr{P}$ of $[P] : \nabla_\sigma^\infty$ admits, for any Egerváry ordering $\prec$ associated to the cover $\mu, \nu$ and such that $x_{\sigma(i)}^{(\nu_{\sigma(i)})} \succ x_{\sigma(i')}^{(\nu_{\sigma(i')})}$ if $i < i'$, a characteristic set $\mathscr{A}$ with main derivatives $x_{\sigma(i)}^{(\mu_i+\nu_{\sigma(i)})}$, that may be extracted from a characteristic set of the prime component $\mathscr{P} \cap \mathrm{R}_{\mu_n}$ of $\mathscr{I}_{\sigma,\mu_n}$, for the same ordering.*

*iii) More precisely, the elements $A$ of $\mathscr{A}$ such that $\mathrm{ord}^\nu A \le \mu_i$ may be chosen in the characteristic set of the prime component $\mathscr{P} \cap \mathrm{R}_{\mu_i}$ of $\mathscr{I}_{\sigma,k}$, for the same ordering.*

Proof. — i) We first prove that, for any integer $k$, the ideal $\mathscr{I}_{\sigma,k}$ is radical. We use Lazard's lemma th. 89 iii) with respect to the main derivatives $x_j^{(\kappa)}$, for $j \in \sigma(F_k)$, $\nu_j + \mu_{\sigma^{-1}(j)} \le \kappa \le \nu_j + k$. We only have to remark that the determinant of the Jacobian matrix of the system defining $\mathcal{I}_{\sigma,k}$ with respect to these derivatives is equal to $\Delta_{\sigma,k}$ in order to conclude. The inclusion (11) is then a consequence of cor. 92 i) and the equality of cor. 92 ii).

As each $\mathcal{I}_{\sigma,k}$ is radical, $[P]:\nabla_\sigma^\infty$, which is equal to $\bigcup_{k\in\mathbf{N}} \mathcal{I}_{\sigma,k}$ by lem. 106 b), is radical too. Then equality (12) is a direct consequence of equality (11).

ii) We can extract from an algebraic characteristic set of $\mathcal{P}\cap \mathrm{R}_k$ for $\prec$ a maximal autoreduced subset with $n$ elements and main derivatives $x_{\sigma(i)}^{(\mu_i+\nu_{\sigma(i)})}$, $1 \le i \le n$.

iii) It is a straightforward consequence of i) inclusion (11). ∎

**Definition 111.** — *The ideals $\mathcal{I}_\sigma$, $\mathcal{I}_{\sigma,k}$ defined in the theorem and the quantities $\Delta_{\sigma,k}$ or $\nabla_\sigma$, will be respectively denoted by ${}_P^\nu\mathcal{I}_\sigma$, ${}_P^\nu\mathcal{I}_{\sigma,k}$, ${}_P^\nu\Delta_{\sigma,k}$ or ${}_P^\nu\nabla_\sigma$ when we need to avoid ambiguities. As the data of $\nu$ defines the cover $\mu,\nu$ up to equivalence and those objects are uniquely defined by $\nu$, "$\mu$" does not appear in the above notations.*

***Remark* 112.** — It is easier in these notations to refer to the permutation $\sigma$, however those objects obviously only depends on the sets $F_k$ and their images $\sigma(F_k)$.

Among the systems possessing such components are the characteristic sets of prime or "regular" ideals[17] (see Boulier *et al.* [6, def. 7]). To this regard, the "lifting of Lazard's lemma" in [6, th. 4] may be seen as a special case of the following corollary that gives a precise meaning in differential algebra to Jacobi's informal claim that one could compute a normal form by differentiating each equation $P_i$ at most $\lambda_i$ time, where $\lambda$ is the minimal canon.

*regular ideals; Boulier*

*lifting of Lazard's lemma*

*minimal canon*

**Corollary 113.** — *Let $P := \{P_i | 1 \le i \le n\} \subset \mathcal{R}\{x_1, \ldots, x_n\}$.*

*i) The ideal $\mathcal{I}_P := [P]:\nabla_P^\infty$ is radical.*

*ii) This ideal is such that, for any minimal cover $\mu,\nu$, $\mathcal{I}_P = \bigcap_{\sigma\in S_n} {}_P^\nu\mathcal{I}_\sigma$ and for any components $\mathcal{P}$ of $\mathcal{I}_P$, there exists an Egerváry ordering $\prec$ associated to $\mu,\nu$, such that $\mathcal{P}$ admits a characteristic set that is included in the char. set of a prime component of $(P_i^{(h)} | 1 \le i \le n,\ 0 \le h \le \ell_i):\nabla_P^\infty$, where $\ell$ is the canon associated to $\mu,\nu$.*

*iii) In particular, with $\mu = \alpha$ and $\nu = \beta$, where $\alpha$, $\beta$ is the Jacobi cover of $A_P$, associated to the minimal canon $\lambda$, any component $\mathcal{P}$ of $\mathcal{I}_P$ admits a char. set for a Jacobi ordering, that is included in the char. set of a prime component of $(P_i^{(h)} | 1 \le i \le n,\ 0 \le h \le \lambda_i):\nabla_P^\infty$, where $\lambda$ is the minimal canon.*

**Proof.** — ii) By prop. 81, $|\mathrm{J}^{\mu,\nu}| = \nabla_P$. Let $\{F_{k_q} | 1 \le q \le r\}$ be the set of all possible sets $F_k$, $0 \le k \le \mu_n$, with $F_{k_q} \subset F_{k_{q+1}}$. For any prime component $\mathcal{P}$ of $\sqrt{\mathcal{I}_P}$, if $\nabla_P \notin \mathcal{P}$, the matrix formed with the rows in $F_{\mu_k}$ of $\mathrm{J}^{\mu,\nu}$ admits a square submatrix of full rank. An easy recurrence shows that one may find sets $G_q$, with $\sharp G_q = \sharp F_{k_q}$ and $G_q \subset G_{q+1}$, such that the square submatrix of $\mathrm{J}^{\mu,\nu}$ formed of rows in $F_{k_q}$ and columns in $G_q$ has full rank. We can find a permutation $\sigma$ such that $\sigma(F_{k_q}) = G_q$. Then, ${}_P^{\nu}\nabla_\sigma \notin \mathscr{P}$ and a char. set of $\mathscr{P}$ for the order $\prec$ of part ii) of the th. can be extracted from a char. set of $\mathscr{I}_{\sigma,\mu_n}$ for the same ordering. This shows that $\sqrt{({}_P^{\nu}\nabla_\sigma | \sigma \in S_n)} \subset \sqrt{\nabla_P}$, so that the two radicals are equal. Then, by prop. 106 ii), we have $\bigcap_{\sigma\in S_n} {}_P^{\nu}\mathscr{I}_\sigma = [P]:\nabla_P^\infty$.

[17] Another notion of "regularity" appears here, *viz.* the possibility to test ideal membership by pseudo-reduction, using its characteristic set. But we will see in subsec. 7.4 that the corresponding "regular components" are also regular according to definition 84 iii).

The assertion iii) is just ii) in the case $\mu = \alpha$ and $\nu = \beta$.

i) By the theorem, ${}_P^{\beta}\mathscr{I}_\sigma$ is radical, so $\mathscr{I}_P = \bigcap_{\sigma\in S_n} {}_P^{\beta}\mathscr{I}_\sigma$ is also radical. ■

The next example shows that the ideals $[P]:\nabla_P^\infty$ or $[P]:\nabla_\sigma^\infty$ may fail to have a
*Boulier* characteristic representation (Boulier [6, def. 8]).

***Example* 114.** — Consider the algebraic system $x_1(x_1-1)=0$, $(x_2-2)(x_2^2-x_1^2)=0$. The ideal $[P]:\nabla^\infty = [P]:\nabla_\sigma^\infty$, for any $\sigma \in S_2$, is radical and is the intersection of the two prime components $[x_1, x_2-2]$ and $[x_1-1, P_2]$. But there is no ordering on the variables for which a characteristic representation could exist.

A simple consequence of the last corollary is that, for any component $\mathscr{P}$ of $[P]:{}_P^{\nu}\nabla_\sigma^\infty$ and any minimal cover $\mu', \nu'$, there exists a permutation $\sigma'$ such that $\mathscr{P}$ is a component of $[P]:{}_P^{\nu'}\nabla_{\sigma'}^\infty$. We can even be a little more precise with the following statement.

**Corollary 115.** — *i) We assume that $[P]:{}_P^{\nu}\nabla_\sigma^\infty$ is non trivial. Let $\{F_{k_q} | 1 \le q \le r\}$ be the set of all possible sets $F_k$, $0 \le k \le \mu_n$. The restriction $\mu^q$, $\nu^q$ of the minimal cover $\mu, \nu$ to the submatrix $A_P^q$ of $A_P$ defined by rows in $F_{k_q}$ and columns in $\sigma(F_{k_q})$ is a minimal cover, which is equivalent to the fact that the corresponding canon contains a maximal transversal set of maxima $a_{i_\kappa^q, j_\kappa^q} + \ell_{i_\kappa^q}$, $1 \le \kappa \le \sharp F_{k_q}$, in the canon $A_P^q + \ell$ associated to $\mu^q$, $\nu^q$.*

*ii) With the same hypotheses and notations, $[P]:{}_P^{\nu}\nabla_\sigma\infty = [P]:{}_P^{\nu'}\nabla_\sigma^\infty$ iff for all $1 \le q \le r$, the $a_{i_\kappa^q, j_\kappa^q} + \ell'_{i_\kappa^q}$ are maximal in the canon $A_P^q + \ell'$ associated to $\mu', \nu'$, which does not depend on the maximal transversal set chosen in $A_P^q + \ell$.*

*iii) There exists a unique minimal canon $\ell'$ such that $[P]:{}_P^{\nu}\nabla_\sigma^\infty = [P]:{}_P^{\nu'}\nabla_\sigma^\infty$, where $\mu', \nu'$ is the minimal cover associated to $\ell'$.*

*iv) For any component $\mathscr{P}$ of $[P]:{}_P^{\nu}\nabla_\sigma^\infty$, there exists a permutation $\sigma'$ such that:*
*a) the sum $\sum_{i=1}^n a_{i,\sigma'(i)}$ is maximal;*
*b) $\mathscr{P}$ is a component of $[P]:{}_P^{\nu'}\nabla_{\sigma'}^\infty$, for any minimal cover $\mu', \nu'$.*

*v) For any component $\mathscr{P}$ of $P$, there exists a char. set $\mathscr{A}$ that may be computed by differentiating $P_i$ at most $\lambda_i$ times, where $\lambda$ is the minimal canon of $A_P$, such that $\mathrm{ord}^{\mathrm{J}}\mathscr{A} \le \mathrm{ord}^{\mathrm{J}} P$ and for any minimal cover $\mu, \nu$, there exists an ordering $\prec$ compatible with $\mathrm{ord}^{\nu}$ such that $\mathscr{A}$ is a char. set of $\mathscr{P}$ for $\prec$.*

**Proof.** — i) As $\prod_{q=1}^r \Delta_{\sigma,k_q}$ divides $\nabla_\sigma$, $\Delta_{\sigma,k_q}$ must be non zero, which implies that the restriction of $\ell$ to $A_q$ is a canon. This is equivalent to the fact that it contains a maximal transversal family.

ii) Maximal transversal families of maxima for a canon correspond to maximal transversal sums and so their set does not depend on the chosen canon. Then, the products appearing in ${}_P^{\nu}\Delta_\sigma$ and ${}_P^{\nu'}\Delta_\sigma$ are the same, so that they are equal, iff the $a_{i_\kappa^q, j_\kappa^q} + \ell'_{i_\kappa^q}$ are maximal in the canon $A_P + \ell'$.

iii) The existence of a minimal canon such that the $a'_{i_k^q, j_k^q + \ell_{i_k^q}}$ are maximal is a direct consequence of rem. 65 and prop. 66. Then, $[P] : {}_P^{\nu}\nabla_\sigma^\infty = [P] : {}_P^{\nu'}\nabla_\sigma^\infty$ by ii).

iv) By cor. 113 ii), $\mathscr{P}$ is a component of $\mathscr{I}_P$. Working with Jacobi's cover, we first prove that there exists a permutation $\sigma'$ such that for all $1 \le i \le n$, $a_{i,\sigma'(i)} = \alpha_i + \beta_{\sigma'(i)}$ and the square submatrix $\mathrm{J}_q$ of J, corresponding to the rows in $[1, k]$ and the columns in $\sigma'([1, k])$, assuming the sequence $\alpha_i$ to be non decreasing, has a determinant that does no belong to $\mathscr{P}$.

This can be done by recurrence. The result is immediate for $n = 1$. If it stands for $n - 1$, let $\mathrm{J}_{\hat{j}}$ denote J deprived of row $n$ and column $j$, then

$$|\mathrm{J}| = \sum_{j=1}^{n} \pm \partial P_n / \partial x_j^{(\alpha_n + \beta_j)} |\mathrm{J}_{\hat{j}}|$$

, so that one may find $j_0$ such that both $\partial P_n / \partial x_{j_0}^{(\alpha_n + \beta_{j_0})}$ and $|\mathrm{J}_{\hat{j}_0}|$ do not belong to $\mathscr{P}$. This means that $a_{n,j_0} = \alpha_n + \beta_{j_0}$ and that, using the recursion hypothesis, one may find a suitable $\sigma' : [1, n-1] \mapsto [1, n] \setminus \{j_0\}$. It is then enough to define $\sigma'(n) := j_0$.

This proves a) and b) in the case of Jacobi's cover. Let $\mu, \nu$ be any minimal cover. Then, in $A_P^q$, we can take $a_{i_k^q, j_k^q} = a_{k,\sigma(k)}$, for all $1 \le k \le r$, so that by ii), $[P] : {}_P^{\alpha}\nabla_{\sigma'}^\infty = [P] : {}_P^{\nu'}\nabla_{\sigma'}^\infty$.

v) By iv), $\mathscr{P}$ admits a char. set $\mathscr{A}$ for some Jacobi ordering $\prec$, with main derivatives $x_{\sigma'(i_1)}^{(\alpha_{i_1} + \beta_{\sigma'(i_1)})}$. By ii) of the theorem, such an $\mathscr{A}$ may be chosen to be included in the algebraic char. set of a component $\mathscr{P} \cap \mathrm{R}_{\alpha_n}$ of $\mathscr{I}_{\sigma,\alpha_n}$, that involves $P_i$, $1 \le i \le n$ up to order $\lambda_i$. And for any cover $\mu, \nu$, $\mathscr{P}$ admits a char. set for an Egerváry ordering $\prec'$, with the same main derivatives. As $\mathscr{A}$ is autoreduced and reduces to 0 all elements of $\mathscr{P}$, it is a char. set of $\mathscr{P}$ for $\prec'$ too. ■

**Definition 116.** — *We call an* Egerváry reduction *the computation of a characteristic set of a prime component in* $[P] : \nabla_P^\infty$ *for an Egerváry ordering using derivatives of the* $P_i$ *of order bounded by the associated canon* $\ell$.

shortest reduction *In the case of a Jacobi ordering, it will be called a* Jacobi reduction *or* shortest reduction. Jacobi reduction

We have of course interest to choose the minimal canon $\ell$, as in rem. 65. The characteristic sets considered in iv) and v) can and should be computed using the Jacobi reduction for better efficiency.

***Example* 117.** — The linear system $P_i := x_i^{(a)} - x_{i+1} = 0$, $1 \le i < n$, $x_n^{(a)} = 0$ admits a single prime component, for which it is a char. set for all Jacobi orderings. The minimal canon is indeed $\lambda_i = 0$.

For an ordering such that $x_1 \ll x_i$, $1 < i \le n$, the only char. set is $\mathscr{A} = \{A_1 := x_1^{(na)}, A_i := x_i - x_1^{(ia)}, 1 < i \le n\}$. It may be computed by an Egerváry reduction with canon $\ell_i = (n - i)a$. We may check that $A_i = \sum_{\iota=1}^{i-1} P_\iota^{(i-\iota)a}$, $1 < i \le n$ and $A_1 = \sum_{\iota=1}^{n} P_\iota^{(n-\iota)a}$.

*change of ordering* *Egerváry reduction* We can now give a last immediate corollary, showing that changes of ordering between characteristic sets computed with Egerváry reductions are reversible. Cases of applications include situations described in sec. 8, *e.g.* in lemma 140.

**Corollary 118.** — *Let $\mathscr{A}$ and $\mathscr{B}$ be two characteristic sets of a prime ideal $\mathscr{P}$, such that $\mathscr{B}$ may be computed from $\mathscr{A}$ by an Egerváry reduction with canon $\ell$, associated to cover $\mu, \nu$, assuming non increasing $\ell_i$. Then, reciprocally, $\mathscr{A}$ may be computed from $\mathscr{B}$ by an Egerváry reduction, with canon $\ell$, using an order on $\mathscr{B}$ compatible with* $\mathrm{ord}^{\nu}$.

Proof. — Assume that the $A_i \in \mathscr{A}$ are ordered by non increasing $\ell_i$, where $\ell$ is the canon associated to the minimal cover $\mu, \nu$, that the main derivative of $A_i$ is $x_i^{(\mu_i+\nu_i)}$ and that $\mathscr{P} = [\mathscr{A}] : H_{\mathscr{A}}^{\infty}$ is equal to $[\mathscr{A}] : {}_{\mathscr{A}}^{\nu}\Delta_{\sigma}^{\infty}$. Then, the main derivative in $B_i$ is $x_{\sigma(i)}^{(\mu_i+\nu_{\sigma(i)})}$ and so $\ell$ is a canon for $\mathscr{B}$. Let $H$ (resp. $E$) be the set of leading derivatives for $\{A_i^{(k)} (\text{resp.} B_i^{(k)}) | 1 \le i \le n,\ \mu_i \le k \le \mu_n\}$, *i.e.* the $x_i^{(k+\nu_i)}$ (resp. $x_{\sigma(i)}^{(k+\nu_{\sigma(i)})}$). As,

$$\left|\partial\eta/\partial v;(\eta,v)\in H\times E\right| = \frac{{}_{\mathscr{A}}^{\nu}\nabla_{\sigma}}{\prod_{i=1}^{n} Sep_{A_i}^{\mu_n-\mu_i}}$$

and

$$\left|\partial\eta/\partial v;(\eta,v)\in H\times E\right| = \frac{{}_{\mathscr{B}}^{\nu}\nabla_{\sigma^{-1}}}{\prod_{i=1}^{n} Sep_{B_i}^{\mu_n-\mu_i}},$$

we have:

$$\left|\partial\eta/\partial v;(\eta,v)\in H\times E\right| \left|\partial v/\partial\eta;(v,\eta)\in E\times H\right| = \frac{{}_{\mathscr{A}}^{\nu}\nabla_{\sigma}}{\prod_{i=1}^{n} Sep_{A_i}^{\mu_n-\mu_i}} \frac{{}_{\mathscr{B}}^{\nu}\nabla_{\sigma^{-1}}}{\prod_{i=1}^{n} Sep_{B_i}^{\mu_n-\mu_i}} = 1,$$

meaning that ${}_{\mathscr{B}}^{\nu}\nabla_{\sigma^{-1}} \notin \mathscr{P}$, and $\mathscr{P}$ is a component of $[\mathscr{B}] : {}_{\mathscr{B}}^{\nu}\nabla_{\sigma^{-1}}^{\infty}$, for which $\mathscr{A}$ is a char. set that may be computed by some Egerváry reduction. ■

***Examples.*** — **119)** We consider the system $x_i' = F_i(x)$, where $F$ is a polynomial function such that any square $s \times s$ submatrix of the Jacobian matrix $\mathrm{J}(F)$, with $s \le n/2$ is invertible. It is quasi-linear and so defines a single prime component. The system is a characteristic set for the minimal canon $\lambda = 0$. It admits $2^n - 1$ classes of minimal covers of zeros and ones, defined by subsets $I \subsetneq [1, n]$ (the class of $[1, n]$ is equivalent to that of $\emptyset$), with $\nu_j = 1$ if $j \in I$ and $\mu_i = 1$ if $i \notin I$. Then $\ell_i = 1$ if $i \in I$ or else $0$. The original system corresponds to $I = \emptyset$ or $I = [1, n]$.

We can find new characteristic sets for the associated Egerváry orderings when $\sharp I \le n/2$. Let $J$ be such that $\sharp J = \sharp I$ and $I \cap J = \emptyset$. We may find a permutation $\sigma$ such that $\sigma$ defines a bijection from $I$ to $J$, leaving $[1, n] \setminus (I \cup J)$ unchanged and $\nabla_{\sigma}$ is equal to the determinant $|\mathrm{J}(F)_{I,J}|$ of the Jacobian matrix restricted to rows in $I$ and columns in $J$. For the corresponding characteristic set, that does not depend on $\sigma$ but only on $I$ and $J$ (see rem. 7.1), the leading derivatives are $x_j$ for $j \in J$, $x_j''$, for $j \in I$ and $x_j'$ for $j \notin I \cup J$.

When $\sharp I > n/2$, a minimal cover exists, but the associated char. sets may also be computed using smaller canons, as in cor. 115 iii). Let indeed in this case $\bar{I}$ (resp. $\bar{J}$) be $[1, n] \setminus I$ (resp. $[1, n] \setminus J$), with $J$ be such that $\sharp I = \sharp J$ and $\bar{I} \cap \bar{J} = \emptyset$. The char. set defined by $I$ and $J$ is the same as the char. set defined by $\bar{I}$ and $\bar{J}$. We notice that for such a system $P$ $\mathrm{J}_P$ is diagonal and that there is a single maximal sum, which is a special case of cor. 115.

*E.g.* for $n=2$, we may consider $x_1'+x_1+x_2$ and $x_2'+x_1+x_2$ and the characteristic sets $\{x_1 + x_2' + x_2, x_2'' + x_2' - x_2 - x_1\}$ for the canon $\ell = (0, 1)$ and $\{x_2 + x_1' + x_1, x_1'' + x_1' - x_1 - x_2\}$ for $\ell = (1, 0)$.

Such systems may be encountered in mechanics: with $p = n/2$, we may imagine that variables $x_1$, ..., $x_p$ are coordinates and $x_{p+1}$, ..., $x_{2p}$ momenta. Then, we can go from Hamilton's equations that are of order 1 in all the $x_i$ to explicit equations of order 2 in the coordinates $x_1$, ..., $x_p$, deduced from Euler-Lagrange equations, and also consider all intermediate cases. A direct computation of these intermediate cases of explicit normal forms can be performed using Routhians, that behave for some coordinates like Lagrangian depending on corresponding velocities and for others like Hamiltonians, depending on momenta. See Landau and Lifshitz [67, § 41 p. 133] for details.

*mechanics* *Euler-Lagrange equations* *Routhian* *Landau* *Lifshitz*

**120)** Consider the system $x_1^{(5)} + x_2'' + x_3''' = 0$, $x_2' = 0$, $x_1''' - x_3' = 0$. We have $\lambda = (0, 1, 2)$, $\alpha = (2, 1, 0)$ and $\beta = (3, 0, 1)$. There are two possible classes of characteristic sets that may be computed using the shortest reduction, *viz.* by differentiating the second equation 1 time and the second 2 times: $\mathscr{A} = \{x_1^{(5)}, x_2', x_3' - x_1'''\}$ and $\mathscr{B} = \{x_3''', x_2', x_1''' - x_3'\}$.

**121)** The system $x_1^{(5)} + x_2'' + x_3''' = 0$, $x_2' + x_3' = 0$, $x_3'' = 0$ admits a single class of char. set for Jacobi orderings, that may be computed using the shortest reduction, *viz.* by differentiating the second and the third equations only 1 time, according to the minimal canon $\lambda = (0, 1, 1)$. It is represented by $\mathscr{A} = \{x_1^{(5)}, x_2' + x_3', x_3''\}$. However, with the same derivatives, we may also compute the following characteristic set, that does not correspond to a Jacobi ordering, but to an Egerváry ordering: $\mathscr{B} = \{x_1^{(5)}, x_2'', x_3' + x_2'\}$ for the canon $\ell = (0, 2, 1)$. By chance, it may be computed with fewer derivatives than the bound given by th. 110.

**122)** We can illustrate cor. 118 with the following example. The system $P$ defined by $x_1'' = 0$, $x_2'' + x_1 = 0$, $x_3'' + x_2 = 0$ admits the 3 canons $\lambda = (0, 0, 0)$, which is the minimal one, $\ell = (0, 2, 0)$ and $\ell' = (0, 2, 4)$. This linear system admits a single component, for which there is a Jacobi reduction associated to $\lambda$ that produces $\mathscr{A} = P$. Using $\ell$, one may compute $\mathscr{B} := \{x_2^{(4)}, x_1 + x_2'', x_3'' + x_2\}$ and using $\ell'$ $\mathscr{C} := \{x_3^{(6)}, x_1 + x_3^{(4)}, x_2 + x_3''\}$. We see that we can go from $\mathscr{A}$ to $\mathscr{B}$ and back using $\ell$, from $\mathscr{A}$ to $\mathscr{C}$ and back using $\ell'$ and from $\mathscr{B}$ to $\mathscr{C}$ and back using $\ell' - \ell$.

*change of ordering*

The corollary is stated using non increasing $\ell_i$. We have kept here for simplicity the original order of the system.

***Remark* 123.** — Jacobi [40, end of § 3 p. 58] claims that the number of possible normal forms of a system, that one may find by the shortest reduction, is equal to the number of monomials in the truncated determinant, or equivalently to the number of transversal maximal sums in the order matrix. Example 121 has already produced some contradiction.

normal form — shortest reduction

Restricting ourselves to normal forms, or classes of characteristic sets, associated to Jacobi orderings does not solve the problem. It is easily seen that the number of normal forms may be smaller than $n!$ for systems such as $x+y+z=0$, $x'+y'+2z'=0$, $x''-y''+z''=0$, for which all 6 possible monomials appear in $\nabla$, but which has only 4 different normal forms: $x=-y-z, z'=0, y''=0$; $y=-x-z, z'=0, x''=0$; $z=-x-y, x'=-y', y''=0$ and $z=-x-y, y'=-x', x''=0$. Furthermore, a system such as $x+y=0$, $x'=0$ has only a single monomial in $\nabla$ but two normal forms: $x=y, y'=0$, $y=x, x'=0$, for the Jacobi ordering associated to the minimal canon $\lambda=(1,0)$.

normal form — characteristic set

The best expression I could find for the number of normal forms associated to a prime component $\mathscr{P}$ of $H_P : \nabla_P^\infty$ for all possible Jacobi orderings, is in fact the general result given by th. 110 ii), *i.e.* the number of permutations $\sigma$ such that ${}^{\beta}\nabla_P^{\sigma} \notin \mathscr{P}$. This is the only assertion of Jacobi that seems impossible to save in any reasonable way.

From a combinatorial standpoint, it means that the maximal number of such characteristic sets is, for a generic system—*i.e.* such that all subdeterminants of $J_P$ that are not identically zero, due to the place of non zero elements, are nonzero—the number of permutations $\sigma$ such that for all $1 \le k \le r$, the minor matrix contained in the $i := \sharp F_k$ first rows and the rows $\sigma(F_k)$ of $A_P + \lambda$ have a set of $i$ transversal maxima—but the set for $i$ need not be included in the one for $i+1$.

transversal maxima

## 7.2 Positive differential dimension

We consider now systems of $s$ equations $P_i \in \mathscr{F}\{x_1, \ldots, x_n\}$ and will extend results to underdetermined systems, *i.e.* when $s < n$. We propose two versions, the first one relies on the canon $A_P^{\boxminus}$ obtained by completing the rectangular order matrix $A_P$ by $n-s$ rows of 0, as done in rem. 3, 36 and 57. With this convention, the canon $\lambda^{\boxminus}$, Jacobi's cover $\alpha$, $\beta$ and $\nabla_P$ have already been defined (def. 76).

underdetermined system — rectangular matrix

A second way is to choose *a priori* a subset of $s$ main variables in order to reduce to the square case.

**Definition 124.** — *A system of $s$ polynomials $P$ being given, we associate to any subset $J \subset [1,n]$ of $s$ integers the system $P_J := \{\phi(P_i) | i \in [1,n]\}$, where $\phi$ is the canonical morphism $\mathscr{F}\{x_j | j \in [1,n]\} \mapsto \mathscr{R}\{x_j | j \in J\}$, with $\mathscr{R} := \mathscr{F}\{x_j | j \notin J\}$ (see rem. 109), the order matrix ${}_J^{\circ}A_P := A_{P_J}$, which is formed of the columns of $A_P$ with index $j \in J$. We define accordingly the associated minimal canons ${}_J^{\circ}\lambda_P$ and minimal cover ${}_J^{\circ}\alpha_P$, ${}_J^{\circ}\beta_P$, the Jacobi number ${}_J^{\circ}\mathcal{O}_P$, the matrix ${}_J^{\circ}J_P$ and the system determinant ${}_J^{\circ}\nabla_P$.*

*We denote by ${}^{\circ}\nabla_P$ the set $\{{}_J^{\circ}\nabla_P | J \subset [1,n],\ \sharp J = s\}$.*

To summarize those notations, $\nabla_P$ means that the cover corresponds to a bigger $n \times n$ order matrix obtained by adding $n - s$ rows of 0 and $^\circ\nabla_P$ is defined with covers of smaller $s \times s$ order matrices obtained by restricting to $s$ columns. It is easily seen that ${}_J\nabla_P$ is equal to ${}_J^\circ\nabla_P$ if ${}_J\mathcal{O}_P = \mathcal{O}_P$; if not, ${}_J\nabla_P$ is equal to 0.

**Theorem 125.** — *With the assumptions of the above definitions, we have the following assertions.*

*i) The ideal* $[P] : {}^\circ\nabla_P^\infty = \bigcap_{I \subset [1,n],\, \sharp I = s} [P] : {}_J^\circ\nabla_P^\infty$ *is radical.*

*ii) The ideal* $[P] : \nabla_P^\infty = \bigcap_{I \subset [1,n],\, \sharp I = s} [P] : {}_J\nabla_P^\infty$ *is radical.*

*iii) We have the inclusion* $[P] : \nabla_P^\infty \supset [P] : {}^\circ\nabla_P^\infty$ *and the prime components of* $[P] : \nabla_P^\infty$ *are those of* $[P] : {}^\circ\nabla_P^\infty$ *of differential order* $\mathcal{O}_P$.

**Proof.** — i) Considering orderings such that all derivatives of the $x_j$ for $j \notin J$ are smaller than those of the $x_j$ for $j \in J$, reduces the problem to the case of a system of $s$ equations in $s$ variables, with polynomial coefficients, for which we can use cor. 113 i) to conclude that $[P] : {}_J^\circ\nabla_P^\infty$ is radical. By prop. 106 ii), we have the decomposition $[P] : {}^\circ\nabla_P^\infty = \bigcap_{J \subset [1,n],\, \sharp J = s} [P] : {}_J^\circ\nabla_P^\infty$.

ii) The case of $[P] : \nabla_P^\infty$ is then a straightforward consequence of i), as ${}_J\nabla_P^\infty$ is equal to ${}_J^\circ\nabla_P^\infty$ or 0.

iii) In the case i), if the restriction to the first $s$ rows and columns of index $j \in J$ of the Jacobi cover $\alpha, \beta$ of $A_P$ is not a minimal cover of ${}_J^\circ A_P$, then ${}_J\nabla_P = 0$, $[P] : {}_J\nabla_P^\infty = [1]$ and there is no corresponding component in $[P] : \nabla_P^\infty$. When ${}_J^\circ\mathcal{O}_P = \mathcal{O}_P$, ${}_J^\circ\nabla_P^\infty$ is equal to ${}_J\nabla_P^\infty$, so that $[P] : {}_J^\circ\nabla_P^\infty = [P] : {}_J\nabla_P^\infty$. ■

***Examples.*** — **126)** Consider the system $P_1 := x_1' + x_2' + x_3 + x_4 = 0$, $P_2 := x_1 - x_2 + 2x_3' - x_4' = 0$. It is prime, so that we have a single component. The choices $J = \{1,3\}$, $\{1,4\}$, $\{2,3\}$ or $\{2,4\}$ provide ${}_J^\circ\mathcal{O}_P = \mathcal{O}_P = 2$, with ${}_J^\circ\lambda_P = (0,0)$ and $P$ is already a characteristic set. The choices $J = \{1,2\}$ and $\{3,4\}$ provide ${}_J^\circ\mathcal{O}_P = 1 < \mathcal{O}_P$ with ${}_J^\circ\lambda_P$ respectively equal to $(0,1)$ and $(1,0)$.

**127)** Let $P_1 := x_2'' + x_3''$ and $P_2 := x_1 + x_2' - x_3'$. The minimal canon is $\lambda = (0,0)$, that corresponds to $J = \{1,2\}$ or $J = \{1,3\}$, with ${}_J^\circ\mathcal{O}_P = 2$. The bound ${}_J^\circ\mathcal{O}_P = \mathcal{O}_P = 3$ is reached for $J = \{2,3\}$, with the canon ${}_J^\circ\lambda_P = \lambda_P^\boxminus = (0,1)$.

**128)** Let $P_1 := \left((x_1')^2 - x_3^2\right)(x_2^2 - x_4^2) = 0$ and $P_2 := x_3x_4$. We have 2 single prime components in $[P] : \nabla_P^\infty = [P] : {}_{\{1,4\}}\nabla_P^\infty$, which are of order 1 and equal to $[x_1' \pm x_3, x_4]$. The ideal $[P] : {}^\circ\nabla_P^\infty$ admits two more components with $J = \{2,3\}$ or $J = \{3,4\}$, which are of order 0 and equal to $[x_2 \pm x_4, x_3]$. There exist 2 other components in $\{P\}$, *viz.* $[x_1', x_3]$ and $[x_2, x_4]$.

In practical situations, we may try to optimize our choice, in order to lower the orders of the main variables or to reduce the order of derivations needed to compute a characteristic set. We know no polynomial time algorithm to compute $J$ with ${}_J^\circ\mathcal{O}_P$ minimal, but Jacobi's algorithm computes a minimal canon in an efficient way.

## 7.3 Minimality of the shortest reduction

normal form
shortest reduction

We need now to examine Jacobi's claim that such reductions in normal form, with the minimal canon $\lambda$ of $A_P$, are the "simplest", meaning that no reduction to normal form may be achieved by using derivatives of $P_i$ of order strictly smaller than $\lambda_i$ [40, § 3 p. 56]. We may first remark that, in the case $\nabla \in \mathcal{P}$, the computation of a characteristic set for $\mathcal{P}$ can require much more derivatives, as in the following examples.

***Examples.* — 129)** Let $P_1 := x_1^{(a)} + x_2^{(a)} + x_1$, $P_2 := x_1^{(a)} + x_2^{(a)} - x_2$. Then $\lambda = (0, 0)$, but the expression of the unique characteristic set $\{x_1, x_2\}$ requires derivatives of $P_1$ and $P_2$ up to order $a$.

**130)** Let $P_1 := \sum_{i=1}^n x_i$, $P_2 := P_1^{(a)} + x_2$, $P_i := P_1^{(a)} + x_{i-1}^{(a)} + x_i$, $2 < i \leq n$. The expression of the unique char. set $\{x_1, \ldots, x_n\}$ involves derivatives of $P_i$ up to order $(n-i)a$, while $\lambda = (a, 0, \ldots, 0)$. For $i = 1$, the order $(n-1)a$ is equal to $\mathcal{O}_P$.

The non vanishing of $\nabla$ is also requested for the minimality, as shown with the next example.

***Example* 131.** — Let $P_1 := x_3 x_1'' + x_2$, $P_2 := x_1$ and $P_3 := x_3$. We have $\lambda = (0, 2, 0)$, but the char. set $\{x_1, x_2, x_3\}$ can be obtained without any strict derivatives.

The following example shows that we need more genericity hypotheses.

***Example* 132.** — Let $P_i := \sum_{j=1}^i x_j^{(ia)}$. Then, $\nabla$ does not vanish, $\lambda_i = (n-i)a$, but we need only derivatives of $P_i$, $1 \leq i \leq n$, up to order $a$ to compute the unique characteristic set $\{x_i^{(ia)} | 1 \leq i \leq n\}$.

Trying to interpret Jacobi's claim, he certainly assumed genericity hypotheses that could have included the non vanishing of the system determinant, considering examples for which the assumption is natural, such as isoperimetric equations (see sec. 1.2). But, replacing the minimal canon $\lambda$ with its analog with respect to some component, *viz.* $\lambda_{\mathcal{P}}$ as defined in def. 93, we may secure a quite general result and even consider a wider notion of characteristic sets.

isoperimetric equations

It is indeed often interesting to consider alternative definitions that may be easier to compute. Among them, weak characteristic sets (see, *e.g.* Golubitsky, Kondratieva, Moreno Maza and Ovchinikov [31]).

weak characteristic set
Golubitsky
Kondratieva, Moreno Maza and Ovchinikov

weak autoreduced set

**Definition 133.** — *A* weak autoreduced set $\mathcal{A} = \{A_1, \ldots, A_r\}$ *is such that* $A_i = I_i v_i^{d_i} + R_i$*, where* $v_i$ *is the main derivative and* $d_i$ *its main degree, and the set of heads* $\{v_i^{d_i} | 1 \leq i \leq r\}$ *is autoreduced. A* weak characteristic set *of an ideal* $\mathcal{I}$ *is a weak autoreduced set* $\mathcal{A} \subset \mathcal{I}$ *of minimal rank.*

weak characteristic set

characteristic set

Obviously, characteristic sets are weak characteristic sets and, for the same ordering, a weak characteristic set has the same rank than any characteristic set. We may remark that the $P_i$ in the ex. 132 form a weak char. set if $x_{i+1} \succ x_i$. The following lemma completes lem. 85 i). Reciprocal implications do not stand.

**Lemma 134.** — *Let $S \in \mathscr{R}\{x\}$, $\ell \in \mathbf{N}^s$, $\mathscr{P}$ be a prime quasi-regular component of $[P]:{}^{\circ}\nabla_P^{\infty}$ and $\mathscr{G}/\mathscr{F}$ the associated differential field extension.* *linearization*

*If the prime ideal $\mathscr{P}$ admits a weak (resp. strong) differential characteristic set that may be extracted from an algebraic weak char. set of the ideal $(P_i^{(k)}|1 \le i \le n,\ 0 \le k \le \ell_i):S^{\infty}$, then the $\mathscr{G}[\delta]$-module $\mathrm{d}P$, as in def. 84, admits a standard basis (resp. reduced standard basis) that may be extracted from a basis of the vector space $\langle \mathrm{d}P_i^{(k)} \mid 1 \le i \le n,\ 0 \le k \le \ell_i\rangle$.* *standard basis* *reduced standard basis*

Proof. — All components of $[P]:{}^{\circ}\nabla_P^{\infty}$ are strongly quasi-regular by lem. 97. Then, by lem. 85 i), we may reduce to the associated module.

*Weak case.* — It is then straightforward that, if $\mathscr{A}$ is a weak char. set of $\mathscr{P}$, then $\mathrm{d}\mathscr{A}$ is a standard basis of $\mathrm{d}\mathscr{P}$ for the same ordering. *standard basis*

*Strong case.* — By i), we know that $\mathrm{d}\mathscr{A}$ is a standard basis. It may not be reduced, but if $A_i$ is reduced with respect to $A_{i'}$ with main derivative $x_j^{(a_j)}$, and $\mathrm{d}A_i$ is not reduced with respect to $\mathrm{d}A_{i'}$, then $\mathrm{ord}_{x_j} A_i = \mathrm{ord}_{x_j} A_{i'}$, so that $\mathrm{d}A_i$ may be reduced by $\mathrm{d}A_{i'}$ without any further differentiation. ▪ *reduced standard basis*

We can now state a general theorem. Assertion i) addresses first the weak case. Considering the strong case in ii), we may lighten a little the hypotheses.

**Theorem 135.** — *Let $P$ be a system of $s$ equations in the differential polynomial ring $\mathscr{R}\{x_1,\ldots,x_n\}$ (see rem. 109) and $\mathscr{P}$ a quasi-regular component of $\{P\}$ of differential dimension $n-s$, defining the field extension $\mathscr{G}/\mathscr{F}$. For any subset $J \subset [1,n]$ with $\sharp J = s$, we define the minimal canon ${}^{\mathscr{P}}_{J}\lambda$ of the square submatrix ${}^{\mathscr{P}}_{J}A$ of $A_{\mathscr{P}}$ defined by the columns of $J$ (see def. 93), its minimal cover ${}^{\mathscr{P}}_{J}\alpha$, ${}^{\mathscr{P}}_{J}\beta$ and ${}^{\mathscr{P}}_{J}\mathrm{J} := \left(\partial P_i/\partial x_j^{({}^{\mathscr{P}}_{J}\alpha_i + {}^{\mathscr{P}}_{J}\beta_j)}\right)$.*

*i) If all submatrices $(s-1)\times(s-1)$ of ${}^{\mathscr{P}}_{J}\mathrm{J}$ have maximal rank modulo $\mathscr{P}$, then:*
*a) the value of ${}^{\mathscr{P}}_{J}\lambda$ does not depend on $J$ and is equal to the minimal canon ${}^{\mathscr{P}}\lambda$ of the matrix $A_{\mathscr{P}}$;*
*b) for all $i_0 \in [1,s]$, the differential ideal $\mathscr{P}$ does not admit a weak characteristic set that may be extracted from the algebraic characteristic set of a component of $(P_i^{(k)}|1 \le i \le n,\ 0 \le k \le \ell_i):{}^{\circ}\nabla_P^{\infty}$, with $\ell_{i_0} < {}^{\mathscr{P}}\lambda_{i_0}$.*

*ii) If, for all $J \subset [1,n]$ with $\sharp J = s$, all submatrices $(s-1)\times(s-1)$ of ${}^{\mathscr{P}}_{J}\mathrm{J}$ containing all rows $i$ with ${}^{\mathscr{P}}_{J}\lambda_i = 0$ have maximal rank modulo $\mathscr{P}$, then, for all $i_0 \in [1,s]$, the differential ideal $\mathscr{P}$ does not admit a strong characteristic set that may be extracted from the characteristic set of a component of $(P_i^{(k)}|1 \le i \le n,\ 0 \le k \le \ell_i):{}^{\circ}\nabla_P^{\infty}$, with $\ell_{i_0} < {}^{\mathscr{P}}\lambda_{i_0}$.*

Proof. — i) a) Let $\Lambda_{\mathscr{P}} := \max_{i=1}^{s} {}^{\mathscr{P}}\lambda_i$, ${}^{\mathscr{P}}\alpha_i := \Lambda_{\mathscr{P}} - {}^{\mathscr{P}}\lambda_i$ and ${}^{\mathscr{P}}\beta_i := \max_{i=1}^{s} {}^{\mathscr{P}}a_{i,j} - {}^{\mathscr{P}}\alpha_i$. The statement is equivalent to saying that for all $J$ the restriction of ${}^{\mathscr{P}}\alpha$, ${}^{\mathscr{P}}\beta$ defines a minimal cover of ${}^{\mathscr{P}}_{J}A$. This is true for some subset $J_0$. If the result does not stand, the property is false for some subset $J$ with $\sharp J = s$. Then, there exists a subset $J_1$ with $J \neq J_1 \subset J_0 \cup J$ and $\sharp J_1 = s$, that satisfies the property and such that $\sharp(J_1 \cap J)$ is maximal. Let then be $j_0 \in J \setminus J_1$ and $i_0 \in [1, s]$ with ${}^{\mathscr{P}}a_{i_0,j_0} = {}^{\mathscr{P}}\alpha_{i_0} + {}^{\mathscr{P}}\beta_{j_0}$. Let be $j_1 \in J_1 \setminus J$, the submatrix of $\mathrm{J}_{\mathscr{P}}$ with columns in $J_1 \setminus \{j_1\}$ and rows in $[1, s] \setminus \{i_0\}$ has a non vanishing determinant, so that there exists $\sigma : [1, s] \setminus \{i_0\} \mapsto J_1 \setminus \{j_0\}$ such that ${}^{\mathscr{P}}a_{i,\sigma(i)} = {}^{\mathscr{P}}\alpha_i + {}^{\mathscr{P}}\beta_{\sigma(i)}$. We can now define a bijection $\sigma : [1, s] \mapsto J_2 := J_1 \setminus \{j_1\} \cup \{j_0\}$ by setting $\sigma(i_0) = j_0$, so that ${}^{\mathscr{P}}a_{i,\sigma(i)} = {}^{\mathscr{P}}\alpha_i + {}^{\mathscr{P}}\beta_{\sigma(i)}$ and ${}^{\mathscr{P}}\alpha$, ${}^{\mathscr{P}}\beta$ is a minimal cover of ${}^{\mathscr{P}}_{J_2}A$ with $\sharp(J_2 \cap J) = \sharp(J_1 \cap J) + 1$, which contradicts our minimality hypothesis.

*linearization* By lem. 85 iii) b), as we are in the quasi-regular case and the differential codimension is $n - s$, $\mathscr{P}$ is strongly quasi-regular and $\mathrm{d}\mathscr{A}$ must depend on all $\mathrm{d}P_i$, so i) b) and ii) stand when ${}^{\mathscr{P}}\lambda_{i_0} = 0$ and we can suppose ${}^{\mathscr{P}}\lambda_{i_0} > 0$.

b) To alleviate notations, we omit the notation $\mathscr{P}$ in the sequel, reducing as we may to the linear case by lem. 85 i). We assume that such a characteristic set $\mathscr{A}$ may be computed for some ordering $\prec$ and look for a contradiction. This implies by lem. 134 that $\mathrm{d}\mathscr{A}$ is in the vector space $E$ generated by the differentials $\mathrm{d}P_i^{(k)}$, for $0 \le k \le \ell_i$.

Without loss of generality, we may assume, up to a renumbering of variables, that the main derivative of $A_i \in \mathscr{A}$ is $x_i^{(\gamma_i)}$ and that $x_{j_1}^{(\beta_{j_1})} \prec x_{j_2}^{(\beta_{j_2})}$ if $1 \le j_2 < j_1 \le s$. We have $\mathrm{d}A_s = \sum_{i=1}^{s} \sum_{\kappa=0}^{k_i} c_{i,\kappa}\, \mathrm{d}P_i^{(\kappa)}$, with $c_{i,k_i} \neq 0$ if $k_i \ge 0$, and $c_{i,\kappa} = 0$ for $\kappa > k_i$ or $\kappa < 0$ by convention. Let $e := \max_{k_i \ge 0} k_i + \alpha_i$. Now, as $k_{i_0} < \lambda_{i_0}$, the sum $\sum_{k_i+\alpha_i=e} c_{i,k_i}\, \mathrm{d}P_i^{(k_i)}$ cannot contain both $i_0$ and $i_1$ with $\lambda_{i_1} = 0 < \lambda_{i_0}$. So, it contains strictly less than $s$ non zero terms. This implies, as the rectangular submatrix of $\mathrm{J}_{\mathscr{P}}$ corresponding to the columns $1 \le j < s$ and rows $i$ with $c_{i,e-\alpha_i} \neq 0$ may be completed to a square matrix of size $s - 1$ that has full rank by hypothesis, that it has full rank too. So, the greatest derivative of $x_s$ appearing in $\mathrm{d}A_s$ is of order at most $e + \beta_s$, but a derivative of some $x_i$, $1 \le i < s$ of order $e + \beta_i$, which is greater, must appear in it. The leading derivative of $\mathrm{d}A_s$ cannot be a derivative of $x_s$, a contradiction.

ii) As in i), we assume, up to a renumbering, that such a char. set $\mathscr{A}$ exists for some ordering $\prec$, the main derivative of $A_i \in \mathscr{A}$ being $x_i^{(\gamma_i)}$, with $x_{j_0}^{(\beta_{j_1})} \prec x_{j_2}^{(\beta_{j_2})}$ for $1 \le j_2 < j_1 \le s$. We look for a contradiction. This implies by lem. 134 ii) that there exists for the same ordering a reduced standard basis $G$ of $\mathrm{d}\mathscr{P}$, in the vector space $E$ generated in by the differentials $\mathrm{d}P_i^{(k)}$, $0 \le k \le \ell_i$.

There exists an element $g_0$ in $G$ with a leading derivative $\mathrm{d}x_{j_0}^{(\gamma_0)}$ with greatest Jacobi order, now defined with respect to the square order matrix ${}^{\mathscr{P}}_{[1,s]}A$. We have ${}^{\mathscr{P}}_{[1,s]}\mathrm{ord}^{\mathrm{J}} g_0 \ge L := \max_{i=1}^{s} \alpha_i$. Indeed, some element in $G$ must depend on $\mathrm{d}P_{i_1}$ with $\alpha_{i_1} = L$. So, for any element that depends on $P_i^{(k)}$, with $k + \alpha_i$ maximal, we have: $k + \alpha_i \ge L$ and the main derivatives in $P_i^{(k)}$ are of Jacobi order at least $L$. Such derivatives cannot depend on $\mathrm{d}P_{i_0}$, so that they cannot all cancel, as the submatrix of ${}^{\mathscr{P}}_{[1,s]}\mathrm{J}$ where row $i_0$ has been suppressed has full rank.

Let $g_0 = \sum_{i=1}^{s} \sum_{\kappa=0}^{k_i} c_{i,\kappa}\, \mathrm{d}P_i^{(\kappa)}$, with $c_{i,k_i} \neq 0$ if $k_i \geq 0$, and $c_{i,\kappa} = 0$ for $\kappa > k_i$ or $\kappa < 0$ by convention. Let $e := \max_{k_i \geq 0} k_i + \alpha_i$. The sum $\sum_{k_i+\alpha_i=e} c_{i,k_i}\, \mathrm{d}P_i^{(k_i)}$ cannot depend on $\mathrm{d}x_j^{(e+\beta_j)}$, with $j \neq j_0$, as $G$ is reduced. Again, as $P_{i_0}$ does not appear in the sum, the submatrix of ${}_{[1,s]}^{\mathscr{P}}\mathrm{J}$ restricted to the rows different from $i_0$, with $\lambda_{i_0} > 0$, and columns different from $j_1$ is of full rank, so that it is impossible, which concludes the proof. ■

Ex. 131 has already illustrated the necessity of using the order with respect to the considered component. It is indeed easily checked that for the unique component $\mathscr{P}$ defined by this system, $\lambda_{\mathscr{P}} = (0,0,0)$.

The next examples illustrate the differences between the criteria i) for weak characteristic sets and ii) for strong ones.

***Examples.* — 136)** We consider the system $P_1 := x_1' + x_2'' + x_3''' = 0$, $P_2 := x_1 - x_2' + x_3''$. The value of ${}_J\lambda = (0,1)$ does not depend on the set $J$ of main derivatives. The bound is met for all weak characteristic set, that are obtained for $J = \{1,2\}$ or $\{2,3\}$, using Jacobi's reduction, and the order for $P_2$ cannot be lowered.

**137)** For the system $P_1 := x_1' + x_2'' + x_3'' = 0$, $P_2 := x_1 + x_2 - x_3$, the value of ${}_J\lambda$ is $(0,1)$ for $J = \{1,2\}$ or $\{1,3\}$, and $(0,2)$ for $J = \{2,3\}$. The hypotheses of i) are not satisfied, but those of ii) are. The Jacobi reduction may be used for all values of $J$. No strong char. set can be computed without considering $P_2'$. But $P$ is already a weak char. set for any ordering $\prec$ such that $x_2'' \succ x_1'$ or $x_3'' \succ x_1'$ and $x_1 \succ x_2, x_3$.

Various phenomena may compensate each other, so that the bound is met, when none of the theorem hypotheses are satisfied, as illustrated by the next example.

***Example* 138.** — Consider the polynomials $P_1 := x_3 + P_2(P_2' + P_3'')$, $P_2 := x_1 + x_1'' + x_2'$ and $P_3 = x_1' + x_2$. It defines a single prime component, with a char. set $\{x_1, x_2, x_3\}$, that may be computed by considering derivatives of $P_i$ up to order $\lambda_i$, with $\lambda = (0,1,2)$. But this is a complete artifact. The vanishing of $P_2$ makes useless to consider derivatives $P_2'$ and $P_3''$ in the expression of $x_3$. So, the condition on the rank of submatrices of $\mathrm{J}_P$ in the th. 135 is not satisfied. But the vanishing of $\nabla$ makes the order lower, so that these derivatives are needed to compute $x_2$ and $x_3$.

The first difficulty in computing characteristic sets for non linear differential systems is to differentiate the equations. Using a classical representation of data, the sizes of the successive derivatives are exponential in the order and it is well possible to saturate the available memory before starting any algebraic elimination.

As Jacobi's method reduces, in the generic case, the number of requested differentiations to the minimum, it suggests many computational applications and

*Pardi ! — Diffalg — Boulier*

easily implemented improvements of existing softwares such as Pardi ! or Diffalg (see Boulier *et al.* [5, 6]).

## 7.4 Regularity

*lifting of Lazard's lemma — quasi-regular component — regular ideal — Boulier — Hubert*

As we have seen, the shortest reduction is a special case of the "lifting of Lazard's lemma", so that the radical ideal $\{P\}$ will split into components, some of them in $[P]:\nabla_P^\infty$ being of maximal order and given by a shortest reduction. Among other possible components, some may be considered "singular" in many ways. We will try to cast more light on the relations between "quasi-regularity", "regularity" (see resp. def. 84 i) and iii) and Boulier's "regular ideals". This attempt of a general definition has no claim to give a definitive answer but tries to underline the necessity to clarify, using differential algebra, a notion that remained, for the general case, in the informal setting inherited from the XVIII$^{e}$ and XIX$^{e}$ centuries[18]. Following the analysis of Ritt and the later work of Hubert [36], the main idea is that such components are envelopes, coming together with reasonable representations of a prime or radical differential ideal using a finite set of polynomials, and are most of the time unwanted. One needs a change of variables or higher order differentiations to get rid of them.

*characteristic set*

The next theorem tries to summarize the situation; Jacobi order plays a special role, mostly in differential dimension 0, but in higher dimension the flexibility of def. 84 iii) allows to fit better the structure of a given system, as shown in iv) for characteristic sets. The definition makes sense for any system but is mostly meaningful for systems with some normalization or reduction property, such as characteristic sets.

**Theorem 139.** — *Let $P$ be a differential system, $r := \mathrm{ord}^{v} P$ (def. 76) and $\mathcal{Q}$ (resp. $\mathcal{S}$) be the intersection of prime regular (resp. singular) components of $\{P\}$ with respect to $v$.*

*i) With $\mathrm{R}_r := \mathcal{F}[x_j^{(k)} | 1 \le j \le n,\ 0 \le k \le r + v_j]$, we have $\{P\} \cap \mathrm{R}_r = \mathcal{Q} \cap \mathrm{R}_r$.*

*ii) The component $\mathcal{P}$ of $\{P\}$ is singular with respect to $v$ iff $\mathcal{P} = \{P\} : S^\infty$ implies $\mathrm{ord}^{v} S > r$.*

*iii) The components of $[P]\ :\ {}^{\circ}\nabla_P^\infty$ are regular.*

*iv) If $\mathcal{A}$ is a characteristic set of a prime ideal $\mathcal{P}$ and all the $A_i \in \mathcal{A}$ are irreducible, then there exists $v$, such that $\mathrm{ord}^{v} \mathcal{A} = 0$. In particular, when $\mathrm{diff.\,dim}\,\mathcal{P} = 0$ is 0, one may choose $v = \beta$. For such a $v$, no regular component $\mathcal{P}_2$ of $\{(\mathcal{A}) : In_{\mathcal{A}}^\infty\}$ contains $Sep_{\mathcal{A}}$, $\mathcal{P}$ is the only regular component of $\{(\mathcal{A}) : In_{\mathcal{A}}^\infty\}$ with respect to $\mathrm{ord}^{v}$, and the intersection $\mathcal{S}$ of singular components of $\{\mathcal{A}\}$ contains $\{\mathcal{A}, Sep_{\mathcal{A}}\}$.*

**Proof.** — i) $\subset$. — Immediate. $\supset$ — Let $\{P\} \cap \mathrm{R}_r = \bigcap_{i=1}^{s} \mathcal{J}_i$, with $\mathcal{J}_i \not\subset \mathcal{J}_{i'}$ when $i \neq i'$. For any prime component $\mathcal{J}_{i_0}$ and $i \neq i_0$, we may find $S_i \in \mathcal{J}_i$ such that $S_i \notin \mathcal{J}_{i_0}$, so that $S := \prod_{i \neq i_0} S_i \notin \mathcal{J}_{i_0}$ and $\mathcal{J}_{i_0} = (\{P\} \cap \mathrm{R}_r) : S^\infty$. Let $\mathcal{A}$ be a diff. char. set included in some algebraic char. set of $\mathcal{J}_{i_0}$ for an Egerváry ordering with respect to $v$, it is a

*Houtain; Ritt — Lagrange; Laplace — Poisson — Hamburger*

[18] See Houtain [34] and the comments of Ritt on the approaches of Lagrange [86, p. 33] or Laplace, Poisson and Hamburger [86, p. 77]

char. set of $\mathscr{P}_{i_0} = [\mathscr{J}_{i_0}] : H_{\mathscr{A}}^\infty = \{P\} : (SH_{\mathscr{A}})^\infty$ ($\supset$ is immediate and $\subset$ stands because $\mathscr{A}$ reduces $P$ to 0). As $\mathrm{ord}^\nu(SH_{\mathscr{A}}) \le r$, $\mathscr{P}$ is regular with respect to $\nu$.

The assertion ii) is a straightforward consequence of the definition.

iii) We have $\mathrm{ord}^{\mathrm{J}}\,{}^\circ\nabla_P \le \mathrm{ord}^{\mathrm{J}} P$, hence the result.

iv) Assuming, as we may, that the main derivative of $A_i$ is $x_i$, we take $\nu_i := \mathrm{ord}_{x_i} A_i$ and for $s < j \le n$, $\nu_j \ge \max_{i=1}^s \mathrm{ord}_{x_j} A_i$. When $s = n$, this means that $\nu = \beta$. Then, $\mathrm{ord}^\nu A = 0$.

The classical theory of char. sets (see e.g. Ritt [86, Chap. IV] implies that $\{(\mathscr{A}) : In_{\mathscr{A}}^\infty\} \cap R_0 \subset \mathscr{P} \cap R_0 = (\mathscr{A}) : In_{\mathscr{A}}^\infty$. We first show that no regular component $\mathscr{P}_2$ can contain $Sep_{\mathscr{A}}$. Indeed, this would mean that some $S$ exists such that $SSep_{\mathscr{A}} \in \{(\mathscr{A}) : In_{\mathscr{A}}^\infty\}$, with $S \notin (\mathscr{A})$ and $\mathrm{ord}^\nu S \le 0$, so $SSep_{\mathscr{A}} \in (\mathscr{A})$, which is impossible. This means that any regular component $\mathscr{P}_2$ belongs to $[(\mathscr{A}) : In_{\mathscr{A}}^\infty] : Sep_{\mathscr{A}}^\infty = \mathscr{P}$ and is so equal to $\mathscr{P}$.

Assume now that a singular component $\mathscr{P}_3$ of $\{\mathscr{A}\}$ does not contain $Sep_{\mathscr{A}}$. There exists $S$ such that $\mathscr{P}_3 = \{\mathscr{A}\} : S^\infty$ (See Boulier *et al.* [6, Sec. 3.1 p. 85]). Let $R$ be the partial reduction of $S$ by $\mathscr{A}$, for some integer $p$, we have $Sep_{\mathscr{A}}^p S = R$ modulo $\{\mathscr{A}\}$, with $\mathrm{ord}^\nu R \le 0$. As $Sep_{\mathscr{A}} \notin \mathscr{P}_3$ and $S \notin \mathscr{P}_3$, we need have $R \notin \mathscr{P}_3$ so that $\mathscr{P}_3 = \{\mathscr{A}\} : R^\infty$ with $\mathrm{ord}^\nu R \le 0$, meaning that $\mathscr{P}_3$ is regular: a contradiction. So, $\mathscr{S} \supset \{\mathscr{A}, Sep_{\mathscr{A}}\}$. ■ Ex. 87 and 104 above illustrate this theorem.

According to def. 84 iii), regular components are those for which one may compute a characteristic set, starting from the initial system $P$ by allowing ourselves addition of two equations, multiplication of an equation by any differential polynomial and factorization, with the restriction that all the intermediate results cannot exceed the initial (Jacobi) order of the system.

The study of singular components requires to consider derivatives of a higher order, for which no bound is known. In practice, the Rosenfeld–Gröbner algorithm (Boulier [6]) provides a regular decomposition $\{P\} = \bigcup_{i=1}^p \mathscr{Q}_i$, where the $\mathscr{Q}_i$ are radical ideals, but we are unable to tell if some of their prime components is included in some other, which amounts to the Ritt problem. No factorization is required, but the consideration of the possible vanishing of the separants provides implicit square-free factorization. *Rosenfeld–Gröbner algorithm*

A prime component $\mathscr{P}_i$, containing some component $\mathscr{P}_j$ in such a decomposition, may correspond to an other kind of singularities: the singularities of the algebraic variety $\mathrm{V}(\mathscr{P}_j)$, and other kind of singularities typical of the differential case and considered by Johnson [52].

The components of $[P] : {}^\circ\nabla_P^\infty$ share many nice properties: they are regular, strongly quasi-regular and one may even compute for them a characteristic set for a Jacobi ordering that does not exceed the Jacobi order of $P$, using derivatives of $P$ that do not exceed this order. (See cor. 115 iv) and v).

# 8 The various normal forms of a system

normal form JACOBI considers in [II/13 b)] (*cf.* [39, p. 9–14]) and [II/23 a) f° 2217 *seq.*] (*cf.* [40, p. 37–43]) the various normal forms that a given system may possess. Systems in normal form include those of the type:

$$x_i^{(a_i)} = f_i(x), \quad 1 \le i \le n,$$

with $\mathrm{ord}_{x_j} f_i < a_j$, for all $1 \le i, j \le n$, which he calls *explicit normal form*. But Jacobi also includes implicit systems $A_i(x) = 0$, with $\mathrm{ord}_{x_j} A_i \le a_j$ and such that characteristic set $|\partial A_i / \partial x_j^{(a_j)}| \neq 0$, so that our characteristic sets enter in his definition of normal forms.

Jacobi claims that, if one cannot reduce a system to an equivalent one, with fewer equations than variables, that is, in our language, if the differential dimension is zero, then one can eliminate all dependent variables, except one, and get an equation of which the order is the order of the system. This is only generically true, and Jacobi was aware of it, for in [II/23 a), f° 2217, note] (*c.f.* [40, p. 37–43]), he introduces the order in some different way, claiming that the reduction to a single equation was sometimes impossible, *e.g.* if each equation $A_i$ depends only of $x_i$.

The order does not depend on the chosen explicit normal form and is equal order to $\sum_{i=1}^{n} a_i$. If we associate to the system a prime differential ideal $\mathscr{P}$, the order is the *algebraic* transcendence degree of the associated differential field extension $\mathscr{G}/\mathscr{F}$. At the time of Jacobi it was referred to as the *number of arbitrary constants appearing in a complete integration*, constants that could be, *e.g.*, initial conditions. Jacobi claims that, in the generic case, the orders of the leading derivatives in a normal form may be arbitrarily chosen, provided that their sum is equal to the order of the system.

Then, he considers ([39, p. 12]) systems possessing fewer possible normal forms, starting with the example of two equations in two variables. His results are expressed in the next lemma in the framework of characteristic sets. We provide an elementary proof that is missing in the manuscript. For the necessary condition expressed by part ii), we treat the general case that requires no extra work. The next section focuses on the elementary special case when the orders of one set of variables are increased and the orders of one set of variables are decreased of the same amounts.

## 8.1 Elementary transformations

In the sequel, we consider a prime ideal $\mathscr{P}$ of differential dimension 0 included in $\mathscr{F}\{x_1, \ldots, x_n\}$. If $\mathscr{A} = \{A_1, \ldots, A_n\}$ (resp. $\mathscr{B} = \{B_1, \ldots, B_n\}$) are char. sets, we assume the leading derivative of $A_i$ (resp. $B_i$) to be $x_i^{(a_i)}$ (resp. $x_i^{(b_i)}$)[19].

[19]The orders $a_i$ or $b_i$ correspond to $\alpha_i$ or $\beta_i$ in Jacobi's notations, in order to reserve these Greek letters to covers. This writing can coexist with the order matrix notation $a_{i,j} = \mathrm{ord}_{x_j} A_i$, so that $a_i$ is

***Lemma* 140.** — *Let $\mathscr{A}$ be a characteristic set of a prime differential ideal $\mathscr{P} \subset \mathscr{F}\{x_1, \ldots, x_n\}$, of differential codimension s.*

*i) If $s = n = 2$, $A_1 \succ A_2$ and $\mathrm{ord}_{x_1} A_2 = b_1 < a_1$, there exists a new characteristic set $\mathscr{B}$ of $\mathscr{P}$ with $\mathrm{ord}_{x_1} B_1 = b_1 := \mathrm{ord}_{x_1} A_2$, and $\mathrm{ord}_{x_2} B_2 = a_2 + a_1 - b_1$.*

*ii) Let $1 \le i_1 < i_2 \le s$, with $A_{i_1} \succ A_{i_2}$, there is no characteristic set $\mathscr{B}$ of $\mathscr{P}$ such that the main derivative of $B_i$ is $x_i^{(a_i)}$ for $i \in [1, s] \setminus \{i_1, i_2\}$ and the main derivatives of $B_{i_k}$ is $x_{i_k}^{(b_{i_k})}$, $k = 1, 2$, with $\mathrm{ord}_{x_{i_1}} B_{i_2} < b_{i_1} < a_{i_1}$.*

Proof. — i) As $x_1^{(a_1)} \succ x_2^{(a_2)}$, the derivative $x_1^{(a_1)}$ does not appear in $A_2$ and $b_1 < a_1$. Consider a new order $\prec_2$ such that $x_1^{(b_1)} \succ_2 x_2^{(a_2)}$ and $x_1^{(b_2-1)} \prec_2 x_2^{(a_2)}$. Then, we may take $B_1 = A_2$ and for $B_2$ the reduction of $A_1$ by $A_2$, which depends on $A_2^{(a_1-b_2)}$, so that it has order $a_2 + a_1 - b_2$ in $x_2$.

ii) If $\mathscr{B}$ is a char. set, then it must reduce $A_{i_2}$ to 0. As $A_i$ and $B_i$ have the same leading derivative, for $i \in I \setminus \{i_1, i_2\}$, the fact that $A_{i_2}$ is irreducible by $A_i$ implies it is irreducible by these $B_i$. Moreover $\mathrm{ord}_{x_{i_1}} A_{i_2} < \mathrm{ord}_{x_{i_1}} B_{i_1}$ and $\mathrm{ord}_{x_{i_2}} A_{i_2} = a_{i_2} < \mathrm{ord}_{x_2} B_{i_2} = a_{i_2} + a_{i_1} - \mathrm{ord}_{x_{i_1}} B_{i_1}$, so that $A_{i_2}$ is irreducible by $B_{i_k}$, $k = 1, 2$. A contradiction that concludes the proof. ■

***Remark* 141.** — In particular, if $\mathrm{ord}_{x_{i_1}} A_{i_2} = -\infty$, *viz.* if no derivative of $x_{i_1}$ appears in $A_{i_2}$, then there is no char. set where the order in $x_{i_1}$ is decreased, the order in $x_{i_2}$ increased, and the other leading derivatives unchanged.

Jacobi comes then ([39, p. 13]) to the case of an arbitrary number of variables and considers the problem of decreasing the order of $m$ variables $x_1, \ldots, x_m$ in a normal form normal form $\mathscr{A}$, when increasing the order of an equal number of variables $x_{m+1}, \ldots, x_{2m}$, the orders of the remaining variables staying unchanged, claiming that, if $\max_{i=m+1}^n \mathrm{ord}_{x_j} A_i = b_j$, $1 \le j \le m$, such a transformation can be achieved, provided that *normal form*

$$|\partial A_i / \partial x_j^{(b_j)}; 1 \le j \le m,\ m < i \le 2m| \neq 0.$$

He further claims that in the new normal form $\mathscr{B}$, the $b_{i+m}$, $1 \le i \le m$, may be chosen to be $a_{i+m} + a_i - b_i$.

The construction sketched by Jacobi, which relies on an elementary sequence of reductions, may sometimes fail and needs to be adapted. We simplify the presentation by using Egerváry reduction (th. 110) and will show that Jacobi's conclusion stands, under an extra hypothesis, up to a permutations of the variables $x_1, \ldots, x_m$. *Egerváry reduction*

***Example* 142.** — Consider the system $x_1''' = x_1''$, $x_2' = x_1''$, $x_3' = x_1 - x_2$, $x_4' = x_2$. Here $m = 2$ and $n = 2m$. We will decrease the order of $x_1$ and $x_2$ and increase the order of $x_3$ and $x_4$. In this case, Jacobi's construction is first to express $x_1$ and $x_2$, using the two last equations: $x_1 = x_3' + x_4'$, $x_2 = x_4'$. Then, we reduce the two first a convenient shorthand for $a_{i,i}$.

equations using these expressions, making derivatives of $x_1$ and $x_2$ disappear; we get: $x_3^{(4)} + x_4^{(4)} = x_3''' + x_4'''$ and $x_4'' = x_3''' + x_4'''$. Jacobi claims that we can now deduce explicit equations of order $4 = a_3 + a_1 - b_1$ in $x_3$ and $2 = a_4 + a_2 - b_2$ in $x_4$, from these two equations, without any further differentiations.

This is obviously not the case. One differential reduction is required to get a new normal form $x_3''' = 0$, $x_4''' = x_4''$, where the order in $x_3$ and $x_4$ is increased. Moreover, by lem. 140 ii), this system admits no normal form with order 4 in $x_3$ (resp. $x_4$) and 2 in $x_4$ (resp. $x_3$). One may check that cond. b) in th. 143 below is not satisfied.

The most general condition we could find to make possible the construction sketched by Jacobi is the existence of a canon, allowing an Egerváry reduction for a suitable permutation of variables, as in the next theorem. This reinterpretation is consistent but the text gives no explicit confirmation.

Of course, a *generic system* would satisfy Jacobi's conclusion, as a generic system admits normal forms with all possible $n$-tuples of orders, but it would not match such a simple process of reduction. In subsection 9.1, we will explicit those genericity conditions and return to the last example (see ex. 158).

**Theorem 143.** — *Let $\mathscr{A} = \{A_1, \ldots, A_n\}$ be a characteristic set of a prime differential ideal $\mathscr{P} \subset \mathscr{F}\{x_1, \ldots, x_n\}$. Let $[1, n] = I_1 \cup I_2 \cup I_3$ be a partition with $\sharp I_1 = \sharp I_2 = m$. For $j \in I_1$, let $b_j := \max_{i \in I_2 \cup I_3} a_{i,j}$. We assume that, for all $j \in I_1$, we have $b_j < a_j$.*

*Let then $\mu, \nu \in \mathbf{Z}^n$ be such that:*

*a) $\mu_i := a_i - b_i$ and $\nu_i := b_i$, for $i \in I_1$;*

*$\mu_i := 0$ and $\nu_i := a_i$, for $i \in I_2 \cup I_3$*

*b) $\forall (i, j) \in I_1^2$ $a_{i,j} \le \mu_i + \nu_j$, then*

*i) $\mu, \nu$ is a minimal cover for the order matrix $A_{\mathscr{A}}$;*

*ii) If $D := |\partial A_i / \partial x_j^{(\beta_j)}; (i, j) \in I_2 \times I_1| \notin \mathscr{P}$, there exists $\sigma \in S_n$ such that $\sigma^2 = \mathrm{Id}$, $\sigma(I_2) = I_1$, $\sigma(I_1) = I_2$, $\sigma_{|I_3} = \mathrm{Id}_{|I_3}$ and ${}_P^\nu\nabla_\sigma \notin \mathscr{P}$;*

*iii) There exists a char. set $\mathscr{B}$ of $\mathscr{P}$ such that $\mathrm{ord}_{x_i} B_i = b_i$ if $i \in I_1$, $\mathrm{ord}_{x_i} B_i = a_i + a_{\sigma_i} - b_{\sigma(i)}$ if $i \in I_2$, $\mathrm{ord}_{x_i} B_i = a_i$ if $i \in I_3$.*

**Proof.** — i) We first notice that $a_{i,i} = \mu_i + \nu_i$, for $i \in [1, n]$, which ensures minimality if $\mu, \nu$ is a cover. For all $j \in I_2 \cup I_3$, $\nu_j = a_j \ge a_{i,j}$, $i \in [1, n]$. We also have $\mu_i \ge 0$ for all $i \in [1, n]$, so $a_{i,j} \le \mu_i + \nu_j$ when $(i, j) \in [1, n] \times (I_2 \cup I_3)$. For $j \in I_1$, $\nu_j = b_j \ge a_{i,j}$, for all $i \in I_2 \cup I_3$, so $a_{i,j} \le \mu_i + \nu_j$ for $(i, j) \in (I_2 \cup I_3) \times I_1$. The condition b) imply that $a_{i,j} \le \mu_i + \nu_j$ in the remaining case $(i, j) \in I_1^2$, so that $\mu, \nu$ is a minimal cover.

ii) Let the leading term of $A_i$ be $In_i (x_i^{(a_i)})^{d_i}$, where $In_i$ denotes the initial of $A_i$, and $\tilde{D} := |\partial A_i / \partial x_j^{(\mu_i + \nu_j)}; (i, j) \in (I_2 \cup I_3) \times (I_1 \cup I_3)|$. Then, the term of greatest degree in the derivatives $x_i^{(a_i)}$, $i \in I_3$ that appears in $\tilde{D}$, is $\prod_{i \in I_3} (d_i In_i (x_i^{(a_i)})^{d_i - 1}) D$. As $In_i \notin \mathscr{P}$, $D \notin \mathscr{P}$ and $(x_i^{(a_i)})^{d_i - 1}$ is irreducible by $\mathscr{A}$, $\tilde{D} \notin \mathscr{P}$.

For any $I \subset I_1$, $J \subset I_2$ with $\sharp J = \sharp I$, we define $\bar{I} := I_1 \setminus I$, $\bar{J} := I_2 \setminus J$ $D_{I,J} := |\partial A_i / \partial x_j^{(\mu_i + \nu_j)}; (i, j) \in \bar{J} \times \bar{I}|$ and $\tilde{D}_{I,J} := |\partial A_i / \partial x_j^{(\mu_i + \nu_j)}; (i, j) \in (I \cup I_2 \cup I_3) \times (J \cup I_1 \cup I_3)|$. We prove in the same way that $\tilde{D}_{I,J} \notin \mathscr{P}$ if $D_{I,J} \notin \mathscr{P}$. Indeed,

the term of greatest degree appearing in $\tilde{D}_{I,J}$, in the derivatives $x_i^{(a_i)}$, $i \in I \cup J \cup I_3$, is $\prod_{i\in I\cup J\cup I_3}(d_i In_i(x_i^{(a_i)})^{d_i-1})D_{I,J}$, which is irreducible by $\mathscr{A}$.

We use the notations and conventions in th. 110 and its proof. We have $F_k = I_2 \cup I_3 \cup F_{1,k}$ where $F_{1,k} := \{i \in I_1 | a_i - b_i \le k\}$. Proceeding as in the proof of cor. 113 ii), let $\{F_{k_q} | 0 \le q \le r\}$ be the set of all possible sets $F_k$, $0 \le k \le \max_{i\in I_1} a_i - b_i$, with $F_{k_q} \subset F_{k_{q+1}}$. In the same way, $G_{q_k} := I_1 \cup I_3 \cup G_{2,q_k}$, with $G_{2,q_k} \subset I_2$. We may choose recursively $G_{2,k_q} \subset G_{2,k_{q+1}} \subset I_2$, with $\sharp G_{2,k_q} = \sharp F_{1,k_q}$, such that $D_{F_{k_q},G_{k_q}} \notin \mathscr{P}$.

We can then find a permutation $\sigma$ such that $\sigma_{|I_3} = \mathrm{Id}_{|I_3}$, $\sigma(I_1) = I_2$, $\sigma(I_2) = I_1$ and $\sigma(F_{1,k_q}) = G_{2,k_q}$, $1 \le q \le r$. Then, $\Delta_{\sigma,k_q} = \tilde{D}_{F_{k_q},G_{k_q}} \notin \mathscr{P}$ $1 \le q \le r$, so that ${}_P^{\nu}\nabla_\sigma \notin \mathscr{P}$.

iii) This is a consequence of th. 110. ▪

***Remarks. —* 144)** Condition b) always stands if the $a_i - b_i$ take a single value, in particular when $\sharp I_1 = \sharp I_2 = 1$.

**145)** By cor. 118, such a change of ordering is also reversible by Egerváry reduction.

The next example illustrates the necessity of choosing a suitable permutation.

***Example* 146. —** We consider the system $\{x_1''' = 0, x_2'' = 0, x_3' = x_2', x_4'' = x_1 + x_2' + x_5, x_5' = x_2\}$. We have $\mu = (3,1,0,0,0)$, $\nu = (0,1,1,2,1)$, that is indeed a minimal cover. We have $F_{1,1} = \{2\}$. As the third equation does not depend on $x_1$, we need to choose $G_{2,1} = \{3\}$, so that $D_{\bar{F}_1,\bar{G}_1} \ne 0$. With $\sigma(1) = 4$ and $\sigma(2) = 3$, we may compute the new system $x_1 = x_4'' - x_3' - x_5$, $x_2' = x_3'$, $x_3'' = 0$, $x_4^{(5)} = 0$, $x_5' = x_2$, using an Egerváry reduction.

What happens if $|\partial A_{m+i}/\partial x_j^{(\beta_j)}; 1 \le i,j \le m| \in \mathscr{P}$? Jacobi concludes with these words: *"Such questions require then a deeper investigation, that I will expose in some other occasion"*. One may guess that Jacobi was thinking of applying the general criterion that will be exposed in the next section 9.

It must be noticed that the requested transformation may sometimes be performed, even in the case when $\{A_{m+1}, \ldots, A_{2m}\}$, considered as a system in $x_1$, ..., $x_m$ alone, does not generate a differential ideal of dimension 0, as in the following example:

***Examples. —* 147)** Consider the explicit normal system of 4 equations in 4 variables

$$x_1''' = x_2'', x_2''' = 0, x_3' = x_1'', x_4 = x_1'.$$

If one wishes to decrease the orders of $x_1$ and $x_2$ and to increase that of $x_3$ and $x_4$, we cannot use the preceding lemma, nor any generalization of it, for the 2 last equations do not depend on $x_2$. However, we can achieve our goal with the following normal form:

$$x_1' = x_4, x_2'' = x_3'', x_3''' = 0, x_4' = x_3'.$$

The next example shows that one can decrease the order of 2 variables, when increasing the order of a single one.

**148)** Consider the system:

$$x_1'' = x_2, x_2' = 0, x_3 = x_1',$$

it is possible to decrease the order of $x_1$ and $x_2$ in the following normal form:

$$x_1' = x_3, x_2 = x_3', x_3'' = 0.$$

Testing the existence of a characteristic set $\mathscr{B}$ with leading derivatives $x_j^{(b_j)}$, for given $b_j$, some characteristic set $\mathscr{A}$ being known, will be the subject of the next section.

## 8.2 All possible shapes of normal forms

*normal form* We denote by orders $\mathscr{A}$ the $n$-tuple $(a_1, \ldots, a_n)$, where $a_j = \mathrm{ord}_{x_j} A_j$, assuming that the main derivative of $A_j$ is a derivative of $x_j$. Let $\mathscr{P}$ be a prime ideal, we denote by orders $\mathscr{P}$ the set $\{\text{orders } \mathscr{A} | \mathscr{A} \text{ a char. set of } \mathscr{P}\}$. We will conclude with some necessary condition on the possible values of orders $\mathscr{P}$, for a prime differential ideal $\mathscr{P}$ of diff. dim. 0 and order $e$. In two variables, it can take any value compatible with the conservation of the order of $\mathscr{P}$, as shown by the next proposition.

*normal form* **Proposition 149.** — *Let $e \in \mathbf{N}$ and $I \subset [0, e]$, there exists a prime differential ideal $\mathscr{P}$ such that* orders $\mathscr{P} = \{(a, e-a) | a \in I\}$.

Proof. — Let $a_1 < a_2 < \cdots < a_r$ be the elements of $I$. Define recursively

$$\begin{array}{rcl@{\qquad}rcl}
{}^1_1A_1 & = & x_1^{(a_1)}, & {}^1_1A_2 & = & 0\\
{}^{i+1}_{\;\;1}A_1 & = & {}^i_1A_2 + {}^i_1A_1^{(a_{i+1}-a_i)}, & {}^{i+1}_{\;\;1}A_2 & = & {}^i_1A_1,\\
{}^r_2A_2 & = & x_2^{(e-a_r)}, & {}^r_2A_1 & = & 0\\
{}^{i-1}_{\;\;2}A_2 & = & {}^i_2A_1 - {}^i_2A_2^{(a_i-a_{i-1})}, & {}^{i-1}_{\;\;2}A_1 & = & {}^i_2A_2.
\end{array}$$

Let then, for $j = 1, 2$ and $1 \le i \le r$, ${}^iA_j = {}^i_1A_j + {}^i_2A_j$ and $\mathscr{A}_i = \{{}^iA_1, {}^iA_2\}$. By construction, for $1 \le i \le r$, ${}^{i+1}A_2 = {}^iA_1$, ${}^{i+1}A_1 = {}^iA_2 + {}^iA_1^{(a_{i+1}-a_i)}$, with $\mathrm{ord}_{x_1} {}^iA_1 = a_i$, $\mathrm{ord}_{x_1} {}^iA_2 = a_{i-1}$, $\mathrm{ord}_{x_2} {}^iA_2 = e - a_i$ $\mathrm{ord}_{x_2} {}^iA_1 = e - a_{i+1}$, so the $\mathscr{A}_i$ are characteristic sets of the same prime diff. ideal $\mathscr{P}$ and by lemma 140, orders $\mathscr{P} = \{(a, e-a) | a \in I\}$.
■

For a greater number of variables, the situation is more complicated. If one tries to visualize the set of possible characteristic sets for a given system in 3 variables, it is convenient to use triangular coordinates, as the sum of the 3 maximal orders in the 3 variables is constant.

***Examples.*** — **150)** The linear prime ideal that describes the 3 functions $x_i = \int C\cos(t + 2i\pi/3)$, where $C$ is an arbitrary constant, admits 6 normal forms, of the 2 following shapes, with $1 \leq i \leq 3$: *normal form*

$$\begin{aligned} x_i''' &= -x_i';\\ x_{i+1[3]}' &= (-x_i' + \sqrt{3}x_i'')/2;\\ x_{i+2[3]}' &= (-x_i' - \sqrt{3}x_i'')/2 \end{aligned} \qquad \text{and} \qquad \begin{aligned} x_i' &= -x_{i+1[3]}' - x_{i+2[3]}';\\ x_{i+1[3]}'' &= (2x_{i+2[3]}' + x_{i+1[3]}')/\sqrt{3};\\ x_{i+2[3]}'' &= -(2x_{i+1[3]}' + x_{i+2[3]}')/\sqrt{3}. \end{aligned}$$

**151)** The system $x_1' - x_2' = 0$, $x_2''' = 0$, $x_3' - x_2' = 0$, admits only 3 normal forms: *normal form*

$$\begin{aligned} x_i''' &= 0;\\ x_j' &= x_i';\\ x_k' &= x_i', \end{aligned} \qquad \text{with } 1 \leq i, j, k \leq 3,\ i \neq j,\ i \neq k \text{ and } j < k.$$

The two examples may be illustrated by such drawings, where the points corresponding to existing normal forms are surrounded by a loop.

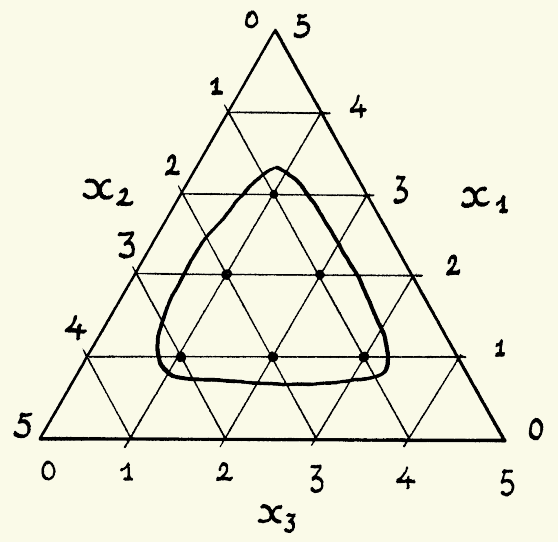

Example 150

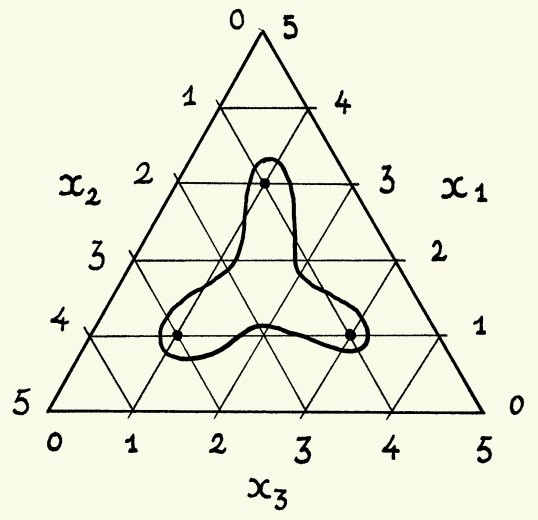

Example 151

Those drawings look very much like these ones, that appear on the margin of manuscript [II 13 b), f° 2206a]. 

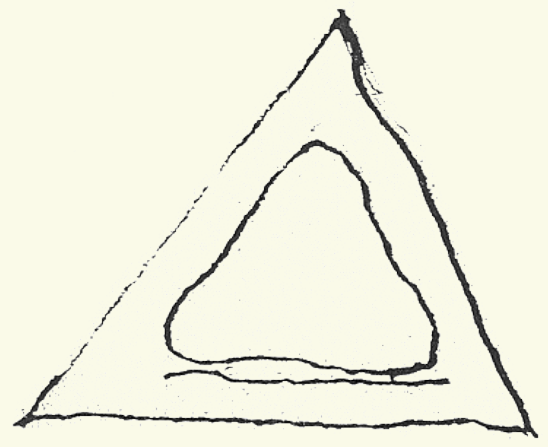

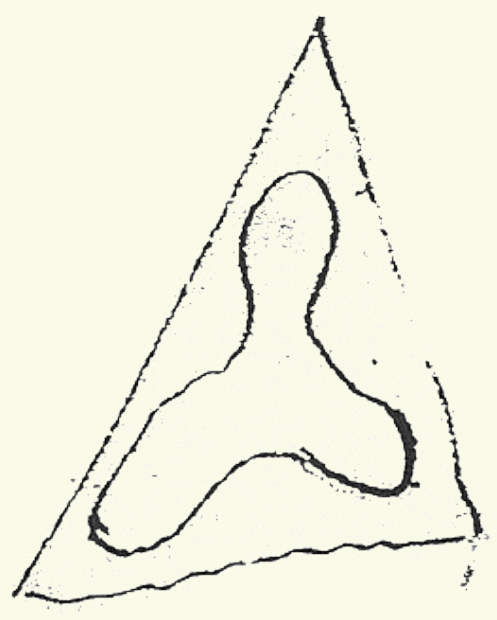

The proposition 152 below shows that, with more that 2 variables, the set orders $\mathscr{P}$ cannot be arbitrary.

**Proposition 152.** — *Let $\mathscr{P} \subset \mathscr{F}\{x_1, x_2, x_3\}$ be a prime differential ideal. Assume that there exist $a_1 > b > c \geq 0$ such that:* *normal form*

*a)* orders *$\mathscr{P}$ contains $(a_1, a_2, a_3)$, $(b, a_2 + a_1 - b, a_3)$ and $(c, a_2, a_3 + a_1 - c)$ and no element $(d_1, d_2, a_3)$ with $a_1 > d_1 > b$ or $(d_1, a_2, d_3)$ with $a_1 > d_1 > c$);*

*b)* orders $\mathscr{P}$ *does not contain any* $(d_1, d_2, a_3)$ *with* $b > d_1 \geq c$*;*

*c)* orders $\mathscr{P}$ *does not contain any* $(a_1, e_2, e_3)$ *with* $a_2 > e_2$ *nor any* $(f_1, f_2, a_3)$ *with* $a_2 > f_2$*.*

*Then,* orders $\mathscr{P}$ *contains* $(c, a_2 + a_1 - b, a_3 + b - c)$*.*

Proof. — We may assume the system to be linear, by lem. 85; which means that we can freely assume $a_{i,j} < a_j$, $i \neq j$. With $\sharp I_1 = 1$, assumption i) b) of th. 143 always stands (rem. 144) and we will always have $D \notin \mathscr{P}$.

So, if $a_{3,1} > a_{2,1}$, th. 143 implies the existence of a char. set with order triplet $(a_{3,1}, a_2, a_3 + a_1 - a_{3,1})$ and and by lem. 140 there is not char. set with orders $(g_1, a_2, g_3)$ with $a_1 > g_1 > a_{3,1}$, so that $a_{3,1} = c$. In the same way, there is no char. set with orders $(h_1, h_2, a_3)$ with $a_1 > h_1 > a_{2,1}$; a contradiction, as $b > c > a_{2,1}$. So $a_{2,1} \geq a_{3,1}$ and there exists a char. set $\mathscr{B}$ with orders $(a_{2,1}, a_2 + a_1 - a_{2,1}, a_3)$ and no char. set with orders $(h_1, h_2, c)$, with $a_1 > h_1 > a_{2,1}$. This implies that $b = a_{2,1}$.

If $\min(a_{1,2}, a_{3,2}) = a_{1,2} > -\infty$ (resp. $= a_{3,2} > -\infty$), then by th. 143, there exists a char. set with orders triplet $(a_1 + a_2 - a_{1,2}, a_{1,2}, a_3)$ (resp. $(a_1, a_{3,2}, a_3 + a_2 - a_{3,2})$). So, condition c) above implies that $a_{1,2} = a_{3,2} = -\infty$.

We can choose $B_1 = A_2$ and $B_3 = A_3$. Proceeding as above, condition b) means, using th. 143 and lem. 140 ii), that $\mathrm{ord}_{x_1} B_2 < c \leq ord_{x_1} B_3$. So, applying again th. 143, we get a new normal form $\mathscr{C}$ with order triplet $(a_{3,1}, a_2 + a_1 - b, a_3 + b - a_{3,1})$.

We can choose $C_1 = A_3$ and $C_3$ to be the reduction of $B_1 = A_2$ by $A_3$, which depends only on $A_2$, $A_3$ and derivatives of $A_3$ up to order $b - c$. Indeed, as $A_3$ does not depend on $x_2$, no reduction of its derivatives by $B_2$ is required.

We have then $\mathrm{ord}_{x_2} C_3 = \mathrm{ord}_{x_2} A_2 = a_2$ and $\mathrm{ord}_{x_2} C_1 = a_{3,2} = -\infty$. We can apply another time th. 143, which implies the existence of a normal form with order triplet $(a_{3,1}, a_2, a_3 + a_1 - a_{3,1})$. As already noticed, $a_{3,1} \geq c$, so that we need have $a_{3,1} = c$ to satisfy a) and orders $\mathscr{C} = (c, a_2 + a_1 - b, a_3 + b - c)$, which achieves the proof. ∎

**153)** Let $\mathscr{A} := \{x_1^{(4)}, x_2' + x_1'', x_3' + x_1\}$. This system admits also the char. set $\mathscr{B} := \{x_1'' + x_2', x_2''', A_3\}$. It also admits the char. set $\mathscr{D} := \{x_1 + x_3', x_2' - x_3''', x_3^{(5)}\}$ and by lem. 140 ii), cond. a) is satisfied. Moreover, as only strict derivatives of $x_2$ and $x_3$ appear in the system, no char. set exists with a leading derivative in $x_2$ or $x_3$ of order less than 1, so condition c) is satisfied. As $B_2$ does not depend on $x_1$, cond. b) stands too. The proposition applies and indeed the char. set $\mathscr{C} := \{A_3, x_2''', x_3''' - x_2'\}$ with a triplet of orders $(0, 3, 3)$ may be exhibited.

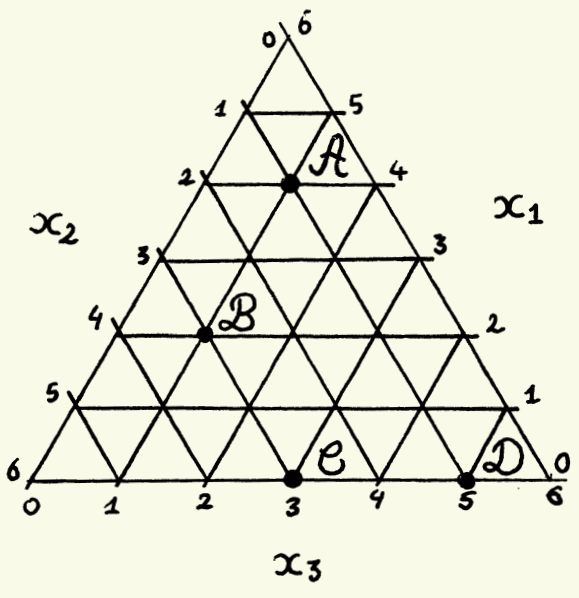

The condition b) and c) in the last theorem are technical and mostly used to fit the hypotheses of th. 143, that are not all necessary. I could build no example showing that one could not dispense with them. Designing a complete set of conditions that may characterize all possible sets of orders $n$-tuples remains an open question, for which th. 156 in the next section offers wider possibilities... *Tam quæstiones altioris indaginis poscuntur.*

# 9 Changes of orderings

CHANGES of orderings on monomials have been considered in the computer algebra literature for standard bases computations (FGLM [24]), or on derivatives (Pardi ! [5]) for characteristic set computations. It may be noticed that the main theoretical works of the XX$^{\text{th}}$ century often restrict to particular orderings, Janet orderings (Janet [48, p. 102][20]), elimination orderings (Ritt [86, p. 3]), but for many applications, one needs to use specific orderings, *e.g.* testing identifiability or observability in control theory requires to eliminate a precise set of indeterminates (see Matera and Sedoglavic [90, 73]).

*FGLM* *Pardi !* *Janet* *Ritt* *Matera and Sedoglavic*

In [II/23 a) f$^{\text{o}}$ 2217–2220] [40, p. 36–43], Jacobi considers, in full generality, the problem of computing a normal form of an ordinary differential system, some normal form being known for a different ordering. The method he gives is quite similar to the tools of contemporary literature and he provides moreover sharp bounds on the requested number of derivations, that may be used to improve the efficiency of many algorithms.

*normal form*

## 9.1 Necessary and sufficient conditions for the existence of a normal form with given orders

Considering a system in explicit normal form, such as

*normal form*

$$x_i^{(e_i)} = F_i(x), \quad 1 \le i \le n, \tag{13}$$

the problem is to compute a new normal form of the shape

$$x_i^{(f_i)} = G_i(x), \quad 1 \le i \le n. \tag{14}$$

In a first step Jacobi divides the indeterminates in three sets. For $i \in I_1$, $f_i > e_i$; for $i \in I_2$, $f_i < e_i$ and for $i \in I_3$, $f_i = e_i$.

Using the derivation operator

$$\delta := \sum_{i=1}^{n} \left( F_i(x) \frac{\partial}{\partial x_i^{(e_i-1)}} + \sum_{k=0}^{e_i-2} x_i^{(k+1)} \frac{\partial}{\partial x_i^{(k)}} \right),$$

Jacobi claims that it is possible to compute the new normal form (14) using the first one (13), completed with the equations $x_i^{(e_i+k)} = \delta^k F_i(x)$, $i \in I_1$, $1 \le k \le f_i - e_i$ iff

$$\left| \frac{\partial\, \delta^k F_i}{\partial x_j^{(a)}} \mid i \in I_1,\ 0 \le k < f_i - e_i; \quad j \in I_2,\ f_j \le a < e_j \right| \ne 0.$$

[20]These orders are already in Riquier [83, p. 195].

We need here a new definition, in order to translate in differential algebra the effect of derivation $\delta$, *i.e.* reducing higher order derivatives as soon as they appear.

**Definition 154.** — *Let $\mathscr{A} = \{A_i | 1 \le i \le n\}$ be a characteristic set of a prime differential ideal in $\mathscr{F}\{x\}$ for an ordering $\prec$. We assume that the $A_i$ appear in an increasing ordering and that a reduction process using $\mathscr{A}$ has been chosen. We denote by $\widetilde{A}'_i$ the reduction of $A'_i$ by the $\widetilde{A}'_j$, $j < i$ and the $A_j$, $1 \le j \le n$, and recursively denote by $\widetilde{A}_i^{(k+1)}$ the reduction of $(\widetilde{A}_i^{(k)})'$ by the $\widetilde{A}'_j$, $1 \le j \le n$, and then by the $A_j$, $1 \le j \le n$. By convention, $\widetilde{A}_i^{(0)} := A_i$.*

***Remark* 155.** — The separant of $\widetilde{A}_i^{(k)}$ is the product of $Sep_i$ by a product of powers of initials and separants of the $A_i$.

For the sake of simplicity, we restrict here to the case of prime ideals. Any finite subset $\mathscr{B}$ of $\{\widetilde{A}_i^{(k)} | 1 \le i \le p, k \in \mathbf{N}\}$ is a characteristic set of the prime algebraic ideal $(\mathscr{B}) : H_{\mathscr{B}}^{\infty} = (\mathscr{B}) : H_{\mathscr{A}}^{\infty}$. In more general situations, splittings may occur that could be considered *à la D5* [18]...

**Theorem 156.** — *Let $\mathscr{P}$ be a prime ideal of differential dimension 0 of $\mathscr{F}\{x_1, \ldots, x_n\}$ and $\mathscr{A} = \{A_1, \ldots, A_n\}$ a characteristic set of $\mathscr{P}$ for some ordering $\prec$, such that the main derivative of $A_i$ is $v_i = x_i^{(e_i)}$.*

*i) Assume that there exists a characteristic set $\mathscr{B} = \{B_1, \ldots, B_n\}$ of $\mathscr{P}$, being such that the main derivative of $B_i$ is $\bar{v}_i = x_i^{(f_i)}$. Let $I_1 := \{i \in [1, n] | f_i < e_i\}$, $I_2 := \{i \in [1, n] | f_i > e_i\}$ and $I_3 := \{i \in [1, n] | f_i = e_i\}$. We have the following propositions:*
*a) $\mathscr{B} \subset \left(\widetilde{A}_i^{(k)} \mid i \in [1, n],\ 0 \le k \le \min(f_i - e_i, 0)\right) : H_{\mathscr{A}}^{\infty}$;*
*b) If $\forall (i, i') \in I_1 \times I_2\ \bar{v}_i \prec \bar{v}_{i'}$, then*

$$\{B_i \mid i \in I_1\} \subset \left(\widetilde{A}_i^{(k)} \mid 1 \le i \le n;\ 0 \le k \le \max(0, f_i - e_i - 1)\right) : H_{\mathscr{A}}^{\infty};$$

*c) if for some $i_0 \in I_2$, $\ell_{i_0} < f_{i_0} - e_{i_0}$, then $B_{i_0} \notin \left(\widetilde{A}_i^{(k)} \mid i \in [1, n],\ 0 \le k \le \ell_i\right) : H_{\mathscr{A}}^{\infty}$.*

*ii) A characteristic set $\mathscr{B}$ satisfying the hypotheses of i) does exist iff*

$$_{\mathscr{A}}^{f}\Delta := \left| \frac{\partial \widetilde{A}_i^{(k)}}{\partial x_j^{(a)}} \mid i \in I_2,\ 0 \le k < f_i - e_i; \quad j \in I_1,\ f_j \le a < e_j \right| \notin \mathscr{P}.$$

**Proof.** — i) a) The char. set $\mathscr{B}$ cannot contain polynomials involving derivatives of each $x_i$ of order higher than $f_i$. So, by cor. 92, it must be included in $\left(\widetilde{A}_i^{(k)} \mid i \in [1, n],\ 0 \le k \le \max(f_i - e_i, 0)\right) : H_{\mathscr{A}}^{\infty}$.

b) If $\forall (i, i') \in I_1 \times I_2\ \bar{v}_i \prec \bar{v}_{i'}$, then $\mathrm{ord}_{x_j} B_i < f_j$, for $j \in I_2$, hence the result, using again cor. 92.

c) If for some $i_0 \in I_2$ $\ell_{i_0} < f_{i_0} - e_{i_0}$, then $x_{i_0}^{(f_0)}$ does not appear in the generators of $\mathscr{C} := \{\widetilde{A}_i^{(k)} | 1 \le i \le s,\ 0 \le k \le \ell_i\}$. So, if $B_{i_0}$ were in $(\mathscr{C}) : H_{\mathscr{A}}^{\infty}$, its initial would be reduced to 0 by $\mathscr{C}$, so would be in $(\mathscr{C}) : H_{\mathscr{A}}^{\infty} \subset \mathscr{P}$, which is impossible as $\mathscr{B}$ is a characteristic set.

ii) ⇒. Using lemma 85 i), the problem is reduced to the existence of a standard basis for $(\mathrm{d}P)_{\mathscr{M}_{\mathscr{P}}}$ with main derivatives $\mathrm{d}x_i^{(f_i)}$. For this, we only need to show that there exists in the vector space $\langle \mathrm{d}\widetilde{A}_i^{(k)} | 1 \le i \le n; 0 \le k \le \min(0, f_i - e_i)\rangle$ elements with the main derivatives $\mathrm{d}x_i^{(f_i)}$, for $i \in I_2 \cup I_3$ and $\mathrm{d}x_i^{(h)}$, $f_i \le h \le e_i$, for $i \in I_1$. This appends iff the determinant

$$\left| \frac{\partial \widetilde{A}_i^{(k)}}{\partial x_j^{(h)}} \;\middle|\; \begin{array}{l} i\in[1,n];\ 0\le k\le \min(f_i-e_i,0);\\ j\in I_1;\ f_j\le h\le e_j\\ \text{or } j\in I_2\cup I_3;\ h=f_j \end{array} \right|,$$

which is the determinant J of th. 89, does not vanish. It is equal to

$$\overset{f}{\underset{\mathscr{A}}{}}\Delta \prod_{i\in I_2} Sep_{\widetilde{A}_i^{(f_i)}} \prod_{i\in I_1\cup I_3} Sep_{A_i}$$

and, as $Sep_{\widetilde{A}_i^{(f_i)}}$ is a product of initials and separants of the $A_i \in \mathscr{A}$ by rem. 155, it does not belong to $\mathscr{P}$. So, a set with the requested main derivatives exists, that must be a standard basis of the module $\mathrm{d}\mathscr{P}$, for any ordering compatible with $\mathrm{ord}^f$, by invariance of the order. This implies the existence of a suitable char. set, by lem. 85.

⇐. If such a char. set exists, the $x_i^{(h)}$ with $e_i \le k < f_i$, for $i \in I_1$, must be independent. As these derivatives may be expressed as functions of the $x_i^{(h)}$ with $f_i \le h < e_i$, the Jacobian determinant of the transformation cannot vanish and is $\overset{f}{\underset{\mathscr{A}}{}}\Delta / \prod_{i\in I_2} \prod_{k=e_i}^{f_i-1} Sep_{\widetilde{A}_i^{(h)}}$, so that $\overset{f}{\underset{\mathscr{A}}{}}\Delta \notin \mathscr{P}$. ■

***Examples.* — 157)** Consider the char. set $\mathscr{A} := \{x_1'' + x_3' + x_4', x_2'' + x_4, x_3'', x_4''\}$. We want to test the existence of a char. set $\mathscr{B}$, with $f_1 = 3$, $f_2 = 4$, $f_3 = 1$ and $f_4 = 0$. We see that $\overset{f}{\underset{\mathscr{A}}{}}\Delta = 1$, so that the existence is granted. It may be computed by using $\widetilde{A}_1 = A_1$, $\widetilde{A}_1' = x_1'''$, $\widetilde{A}_2 = A_2$, $\widetilde{A}_2' = x_2''' + x_4'$ $\widetilde{A}_2'' = x_2^{(4)}$, that provides the solution $\mathscr{B} = \{x_1''', x_2^{(4)}, x_3' + x_4' + x_1'', x_4 + x_2''\}$.

**158)** As an illustration of this powerful result, we can go back to ex. 142 where $\mathscr{A} = \{x_1''' - x_1'', x_2' - x_1'', x_3' - x_1 + x_2, x_4' - x_2\}$. It is indeed impossible to have a new char. set of order 4 in $x_3$ and 2 in $x_4$, as $\widetilde{A}_3''' = x_3'''$, so that $\overset{f}{\underset{\mathscr{A}}{}}\Delta = 0$, whenever $f_3 = 4$. We also have $\widetilde{A}_4''' = x_4''' - x_1''$ and $\widetilde{A}_4^{(4)} = x_4^{(4)} - x_1''$, so that no char. set exists with $f_4 = 4$.

This theorem justifies also the informal claim that a *generic* system admits all order $n$-tuples $f$, provided that the system order is preserved, *i.e.* $\sum_{i=1}^n f_i = \sum_{i=1}^n e_i$ (see subsec. 8.2). The Zariski open sets where a given $n$-tuple admits a char. set is made explicit by the condition $\overset{f}{\underset{\mathscr{A}}{}}\Delta \notin \mathscr{P}$.

## 9.2 Bounding the order of differentiation of the initial system

From a computational standpoint, Jacobi was not fully satisfied with the last theorem. Claiming that it was sometimes more efficient to use derivatives of the $A_i$ instead of the $\widetilde{A}_i^{(k)}$ obtained by substitutions. This strongly suggests a practical experience of computing changes of ordering, although no explicit example is found in his manuscripts. We have already noticed that differentiation, introducing new derivatives, produces in the non linear case an exponential growth of the equations in dense representation. The situation becomes worse if substitutions are done at the same time, for then the degree can increase too. The best known *Bézout* bounds for the required eliminations imply to use Bézout's theorem, and the degrees will be the smallest using the $A_i^{(k)}$ instead of the $\widetilde{A}_i^{(k)}$. We see that Jacobi's intuition of the complexity issues meets here again contemporary research, such *D'Alfonso* as D'Alfonso *et al.* [14, 15, 17, 16], in the spirit of the Kronecker algorithm [30]. *Kronecker algorithm*

This problem is considered in § 18 of [II/23 a)] [40, p. 40–43]. The end of this manuscript seems lost and the sentence at the end of f° 2220 remains unachieved, but we can understand the general idea.

We use, as above, the notation $a_{i,j} := \mathrm{ord}_{x_j} A_i$ (*cf.* def. 71). With the notations and hypotheses of th. 156, one needs to differentiate equation $A_i$ $f_i - e_i$ times if $i \in I_2$. Then, generically, $A_i$, $1 \le i \le n$, must be differentiated $\ell_i$ times in order to compute the reduced derivatives $\widetilde{A}_i^{(h)}$ for $i \in I_2$, up to order $h = f_i - e_i$, where $\ell_i$ is such that: $\ell_i \ge f_i - e_i$ for $i \in I_2$ and $a_{j,j} + \ell_j \ge \max_i a_{i,j} + \ell_i$, so that the necessary reductions could be performed. The minimal solution of this problem is obtained by computing the unique minimal canon of the matrix $a_{i,j} + \max(f_i - e_i, 0)$, using the methods of subsection 3.6.

***Example* 159.** — Consider the char. set $\mathscr{A} := \{x_1'' + x_2 + x_3' + x_4, x_2''', x_3'' - x_4', x_4''\}$. We may compute the new char. set $\mathscr{B} := \{x_1^{(5)}, x_2 + x_3' + x_4 + x_1'', x_3'' - x_4; x_4''\}$ using derivatives of the $A_i$ up to $\ell_i$, with $\ell = (3, 0, 2, 1)$, which is the minimal canon of matrix

$$\left(a_{i,j} + \max(f_i - e_i, 0)\right) = \begin{pmatrix} 5 & 3 & 4 & 3 \\ -\infty & 3 & -\infty & -\infty \\ -\infty & -\infty & 2 & 1 \\ -\infty & -\infty & -\infty & 2 \end{pmatrix}.$$

It is easily checked that these orders cannot be lowered. But if we change $A_3$ to $x_3'' + x_4'$, $\widetilde{A}_1''' = x_1^{(5)} = A_1''' - A_2 - A_3''$, so that no strict derivative of $A_4$ is needed.

***Remark* 160.** — This bound is sharp. It is here convenient to complete the system equations with left members: $A_i = T_i$. Furthermore, we may express explicitly the dependency of $x_i^{(e_i+\ell_i)}$, with respect to $T_j^{(\ell_j)}$. The computation of $\widetilde{A}_i^{(\ell_i)}$ will actually require to differentiate $A_j$ up to order $\ell_j > 0$ if $\partial x_i^{(e_i+\ell_i)}/\partial T_j^{(\ell_j)} \neq 0$. To compute this value, we use the path relation in the canon $a_{i,j} + \ell_i$, as in alg. 9 d) p. 16. Let $\Pi$ denote the set of paths $\pi = \{\pi_0 = j, \ldots, \pi_r = i\}$, from row $j$ to row $i$ ($r$ may depend on $\pi$). An easy recurrence on $r$ shows that, denoting by $r(\pi)$ the length of path $\pi$, we have:

$$\frac{\partial x_i^{(e_i+\ell_i)}}{\partial T_j^{(\ell_j)}} = \sum_{\pi\in\Pi} \frac{\prod_{h=0}^{r(\pi)-1} -\partial A_{\pi_{h+1}}^{(\ell_{\pi_{h+1}})} / \partial x_{\pi_h}^{(a_{\pi_{h+1},\pi_h}+\ell_{\pi_{h+1}})}}{\prod_{h=0}^{r(\pi)} \partial A_{\pi_h} / \partial x_h^{(a_{\pi_h,\pi_h})}}.$$

***Example* 161.** — We go back to ex. 159. We have $\widetilde{A}_1''' = A_1''' - A_2 - A_3'' - 2A_4'$; indeed, there is a single path of length 1 from rows 2 and 3 to row 1, but two paths from row 4 to row 1, one of length 1: $(4,1)$ and one of length 2: $(4,3,1)$, both with with the same coefficient: 1. Taking $x_3'' + x_4$ for $A_3$, the coefficients for the 2 paths cancel and $\widetilde{A}_1''' = A_1''' - A_2 - A_3''$.

***Remark* 162.** — Using explicit normal forms, we may assume that the leading derivatives $x_i^{(e_i)}$ do not appear in the right members $F_i(x)$. It is no longer the case with characteristic sets. All we know is that the leading derivative of $A_i$ may only appear in $A_j$ with a strictly smaller degree. But we may, without changing the main derivatives and by purely algebraic computations, reduce to a case where the $A_i$ with $f_i > e_i$ do not depend on the main derivatives of the $x_j$ for which $f_j \leq e_j$. This may be an advantage for some systems. *explicit normal form*

The next example illustrates this situation.

***Example* 163.** — Let $\mathscr{A} = \{x_1^{(e_1)} + x_2^{(e_2)}, (x_2^{(e_2)})^2 + x_2\}$. On may also choose a char. set with the same leading derivatives, but different degrees: $\mathscr{B} = \{(x_1^{(e_1)})^2 + x_2, x_2^{(e_2)} + x_1^{(e_1)}\}$. In order to compute $\mathscr{C} := \{[(x_1^{(e_1)})^2]^{(e_2)} - x_1^{(e_1)}, x_2 + (x_1^{(e_1)})^2\}$, one must differentiate $B_1$ $e_2$ times, but no strict derivative of $B_2$ is required. Working with $\mathscr{A}$, we need to differentiate both $A_1$ an $A_2$ $e_2$ times.

## 10 Resolvents

JACOBI's treatment of resolvent computations deserves a special interest. Together with algorithmic issues in sec. 2 and 3 and his proof of the bound itself (th. 94), it is one of the only places in his manuscripts where he gives the precise framework of a proof, that involves a wide varieties of the tools that he has introduced: the bound on the order (sec. 6), the shortest reduction (sec. 7), changes of orderings (sec. 9) and combinatorics issues related to minimal canons, subject to lower values (subsec. 3.6).

*Cluzeau and Hubert* For a modern approach of the question, one may refer to Cluzeau and Hubert [12]. Following Ritt [86, chap. II § 22], we define resolvents in the following way. We also introduce notions of *weak* and *local resolvent*, closer to Jacobi's approach.

**DEFINITION 164.** — *Let $\mathscr{P}$ be a prime differential ideal. We call a* differential resolvent[21] *of $\mathscr{P}$ the data of two differential polynomials $P$ and $Q$, together with a* *resolvent*

[21]As we are only concerned here with the differential case, we omit "differential" in the sequel.

*characteristic set $\mathscr{A}$ of the prime differential ideal $[\mathscr{P}, Qw - P] : Q^\infty$ (in $\mathscr{F}\{x_1, \ldots, x_n, w\}$), such that $A_{n+1}(w)$ depends only on $w$ and, for $1 \le i \le n$, $A_i$ is linear in its main derivative $x_i$. We may choose for $A_{n+1}(w)$ a prime polynomial that is unique up to* *resolvent polynomial* *multiplication by a non zero element of $\mathscr{F}$. We call it the* resolvent polynomial *of $\mathscr{P}$ for $P/Q$.*

*local resolvent* *We call a local resolvent a char. set of $\mathscr{P}$ where $A_{n+1}$ depends only on $w$, but $A_i$* *weak resolvent* *is only requested to be of order 0 in $x_i$ and a weak resolvent a char. set $\mathscr{B}$ such that $B_i$ is of order 0 in $x_i$, for an ordering $\prec$ such that $w^{(\mathcal{O})} \succ x_i \succ w^{(\mathcal{O}-1)}$, $1 \le i \le n$, so that $A_{n+1}$ may depend on the $x_i$, $1 \le i \le n$.*

*In the sequel, $\mathscr{A}$ stands for a local resolvent and $\mathscr{B}$ for a weak resolvent. Considering only the case $w = x_{j_0}$, we will use $A_{j_0} = A_{n+1}(x_{j_0})$ (resp. $B_{j_0} = B_{n+1}(x_{j_0})$) and avoid to introduce a new useless letter $w$.*

Resolvents are local resolvents, but with local resolvents, assuming $x_i \succ x_{i'}$, $A_i$ may sometimes also depend of $x_{i'}$. The weak notion may be convenient, even when a resolvent (with $A_j$ linear in $x_j$) exists, as it allows to express $x_j$ as an implicit algebraic function involving only derivatives of $w$ up to order $\mathcal{O} - 1$ (see below th. 175 iii) that is a special case of th. 156 i) b). In a local setting, that fit best Jacobi's statements, we may always define algebraic functions $F_i$, such that $w^{(\mathcal{O})} = F_0(w)$ and $x_i = F_i(w)$, where the $F_i$ are of order strictly less than $\mathcal{O}$.

***Remarks.* — 165)** Considering the resultant of $A_0(w)$ and $A_j(x_j, w)$, we see that *weak resolvent* the existence of a resolvent implies that of a weak resolvent. The converse obviously does not stand, but the existence of a local and of a weak resolvent are equivalent, which reduces in practice to a change of ordering that involves only algebraic computations, as the main derivatives remain the same.

*weak and local resolvents* **166)** Resolvents, local resolvents and weak resolvents may coïncide. It is always the case for linear ideals.

***Example* 167. —** Consider the differential system $x_1'' - x_2^2/2 = 0$, $x_2 x_2' - x_3^2/2 = 0$, $x_3 x_3' - x_2 - x_3 = 0$. With $w = x_1$, we have a resolvent char. set $((x_1^{(4)})^2/2 - x_1''' - x_1'')^2 - 4x_1'' x_1'''$, $x_1^{(4)} x_2 - (x_1^{(4)})^2/2 + x_1''' - x_1''$, $x_1^{(4)} x_3 - (x_1^{(4)})^2/2 - x_1''' + x_1''$ and a weak resolvent char. set $x_1^{(4)} - x_2 - x_3$, $x_2^2 - 2x_1'' = 0$, $x_3^2 - 2x_1'''$.

*differential resultant* The resolvent polynomial is also a special case of *differential resultant.* On *Carrà Ferro* this topic, see Ritt [84, chap. III § 34 p. 47], Carrà Ferro [8] or Li *et al.* [70], that *Li* contains complexity results relying on Jacobi's bound.

**Definition 168. —** *If $A$ is the order matrix of a system $P$ of $n$ equations in $n$ variables, let $I \subset [1, n]$, we denote $[1, n] \setminus I$ by $\bar{I}$. For $I, J \subset [1, n]$, we denote by ${}_{I,J}$ the matrix $(a_{\hat{\imath},\hat{\jmath}} | \hat{\imath} \in I, \hat{\jmath} \in J)$ and by $\mathcal{O}_{I,J}$ its tropical determinant.*

*We will also denote $A_{\bar{I},\bar{J}}$ by $\bar{A}_{I,J}$ and by $\bar{\mathcal{O}}_{I,J}$ its tropical determinant. We denote* *tropical determinant* *by $\bar{\mathcal{O}}_{i,j}$ the tropical determinant of the matrix $\bar{A}_{i,j} := (a_{\hat{\imath},\hat{\jmath}} | \hat{\imath} \ne i, \hat{\jmath} \ne j)$.*

*shortest reduction* Jacobi implicitly assumes that he considers the situation where the shortest reduction (def. 116) may be used, which corresponds in our setting to components $\mathcal{P} \in [P]:\nabla_P^\infty$. So, we know that the order of $\mathcal{P}$ is precisely Jacobi's bound $\mathcal{O}$, that must be the order of a resolvent polynomial. *resolvent polynomial*

Jacobi [40, § 4 p. 58–63] also restricts to the case where a resolvent[22] exists when choosing $w = x_{j_0}$. In order to compute it, he claims that one needs to differentiate each equation $P_{i_0}$ up to an order that corresponds to the tropical determinant[23] $\bar{\mathcal{O}}_{i_0,j_0}$ of the order submatrix $\bar{A}_{i_0,j_0}$.

The next theorem is an attempt to give a precise meaning to this statement in differential algebra, when using the strong bound. Then, Jacobi's bound stands only for the computation of the resolvent polynomial. We will see that the computation of the parametrization $A_j(x_j, w)$ may require higher order of derivation of $P_i$, unless when $\bar{\mathcal{O}}_{i,j_0} \geq \lambda_i$, which stands in the case of an order matrix $A_P$ of non negative integers, according to lem. 12. A convention equivalent to the weak bound was perhaps implicitly assumed by Jacobi. His proof relies on a combinatorial argument that is used to prove assertion i) of the next lemma. *weak bound*

**Lemma 169.** — *We assume that there exists a component $\mathcal{P}$ of $[P]:\nabla_P^\infty$ with a weak resolvent for $w = x_{j_0}$. Up to a permutation, we may also assume that $\sum_{i=1}^n a_{i,i} = \mathcal{O}_P$. Let $\ell$ be the minimal canon of $A_P$, subject to the condition[24] $\ell_{j_0} + a_{j_0,j_0} \geq \mathcal{O}_P + D$, with $D \in \mathbf{N}$.*

*i) For any $i_0 \in [1,n]$, we have c) $\Longrightarrow$ b) $\Longrightarrow$ a). When*

$$D > (n-1) \max_{(i,j)\in[1,n]^2} a_{i,j}, \tag{15}$$

*we have furthermore a) $\Longrightarrow$ c), so that when equation (15) stands, these statements are equivalent:*

*a) $\bar{\mathcal{O}}_{i_0,j_0} \neq -\infty$;*

*b) $\ell_{i_0} = \bar{\mathcal{O}}_{i_0,j_0} + \ell_{j_0} - \bar{\mathcal{O}}_{j_0,j_0}$;*

*c) There is a path from row $i_0$ to row $j_0$ by the reflexive transitive closure of the path relation, as defined in step c) of Jacobi's algorithm (alg. 9).*

---

[22] The word resolvent was not used by Jacobi, but he evokes the notion as something well known in the the mathematical folklore of his time: *"It is usual that this type of normal forms be considered before others by mathematicians".* *normal form*

[23] Nanson [76] and Jordan [54] proposed independently heuristic methods for proving Jacobi's bound, that rely on resolvent computations. The first considers the case $n = 3$ and the second the case $n = 4$, recursively using formula $\mathcal{O} = \max_i a_{i,j_0} + \bar{\mathcal{O}}_{i,j_0}$ and the bound $\bar{\mathcal{O}}_{i,j_0}$ on the order of differentiation of each equation $P_i$, which is guessed using informal considerations on the number of derivatives of the equation $P_i$, and the number of derivatives of the $x_j$, $j \neq j_0$ that the computation of a resolvent requests. Their relation is expressed by the formula *Nanson; Jordan*

$$\sum_{i=1}^n (\bar{\mathcal{O}}_{i,j_0} + 1) = 1 + \sum_{j \neq j_0} \left( \max_{i=1}^n (\bar{\mathcal{O}}_{i,j_0} + a_{i,j}) + 1 \right),$$

so that there are exactly one more algebraic equations than derivatives of the $x_j$, $j \neq j_0$, involved in them, making their potential elimination possible. The last formula is proved in cor. 173.

[24] It may be computed as in subsec. 3.6. By prop. 54, the path relation does not depend on the choice of the permutation $\sigma$. The canon itself is also independent of this choice: if $a_{j_0,j_0} + \ell_{j_0}$ is maximal, then for any permutation $\sigma$ such that $\mathcal{O} = \sum_{i=1}^n a_{i,\sigma(i)}$, the quantity $a_{\sigma^{-1}(j_0),j_0} + \ell_{\sigma^{-1}(j_0)}$ must be maximal too.

*ii) Let $I \subset [1,n]$ be the set of rows $i$ such that, for $D > (n-1)\max_{i,j} a_{i,j}$, there is a path to row $j_0$ from row $i$ by the path relation and $\bar{I} := [1,n] \setminus I$. For any $i_0 \in \bar{I}$, $\ell_{i_0} = \lambda_{i_0}$, where $\lambda$ is the minimal canon of $A_P$. The submatrix $A_{I,\bar{I}}$ only contains $-\infty$ elements and $\mathcal{O}_{\bar{I},\bar{I}} = 0$.*

*iii) For $i \in \bar{I}$, $a_{i,i} = 0$. If $i, j \in \bar{I}$ and $\ell_i < \ell_j$, then $a_{j,i} = -\infty$.*

*iv) If $D = 0$ and $\bar{I} = \emptyset$, then $\forall 1 \le i \le n$, $\ell_i = \bar{\mathcal{O}}_{i,j_0}$.*

Proof. — i) b) $\Rightarrow$ a). — The quantity $\bar{\mathcal{O}}_{i_0,j_0} = \ell_{i_0} - \ell_{j_0} + \bar{\mathcal{O}}_{j_0,j_0} = \ell_{i_0} - \ell_{j_0} + \sum_{i \neq j_0} a_{i,i}$ is an integer. c) $\Rightarrow$ b). — We proceed as indicated by Jacobi[25]. Let $\iota : k \in [0,r] \mapsto \iota_k \in [1,n]$ be an injection such that there is an elementary path from row $\iota_k$ to row $\iota_{k+1}$, with $\iota_0 = i_0$ and $\iota_r = j_0$. Then,

$$\sum_{k=0}^{r} \left(a_{\iota_k,\iota_k} + \ell_{\iota_k}\right) = a_{j_0,j_0} + \ell_{j_0} + \sum_{k=0}^{r-1} \left(a_{\iota_{k+1},\iota_k} + \ell_{\iota_{k+1}}\right). \tag{16}$$

We remark that the $a_{i,i} + \ell_i$, $i \notin \mathrm{Im}(\iota)$ and the $a_{\iota_{k+1},\iota_k} + \ell_{\iota_{k+1}}$, $0 \le k < r$ form a transversal family of $n-1$ maximal elements in $\bar{A}_{i_0,j_0} + \ell$, so that

$$\sum_{i \notin \mathrm{Im}(\iota)} a_{i,i} + \sum_{k=0}^{r-1} a_{\iota_{k+1},\iota_k} = \bar{\mathcal{O}}_{i_0,j_0}.$$

Adding $\sum_{i \notin \mathrm{Im}(\iota)} \left(a_{i,i} + \ell_i\right)$ to both sides of (16), we get

$$\sum_{i=1}^{n} a_{i,i} (= \mathcal{O}_P) + \sum_{i=1}^{n} \ell_i = \bar{\mathcal{O}}_{i_0,j_0} + \sum_{i \neq i_0} \ell_i + \ell_{j_0} + a_{j_0,j_0},$$

which implies $\ell_{i_0} = \bar{\mathcal{O}}_{i_0,j_0} + \ell_{j_0} + a_{j_0,j_0} - \bar{\mathcal{O}}_P = \bar{\mathcal{O}}_{i_0,j_0} + \ell_{j_0} - \bar{\mathcal{O}}_{j_0,j_0}$.

The proof of i) a)$\Rightarrow$ c) requires ii) and is postponed.

ii) We first prove that $a_{i_1,i_0} = -\infty$ for all $(i_1, i_0) \in I \times \bar{I}$. By lem. 51 i), there is a path from any row $i_0 \in \bar{I}$ to a row $i$ with $\ell_i = 0$. Then, $\ell_{i_0} = \lambda_{i_0}$ and by lem. 51 ii), $\lambda_{i_0} \le \rho \max_{i,j} a_{i,j}$, where $\rho$ is the length of the path. As there is no path from $i_0$ to $j_0$, there is no path from $i_0$ to any row $i_1 \in I$, so that $\rho \le \sharp\bar{I} - 1 \le n-2$ and $a_{i_0,i_0} + \lambda_{i_0} \le (n-1)\max_{i,j} a_{i,j} < D$. Let $\iota$ be a path from $i_1$ to $j_0$, with $\iota_1 = i_1$ and $\iota_r = j_0$. Proceeding as in the proof of lem. 51 ii), we see that $\ell_{i_1} = \ell_{j_0} + \sum_{k=1}^{r-1} \left(a_{\iota_{k+1},\iota_k} - a_{\iota_k,\iota_k}\right) \ge \mathcal{O}_P - a_{j_0,j_0} + D - \sum_{k=1}^{r-1} a_{\iota_k,\iota_k} \ge D$. So, $a_{i_0,i_0} + \ell_{i_0} \ge a_{i_1,i_0} + \ell_{i_1}$ implies $a_{i_1,i_0} \le \left(a_{i_0,i_0} + \ell_{i_0}\right)(< D) - \ell_{i_1}(\ge D) < 0$ and $a_{i_1,i_0} = -\infty$.

Let $\Delta$ be the determinant of the submatrix of $J_P$ corresponding to the rows and columns in $\bar{I}$: $\Delta$ is a factor of $\nabla_P$ and does not vanish modulo $\mathcal{P}$. This means that the $x_{i_0}$, $i_0 \in \bar{I}$ belong to an extension of $\mathcal{F}\langle x_{j_0}\rangle$ of order $\mathcal{O}_{\bar{I},\bar{I}}$. As a weak resolvent exists, we know that $x_{i_0}$ is algebraic over $\mathscr{F}\langle x_{j_0}\rangle$, so that $\mathcal{O}_{\bar I,\bar I} = 0$.

[25]One may look at formula (3) p. 16 for an illustration of the situation. The notion of increasing path is used here in the reverse way: one deduces a maximal transversal sum for $\bar{A}_{i_0,j_0}$ from a maximal transversal sum for $A_P$.

$$\begin{array}{cc} & \begin{array}{cc} \overbrace{\qquad\qquad}^{I} & \overbrace{\qquad\qquad}^{\bar I} \end{array} \\ \begin{array}{c} I \left\{ \vphantom{\begin{matrix}a\\b\\c\end{matrix}} \right. \\ \bar I \left\{ \vphantom{\begin{matrix}a\\b\\c\end{matrix}} \right. \end{array} & \left( \begin{array}{c|c} \qquad\qquad & -\infty \\ & \\ \hline & \begin{matrix} 0 & & \\ & \ddots & \\ & & 0 \end{matrix} \end{array} \right) \end{array}$$

a) $\Rightarrow$ c). — If $i_0 \notin I$, then by ii), we know that the non $-\infty$ elements in $A_{i_0,j_0}$ are contained in $\sharp\bar I - 1$ rows and $\sharp I - 1$ columns. So, according to th. 17, $\mathcal{O}(A_{i_0,j_0}) = -\infty$. This achieves the proof of i).

iii) As $\mathcal{O}_{\bar I,\bar I} = \sum_{i\in\bar I} a_{i,i} = 0$, $a_{i,i} = 0$ for any $i \in \bar I$. As $a_{i,i} + \ell_i$ is maximal, if $\ell_i < \ell_j$, we need have $a_{j,i} = -\infty$.

iv) By i), we know that, when $I = [1, n]$, $\ell_i = \bar{\mathcal{O}}_{i,j_0} + D$ is the minimal canon for $A_P$, such that $\ell_{j_0} \geq \mathcal{O}_P - a_{j_0,j_0} + D$, for $D$ great enough. This implies that $\ell_i = \bar{\mathcal{O}}_{i,j_0}$ is the minimal canon for $A_P$, such that $\ell_{j_0} \geq \mathcal{O}_P - a_{j_0,j_0} = \bar{\mathcal{O}}_{j_0,j_0}$. ■

***Remarks.* — 170)** In the case $\bar I = \emptyset$, this shows that we can compute the $n$ tropical subdeterminants $\bar{\mathcal{O}}_{i,j_0}$, $1 \leq i \leq n$, with total cost $O(n^3)$, by using Jacobi's algorithm to compute the canon $\ell$ as in section 3.6, which is faster than applying the algorithm $n$ times.

**171)** In the case $\bar I = \emptyset$ and $D = 0$, this canon $\ell$ is explicitly considered by Jacobi, but not the associated reduction. Computing a shortest reduction (th. 110) with *Egerváry* respect to this canon may be an efficient step in the computation of a resolvent, *resolvent* that is implicit in Jacobi's proof [40, § 4 p. 59].

**172)** In the case $\bar I \neq \emptyset$, we have a block triangular decomposition, as in subsec. 4.3: *block decomposition* the equations $P_i$, $i \in I$ do not depend on the $x_i$ $i \in \bar I$ and the $P_i$, $i \in \bar I$ with $\ell_i \leq \ell_{i_0}$ do not depend on the $x_j$ with $\ell_j > \ell_{i_0}$.

As already noticed, the proof of assertion i) c)$\Rightarrow$ b) follows Jacobi's construction. Jacobi [40, p. 61] denotes the term $a_{\iota_k,\iota_k}$ by $S^{(k)}$ and $a_{\iota_{k+1},\iota_k}$ by $\overline{S^{(k)}, S^{(k+1)}}$. The same construction may be used to prove the following tropical identity, of which the relation of footnote 23 is an easy consequence.

**Proposition 173.** — *Let $A$ be a square $n \times n$ matrix, then:*

$$\sum_{i=1}^{n} \bar{\mathcal{O}}_{i,j_0} = \sum_{j\neq j_0} \left( \max_{i=1}^{n} (\bar{\mathcal{O}}_{i,j_0} + a_{i,j}) \right).$$

Proof. — Assume for simplicity that $j_0 = n$ and that the $a_{i,i}$, $1 \leq i \leq n$, are transversal maxima. Then $\max_{i=1}^{n}(\bar{\mathcal{O}}_{i,n} + a_{i,j})$ is equal to $\mathcal{O}_{A_{[j]}}$, where $A_{[j]}$ is obtained by replacing column $n$ by column $j$ in $A$. Let $B = A + \ell$ be the minimal canon for $A$, subject to the condition $\ell_n \geq \bar{\mathcal{O}}_{n,n}$. As in the proof of lem. 169 i) c)⇒ b), this implies that there is a path from any row $i_0$ in $B$ to row $n$. Let $\iota : [0, r] \mapsto [1, n]$ be this path, with $b_{\iota_{k+1},\iota_k} = b_{\iota_k,\iota_k}$, $0 \leq k < r$, $\iota_0 = i_0$ and $\iota_r = n$. Again, as in the proof of lem. 169 i) c)⇒ b), we see that the $a_{i,i}$, $i \notin \mathrm{Im}(\iota)$, and $a_{\iota_{k+1},\iota_k}$, $0 \leq k < r$, form a family of $n-1$ transversal maxima in $\bar{B}_{j,n}$, that may be completed with the element $a_{j,j}$ in column $n$, to form a family of $n$ transversal maxima in $A_{[j]}$. So $\mathcal{O}_{A[j]} = \bar{\mathcal{O}}_{j,n} + a_{j,j}$, which implies

$$\sum_{j=1}^{n-1} \mathcal{O}_{A_{[j]}} = \sum_{j=1}^{n-1} \bar{\mathcal{O}}_{j,n} + \sum_{j=1}^{n-1} a_{j,j} = \sum_{j=1}^{n} \bar{\mathcal{O}}_{j,n},$$

hence the proposition. ∎

The next lemma is a ready-to-use specialization of cor. 92, adapted to common situations in differential algebra, that we will use in the proof of th. 175 bellow.

**Lemma 174.** — *Let $P_i \in \mathcal{F}\{x_1, \ldots, x_n\}$, $1 \leq i \leq n$ be differential polynomials, such that the subset $P_I := \{P_i | i \in I\}$ does not depend on the $x_j$ for $j \in \bar{I}$.*

*Then, a) $\nabla_{P_I}$ divides $\nabla_P$; b) for any prime component $\mathcal{P}$ of $[P] : \nabla_P^\infty$, $\mathcal{P} \cap \mathcal{F}\{x_j | j \in I\}$ is a prime component of the radical ideal $[P_I] : \nabla_P^\infty \supset [P_I] : \nabla_{P_I}^\infty$.*

Proof. — To prove i) a), we proceed as in the proof of lem. 169 ii): $a_{i,j} = -\infty$ for $(i, j) \in I \times \bar{I}$, so that $\mathrm{J}_P$ has a block decomposition and $|\mathrm{J}_P| = |\mathrm{J}_{I,I}| \cdot |\mathrm{J}_{\bar{I},\bar{I}}|$. For b), we may then use th. 110 ii) with a minimal cover $\mu, \nu$, such that $\min_{i\in\bar{I}} \mu_i > \max_{i\in I} \mu_i$. *E.g.*, we may use the Jacobi cover $\alpha, \beta$ and define $\mu_i = \alpha_i + h$, for $i \in \bar{I}$, with $h = \max_{i\in I} \mu_i - \min_{i\in\bar{I}} \mu_i + 1$, $\nu_j = \beta_j - h$ for $j \in \bar{I}$, keeping them unchanged in the remaining cases. By cor. 113, we can choose a permutation $\sigma$ such that $\mathcal{P}$ is a component of $[P] : \nabla_\sigma^\infty$ and the proof of th. 110 ii) shows the existence of a char. set $\mathcal{A}$ such that its first $\sharp I$ elements, that belong to a prime component of $[P_I] : \nabla_P^\infty$, also belong to $\mathcal{P} \cap \mathcal{F}\{x_j | j \in I\}$ and reduce all $Q \in \mathcal{P} \cap \mathcal{F}\{x_j | j \in I\}$ to 0, so that they form a char. set of this prime ideal. Hence the lemma. ∎

We can now state the main result of this section.

**Theorem 175.** — *Let $P := \{P_1, \ldots, P_n\} \subset \mathcal{F}\{x_1, \ldots, x_n\}$ and $a_{i,j} := \mathrm{ord}_{x_j} P_i$, with the convention of def. 71. Assume that a weak resolvent char. set $\mathcal{A}$ exists for a prime component $\mathcal{P}$ of $[P] : \nabla_P^\infty$ (cf. def. 76 and cor. 113), for $w = x_{j_0}$. Let $\mu, \nu$ be the minimal cover associated to $\ell$, defined using the notations of lem. 169 with $D = 0$. We further define $\mathrm{R}_k^{\tilde{\nu}} := \mathcal{F}[x_j^{(\kappa)} | 1 \leq i \leq n,\ 0 \leq \kappa \leq \tilde{\nu}_j + k]$, where $\tilde{\mu}, \tilde{\nu}$ is the minimal cover associated to any canon $\tilde{\ell}$.*

weak or local resolvent *i) Let $L := \max_{i=1}^n \ell_i$, a or weak resolvent of $\mathcal{P}$ belongs to the prime component $\mathcal{P} \cap \mathrm{R}_L$ of $[P_i^{(k)} | 1 \leq i \leq n,\ 0 \leq k \leq \ell_i] : \nabla_P^\infty$.*

*ii) If $i_0 \in I$, the polynomials $A_{i_0}$, or $B_{i_0}$, of a local resolvent $\mathcal{A}$, or a weak resolvent $\mathcal{B}$ of $\mathcal{P}$, belong to the component $\mathcal{P} \cap \mathrm{R}_{\bar{L}}^{\bar{\nu}}$ of $[P_i^{(k)} | i \in I,\ 0 \leq k \leq \bar{\mathcal{O}}_{i,j_0}] : \nabla_P^\infty$, where $\bar{\mu}$, $\bar{\nu}$ is the minimal cover associated to the canon $\bar{\ell}_i := \bar{\mathcal{O}}_{i,j_0}$ of $A_{I,I}$ and $\bar{L} := \max_{i\in I} \bar{\ell}_i$.*

*iii) If $i_0 \neq j_0$, $i_0 \in I$ and if $\mathscr{B}$ is a weak resolvent with $B_{i_0} = B_{i_0}(w, x_h; h \in H)$ and $\bar{\mathcal{O}}_{h,j_0} + a_{h,i_0} > 0$ for all $h \in H$, then $B_{i_0}$ belongs to the component $\mathscr{P} \cap \mathrm{R}^{\bar{\nu}}_{\bar{L}-1}$ of $[P_i^{(k)} | i \in I,\ 0 \le k < \bar{\mathcal{O}}_{i,j_0}] : \nabla_P^\infty$.*

*iv) For any $i_0 \in \bar{I}$, let $J_{i_0} := \{j \in \bar{I} | \ell_j \ge \ell_{i_0}\} \cup I$, the polynomials $A_{i_0}$, or $B_{i_0}$, of a local resolvent $\mathscr{A}$, or a weak resolvent $\mathscr{B}$ of $\mathscr{P}$, belong to the component $\mathscr{P} \cap \mathrm{R}^{\hat{\nu}}_{\hat{L}}$ of $[P_i^{(k)} | i \in J_{i_0},\ 0 \le k \le \hat{\ell}] : \nabla_P^\infty$, where $\hat{\ell}$ is the canon of $A_{J_{i_0}, J_{i_0}}$ defined by $\ell_i - \ell_{i_0}$, for $i \in \bar{I}$, or $\max(\ell_i - \ell_{i_0}, \bar{\mathcal{O}}_{i,j_0})$, for $i \in I$, $\hat{\mu}$, $\hat{\nu}$ denotes the associated cover and $\hat{L} := \max_{i \in J_{i_0}} \hat{\ell}_i$.*

*v) For any $i_0 \in I$ and any $0 \le k < a_{i_0,i_0} + \bar{\mathcal{O}}_{i_0,j_0}$, $x_{i_0}^{(k)}$ is algebraic over the field $\mathscr{F}(x_{j_0}^{(\kappa)} | 0 \le \kappa < \mathcal{O})$ and there exists a relation $Q(x_{i_0}^{(k)}, x_{j_0}) = 0$ that belongs to the component $\mathscr{P} \cap \mathrm{R}^{\bar{\nu}}_{\bar{L}-1}$ of $[P_i^{(k)} | i \in I,\ 0 \le k < \bar{\mathcal{O}}_{i,j_0}] : \nabla_P^\infty$.*

Proof. — i) We have $\mathcal{O} = \max_{i=1}^n a_{i,j_0} + \bar{\mathcal{O}}_{i,j_0}$, so that $\bar{\mathcal{O}}_{i,j_0} \le \mathcal{O}$. By lem. 169 i), $\ell_{j_0} \ge \bar{\mathcal{O}}_{j_0,j_0}$, so that $\mathrm{ord}^\nu x_{j_0}^{(\mathcal{O})} = \mu_{j_0} + \mathcal{O}_{j_0,j_0} \le \mu_{j_0} + \ell_{j_0} = L$, as $\mu_i = L - \ell_i$ (prop. 20). For any $1 \le i \le n$, $\mathrm{ord}^\nu x_i = \mu_i - a_{i,i} \le L$, so for any $Q$ that belongs to some weak or local resolvent, $\mathrm{ord}^\nu_{x_{j_0}} Q \le L$ and $\mathrm{ord}^\nu_{x_i} Q \le L$, for $i \neq j_0$. This means that $Q \in \mathscr{P} \cap \mathrm{R}^\nu_L$, which is a component of $[P_i^{(k)} | 1 \le i \le n,\ 0 \le k \le \ell_i] : \nabla_P^\infty$ by th. 110 i) equ. (11).

ii) By lem. 174 i) $\mathscr{P} \cap \mathscr{F}\{x_j | j \in I\}$ is a component of $[P_I] : \nabla_{P_I}^\infty$. By lem. 169 iv), the minimal cover for $P_I$ subject to the condition $\ell_{j_0} \ge \mathcal{O}_{j_0,j_0} = \mathcal{O} - a_{j_0,j_0}$ is $\bar{\ell}_i = \bar{\mathcal{O}}_{i,j_0}$. The result is then a straightforward consequence of i)

iii) If $i_0 \neq j_0$ and $\mathscr{B}$ is a weak resolvent, we have $\mathrm{ord}_{x_{j_0}} B_{i_0} \le \mathcal{O} - 1$. We may as in ii) reduce to the case $I = [1, n]$. Then, for any $h \in H$, $\mathrm{ord}^{\bar{\nu}} x_h = \bar{\mu}_h - a_{h,h} = \bar{L} - \bar{\mathcal{O}}_{h,j_0} - a_{h,h}$, so if $\bar{\mathcal{O}}_{h,j_0} + a_{h,h} > 0$, $\mathrm{ord}^{\bar{\nu}} B_{i_0} < \bar{L}$. This implies that $B_{i_0}$ belongs to $\mathscr{P} \cap \mathrm{R}^{\bar{\nu}}_{\bar{L}-1}$, which is a component of $[P_i^{(k)} | i \in I,\ 0 \le k \le \bar{\mathcal{O}}_{i,j_0} - 1] : \nabla_P^\infty$ by th. 110 i) equ. (11).

iv) By lem. 169 iii), the variables $x_j$, with $\ell_j < \ell_{i_0}$, do not appear the equations $P_i$, $i \in J_{i_0}$. By lem. 174, we can restrict to this subsystem, for which $\hat{\ell}$ is a canon, subject to $\hat{\ell}_{j_0} \ge \bar{\mathcal{O}}_{j_0,j_0}$. Indeed, if $\hat{\ell}_i = \ell_i - \ell_{i_0}$, then $a_{i,i} + \hat{\ell}_i \ge a_{i',i} + \hat{\ell}_{i'}$ for all $i'$ such that $\hat{\ell}_{i'} = \ell_{i'} - \ell_{i_0}$, as $\ell$ is a canon. If $\hat{\ell}_{i'} = \bar{\mathcal{O}}_{i',j_0}$, then $a_{i,i} + \hat{\ell}_i \ge a_{i,i} + \bar{\mathcal{O}}_{i,j_0} \ge a_{i',i} + \bar{\mathcal{O}}_{i',j_0}$, as $\bar{\mathcal{O}}$ is a canon of the order matrix restricted to rows and columns in $I$. In the same way, if $\hat{\ell}_i = \bar{\mathcal{O}}_{i,j_0}$, then $a_{i,i} + \hat{\ell}_i \ge a_{i',i} + \ell_{i'}$ when $\ell_{i'} = \bar{\mathcal{O}}_{i,j_0}$ and if $\hat{\ell}_{i'} = \ell_i - \ell_{i_0}$, then $a_{i,i} + \hat{\ell}_i \ge a_{i,i} + \ell_i - \ell_{i_0} \ge a_{i',i} + \ell_{i'} - \ell_{i_0}$, as $\ell$ is a canon. So $\hat{\ell}$ is greater than the canon $\ell$ for the system $P_i$, $i \in J_0$ and we only have to use i).

v) Using ii), we can reduce to the case $\bar{I} = \emptyset$. Then, we may increase the set of variables with $x_{n+1}$, satisfying $x_{n+1} - x_{i_0}^{(k)} = 0$, so that $n + 1 \in \bar{I}$. As $k < a_{i_0,i_0} + \bar{\mathcal{O}}_{i_0,j_0}$, $\ell_i = \bar{\mathcal{O}}_{i,j_0}$ for $i \in I$ and $\ell_{n+1} = 0$ is the minimal canon of the order matrix, subject to the condition $\ell_{j_0} \ge \bar{\mathcal{O}}_{j_0,j_O}$. The result is then again a consequence of i). ■

The next examples illustrate some of the possible cases.

***Examples.* — 176)** Consider the system $x_1' - (x_2 - x_3) = 0$, $x_2' - (x_1 + x_3) = 0$, $x_3' - (x_1 - x_2) = 0$. With $j_0 = 1$, $\bar{I} = \emptyset$ and we have $\ell_1 = \bar{\mathcal{O}}_{1,1} = 2$, $\ell_2 = \bar{\mathcal{O}}_{2,1} = \ell_3 = \bar{\mathcal{O}}_{3,1} = 1$. Indeed, a resolvent representation may be computed using derivatives of the equations $P_i$ up to $\ell_i$. Moreover, we may check that $A_2 = x_2 - (x_1'' + x_1')/2 = (-P_1' - P_2 + P_3 - P_1)/2$ and $A_3 = x_3 - (x_1'' - x_1')/2 = (-P_1' - P_2 + P_3 + P_1)/2$ are expressed using derivatives of the $P_i$ up to order $\ell_i - 1$, according to assertion ii) of the last theorem. We also have $A_1 = P_1'' + P_2' - P_3' + A_2' + A_3'$.

**177)** We now look at the system $x_1'' - x_2 = 0$, $x_2' - x_1 = 0$ $x_3 - x_1' = 0$, $x_4 - x_3''' = 0$, $x_5 - x_4' = 0$, with $j_0 = 1$. Then, $\bar{I} = \{3, 4, 5\}$; we may check that $\bar{\mathcal{O}}_{\bar{I},\bar{I}} = 0$. In this case, $\ell = (3, 2, 4, 1, 0)$.

i) For the restricted system $P_I$, $\hat{\ell} = (1, 0)$, which is coherent with $A_1 = P_1' + P_2$.

ii) The computation of $A_2$ requires less derivatives as $A_2 = P_1$, which is coherent with the bound $\hat{\ell}_i - 1$: in the linear case, resolvent are also weak resolvents.

iii) To compute $A_3 = P_3$ no strict derivatives is required, for $A_4 = x_4 - x_1' = P_4 + P_3''' + P_2' + P_1''$, the bound $\ell_i - \ell_4$ of iv) is reached and for $A_5 = x_5 - x_2 = P_5 + P_4' + P_3^{(4)} + P_2'' + P_1''' + P_1$, the bound $\ell_i$ on the orders is reached too.

**178)** We go back to ex. 167. We have $\mathcal{O}_{1,1} = 2$, $\mathcal{O}_{2,1} = 1$, $\mathcal{O}_{3,1} = 0$, that corresponds to the orders of derivations requested to compute a resolvent char. set in ii). The computation of $B_2 = -2P_1$ and $B_3 = -2(P_1' + P_2)$ in a weak resolvent only requires $P_1$, $P_1'$ and $P_2$, which is coherent with the bound $\mathcal{O}_{i,j_0} - 1$ in iii).

**179)** With the system $x_1' - x_2 = 0$, $x_2 - x_1 = 0$, $\bar{\mathcal{O}}_{2,1} + a_{2,2} = 0$ and the condition in iii) is not satisfied; indeed $A_2 = P_2$.

**180)** We conclude these examples with $x_1^{(e_1)} - x_2'$, $x_2^{(e_2)} - x_2$, $x_3 - x_1^{(f_3)}$ and $x_4 - x_3^{(f_4)}$, for $w = x_1$. In this case, $I = \{1, 2\}$ and $\bar{I} = \{3, 4\}$, $\bar{\mathcal{O}}_{1,1} = e_2$, $\bar{\mathcal{O}}_{2,1} = 1$. We have $A_1 = x_1^{(e_1+e_2)} - x_1^{(e_1)} = P_1^{(e_2)} + P_2' - P_1$. If $f_3 + f_4 \le e_1 + e_2$, $\ell = (e_2, 1, f_4, 0)$, or else $\ell = (f_3 + f_4 - e_1, f_3 + f_4 - e_1 - e_2 + 1, f_4, 0)$.

We have $A_3 = x_3 - x_1^{(e_1+a-1)} = P_3 + P_1^{(f_3-e_1)} + \sum_{k=0}^{p-1} P_2^{(ke_2+a)} - P_1^{(a-1)}$, where $a, p \in \mathbf{N}$ are such that $pe_2 + a = f_3 - e_1$ and $0 \le a < e_2$, so that the bound of iv) is reached. We also get $A_4 = x_4 - x_1^{(e_1+a-1)} = P_4 + P_3^{(f_4)} + P_1^{(f_3+f_4-e_1)} + \sum_{k=0}^{p-1} P_2^{(ke_2+a)} - P_1^{(a-1)}$, where $a, p \in \mathbf{N}$ are such that $pe_2 + a = f_3 + f_4 - e_1$ and $0 \le a < e_2$

Choosing $x_1^{(e_1)} - \sum_{k=0}^{g} c_k x_2^{(k)}$ for $P_1$, where the $c_k$ are constants such that $c_g \neq 0$ and $\sum_{k=0}^{g} c_k \omega^k \neq 0$ for any $e_2^{\text{th}}$ root of unity $\omega$, the system admits a resolvent for $w = x_1$ and the bounds $\bar{\mathcal{O}}_{2,1} = g$ of ii) or $\bar{\mathcal{O}}_{2,1} = g - 1$ of iii) are reached for all $g \ge 0$, as well as all other bounds deduced from the canon $\ell = (e_2, g, f_4, 0)$, if $f_3 + f_4 \le e_1 + e_2$, or else $\ell = (f_3 + f_4 - e_1, f_3 + f_4 - e_1 - e_2 + g, f_4, 0)$.

All the bounds given in th. 175 are sharp. Example 117, with $w = x_1$ shows that the bound given in ii), *i.e.* in the case when $I = [1, n]$, is reached for the computation of the resolvent polynomial. Some situations with $\sharp I = \sharp \bar{I} = 2$ are covered in the last example. The construction of suitable classes of examples in all remaining cases is left to the reader.

The way in which Jacobi presents his results is precise enough to suggest a complete proof for the theorem, that seems to assume first the construction of a normal form for a Jacobi ordering, using the minimal canon. In our presentation, this intermediate step is bypassed by lemma 174 and th. 110. Then, the computation of the resolvent representation follows the ideas exposed in sec. 9. The orders of the maximal derivatives required are indeed evaluated as in subsec. 9.2, the basic idea provided by Jacobi being reproduced in the proof of lem. 169 ii). *normal form*

We see that it is possible to prove Jacobi's claim for the resolvent polynomial, even with the strong bound. Considering the resolvent char. set, things are more complicated: a slightly better bound is obtained for weak resolvents when $\bar{\mathcal{O}}_{i,j_0} \neq -\infty$, but things get worse in the general case. We do not know how Jacobi considered the case when a variable did not appear in some equation. We can only note that his statements are consistent with the weak bound and keep in mind that his approach was to work first with generic situations, rather than with the general case, so that he possibly neglected this issue. *resolvent polynomial* *weak bound*

To conclude this section, we may notice the following interpretation of the bound $\bar{\mathcal{O}}_{i_0,j_0}$: reducing to the linear case, we consider the contribution of $x_{j_0}$ in equation $P_i$ as some right hand term and rewrite it $L_i(x_j | j \neq j_0) = F_i(x_{j_0})$. Then $\bar{\mathcal{O}}_{i_0,j_0}$ corresponds to the maximal order of the system $L_i$, $i \neq i_0$ in the variables $x_j$, $j \neq j_0$ (see th. 175 v). So, $L_{i_0}$ and its derivatives up to order $\bar{\mathcal{O}}_{i_0,j_0}$ must be dependent, producing a relation of order $\mathcal{O}$ between second members. In the case $I = [1, n]$, this informal approach can become a rigorous proof for ii) iii) and v) in th. 175.

## Conclusion

The first draft of this paper was intended to follow the translations [39, 40]. During the final revision, I tried to keep length and technicality to acceptable limits, with a hope that a common presentation of an important corpus that emerged around 1840 could justify some exception to our contemporary standards. The plan and internal mathematical structure of the abandoned *Phoronomia* is lost and the available short fragments, remained scattered in a few posthumous publications, or some manuscripts from the Berlin archives, provide limited information. It is still possible to restore some obvious links between concepts or tools kept alive or slowly independently rediscovered by various mathematical communities: Jacobi's bound (Ritt [84], Kondratieva *et al.* [60]), shortest reduction (Shaleninov [91], Pryce [81]), changes of ordering (Boulier *et al.* [5]), differential resolvent ( Cluzeau and Hubert [12]). Jacobi's contributions to combinatorics also make him appear as a pioneer of graph theory, operational research or shortest path algorithms. Later developments of the subject ( Kőnig [62], Egerváry [23], Kuhn [64]) cast a new light on Jacobi's work. We have seen how Egerváry's covers could be used to design differential orderings and suggest original differential reductions. There is indeed a strong link between combinatorics and differential tools and the pioneering contribution of Jacobi to the tropical determinant is perhaps the first step of an recently emerging *tropical differential algebra* (see Aroca *et al.* [2] and the references therein).

We have seen that the corpus of results contained in Jacobi's posthumous manuscripts provide a large set of applicable methods for the resolution of ordinary differential systems. From the automatization of easy ideas, such as looking for block decompositions, to more sophisticated tools, allowing to produce simpler normal form reductions or better ways to perform changes of ordering, they can improve in many ways the existent computer algebra algorithms.

Jacobi's bound by itself is able to replace advantageously Ritt's analog of the Bézout bound, in all situations where it is proved, *i.e.* at this time quasi-regular components (Kondratieva *et al.* [60, 61]) or systems of two equations (Ritt [86, Ch. VII § 6 p. 136]). Its interest to produce sharper complexity bounds, upper and lower, is obvious. The recent work of Li *et al.* [70] on differential resultants is a nice illustration. A difference algebra analog also exists, that has been proved by Hrushowski [35].

Using these tools as widely as possible is a promising task, and generalizing them to arbitrary systems a challenging goal. This paper has no claim to exhaustivity and the proofs that were suggested by Jacobi's arguments are conjectural reconstructions, that are meant to encourage the reading of the original works.

# Primary material

## Manuscripts

References from *Nachlaß Karl Gustav Jacob Jacobi (1804-1851)*, Archiv der Berlin-Brandenburgischen Akademie der Wissenschaften, 28.11.1958, 34 p.

[II/13 a)] *Letter from* Sigismund Cohn *to* C.W. Borchardt. Hirschberg, August, 25th 1859, in German, 3 p. Quoted p. 4.

[II/13 b)] Carl Gustav Jacob Jacobi, Manuscript *De ordine systematis æquationum differentialium canonici variisque formis quas inducere potest*, folios 2186–96, 2200–2206. 35 p. Basis of Cohn's transcription. English translation in [39]. Quoted p. 11, 54, 74.

[II/23 a)] *Reduction simultaner Differentialgleichungen in ihre canonische Form und Multiplicator derselben.*, manuscript by Jacobi, pages: 2214–2237. Five different fragments: 2214–2216; 2217–2220 (§ 17-18); 2221–2225 (§ 17); 2226–2229; 2230–2232, 2235, 2237, 2236, 2238 (numbered from 1 to 13). English partial translation in [40]. Quoted p. 74.

[II/23 b)] *De aequationum differentialium systemate non normali ad formam normalem revocando*, manuscript by Jacobi p. 2238, 2239–2241, 2242–2251. 25 p. Envelop by Borchardt. The basis of [40]. Quoted p. 54.

## Journals and complete works

[Crelle 27] *Journal für die reine und angewandte Mathematik*, XXVII, Berlin, 1844. Contains [41, I].

[Crelle 29] *Journal für die reine und angewandte Mathematik*, XXIX, Berlin, 1845. Contains [41, II].

[Crelle 64] *Journal für die reine und angewandte Mathematik*, LXIV, (4), p. 297-320, Berlin, Druck und Verlag von Georg Reimer, 1865. Quoted p. 54.

[GW IV] *C.G.J. Jacobi's gesammelte Werke*, vierter Band, herausgegeben von K. Weierstrass, Berlin, Druck und Verlag von Georg Reimer, 1886. Contains [41, 43, 44, 45]. Quoted p. 10, 10.

[GW V] *C.G.J. Jacobi's gesammelte Werke*, fünfter Band, herausgegeben von K. Weierstrass, Berlin, Druck und Verlag von Georg Reimer, 1890. Contains [39, 40].

[VD] *Vorlesungen über Dynamik von C.G.J. Jacobi nebstes fünf hinterlassenen Abhandlungen desselben*, herausegegeben von A. Clebsch, Berlin, Druck und Verlag von Georg Reimer, 1866. Contains [40].

# References

[1] Adelson-Velsky (Georgii) and Landis (Evgenii), "An algorithm for the organization of information", *Doklady Akademiya Nauk SSSR* 146, 263–266, 1962; English translation, *Soviet Math.* 3, 1259-1263, 1962. Quoted p. 32.

[2] Aroca (Fuensanta), Garay (Cristhian) and Toghani (Zeinab), "The fundamental theorem of tropical differential algebraic geometry", *Pacific Journal of Math.*, 283, (2), 257–270, 2016. Quoted p. 94.

[3] Bass (Hyman), Connell (Edwin H.) and Wright (David), "The Jacobian conjecture: reduction of degree and formal expansion of the inverse", *Bull. Amer. Math. Soc.* (New Series), 7, (2), 287–330, 1982. Quoted p. 54.

[4] Bellman (Richard), "On a routing problem", *Quarterly of Applied Mathematics*, 16, 87–90, 1958. Quoted p. 4 and 38.

[5] Boulier (François), Lemaire (François) and Moreno Maza (Marc), "PARDI !", *ISSAC 2001*, ACM Press, 38–47, 2001. Quoted p. 58, 72, 81, and 94.

[6] Boulier (François), Lazard (Daniel), Ollivier (François) and Petitot (Michel), "Computing representations for radicals of finitely generated differential ideals", Special issue Jacobi's Legacy of *AAECC* 20, (1), 73–121, 2009. Quoted p. 49, 51, 61, 62, 72, and 73.

[7] Burkard (Rainer), Dell'Amico (Mauro) and Martello (Silvano), *Assignment problems*, SIAM, Philadelphia, 2012. Quoted p. 17 and 41.

[8] Carrá Ferro (Giuseppa), "A resultant theory for ordinary algebraic differential equations", In: Mora T., Mattson H. (eds) *Applied Algebra, Algebraic Algorithms and Error-Correcting Codes*, LNCS vol 1255, 55–65, Springer, Berlin, 1997. Quoted p. 86.

[9] Castro-Jiménez (Francisco-Jesus), "Calculs effectifs pour les idéaux d'opérateurs différentiels", *Actas de la II Conferencia Internacional de Geometría Algebraica*, La Rábida, Travaux en Cours 24, Hermann, Paris, 1987. Quoted p. 50 and 55.

[10] Castro-Jiménez, (Francisco-Jesus) and Granger, (Michel), "Explicit calculations in rings of differential operators", *Séminaires & Congrès*, 8, 89–128, SMF, 2004. Quoted p. 50 and 55.

[11] Chrystal (George), "A fundamental theorem regarding the equivalence of systems of differential equations, and its application to the determination of the order and the systematic solution of a determinate system of such equations", *Transactions of the Royal Society of Edinburgh*, vol. 38, n$^{o}$ 1, 163–178, 1897. Quoted p. 56.

[12] Cluzeau (Thomas) and Hubert (Évelyne), "Resolvent Representation for Regular Differential Ideals", *AAECC*, 13, (5), 395–425, 2003. Quoted p. 85 and 94.

[13] Cohn (Richard M.), "Order and dimension", *Proc. Amer. Math. Soc.* 87, (1), 1–6, 1983. Quoted p. 4, 51, and 55.

[14] D'Alfonso (Lisi), Jeronimo (Gabriela) and Solernó (Pablo), "On the complexity of the resolvent representation of some prime differential ideals", *Journal of Complexity*, 22, (3), 396–430, 2006. Quoted p. 84.

[15] D'Alfonso (Lisi), Jeronimo (Gabriela) and Solernó (Pablo), "A linear algebra approach to the differentiation index of generic DAE systems", *AAECC*, 19, (6), 441–473, 2008. Quoted p. 84.

[16] D'Alfonso (Lisi), Gabriela Jeronimo (Gabriela), Ollivier (François), Sedoglavic (Alexandre) and Solernó (Pablo), "A geometric index reduction method for implicit systems of differential algebraic equations", *J. Symb. Comput.* 46 (10) 1114–1138, 2011. Quoted p. 84.

[17] D'Alfonso (María Elisabet), "Métodos simbólicos para sistemas de ecuaciones álgebro-diferenciales", PhD thesis, University of Buenos Aires, 2006. Quoted p. 84.

[18] Della Dora (Jean), Dicrescenzo (Claire) and Duval (Dominique), "About a new method for computing in algebraic number fields", *Eurocal'85*, Springer, LNCS 204, 289-290, 1985. Quoted p. 82.

[19] Dijkstra (Edsger Wybe), "A note on two problems in connexion with graphs", *Numerische Mathematik*, 1, 269–271, Springer, 1959. Quoted p. 4 and 38.

[20] Dinic, (E.A.) [Dinic (Yefim)] and Kronrod (M.A.), "An algorithm for the solution of the assignment problem", *Soviet Math. Dokl.* 69 (6), 1324–1326, 1969. Quoted p. 31.

[21] Duan, (Ran) and Pettie, (Seth), "Approximating Maximum Weight Matching in Near-Linear Time", *Proceedings of the* $51^{rst}$ *Annual IEEE Symposium on Foundations of Computer Science (FOCS)*, 673–682, 2010. Quoted p. 42.

[22] Edmonds, (J.) and Karp, (R.M.), "Theoretical improvements in algorithmic efficiency for network flow problems", *J. ACM*, 19, 248–264, 1972. Quoted p. 31.

[23] Egerváry (Jenő), "Matrixok kombinatorius tulajdonságairól" [In Hungarian: On combinatorial properties of matrices], *Matematikai és Fizikai Lapok*, vol. 38, 1931, 16–28; translated by H.W. Kuhn as Paper 4, Issue 11 of *Logistik Papers*, Georges Washington University Research Project, 1955. Quoted p. 4, 12, 21, 22, and 94.

[24] Faugère (Jean-Charles), Gianni (Patricia), Lazard (Daniel) and Mora (Teo), "Efficient Computation of Zero-dimensional Gröbner Bases by Change of Ordering", *Journal of Symbolic Computation*, 16, (4), 329–344, 1993. Quoted p. 81.

[25] Ford, (Lester Randolph Jr.), "Network Flow Theory", Paper P-923, RAND Corporation, Santa Monica, California, 1956. Quoted p. 4 and 38.

[26] Fredman (Michael Lawrence) and Tarjan (Robert Endre), "Fibonacci Heaps and Their Uses in Improved Network Optimization Algorithms", *Journal of the ACM*, 34, (3), 596–615, 1987. Quoted p. 38 and 42.

[27] Frobenius (Ferdinand Georg), "Über Matrizen aus nicht negativen Elementen", *Sitzungsberichte der Königlich Preußischen Akademie der Wissenschaften zu Berlin*, 456–477, 1912. Quoted p. 20 and 30.

[28] Frobenius (Ferdinand Georg), "Über zerlegbare Determinanten", *Sitzungsberichte der Königlich Preußischen Akademie der Wissenschaften zu Berlin*, 274–277, 1917. Quoted p. 20.

[29] Gabow (Harold N.) and Tarjan (Robert Endre), "Faster scaling algorithms for network problems", *SIAM J. Comput.*, 18, 1013–1036, 1989. Quoted p. 41.

[30] Giusti (Marc), Lecerf (Grégoire) and Salvy (Bruno), "A Gröbner free alternative for polynomial system solving", *Journal of Complexity*, 17, (1), 154–211, 2001. Quoted p. 52 and 84.

[31] Golubitsky (Oleg), Kondratieva (Marina), Moreno Maza (Marc) and Ovchinnikov (Alexei), "Bounds and algebraic algorithms in differential algebra: the ordinary case", in Decker, W., Dewar, M., Kaltofen, E., Watt, S.M. (Eds.), *Challenges in Symbolic Computation Software*, N° 06271 in Dagstuhl Seminar Proceedings, Internationales Begegnungs und Forschungszentrum für Informatik (IBFI), Schloss Dagstuhl, Germany, 2007. Quoted p. 68.

[32] Grier (David Alan), *When computers were human*, Princeton University Press, 2005. Quoted p. 5.

[33] Hopcroft (John E.) et Karp (Richard M.), "An $n^{5/2}$ algorithm for maximum matchings in bipartite graphs", *SIAM Journal on Computing*, 2, (4), 225–231, 1973. Quoted p. 16, 21, 26, and 29.

[34] Houtain (Louis), "Des solutions singulières des équations différentielles", *Annales des universités de Belgique*[26], années 1851–1854, 973–1323. Quoted p. 49 and 72.

[26] Reference established from *Bull. Amer. Math. Soc.* 12 (1906), p. 212. In the copy I could consult, there is no publisher indication, and the handwritten date 1852.

[35] Hrushovski (Ehud), *The elementary theory of Frobenius automorphisms*, preprint, 1996, revised 2004, ..., 2014. https://arxiv.org/abs/math/0406514 Quoted p. 94.

[36] Hubert (Évelyne), "Essential components of an algebraic differential equation", *J. Symbolic Computation*, 21, 1–24, 1999. Quoted p. 49 and 72.

[37] Ibarra (Oscar H.) and Moran (Shlomo), "Deterministic and probabilistic algorithms for maximum bipartite matching via fast matrix multiplication", *Inform. Process. Lett.*, 13, (1), 12–15, 1981. Quoted p. 30.

[38] Jacobi (Carl Gustav Jacob), "De eliminatione variabilis e duabus aequationibus algebraicis", *Journal für die reine und angewandte Mathematik*, 15, 101-124, 1836. Quoted p. 5.

[39] Jacobi (Carl Gustav Jacob), "De investigando ordine systematis aequationum differentialum vulgarium cujuscunque ", [Crelle 65], 297-320, [GW V], 193-216. English translation, *AAECC*, 20, (1), 7–32, 2009. Quoted p. 4, 5, 10, 11, 14, 15, 16, 18, 26, 31, 35, 44, 54, 74, 75, 94, and 95.

[40] Jacobi (Carl Gustav Jacob), "De aequationum differentialum systemate non normali ad formam normalem revocando", [VD], 550–578 and [GW V], 485-513. English translation, *AAECC*, 20, (1), 33–64, 2009. Quoted p. 4, 9, 15, 19, 58, 66, 68, 74, 81, 84, 87, 89, 94, and 95.

[41] Jacobi (Carl Gustav Jacob), "Theoria novi multiplicatoris systemati aequationum differentialium vulgarium applicandi", first published in [Crelle 27] Heft III 199–268 (part I) and [Crelle 29] Heft III 213–279 and Heft IV 333–376 (part II), reproduced in [GW IV], 317–509. Quoted p. 4, 9, and 95.

[42] Jacobi (Carl Gustav Jacob), "Lettre adressée à M. le président de l'Académie des sciences", *Journal de math. pures et appliquées*, V, 1840, 350-351. Quoted p. 5 and 9.

[43] Jacobi (Carl Gustav Jacob), "Sur un nouveau principe de la mécanique analytique", *C. R. Acad. Sci. Paris*, XV, 1842, 202–205; [GW IV], 289–294. Quoted p. 7, 9, and 95.

[44] Jacobi (Carl Gustav Jacob), "Sul principo dell'ultimo multiplicatore e suo uso como nuovo principio generale di meccanica", *Giornale arcadico*, XCIX, 129–146, 1844; [GW IV] 511-522. Quoted p. 7 and 95.

[45] Jacobi (Carl Gustav Jacob), "Sur un théorème de Poisson", *C. R. Acad. Sci. Paris*, XI, 1840, p. 529; [GW IV], 143–146. Quoted p. 9 and 95.

[46] Jacobi (Carl Gustav Jacob), *Canon arithmeticus, sive talulæ quibus exhibentur pro singulis numeris primis...*, Berolini, typis academicis, 1839. Quoted p. 5.

[47] *Briefwechsel zwischen C.G.J. Jacobi und M.H. Jacobi*, herausgegeben von W. Ahrens, Abhandlungen zur Geschichte der Mathematischen Wissenschaften, XXII. Heft, Leipzig, Druck und Verlag von B.G. Teubner, 1907. Quoted p. 4 and 5.

[48] Janet (Maurice), "Sur les systèmes d'équations aux dérivées partielles", *Journal de mathématiques pures et appliquées*, 8e série, tome 3, 65–152, Gauthier-Villars, Paris, 1920. Quoted p. 81.

[49] Johnson (Donald B.), "Efficient algorithms for shortest paths in sparse networks", *Journal of the ACM*, 24, (1), 1–13, 1977. Quoted p. 38.

[50] Johnson (Joseph), "Differential dimension polynomials and a fundamental theorem on differential modules", *Am. J. Math.*, 91, 251-257, 1969. Quoted p. 49.

[51] Johnson (Joseph), "Kähler Differentials and Differential Algebra", *The Annals of Mathematics*, 2nd Ser., 89, (1), 92-98, 1969. Quoted p. 49.

[52] Johnson (Joseph), " A notion of regularity for differential local algebras", *Contributions to algebra*, A Collection of Papers Dedicated to Ellis Kolchin, Edited by Hyman Bass, Phyllis J. Cassidy and Jerald Kovacic, Elsevier, 211–232, 1977. Quoted p. 49 and 73.

[53] Johnson (Joseph), "Systems of *n* partial differential equations in *n* unknown functions: the conjecture of M. Janet", *Trans. of the AMS*, vol. 242, 329–334, 1978. Quoted p. 4 and 49.

[54] Jordan (Camille), "Sur l'ordre d'un système d'équations différentielles ", *Annales de la société scientifique de Bruxelles*, vol. 7, B., 127–130, 1883. Quoted p. 87.

[55] Jüttner (Alpár), "On the efficiency of Egerváry's perfect matching algorithm", TR-2004-13, Egerváry Research Group, Budapest, 2004. ISSN 1587–4451. Quoted p. 22.

[56] Keller (Ott-Heinrich), "Ganze Cremona-Transformationen" *Monatsh. Math. Phys.*, 47, 299–306, 1939. Quoted p. 54.

[57] Knuth (Donald E.), *The art of computer programming III. Sorting and searching*, Addison-Wesley, Reading, 1973. Quoted p. 32.

[58] Koenigsberger (Leo), *Carl Gustav Jacob Jacobi. Festschrift zur Feier der hundertsten Wiederkehr seines Geburstages*, B.G. Teubner, Leipzig, 1904, xviii, 554 p. Quoted p. 5.

[59] Kolchin (Ellis Robert), *Differential algebra and algebraic groups*, Academic Press, New York, 1973. Quoted p. 49 and 50.

[60] Kondratieva (Marina Vladimirovna), Mikhalev (Aleksandr Vasil'evich), Pankratiev (Evgeniĭ Vasil'evich), "Jacobi's bound for systems of differential polynomials" (in Russian), Algebra. M.: MGU, s. 79-85, 1982. Quoted p. 4, 49, and 94.

[61] Kondratieva (Marina Vladimirovna), Mikhalev (Aleksandr Vasil'evich), Pankratiev (Evgeniĭ Vasil'evich), "Jacobi's bound for independent systems of algebraic partial differential equations", *AAECC*, 20, (1), 65–71, 2009. Quoted p. 4, 20, 49, 50, 54, and 94.

[62] Kőnig (Dénes), "Gráfok és mátrixok", *Matematikai és Fizikai Lapok* 38: 116–119, 1931. Quoted p. 20 and 94.

[63] Kőnig (Dénes), *Theorie der endlichen und unendlichen Graphen*, Math. in Monogr. und Lehrb. XVI., Akad. Verlagsges., Leipzig, 1936; AMS Chelsey Publishing, Rhode Island, 2001. Quoted p. 20.

[64] Kuhn (Harold W.), "The Hungarian method for the assignment problem", *Naval Research Logistics Quarterly*, 2, 83–97, 1955. Quoted p. 4, 12, 21, and 94.

[65] Kuhn (Harold W.), "A tale of three eras: The discovery and rediscovery of the Hungarian Method", *European Journal of Operational Research*, 219, 641–651, 2012. Quoted p. 12 and 19.

[66] Munkres (James), "Algorithms for the assignment and transportation problems", *J. Soc. Industr. Appl. Math.*, 5 (1957), 32–38. Quoted p. 17.

[67] Landau (Lev Davidovich) and Lifshitz (Evgeny Mikhailovich), *Theoretical physics*, vol. 1, Mechanics, $3^{rd}$ ed., Pergamon Press, 1976. Quoted p. 65.

[68] Lando (Barbara A.), "Jacobi's bound for the order of systems of first order differential equations", *Trans. Amer. Math. Soc.*, 152, 119–135, 1970. Quoted p. 43 and 45.

[69] Le Gall, (François), "Powers of tensors and fast matrix multiplication", *Proceedings of ISSAC'2014*, ACM Press, 296–303, 2014. Quoted p. 41.

[70] Li (Wei), Yuan (Chun-Ming) and Gao (Xiao-Shan), "Sparse Differential Resultant for Laurent Differential Polynomials", *Foundations of Computational Mathematics*, 15, (2), 451–517, 2015. Quoted p. 86 and 94.

[71] Maclagan (Diane) and Sturmfels (Bernd), *Introduction to Tropical Geometry*, Graduate Studies in Mathematics 161, AMS, Providence, 363 p., 2015. Quoted p. 4 and 36.

[72] Martello, (Silvano), "Jenő Egerváry: from the origins of the Hungarian algorithm to satellite communication", *Cent. Eur. J. Oper. Res.*, 18, (1), 47–58, 2010. Quoted p. 4.

[73] textscMatera, (Guillermo) and Sedoglavic, (Alexandre), "Fast computation of discrete invariants associated to a differential rational mapping", *Journal of Symbolic Computation* 36, (3–4), 473–499, 2003. Quoted p. 81.

[74] Morrison (Sally), "The Differential Ideal $[P]\ :\ M^{\infty}$", *Journal of Symbolic Computation*, 28, 631–656, 1999. Quoted p. 51.

[75] Monge (Gaspard), "Mémoire sur la théorie des déblais et des remblais", *Histoire de l'Académie royale des Sciences*, [année 1781. Avec les Mémoires de Mathématique & de Physique pour la même Année] (2e partie) (1784) [Histoire: 34–38, Mémoires :] 666–704. Quoted p. 11.

[76] Nanson (Edward John), "On the number of arbitrary constants in the complete solution of ordinary simultaneous differential equations", *Messenger of mathematics*, vol. 6, 77–81, 1876. Quoted p. 87.

[77] Nucci (Maria Clara) and Leach (Peter G.L.), "Jacobi's last multiplier and symmetries for the Kepler problem plus a lineal story", *J. Phys. A: Math. Gen.*, 37, 7743–7753, 2004. Quoted p. 7.

[78] Nucci (Maria Clara) and Leach (Peter G.L.), "The Jacobi last multiplier and its applications in mechanics", *Physica Scripta*, 78, (6), 2008. Quoted p. 7.

[79] Ollivier, (François) and Sadik (Brahim), "La borne de Jacobi pour une diffiété définie par un système quasi régulier" (Jacobi's bound for a diffiety defined by a quasi-regular system), *Comptes rendus Mathematique*, 345, (3), 139–144, 2007. Quoted p. 49, 50, 54, and 57.

[80] Ollivier, (François), "Jacobi's bound and normal forms computations. A historical survey", *Differential Algebra And Related Topics*, World Scientific Publishing Company, 2008. ISBN 978-981-283-371-6 Quoted p. 5.

[81] Pryce (John D.), "A simple structural analysis method for DAEs", *BIT: Numerical Mathematics*, vol. 41, no. 2, 364–394, Swets & Zeitlinger, 2001. Quoted p. 4, 58, and 94.

[82] Rabinowitsch, (J. L.) [Rainich (George Yuri)], "Zum Hilbertschen Nullstellensatz", *Math. Ann.*, 102, p. 520, 1929. Quoted p. 52 and 59.

[83] Riquier (Charles), "Les systèmes d'équations aux dérivées partielles", Gauthier-Villars, Paris, 1910. Quoted p. 81.

[84] Ritt, (Joseph Fels), *Differential equations from the algebraic standpoint*, Amer. Math. Soc. Colloq. Publ., vol. xiv, A.M.S., New-York, 1932. Quoted p. 86 and 94.

[85] Ritt, (Joseph Fels), "Jacobi's problem on the order of a system of differential equations", *Annals of Mathematics*, vol. 36, 1935, 303–312. Quoted p. 4, 5, 43, and 55.

[86] Ritt, (Joseph Fels), *Differential Algebra*, Amer. Math. Soc. Colloq. Publ., vol. xxxiii, A.M.S., New-York, 1950. Quoted p. 7, 49, 50, 72, 73, 81, 85, and 94.

[87] Saltykow, (Nikolaï Nikolaïewitch), "L'Œuvre de Jacobi dans le domaine des équations différentielles du premier ordre", *Bulletin des sciences mathématiques*, Gauthier-Villars, Paris, 213–228, 1939. Quoted p. 5.

[88] Saltykow, (Nikolaï Nikolaïewitch), "Ordre d'un système d'équations différentielles ordinaires de la forme générale", *Bulletin (Académie serbe des sciences et des arts. Classe des sciences mathématiques et naturelles. Sciences mathématiques)*, n° 1, 141–148, 1952. Quoted p. 4.

[89] Sankowski (Piotr), “Maximum weight bipartite matching in matrix multiplication time”, *Theoretical Computer Science*, 410, 4480–4488, 2009. Quoted p. 41.

[90] Sedoglavic (Alexandre), “A probabilistic algorithm to test local algebraic observability in polynomial time”, *Journal of Symbolic Computation*, 33, (5), 735–755, 2002. Quoted p. 81.

[91] Shaleninov (A.A.), “Removal of topological degeneracy in systems of differential evolution equations”, *Automation and Remote Control*, 51, 1599–1605, 1991 (English version). Translation from *Avtomatika i Telemekhanika*, Issue 11, 163–170, 1990. Quoted p. 4, 58, and 94.

[92] Schrijver (Alexander), *Combinatorial Optimization. Polyhedra and Efficiency*, Springer, 2003. Quoted p. 20.

[93] Schrijver (Alexander), “On the history of combinatorial optimization (till 1960)”, *Handbook of Discrete Optimization*, K. Aardal, G.L. Nemhauser, R. Weismantel, eds., Elsevier, Amsterdam, 1–68, 2005. Quoted p. 11, 12, 20, and 41.

[94] Schrijver (Alexander), “On the history of the shortest path problem”, *Documenta Math.*, extra volume: Optimization Stories, 155–167, 2012. Quoted p. 38.

[95] Tomizawa, (N.), “On some techniques useful for solution of transportation network problems”, *Networks*, 1, 173–194, 1971. Quoted p. 31.

[96] Volevich, (Leonid Romanovich), “On general systems of differential equations”, *Dokl. Akad. Nauk SSSR*, 132, (1) 20–23, 1960. English translation *Soviet. Math.*, 1, 458–465, 1960. Quoted p. 4.

# Index

# Notations

Notations are listed by alphabetical order, using TEX spelling for non roman letters and symbols, *e.g.* look for *nabla* for $\nabla$, *lambda* for $\lambda$ or *odot* for $A \odot B$.

$\mathscr{A}$ A char. set with elements $A_i$. In sec. 10, $\mathscr{A}$ is a resolvent or a local resolvent char. set, def. 164.

$A_P$ Order matrix of the system $P$, def. 71.

$A^{\boxminus}$ Rectangular matrix completed with rows or columns of 0 to make it square, rem. 75.

$a_{i,j}$, def. 71.

$a_i$, $b_i$, footnote 19.

$\alpha, \beta$ Jacobi's minimal cover, def. 22.

$\alpha_P$, $\beta_P$ def. 76.

${}^{\circ}_{J}\alpha_P$, ${}^{\circ}_{J}\beta_P$, def. 124.

$\mathscr{B}$, $\mathscr{C}$ Char. set with elements $B_i$ or $C_i$, *e.g.* in sec. 8.1. In sec. 10, $\mathscr{B}$ is a weak resolvent char. set, def. 164.

$\mathscr{D}$, $\mathscr{D}_{\mathscr{P}}$ Ring of differential operators, def. 84.

$\mathrm{d}P$ Linearized system, p. 49.

$\mathrm{d}_{\mathscr{P}}$, def. 84.

$\Delta_{\sigma,k}$, th. 110.

${}^{v}_{P,}\Delta_{\sigma,k}$, def. 111.

$F_k$, th. 110.

$\mathscr{F}$ A differential field of char. 0.

$\mathscr{G}/\mathscr{F}$, $\mathscr{G}_{\mathscr{P}}$ Differential field extension defined by some prime differential ideal $\mathscr{P}$, def. 84.

$x_i \gg x_j$, p. 7.

$H_{\mathscr{A}}$ Product of initial and separants of a char. set $\mathscr{A}$, p. 7.

$\mathscr{I}_\sigma$, $\mathscr{I}_{\sigma,k}$, th. 110.

${}^{v}_{P}\mathscr{I}_\sigma$, def. 111.

$In_i = In_{A_i}$ Initial of $A_i$, $In_{\mathscr{A}}$ Product of initialss of $\mathscr{A}$, p. 7.

$\mathrm{J}$, def. 76.

${}^{\circ}_{J}\mathrm{J}_P$, def. 124.

$\kappa$, def. 96.

$M$, $\bar{M}$, $\hat{M}$, th. 175.

$\lambda$ Minimal canon, algo. 9.

$\lambda_P$, $\lambda_P^{\boxminus}$ def. 76.

$\lambda^{\boxminus}$ Minimal canon of $A^{\boxminus}$, def. 76.

${}^{\circ}_{J}\lambda_P$, def. 124.

$x_i \ll x_j$, p. 7.

$\mathscr{M}$, def. 84 and def. 96.

$\mathscr{M}_0$, def. 96.

$\mathrm{M}$ A totally ordered commutative monoïd, def. 1

$\mu, \nu$ A minimal cover, def. 19.

$n$ Number of variables.

$\nabla$ System determinant, def. 76.

${}^{v}_{P}\nabla_\sigma$, def. 111.

${}^{\circ}_{J}\nabla_P$, def. 124.

$\mathscr{O}_A$, Jacobi's number (tropical determinant) of the matrix $A$, def. 1

${}^{\circ}_{J}\mathscr{O}_P$, def. 124.

$\bar{\mathscr{O}}_{i,j}$, $\mathscr{O}_{I,J}$ Def. 168.

$A \odot B$, Tropical matrix multiplication, p. 36.

$\mathrm{ord}^{v} P$, $\mathrm{ord}^{\mathrm{J}}$, def. 76.

$\mathrm{ord}_{\mathscr{P}}$ Order with respect to a prime component, def. 96.

orders, p. 78.

$P$ A differential system in $\mathscr{F}\{x_1, \ldots, x_n\}$

$\mathscr{P}$ Prime differential ideal.

$\pi_A$ Path relation, def. 55.

$\mathscr{Q}$ Radical differential ideal, proof of th. 89.

$\mathscr{R}$, rem. 109.

$\mathrm{R}_k$, $\mathrm{R}_k^{v}$, th. 110, th. 139.

$s$ Number of equations in the system $P$.

$S_{s,n}$, def. 1.

$Sep_i = Sep_{A_i}$ Separant of $A_i$, $Sep_{\mathscr{A}}$ Product of separants of $\mathscr{A}$, p. 7.

$\sigma$ An injection $\mathbf{N}^s \mapsto \mathbf{N}^n$, def. 1.

$\Theta$ Commutative monoïd generated by the derivations, def. 84.

$|A|_T$, Tropical determinant, p. 36.

$\mathrm{Y}$ Set of derivatives of the $x_i$, def. 96 and lem. 85

$\upsilon_A$, Main derivative of $A$, lem. 85.

$x$, $y$, $z$ Variables.